\documentclass[12pt]{amsbook}
\usepackage{amsmath,amsfonts,amssymb,extarrows,geometry,mathdots,mathtools}
\usepackage[T1]{fontenc}
\usepackage[latin1]{inputenc}
\usepackage[english]{babel}
\usepackage[font=footnotesize,labelfont=bf,labelsep=colon]{caption}

\usepackage{fbox}

\usepackage{chngcntr}
\counterwithout{figure}{chapter}

\usepackage{pgfplots}
\pgfplotsset{compat = newest}
\usepackage{xcolor}
\usepackage{xparse}
\usepackage{epigraph}
\usepackage{csquotes}

\usepackage{tikz}
\usetikzlibrary{angles,quotes,arrows,arrows.meta,chains,matrix,calc,
intersections,positioning,fit,decorations,
decorations.pathmorphing,decorations.text}
\usepackage{tikz-3dplot}
\usepackage{tikz-cd}

\usepackage{enotez}
\let\footnote=\endnote
\setenotez{list-name = References and Notes for Further Study}

\newlength{\layerwd}
\newcounter{outermost}

\NewDocumentEnvironment{onion}{sm}{
    \begin{tikzpicture}
        \setlength{\layerwd}{#2}%
        \setcounter{outermost}{0}
}{%
    \foreach \A in {0,...,\theoutermost}{\draw[thick] (0,0) circle 
(\A*\layerwd+\layerwd);}
    \end{tikzpicture}
}
\NewDocumentCommand{\annulus}{sO{lightgray}mmmo}{%
    \filldraw[thick,fill=#2] (#4:#3*\layerwd) 
        arc [radius=#3*\layerwd, start angle=#4, delta angle=#5-#4] 
        -- (#5:#3*\layerwd+\layerwd) 
        arc [radius=#3*\layerwd+\layerwd, start angle=#5, delta angle=#4-#5] 
        -- cycle; 
    \pgfmathsetmacro{\tmp}{(#5-#4)/2 +#4} 
    \IfNoValueF{#6}{%
        \IfBooleanTF{#1}
        {%
            \begingroup
            
                \pgfmathsetmacro{\rpTF}{ifthenelse(\tmp>180,"false","true")}
                \def\\{\space} 
                \path[rotate=\tmp-180,postaction={
                    decorate,
                    decoration={
                        text along path,
                        raise=-3pt,
                        text align={align=center},
                        reverse path=\rpTF,
                        text=#6
                    }
                }] (0,0) circle (#3*\layerwd+0.5*\layerwd);
            \endgroup
        }
        {%
            \pgfmathsetmacro{\rpTF}{ifthenelse(\tmp>180,\tmp+90,\tmp-90)}
            \node[inner sep=0pt, 
            text width=#3*\layerwd*3+\layerwd,
            align=center,
            rotate=\rpTF,
            font=\footnotesize] at (\tmp:#3*\layerwd+0.5*\layerwd)
            {#6};
        }%
    }%
    \ifnum\theoutermost<#3\setcounter{outermost}{#3}\fi
}

\begin{document}

\title{THE CODE OF MATHEMATICS
 \\ From Truth, Proof, and Computability to Homotopy Type Theory}
 
\author{Stefan M\"uller-Stach}

\frontmatter

\maketitle

\renewcommand{\figurename}{Fig.}


\chapter*{Prefaces}

\ \\
{\bf Preface to the first edition}

\setlength{\epigraphwidth}{8cm}
\setlength{\epigraphrule}{0pt}
\epigraph{La musique est une math\'ematique myst\'erieuse dont les \'el\'ements 
participent de l'infini.}{Claude Debussy (Musica, May 1903).}

Our Earth is part of a gigantic universe\index{Universe}, the fundamental 
principles of which are explained by physics. Despite all the mysteries it still 
holds, we know a great deal about it. In mathematics, there are 
universes\index{Universe} of a completely different nature. They form infinite 
reservoirs of objects that are not found in reality. Together with the precise 
mathematical calculus, they unfold from them
a world of fascinating beauty and diversity.

The connection between mathematics and physics goes much deeper than this 
analogy. Numerous sciences cannot express their theories without 
the aid of mathematical structures and mathematical 
calculus forms the basis of digitisation and other technologies. 
Galileo Galilei\index{Galilei, Galileo} once called mathematics the language 
of nature. Eugene Wigner\index{Wigner, Eugene} spoke of the
\begin{quote} 
Unreasonable effectiveness of mathematics in the natural sciences.\footnote{See 
\cite{wigner}.}
\end{quote}
It was a dream far beyond that of Gottfried 
Wilhelm Leibniz\index{Leibniz, Gottfried Wilhelm} to construct a universal 
scientific language -- the lingua universalis\index{Lingua universalis} -- 
which generalises mathematics. As Umberto Eco\index{Eco, Umberto} and others 
have noted, this idea has never been fully realised. Its systematic limits 
remained -- despite important milestones in computer science -- 
unclear.\footnote{See \cite[Ch. 15]{eco1993}.}

In May 2018, I gave a lecture titled \enquote{Truth, Proof, Thought, Identity} 
at the Mainz Studium Generale. In it, I explained the basics 
of homotopy type theory\index{Homotopy type theory}, which builds on 
dependent type theory\index{Type theory} and can better grasp the 
concepts of equality\index{Equality}\index{Identity}, 
isomorphism\index{Isomorphism} and equivalence\index{Equivalence} 
through built-in topological\index{Topology} concepts than the traditional 
approach via set theory\index{Set theory}. Such ideas are related to 
fundamental questions of computer science, philosophy and physics. 

Since type theory\index{Type theory} resembles the code of a programming 
language, this approach enables the machine-checkable verification of proofs 
and algorithms. This development has the potential to possibly change the 
research, learning and publication culture far beyond mathematics, as 
Leibniz\index{Leibniz, Gottfried Wilhelm} had envisaged.

This book presents a detailed version of the lecture for an interested 
audience. To achieve good access, a natural structure of mathematics is 
outlined, which deviates from the usual curriculum. Since in particular the 
later chapters have a higher level of difficulty, numerous notes and references 
to further literature should facilitate reading and compensate for the missing 
details in the popular scientific presentation. The selected aspects span a 
historical arc from antiquity to the present day. This once again illustrates 
that scientific knowledge gain usually takes place on long time scales. 

This text is based on a machine English translation of the German first 
edition. The post editing work and all illustrations were done by myself. 
German quotes and titles of publications in the main text have been translated 
with references pointing to the original source. I would like to thank the 
following people for remarks: Lambert Alff, Nikoo Azarm, Carl-Heinz Barner, 
Manuel Blickle, Elke Brendel, Annika Denkert, Petra Gerster, J\"urgen Jost, 
Thomas Metzinger, Sieglinde M\"uller-Stach, Verena N\"orthen, Andreas R\"odder, 
Silvana R\"odder, Andreas R\"udinger, Peter Stoll, Thomas Streicher, Christian 
Tapp, Ulrich Volp and Rainer Wieland. 

\ \\
Stefan M\"uller-Stach, Mainz and Neustadt-Haardt, March 2024 

\ \\
{\bf Preface to the second edition} 

\ \\
The second English edition of this book contains, besides corrections of 
errors, numerous extensions to provide better access to the material. Hints 
were provided by Markus Blumenstock, Felix Cherubini, Verena N\"orthen, and 
Andreas R\"udinger. 

In Chapter 1 we added additional explanations in the sections on Platonism and 
realism. In Chapter 7 we have included additional material on homotopical 
categories and on recent developments in model-independent category theory. In 
Chapter 8 we extended the presentation of dependent type theory and its 
interpretations in homotopical categories and included an additional 
introduction to categories with universes and contextual categories. From 
Chapter 9 the corresponding sections were removed, so that it hardly requires 
expert knowledge anymore now.

Besides the machine translation of the first English edition, the author has 
not used any AI-tools for writing new text parts in this second English 
edition. The second English and German editions were created at the same time. 

\ \\
Stefan M\"uller-Stach, Mainz and Neustadt-Haardt, June 2026 

\chapter*{Prelude}

As a guide, we want to anticipate some themes and questions that we 
will encounter. Two fundamental questions are: 
\begin{quote}
What is a thought?  
\end{quote}
and 
\begin{quote}
Where are the abstract concepts of our thinking located? 
\end{quote} 
Despite numerous insights from philosophy and the life sciences, there are 
still no simple answers to these two questions about the nature of thoughts and 
the location of our thinking and consciousness. Leibniz\index{Leibniz, 
Gottfried Wilhelm} pointed out in his mill example\index{Mill 
example}\footnote{See \S 17 of Leibniz's \enquote{Monadology}\index{Leibniz, 
Gottfried Wilhelm} at www.projekt-gutenberg.org.} that any journey through our 
brain would only reveal a machine -- a mill\index{Mill example} from the inside 
-- without perceptions, consciousness or other qualia\index{Qualia} and emergent 
phenomena of life being recognisable.

In the philosophy of Plato\index{Plato} and its tradition, abstract objects and 
concepts were located in a world of ideas\index{Platonic idealism}\index{Plato} 
outside of physical reality. The Platonic world of ideas has a natural 
counter-position in nominalism\index{Nominalism}. An introductory question to 
this complex of topics is: 
\begin{quote}
What is a number? 
\end{quote} 
Numbers are universally known objects that do not exist in reality themselves, 
but only in the form of counts. We handle small numbers confidently, but very 
large numbers completely elude our imagination. The same applies to ideal 
geometric shapes, such as circles, which exist in their mathematical purity only 
approximately in reality. Such figures are special cases of a general concept of 
space. Algebraic structures like numbers can be assigned suitable spaces, so 
that we can consider the concept of number as a part of the concept of space. 
Thus, the nature of numbers and geometric objects is naturally connected with 
the following question:
\begin{quote}
What is the most general concept of space?
\end{quote}
We approach this question by dealing in detail with examples of topological 
spaces\index{Topology} and their generalisations. 

Since we are already deeply involved in mathematics with such discussions, we 
deal with its working methods and its influence:
\begin{quote} 
What is the working method of mathematics?
\end{quote}
and 
\begin{quote}
What role does it play in our culture? 
\end{quote}
In doing so, we exemplarily address some basic mathematical concepts and 
describe various useful algorithms\index{Algorithm} and exciting applications.

The foundation of mathematics is its syntax\index{Syntax}. This includes 
underlying deductive systems\index{Deductive system}, in which proofs can be 
conducted using a formal language. The following question can be asked:
\begin{quote}
Is provability\index{Provability} a kind of computability\index{Computability}? 
\end{quote} 
This is indeed correct. Every proof is a computation\index{Computability} in a 
deductive system\index{Deductive system}. In a certain sense, conversely, every 
computation is a proof for the statement that asserts the result of the 
computation. The question can therefore be answered in the 
affirmative in a very general sense. 

However, it remains unclear where the boundaries of 
computability\index{Computability} and provability\index{Provability} lie. In 
this context, a difficult question is the decision problem\index{Decision 
problem}\index{Entscheidungsproblem} (German Entscheidungsproblem) of 
Hilbert\index{Hilbert, David} and Ackermann\index{Ackermann, Wilhelm}: 
\begin{quote}
Can it be decided which propositions in an axiomatic theory are 
provable\index{Provability} by deductive methods? 
\end{quote} 
Propositions correspond to formulas without free variables in the formal 
language\index{Formal language} of a deductive system\index{Deductive system}. 
Sometimes they are also called sentences or theorems and often we refer to 
them informally as statements. A decision would be an 
algorithm\index{Algorithm} that always terminates and delivers the result $1$ 
if the proposition is provable\index{Provability} and $0$ otherwise. However, 
the answer to this question is a clear no, as Alan Turing\index{Turing, Alan} 
and Alonzo Church\index{Church, Alonzo} showed in their famous 
works\footnote{See \cite{church1936,turing1936}.} of 1936. In addition to Emil 
Post\index{Post, Emil} and Stephen Kleene\index{Kleene, Stephen}, they provided 
a precise concept of computability\index{Computability} and used a trick by 
G\"odel\index{Goed@G\"odel, Kurt}, in which mathematical propositions 
or Turing machines\index{Turing machine} are assigned a natural number as 
G\"odel number\index{Goed@G\"odel, Kurt}. This leads to the 
decision problem\index{Decision problem}\index{Entscheidungsproblem} and other 
undecidable problems\index{Undecidability} being reduced to a yes/no decision 
for natural numbers, which can be technically related to the 
undecidability\index{Undecidability} of the halting problem\index{Halting 
problem} for Turing machines\index{Turing machine}.

An open question is whether there is a concept of 
computability\index{Computability} beyond Turing computability\index{Turing 
machine}, often called hypercomputing\index{Hypercomputing}, and whether the 
human mind is superior to a computer or not:
\begin{quote}
What distinguishes human intelligence from a computer? 
\end{quote}
Alan Turing\index{Turing, Alan} himself thought intensively 
about such questions and developed the 
Turing test\index{Turing test}\footnote{See \cite{turing1950}.} in this context.

Part of the book is dedicated to the concepts of truth\index{Truth} and 
semantics\index{Semantics}: 
\begin{quote}
How can the concept of truth\index{Truth} be defined? 
\end{quote} 
and 
\begin{quote}
What role does semantics\index{Semantics} play in mathematics? 
\end{quote}
The first question was already asked in antiquity by Aristotle\index{Aristotle} 
and others and has led over many centuries to a multitude of 
attempts at explanation, of which the 
correspondence theory\index{Correspondence theory of truth} of 
truth\index{Truth} was the most widespread.

A possible verification of the concept of truth\index{Truth} is aimed at by 
Leibniz's\index{Leibniz, Gottfried Wilhelm} search for a universal 
scientific language, the lingua universalis\index{Lingua universalis}. In 
such a language, the verification of the truth\index{Truth} 
of statements in a calculus adapted to the subject of investigation would be 
possible by proof in a syntactic\index{Syntax} way. This approach is 
fundamentally pragmatic, as the concept of 
provability\index{Provability} in deductive systems\index{Deductive 
system} is well understood. The ideas of Leibniz\index{Leibniz, Gottfried 
Wilhelm} led to the emergence of mathematical logic.

Building on this, Alfred Tarski\index{Tarski, Alfred} provided a 
language-analytical definition of truth\index{Truth}\footnote{Although 
Tarski\index{Tarski, Alfred} claims to want to further develop the 
correspondence theory of truth, his approach is of a different nature, see 
\cite{tarski1935,tarski1969}. We cannot fully represent the extensive 
literature on the concept of truth. See, however,  
\cite{halbach1996,reclam,schrenk}.} in the 20th century, which works well at 
least for statements in deductive systems\index{Deductive system} with formal 
object languages\index{Object language}\index{Formal language} $L$, as they 
occur, for example, in mathematics. The 
self-referential everyday language of daily life, which carries its 
truth predicate\index{Truth predicate} within itself, 
does not meet these conditions due to the existence\index{Existence} of 
self-referential antinomies\index{Antinomy}. Tarski\index{Tarski, Alfred} used 
a metalanguage\index{Metalanguage} $M$, which is richer than $L$, to realise 
a truth predicate\index{Truth predicate} for statements in $L$ 
within $M$. The chosen mathematical model\index{Model theory} associated with 
$M$ is referred to as the semantics\index{Semantics} of $L$ and the possible 
transitions from $L$ to $M$ are referred to as 
interpretations\index{Interpretation}. In many cases, $M$ is 
given by the metalanguage\index{Metalanguage} of the axiomatic 
set theory\index{Set theory} with the 
Zermelo-Fraenkel axioms\index{Zermelo-Fraenkel axioms} and possibly 
additional axioms about Grothendieck universes\index{Universe}.

The most important example of Tarski's\index{Tarski, Alfred} method is the 
formal object language\index{Object language} $L_\mathrm{ar}$ of 
Dedekind-Peano arithmetic\index{Dedekind-Peano arithmetic}. In addition to 
the set-theoretic\index{Set theory} interpretation\index{Interpretation} in 
the standard model\index{Model theory} 
\[
\mathbb{N}=\{0,1,2,\ldots\},
\]
$L_\mathrm{ar}$ allows various interpretations\index{Interpretation} in 
non-standard models\index{Non-standard model} with differing properties. 
A statement in $L_\mathrm{ar}$ is called true\index{Truth} if 
its interpretation\index{Interpretation} in the respective 
set-theoretic\index{Set theory} model\index{Model theory} is fulfilled. This 
is proven with the help of the axioms of set theory\index{Set theory}. 
Tarski\index{Tarski, Alfred} has shown that the 
truth predicate\index{Truth predicate} for 
Dedekind-Peano arithmetic\index{Dedekind-Peano arithmetic} is not definable 
within $L_\mathrm{ar}$. In his theorem on the 
undefinability of truth\index{Truth}, he even proved that the 
G\"odel numbers\index{Goed@G\"odel, Kurt} of the true theorems in the 
standard model\index{Model theory} form an undecidable\index{Decidability} set. 

With Tarski's\index{Tarski, Alfred} concept of truth\index{Truth}, we can 
approach a new question: 
\begin{quote}
Is truth\index{Truth} the same as provability\index{Provability}? 
\end{quote} 
Since the truth\index{Truth} of a statement in $L$ is defined by 
provability\index{Provability} in a selected 
interpretation\index{Interpretation} with the help of a 
metalanguage\index{Metalanguage} $M$, 
this question seems to be clarified in a certain way, albeit not 
within $L$. Indeed, there are statements in certain deductive 
systems\index{Deductive system} $L$ that are neither provable\index{Provability} 
nor refutable in $L$ and whose validity can only be clarified in a selected 
metalanguage\index{Metalanguage}. This was proven by Kurt 
G\"odel\index{Goed@G\"odel, Kurt} in his famous 
incompleteness theorem\index{Incompleteness theorem}.\footnote{See  
\cite[Vol. I, 1931]{goedel}.} The 
continuum hypothesis\index{Continuum hypothesis} falls into a similar 
category, as it cannot be decided with the usual Zermelo-Fraenkel 
axioms\index{Zermelo-Fraenkel axioms} of set theory\index{Set theory}. Such 
statements are called undecidable\index{Undecidability}, even though the 
concept of undecidability\index{Undecidability} is used in a different way than 
in the decision problem\index{Decision problem}\index{Entscheidungsproblem}. 
The answer to our question is thus: 
\begin{quote}
Truth\index{Truth} depends on an interpretation\index{Interpretation}. 
The proof is carried out in the corresponding 
metalanguage\index{Metalanguage}. 
\end{quote} 
Because of these observations, it makes sense to consider hierarchies of 
extensions of formal languages\index{Formal language}.  

Mathematical theories are usually formulated in the language of predicate 
logic and set theory\index{Set theory}. However, this 
is not the only possibility. The search for alternatives to traditional 
set theory\index{Set theory} has been approached on the one hand in the form 
of (dependent) type theory\index{Type theory} through the work of the logician 
Per Martin-L\"of\index{Martin-L\"of, Per} and on the other hand through 
higher categories\index{Category theory} and 
infinity categories\index{Infinity category}, building on ideas from 
Alexander Grothendieck\index{Grothendieck, Alexander}.\footnote{See  
\cite{grothendieck,martin-loef1984}.} This leads us to an overarching theme:
\begin{quote} 
Which form of the foundations of mathematics is most suitable for structuralist 
thinking? 
\end{quote}
By foundations of mathematics, we understand thought systems that  
can describe large parts of mathematics. Vladimir 
Voevodsky\index{Voevodsky, Vladimir} and others have significantly further 
developed the theory of Martin-L\"of\index{Martin-L\"of, Per}. This research 
direction, also known as homotopy type theory\index{Homotopy type theory}, has 
incorporated topological\index{Topology} and homotopy 
theoretical\index{Homotopy theory} aspects into its basic concepts. Through 
this artifice, type theory\index{Type theory} has become a mature theory and 
can handle the concept of equality\index{Equality}\index{Identity} and its 
variants isomorphism\index{Isomorphism}, symmetry and 
equivalence\index{Equivalence} more adequately than, for example, set 
theory. The common concept of equality\index{Equality}\index{Identity} has 
proven to be too rigid over time. The question behind this is:
\begin{quote} 
What does the equation $A=A$ mean? 
\end{quote}
This question, which is also significant in philosophy, leads in type 
theory\index{Type theory} to the definition of an identity 
type\index{Identity type}, which is the starting point for new forms of 
equality\index{Equality}\index{Identity}. 

Dependent type theory\index{Type theory} forms through its deductive 
system\index{Deductive system} the basis of new software tools, which allow the 
provable\index{Provability} verification of proofs and algorithms and in the 
future will enable both powerful assistance systems for research and teaching 
as well as likely change publication practices. Large parts of mathematics have 
already been formalised in such systems. 

Through type theory\index{Type theory} and category 
theory\index{Category theory}, new concepts emerge for the foundations of 
mathematics. Often a given dependent type theory\index{Type theory} is viewed 
as an object language\index{Object language} and a (higher) 
categorical\index{Category theory} semantic\index{Semantics} 
interpretation\index{Interpretation} is considered. 
Through such interpretations\index{Interpretation}, the 
traditional mixing of syntax\index{Syntax} and semantics\index{Semantics} 
in set theory\index{Set theory} can be resolved. Going even further, 
these three foundations of mathematics can each be viewed as 
an object language\index{Object language} or metalanguage\index{Metalanguage} 
and are mutually interpretable. From such an 
overarching perspective, a certain disenchantment of the concept 
of semantics\index{Semantics} results, which is reduced to the 
syntax\index{Syntax} of a chosen metalanguage\index{Metalanguage}. On the 
other hand, the mathematical concept of truth\index{Truth} thereby receives a 
solid foundation that can no longer be relativised. 

As a result of our considerations, a new perspective on mathematics emerges. It 
opens up a more structuralist understanding than before, and significantly 
expands the possibilities for the verification of proofs\index{Provability} and 
algorithms\index{Algorithm}. This development goes beyond mathematics and makes 
a relevant contribution to Leibniz's\index{Leibniz, Gottfried Wilhelm} goal of 
the lingua universalis\index{Lingua universalis}, and thus to the way in which 
we work in science in the future.

\tableofcontents

\mainmatter

\chapter{Fundamental Questions}

What is mathematics and what are its subjects? This question is not easy to 
answer. From observing the world, we have always drawn inspiration 
for mathematical concepts. Essentially, however, it is an a priori science, as 
it is based neither on experience nor on other assumptions. Its simplest 
features seem to be evolutionarily ingrained in our brain. In 
addition to the use of words, also symbols are necessary for this kind of human 
thinking. In his \enquote{Dialogus} from 1677, Leibniz\index{Leibniz, Gottfried 
Wilhelm} wrote fittingly:
\begin{quote}
Cogitationes fieri possunt sine vocabulis ... At non sine aliis 
signis.\footnote{See \cite[Ser. VI, Vol. 4A, No. 8]{leibniz}.}
\end{quote}

Based on logic,\footnote{Logical concepts are part of mathematics.} which has 
served as a tool for proofs since antiquity, mathematical concepts are 
formulated in deductive systems\index{Deductive system}. Their formal 
languages\index{Formal language} must be rich enough to describe all 
mathematical objects that underlie the individual concepts being considered. 
Such systems are called foundations of mathematics. Axiomatic set 
theory\index{Set theory} is the most widespread form among them. Two other 
forms of foundations are category theory\index{Category theory} and type 
theory\index{Type theory}. 

When the first arithmetic and geometric concepts emerged in our cultures, 
the question already arose as to what form of existence\index{Existence} the 
underlying objects possess and in what way they can be identified with each 
other. The Platonic world of ideas\index{Platonic idealism} postulated an 
existence\index{Existence} at a specially designated place. It is closely 
related to the concept of equality\index{Equality}\index{Identity} among 
abstract objects. The aim of better understanding the role of different forms 
of equality\index{Equality} and even more general equivalence 
concepts\index{Equivalence} has triggered promising developments in research in 
recent mathematics.

\section{Mathematical Objects and Their Identification}

We will first consider abstract objects and illustrate the mathematical 
way of thinking using fundamental arithmetic and geometric concepts 
to lead from there to deeper questions about the foundations of 
mathematics.

\medskip
If you believe you know what concept the natural numbers 
\[
0,1,2,3,4,5,\ldots 
\]
and their arithmetic as a whole form, what the concept of a number is or what 
kind of object the single number
\[
5
\]
represents, then you may be mistaken. Even mathematically trained 
people have always had difficulty explaining this. It was not until the 19th 
century that Georg Cantor\index{Cantor, Georg}, Richard 
Dedekind\index{Dedekind, Richard} and Gottlob Frege\index{Frege, Gottlob} began 
to define the natural numbers precisely and to prove their 
desired properties. 

The notation for the number $5$ varied in different world languages and 
throughout history. This is not a problem, because we all -- even small 
children, some animal species and people in cultures without significant 
schooling -- have a sense for small numbers and elementary 
mathematical problems.\footnote{See the neuroscientific research 
by Stanislas Dehaene\index{Dehaene, Stanislas}.} Any intelligent 
extraterrestrial being that would visit us would probably 
also have a concept of it, because it would be able to make a connection when 
looking at our hands with their $5$ fingers and would probably have a 
designation for this number, depending on what form of 
communication it uses. But is this the same idea for all people and are all 
avatars of the number $5$ somehow the same? In other words, does the 
equation\index{Equality}\index{Identity}
\[
5=5 
\]
always hold? Our intuition tells us that this is correct and all variants of 
the number $5$ can be identified. This impression probably comes from the fact 
that we all mean the same thing in our daily dealings with it and can 
communicate perfectly well about small numbers. Cognitive studies on early 
childhood suggest that a sense for small numbers of objects is innate or can at 
least be learned quickly, similar to language acquisition. However, we do not 
understand why these facilities are present in the brain and how they are laid 
out there.

To make matters worse, there are complicated additive and multiplicative 
relationships between the numbers that occur within arithmetic. The equations 
\[
5=1+1+1+1+1 
\]
or 
\[
5=3+2
\]
are to be regarded as non-trivial statements about the number $5$. Here, a 
big problem apparently arises, namely the question of the identifiability of 
mathematical objects. The terms equality\index{Equality}, 
identity\index{Identity}, isomorphism\index{Isomorphism} and 
equivalence\index{Equivalence} are commonly used in mathematics for this. 
The equation 
\[
5=5 
\]
and especially the different representations of the number $5$ as sums of other 
smaller numbers touch deeper levels of mathematics. Richard 
Dedekind\index{Dedekind, Richard} and Gottlob Frege\index{Frege, Gottlob} 
recognised that natural numbers have infinitely many different realisations and 
equations like $5=3+2$ are in reality subtle assertions of 
equivalences\index{Equivalence} between different objects, which need to be 
mathematically explained and proven. Frege\index{Frege, Gottlob} attempted to 
solve the definition of natural numbers through his method of abstraction. 
Georg 
Cantor\index{Cantor, Georg} developed further concepts of transfinite 
numbers.\footnote{Contained in \cite{cantor,dedekind,frege1884,frege1893}.}

The strongest form of equality\index{Equality}\index{Identity}, i.e., sameness, 
such as $5=5$, we will call judgemental (or definitional)
equality\index{Equality}\index{Identity}. Other forms of 
equality\index{Equality}\index{Identity}, such as $5=3+2$, are called 
propositional equality\index{Equality}\index{Identity}. Such different variants 
of equality\index{Equality}\index{Identity} become a subject of research in 
modern mathematics. This leads -- besides a more precise definition of the 
concept of equality\index{Equality}\index{Identity} -- to a better 
understanding of generalisations such as the concepts of 
isomorphism\index{Isomorphism}, equivalence\index{Equivalence} and 
univalence\index{Univalence}. 

In everyday life and in other sciences, the concept of 
equality\index{Equality}\index{Identity} is also important and has very 
pragmatic aspects. However, there are also gradations there. One aspect of this 
is the recognisability of objects and people. In particular, the word 
identity\index{Identity}, and less often the word equality\index{Equality}, is 
used in connection with the identification\index{Identity}\index{Equality} of 
people by means of IDs or fingerprints. A person changes a lot in life. Parts 
of our bodies have a shorter lifespan than we do. When we go to the 
hairdresser, or lose some hair or skin cells in some other way, no one will 
doubt that we have remained equal as a person and have retained our 
identity\index{Identity}\index{Equality}, for a human being is more than a 
collection of individual objects. In addition to the \enquote{external}, a 
person changes over time also as a personality and is longer the same person in 
characteristics and appearance. After abrupt changes, the 
externally identical person may even no longer be equal as a person. The 
concept of equality\index{Identity}\index{Equality} in humans 
thus also touches on the psychological and social aspects of life and 
is, as in mathematics, to be distinguished from the concept of the same, which 
means complete agreement.

A composed piece of music which is provided with a fixed score  
will be played slightly differently at each performance. Nevertheless, it 
is equal as a piece of music that is being listened to. So there is a kind of 
concept of equality\index{Equality}\index{Identity} for pieces of music and 
musical performances which is not based on complete agreement. This is even 
more liberal in jazz compositions, where improvisation is done over given 
passages.\footnote{See \cite{feige} on the philosophy of jazz.} Here too, we 
never hear the same performance, but the concept of 
equality\index{Equality}\index{Identity} is 
more generously defined. Everyone immediately understands this example and 
recognises the same pieces of music without any problems. The same applies to 
other works of art. Every book, every piece of music and picture, which 
develops in its creation, remains equal from a certain stage, even if 
certain parts are still changed. Once an 
identity\index{Identity}\index{Equality} of the artwork has emerged, 
it remains preserved in the future.

\section{The Concept of Equality\index{Equality}\index{Identity} in 
Frege\index{Frege, Gottlob} and Leibniz\index{Leibniz, Gottfried Wilhelm}}

Frege\index{Frege, Gottlob} dealt with the concept of 
equality\index{Equality}\index{Identity} in his influential essay \enquote{On 
sense and meaning}.\footnote{Contained in \cite{frege1962}.} In it, he found a 
wonderful geometric example. He used triangle geometry and considered the three 
medians in a triangle. In his example, the assertion is that the 
three intersection points of any two pairs of medians in the triangle coincide. 
This is equivalent to an equation\index{Equality}\index{Identity} $A=B$ between 
two such points (see Fig.~\ref{frege_example}). 

\begin{figure}[ht!]
\begin{tikzpicture} 
\pgfmathsetmacro{\xA}{0}
\pgfmathsetmacro{\yA}{0}
\pgfmathsetmacro{\xB}{5}
\pgfmathsetmacro{\yB}{0}
\pgfmathsetmacro{\xC}{3}
\pgfmathsetmacro{\yC}{2}

\coordinate (A) at (\xA,\yA);
\coordinate (B) at (\xB,\yB);
\coordinate (C) at (\xC,\yC);

\pgfmathsetmacro{\xW}{(\xA + \xB + \xC)/3}
\pgfmathsetmacro{\yW}{(\yA + \yB + \yC)/3}
\coordinate[label=above left:$A$] (W) at (\xW,\yW);

\draw[path picture={
    \foreach \p in {A,C}\draw[color=black,line 
width=0.4mm,shorten >=-30cm](\p)--(W);
    \draw[dashed,shorten >=-30cm](B)--(W);
  }]
  (A)--(B)--(C)--cycle;
\node at (6,1) {$A=B$};

\pgfmathsetmacro{\xA}{7}
\pgfmathsetmacro{\yA}{0}
\pgfmathsetmacro{\xB}{12}
\pgfmathsetmacro{\yB}{0}
\pgfmathsetmacro{\xC}{10}
\pgfmathsetmacro{\yC}{2}

\coordinate (A) at (\xA,\yA);
\coordinate (B) at (\xB,\yB);
\coordinate (C) at (\xC,\yC);
\pgfmathsetmacro{\xW}{(\xA + \xB + \xC)/3}
\pgfmathsetmacro{\yW}{(\yA + \yB + \yC)/3}
\coordinate[label=above right:$B$]  (W) at (\xW,\yW);

\draw[path picture={
    \foreach \p in {B,C}\draw[color=black,line width=0.4mm,shorten 
>=-30cm](\p)--(W);
    \draw[dashed,shorten >=-30cm](A)--(W);
  }]
  (A)--(B)--(C)--cycle;
\end{tikzpicture}
\caption{\label{frege_example}In Frege's\index{Frege, Gottlob} example, the 
intersection points $A$ and $B$ of two different pairs of medians in the 
triangle coincide.}
\end{figure}

The proof for this arises from the fact that the three medians of a triangle 
always intersect at a common point, namely the centroid. Both $A$ and $B$ are 
therefore identical with this point and thus $A$ and $B$ are identical 
themselves. The equation\index{Identity}\index{Equality} $A=B$ therefore 
expresses the correctness of the geometric theorem of the centroid and is not a 
random identity\index{Identity}\index{Equality} of two a priori different 
points.

By the way, you may have noticed that in the above proof the two triangles in 
the figure are supposed to symbolise the same triangle, although they are only 
equal. To better see the respective intersection points, two 
identical\index{Equality}\index{Identity} copies of the triangle were drawn, a 
common trick. Our thinking thus uses variants of the concept of 
equality\index{Equality}\index{Identity} without hesitation.

In the same essay, Frege\index{Frege, Gottlob} pointed out that in this example 
the two points $A$ and $B$ coincide and thus have the same meaning (i.e., 
reference or denotation), but that the respective sense, which lies in the 
definition of the points, is different. The difference between sense and 
meaning of signs, in his opinion, is expressed by the 
equation\index{Identity}\index{Equality}\footnote{See the articles on sense 
and meaning in \cite{reclam,schrenk}.} 
\[
A=B. 
\]
In his text, Frege\index{Frege, Gottlob} gives further examples of this. One of 
them comes from the two designations morning star and evening star for the 
planet Venus. The equation 
\[
\text{morning star }=\text{ evening star} 
\]
expresses the coinciding meaning of both as planet Venus, while the sense on 
both sides of the equation differs by the time of day during observation. This 
distinction leads to a differentiation of the concept of 
equality\index{Equality}\index{Identity} in sense and meaning, or equivalently 
in intensional\index{Intensionality} and extensional\index{Extensionality} 
ways. Such investigations form the beginning of a philosophy of language with 
Frege\index{Frege, Gottlob}, which in the course of his text leads to the 
consideration of thoughts and statements and their truth\index{Truth}.

Even before Frege\index{Frege, Gottlob}, Gottfried Wilhelm 
Leibniz\index{Leibniz, Gottfried Wilhelm} introduced a concept of 
equality\index{Identity}\index{Equality} that still plays an important role in 
mathematics and elsewhere today. In his own words, he said: 
\begin{quote} 
Eadem sunt quorum unum potest substitui alteri salva veritate.\footnote{See 
\cite[Ser. VI, Vol. 4A, No. 178]{leibniz} and \cite{lolli2017}.}
\end{quote}
Translated into today's language, the equality\index{Identity}\index{Equality} 
$A=B$ applies to two objects if in all investigations the object $A$ can be 
replaced by $B$ while maintaining the truth\index{Truth}. This is therefore a 
substitutional concept of equality\index{Identity}\index{Equality}, which is 
weaker than the concept of equality\index{Identity}\index{Equality} of the
indistinguishability and additionally presupposes a concept of 
truth\index{Truth}. It becomes a rule of inference in formal mathematical 
calculi -- called Leibniz's invariance rule\index{Leibniz's invariance rule} -- 
in the sense that identical\index{Identity}\index{Equality} arguments in 
relations are interchangeable. This will encounter us later in our 
considerations of type theory\index{Type theory} by Per 
Martin-L\"of\index{Martin-L\"of, Per}.

Leibniz\index{Leibniz, Gottfried Wilhelm} also defined a broader concept of 
similarity or equivalence\index{Equivalence}, which was extremely clear and 
forward-looking for its time. He distinguished between the similarity and 
congruence of two geometric figures:
\begin{quote}
Similia sunt quae singulatim discerni non possunt.\footnote{See \cite[Ser. VII, 
Mathesis, No. 57]{leibniz}.}
\end{quote}
Translated into modern language, two figures are therefore similar if they are 
indistinguishable when viewed separately. Leibniz\index{Leibniz, Gottfried 
Wilhelm} thus differentiated between the size and shape of a figure. When two 
figures are viewed separately, their size recedes into the background and only 
their shape remains. This means that in geometric figures there is not only the 
strong concept of equality\index{Equality}\index{Identity}, i.e., congruence, 
but also the weaker concept of similarity or equivalence\index{Equivalence}, 
which is just as significant. This example shows that Leibniz\index{Leibniz, 
Gottfried Wilhelm} played a role in the discovery of the concept of 
equivalence\index{Equivalence}.\footnote{See \cite{derisi2007,jost2019} for 
Leibniz's\index{Leibniz, Gottfried Wilhelm} ideas in general.}

The examples discussed so far are very illustrative and unproblematic. So what 
is the actual problem of equality\index{Identity}\index{Equality} or 
equivalence\index{Equivalence}? Mathematical objects such as vector spaces, 
groups or manifolds are not unique\index{Non-uniqueness problem}. 
Isomorphisms\index{Isomorphism} and equivalences\index{Equivalence} identify 
them in mathematics and in 
applications in the form of structure-preserving transformations. For example, 
space-time theories in physics are modelled on manifolds, but the physical 
reality is of course independent of the chosen modelling. The fundamental 
problem is that it is not easy to decide for two given objects whether they are 
isomorphic\index{Isomorphism} or equivalent\index{Equivalence}. A well-known 
example of this is the homeomorphism problem\index{Homeomorphism} for two 
topological\index{Topology} spaces, which is undecidable\index{Undecidability} 
in the sense of computability theory\index{Computability}. Even if 
isomorphisms\index{Isomorphism} and equivalences\index{Equivalence} from 
previous considerations are available, they may not be able to be specified 
concretely or this information may be lost in the course of proofs due to lack 
of suitable bookkeeping. Recent mathematical research is trying to remedy such 
difficulties.

\section{The Concept of Equality\index{Equality}\index{Identity} in Philosophy}

We want to take a small detour into the field of philosophy. The problem of 
equality\index{Equality}\index{Identity} is not only a question of mathematics, 
but also a deep problem of philosophy. Martin Heidegger\index{Heidegger, 
Martin} gave a lecture titled \enquote{The theorem of the identity} in Freiburg 
(June 27th 1957), in which the question of 
identity\index{Identity}\index{Equality} is brought into connection with 
existence\index{Existence}, i.e., being:\footnote{See the articles on being in 
\cite{reclam,schrenk}.} 
\begin{quote}
The principle of identity is commonly formulated as: $A=A$. The 
principle is considered the supreme law of thought. We try to think about this 
principle for a while. Because we want to learn from the principle what 
identity is ... What the principle of identity, heard from its basic tone, 
states, is exactly what all Western thinking thinks, namely this: 
The unity of identity forms a basic feature in the being of 
beings.\footnote{Machine translation of \cite{heidegger}.}
\end{quote}
This essay by Heidegger\index{Heidegger, Martin} is not easy to 
understand and delves deeply into the world of metaphysics\index{Metaphysics}. 
At the end of his text, Heidegger\index{Heidegger, Martin} refers to 
Parmenides\index{Parmenides} and quotes him with the words:
\begin{quote}
The same is thinking as well as being.\footnote{Machine translation of 
\cite{heidegger}.}
\end{quote}
So, Heidegger\index{Heidegger, Martin} is concerned with the 
identity\index{Identity}\index{Equality} of the human being with his being and 
thinking. Behind the concepts of being and self is the question of the nature 
of human consciousness.  

There is a remarkable final remark in this text which classifies 
computability\index{Computability} in contrast to thinking:
\begin{quote}
Today, the thinking machine calculates thousands of relationships in a second. 
They are essenceless despite their usefulness.\footnote{Machine translation of 
\cite{heidegger}.}
\end{quote}
Behind this is the question of whether human thinking fundamentally differs 
from computability\index{Computability}. Heidegger\index{Heidegger, Martin} was 
convinced that our thinking has a completely different essence than the 
calculations of the then emerging computers. 

In addition to Heidegger\index{Heidegger, Martin} and 
Leibniz\index{Leibniz, Gottfried Wilhelm}, Johann Gottlieb 
Fichte\index{Fichte, Johann Gottlieb} dealt with the question of 
equality\index{Identity}\index{Equality}. He pointed out the deep content of 
the 
equation $A=A$ and emphasised, similar to Henri Poincar\'e\index{Poincar\'e, 
Henri}, that this statement presupposes a judgement\index{Judgement} of a 
subject.\footnote{See \cite{fichte1794,tasic2012}.} 

The English philosophers John Locke\index{Locke, John} and David 
Hume\index{Hume, David} also thought about 
equality\index{Equality}\index{Identity}. Thus, Locke\index{Locke, John} writes 
in his book \enquote{An essay concerning humane understanding} from 1690:
\begin{quote}
Another occasion the mind often takes of comparing, is the very being of 
things, when, considering anything as existing at any determined time and 
place, we compare it with itself existing at another time, and thereon form the 
ideas of identity and diversity.\footnote{See \cite[Book II, 
Ch. 27]{locke1690}.}
\end{quote}
Hume\index{Hume, David} encountered similar questions in his essay \enquote{Of 
personal identity} from 1739:
\begin{quote}
There are some philosophers, who imagine we are every moment intimately 
conscious of what we call our self; that we feel its existence and its 
continuance in existence; and are certain, beyond the evidence of a 
demonstration, both of its perfect identity and simplicity ... Unluckily all 
these positive assertions are contrary to that very experience, which is 
pleaded from them, nor have we any idea of self after the manner it is here 
explained ... The mind is a kind of theatre, where several perceptions 
successively make their appearance; pass, re-pass, glide away, and mingle in 
an infinite variety of postures and situations. There is properly no simplicity 
in it at one time, nor identity in different; whatever natural propension we 
may have to imagine that simplicity and identity. The comparison of the theatre
must not mislead us. These are only the successive perceptions that constitute 
our mind; nor do we have the most distant notion of the place where these 
scenes are represented, or of the materials of which it is composed ... We have 
a distinct idea of an object that remains invariable and uninterrupted through 
a supposed variation of time; and this idea we call that of identity or 
sameness ... that all the nice and subtle questions concerning personal identity 
can never possibly be decided, and are to be regarded rather as grammatical 
than as philosophical difficulties. Identity depends on the relations of ideas; 
and these relations produce identity.\footnote{See 
\cite[Book I, Sec. 6]{hume}.} 
\end{quote}
These sentences show what a critical mind David Hume\index{Hume, David} was. 
He examined concepts in the finest detail and made detailed judgements, which 
were always committed to scepticism. Thus, the first sentence questions the 
whole idea of consciousness and separates it from the self of humans. 
Hume\index{Hume, David} noted that man's assessment of his self might be more 
of a narrative. Indeed, the puzzle of the self and consciousness remains 
unsolved to this day. 

Hume\index{Hume, David} contributed greatly to the philosophy of the 
Enlightenment and spent a long time in the Parisian salon of Baron Paul Henri 
Thiry d'Holbach\index{Holb@d'Holbach, Paul Henri Thiry}, the son of a 
wine-growing family from Edesheim in the German region of Palatinate. He 
brought his philosophical views into the discourse and had close personal 
contact with the people present there. Among them was the unique Denis 
Diderot\index{Diderot, Denis}, who, together with the mathematician Jean 
d'Alembert\index{Alem@d'Alembert, Jean}, was at that time editor of the famous 
Encyclopedia. This salon\footnote{The book \cite{blom2011} describes the lively 
goings-on in this salon. Famous is an anecdote, after which Hume\index{Hume, 
David} made the remark on his first arrival in the salon that he did not know 
any atheists and d'Holbach\index{Holb@d'Holbach, Paul Henri Thiry} replied in 
essence that he could name him $15$ of the $18$ present immediately and $3$ 
more 
who were not quite decided yet.} was, like the Encyclopedia, committed to the 
spirit of the Enlightenment. Hume\index{Hume, David} will encounter us again in 
connection with the concept of numbers by 
Frege\index{Frege, Gottlob}. 

\section{The Platonic World of Ideas\index{Platonic idealism}}

Let's assume for a moment that we can identify all avatars of the number $5$. 
Where then does this one idea of the number $5$ live, if there is only one? 
Starting with Plato\index{Plato}, many people, especially 
Leibniz\index{Leibniz, Gottfried Wilhelm} and Frege\index{Frege, Gottlob}, have 
asked similar questions also for geometric shapes. 

\begin{figure}[ht!]
\begin{tikzpicture}
\draw (0,0) circle (30pt);
\end{tikzpicture}
\caption{A circle figure is always an approximation since a true mathematical 
circle would be invisible.}
\end{figure}

Consider a circle. With the depiction of geometric figures, different questions 
and problems are associated. A real circle would basically have no extension 
and zooming in should not produce a pixelated circle ring as in this 
illustration. Moreover, a real circle would actually be invisible. It can 
therefore be said with justification that a real circle has never been drawn or 
seen anywhere in the world. What we draw or see is always just an image of our 
mental representation of an ideal circle. All these approximations of a 
circle seem to us equivalent\index{Equivalence} to this ideal representation. 

The non-existence of numbers and real circles in reality brings 
us to the question of where such abstract objects can be located. This 
leads us to the Platonic world of ideas\index{Platonic idealism}, 
also called Platonism or Platonic realism.\footnote{See
\cite{balaguer1998,cassou2011} on the Platonic world of ideas.} It attributes 
a kind of existence\index{Existence}, i.e., a being, outside of our thoughts 
and outside of reality to abstract concepts such as numbers or geometric 
objects as well as the theorems that apply under these objects. The Platonic 
view\index{Platonic idealism} goes back to Plato\index{Plato}, a student of 
Socrates\index{Socrates}, who liked to record his thoughts in dialogues. 
Plato\index{Plato} and his school considered this world of ideas\index{Platonic 
idealism} -- in contrast to Aristotle\index{Aristotle} -- as the superior world 
compared to physical reality. What kind of existence\index{Existence} lies 
behind the Platonic world of ideas\index{Platonic idealism} remains unclear 
from the surviving texts and it seems that Plato\index{Plato} occasionally 
adjusted his position on this point throughout his life. 

\begin{figure}[ht!]
\begin{tikzpicture}[scale=.9,auto=center,every 
node/.style={circle,fill=lightgray}] 

  \node (a1) at (1,5)  {Reality}; 
  \node (a2) at (3,7)  {\hskip0.2cm Ideas \hskip0.2cm}; 
  \node (a3) at (5,5)  {\hskip0.22cm Mind \hskip0.22cm}; 
  
  \draw (a3) -- (a1);  
 
  \draw (a1) -- (a2);  
  \draw (a2) -- (a3);  

\end{tikzpicture}
\caption{\label{fig:platon}A picture of the Platonic world of 
ideas\index{Platonic idealism}. It shows three places: reality, Platonic world 
of ideas, and our mind.}
\end{figure}

To explain the Platonic world of ideas\index{Platonic idealism}, a three-part 
image is suitable (see Fig.~\ref{fig:platon}). This model serves the careful 
mental separation of the three places reality, world of ideas and mind. The 
mind combines the possible human thoughts into one place. There are direct 
connections between mind, reality and ideas.\footnote{Ren\'e 
Descartes\index{Descartes, Ren\'e} referred to such mutual relationships as 
dualisms. They play a role in the body-soul and in the mind-matter problem.} 

Based on this image, there are several attitudes towards these three places. No 
one will doubt that our mind exists, because otherwise our thinking would not
take place. The concept of reality is different. The 
immaterialism\index{Immaterialism} of George Berkeley\index{Berkeley, George} 
claims that the material reality has no own existence\index{Existence} 
and the world is only perceptible through perception.\footnote{There are 
further intermediate forms like conceptualism, which was represented by William 
of Ockham\index{Ockham, William of}. For immaterialism see \cite{berkeley}.}

If we accept mind and reality, one attitude can be that the world of ideas is  
part of its own third place and strictly distinguishes itself from reality 
and mind. This picture is exactly the Platonic world of ideas\index{Platonic 
idealism}. Secondly, there is the counterposition that there is no own place 
for the ideas and they are either part of reality or our thoughts. This view is 
called nominalism\index{Nominalism}.\footnote{See \cite{linnebo,field1980} for 
nominalism. In recent times, Hartry Field\index{Field, Hartry}, Paul 
Benacerraf\index{Benacerraf, Paul} and John Searle\index{Searle, John} held 
nominalistic views.} It means, that ideas appear in two forms, either in the 
form of natural laws and facts in physical reality or as human or divine 
thoughts.\footnote{If we reject the Platonic world of ideas and at the same 
time in nominalism only acknowledge the world of ideas in the form of human or 
divine thoughts, then we obtain a pantheistic ontological proof of God, because 
the laws of nature were already present at the Big Bang, when no humans 
existed.}

In mathematics, the Platonic view\index{Platonic idealism} is widely held. This 
may be because mathematical theorems do not seem to be tautologies, but 
precious treasures that must be discovered. Indeed, many theorems in 
mathematics 
have a surprising statement and their proofs are hard to find. Examples of such 
theorems include Wiles' theorem\index{Wiles, Andrew} or the $4$-squares 
theorem. When they are found, or rather rediscovered, it feels like the 
discovery of a continent or the ascent of a high mountain.

Abstract thinking takes place in nominalism\index{Nominalism} with signs and 
words and therefore with names, which explains the designation of this current 
in the philosophy of mathematics. This view already played a role in medieval 
thinking with Roscellinus\index{Roscellinus} or William of Ockham\index{Ockham, 
William of}. The economy of thought, expressed by the principle of Ockham's 
razor\index{Ockham's Razor}, suggests a preference for 
nominalism\index{Nominalism}, rejecting the unnecessary realism\index{Realism} 
of objects in the Platonic world of ideas\index{Platonic idealism}, which is 
also called Plato's beard as a counter-design to Ockham's razor\index{Ockham's 
Razor}.

The question of whether the Platonic\index{Platonic idealism} or the 
nominalistic view is to be preferred is related to the dispute over 
universals\index{Universals}\footnote{See the article on universals in 
\cite{reclam}.} in philosophy. This is usually only applied to 
universals\index{Universals}, i.e., to universal superordinate terms that name 
totalities of similar objects. A well-known theological example is the 
Trinity. A mathematical example of universals\index{Universals} comes from 
Frege\index{Frege, Gottlob}. He defined numbers like $5$ as the class of all 
finite collections consisting of $5$ elements. More precisely, he identified 
all 
sets with $5$ elements to a new object that represents this class. However, 
this 
approach got him into trouble through Russell's paradox\index{Paradox}, as 
we will see. It is astonishing that the dispute over 
universals\index{Universals}, the Platonic world of ideas and 
nominalism\index{Nominalism} not only appeared in medieval 
scholasticism\index{Scholasticism}, but played a role throughout church 
history. 
Anselm of Canterbury\index{Canterbury, Anselm of}, who co-founded 
scholasticism\index{Scholasticism} in the 11th century and whose name stands 
for 
an ontological\index{Ontology} proof of God\index{Proof of God}, was a 
representative of the Platonic ideas\index{Platonic idealism}. For him, it was 
plausible that the existence\index{Existence} of a divine being could be 
compared with the existence\index{Existence} of abstract objects or 
universals\index{Universals}. His proof of God\index{Proof of God} was based on 
the maximality of God's positive properties, which would not be fulfilled 
without his existence\index{Existence}. Kurt G\"odel\index{Goed@G\"odel, Kurt}, 
who confessed to the Platonic ideas\index{Platonic idealism}, sketched in his 
last years an ontological\index{Ontology} proof of God\index{Proof of God} with 
the help of a suitable modal logic, which is similar to that 
of Canterbury\index{Canterbury, Anselm of} in the assumption of 
maximality.\footnote{All proofs of God are -- unsurprisingly -- based on such 
axioms. See \cite[Vol. III, 1970]{goedel}.} 

\section{Doubts About the Platonic World of Ideas\index{Platonic idealism}}

Until today, no one can decide whether the Platonic world of 
ideas\index{Platonic idealism} or nominalism\index{Nominalism} is the correct 
point of view. Although this historical debate plays only a subordinate 
role in our thinking about abstract objects, as many sections in the present 
book show, we would like to look more closely at different aspects of this 
issue.

There are legitimate doubts about the Platonic world of ideas\index{Platonic 
idealism}. They stem from the unresolved question of what kind of 
existence\index{Existence} is attributed to ideas at a fictitious third place. 
Such considerations belong to the field of metaphysics\index{Metaphysics}. For 
most people, it is clear that material objects have an 
existence\index{Existence} in reality, i.e., in physical reality. In 
philosophy, even this is not self-evident and has led to the question of 
realism\index{Realism},\footnote{See \cite[Ch. IV]{schrenk} and 
\cite{clarke2020,willaschek2015}.} which asks about the 
existence\index{Existence} 
of an external world independent of our consciousness. There are different and 
contradictory attitudes to this. The sceptical view of 
antirealism\index{Antirealism} considers the world around us as a kind of 
fiction that arises in our consciousness as an internal construction. A 
mathematical antirealist positions which we will get to know in this book is 
intuitionism\index{Intuitionism} in its original form. 
Realism\index{Realism}, with its different manifestations, offers a 
counter-position. Among the realistic\index{Realism} positions there are 
variants accepting the existence\index{Existence} of thoughts and facts 
independent of the human mind, however, not in the extremal form of an external 
Platonic world of ideas\index{Platonic idealism}. This discourse thus touches 
upon fragmentations of the philosophical question of 
existence\index{Existence}.\footnote{Kant\index{Kant, Immanuel}, 
Frege\index{Frege, Gottlob} and Russell\index{Russell, Bertrand} had doubts 
whether existence is a predicate or rather a property of properties. See 
Frege's \enquote{Dialogue with P\"unjer on existence} in \cite[p. 
60]{frege1983}, \cite{mcdaniel}, \cite{schrenk}, and the article on existence 
in \cite{reclam}.} 

There always have been doubts in which way abstract objects and ideas could  
exist in a stronger Platonic sense\index{Platonic idealism}. The philosopher 
Paul Benacerraf\index{Benacerraf, Paul} contributed to this debate by 
identifying a few problems in the Platonic world of ideas\index{Platonic 
idealism}. On the one hand, he detected an identification problem arising 
from the non-uniqueness\index{Non-uniqueness problem} of abstract objects. 
In a similar way as we have illustrated for the natural number $5$, he saw that 
it is inevitable that every natural number necessarily has (infinitely) many 
mutually distinct avatars of which none is preferred. This contradicts the idea 
of a Platonic world of ideas\index{Platonic idealism} in a fundamental way. On 
the other hand, Benacerraf\index{Benacerraf, Paul} detected an 
epistemological\index{Epistemology} access problem caused by the 
inaccessibility of Platonic\index{Platonic idealism} mathematical objects 
for our human senses. In his words, mathematics cannot simultaneously have a 
reasonable semantics\index{Semantics} (usually set-theoretic) and a suitable 
ontology\index{Ontology}, since the truth\index{Truth} conditions for 
mathematical statements are not compatible in both approaches. The objections 
raised by Benacerraf\index{Benacerraf, Paul} are often called  
Benacerraf's\index{Benacerraf, Paul} dilemma\index{Benacerraf's 
dilemma}.\footnote{See \cite{benacerraf1965}, \cite{benacerraf1973} and  
\cite{field1980} as well as the books \cite{field1989}, \cite{reclam}, 
\cite{linnebo}, \cite{schrenk}, and \cite{willaschek2015} on realism and 
Benacerraf's dilemma. See \cite{clarke2020} and \cite{nagel} for the analogy 
between morality and mathematics in this respect.}

Hartry Field\index{Field, Hartry} has described Benacerraf's 
dilemma\index{Benacerraf's dilemma} more vividly as follows: 
\begin{quote}
We start out by assuming the existence of mathematical entities that obey the 
standard mathematical theories; we grant also that there may be positive 
reasons for believing in those entities. [...] But Benacerraf's challenge 
[...] is to [...] explain how our beliefs about these remote entities can so 
well reflect the facts about them.\footnote{See \cite{field1989}.} 
\end{quote}

A possible solution of the non-uniqueness\index{Non-uniqueness problem} 
problem of mathematical objects would be branched Platonic parallel worlds in 
the sense of possible worlds\index{Theory of possible worlds} according to 
Leibniz\index{Leibniz, Gottfried Wilhelm} which house the respective avatars of 
numbers, geometric figures and other objects. In the further course of this 
book, we will get to know a new structuralist\index{Structuralism} 
approach that allows us to identify\index{Equality}\index{Identity} objects 
with the same structure. Such approaches defuse Benacerraf's 
dilemma\index{Benacerraf's dilemma} but do not yet justify the Platonic world 
of ideas\index{Platonic idealism} in a satisfactory way. 

Regardless of one's attitude towards realism\index{Realism}, it appears 
reasonable to embed abstract objects in mathematics (and other sciences) into a 
suitable formal scientific language in which statements can be checked via 
mathematical proofs. Rudolf Carnap\index{Carnap, Rudolf} dealt with this 
question intensively and came to the following conclusion: 
\begin{quote}
It is the purpose of this article to clarify this controversial issue. The 
nature and implications of the acceptance of a language referring to abstract 
entities will first be discussed in general; it will be shown that using such a 
language does not imply embracing a Platonic ontology but is perfectly 
compatible with empiricism and strictly scientific thinking. Then the special 
question of the role of abstract entities in semantics will be discussed. It is 
hoped that the clarification of the issue will be useful to those who would 
like to accept abstract entities in their work in mathematics, physics, 
semantics, or any other field; it may help them to overcome nominalistic 
scruples.\footnote{See \cite{carnap1950}.} 
\end{quote}

\section{The Language of Set Theory\index{Set theory}}

An informal concept of collections of objects has always been present in human 
minds. In 1851, Bernhard Bolzano\index{Bolzano, Bernard} described a naive 
concept of sets\index{Set theory} in his book \enquote{Paradoxes of the 
infinite}\footnote{See \cite{bolzano1851}.} which was common in his time. 
Between 1872 and 1882, Georg Cantor\index{Cantor, Georg} and Richard 
Dedekind\index{Dedekind, Richard} exchanged numerous letters. Both had 
independently developed the basic concepts of set theory\index{Set theory} in a 
rigorous manner and openly shared their knowledge with each other.\footnote{See 
\cite{cavailles}.}

Dedekind\index{Dedekind, Richard} developed set theory\index{Set theory} in his 
book \enquote{What are and what should the numbers be?}\footnote{Contained in 
\cite{dedekind}.} as part of his investigations into natural numbers. 
Cantor\index{Cantor, Georg} had already needed aspects of set theory\index{Set 
theory} in his previous investigations into trigonometric series. Later, he 
wrote a series of works titled \enquote{On infinite linear 
point-manifolds}.\footnote{Contained in \cite{cantor}.} As David 
Hilbert\index{Hilbert, David} once said, Cantor\index{Cantor, Georg} had 
created a paradise within mathematics with the discovery of the world of 
infinite sets\index{Set theory} and their ascending transfinite hierarchy.

Both used -- like Bolzano\index{Bolzano, Bernard} -- a concept of 
sets\index{Set theory} in which sets were thought of as collections of 
arbitrary elements. This idea is often referred to today as naive or material 
set theory\index{Set theory}, because in it sets are seen in such a way that 
all elements $x \in A$ of a set $A$ are known and it can be checked for each 
element $x$ whether $x$ is an element of $A$. Two sets are equal if they have 
the same elements. An important feature of this concept is the fact that an 
element $x$ can occur in several sets. Thus, a naive set of the form
\[
A=\{\text{Jan's dog}, \text{ Heike's cat}\} 
\]
can exist, where the description stands for two specific animals and their 
owners. Set theory\index{Set theory} has various operations, such as the 
mappings 
\[
f \colon A \xlongrightarrow{~~~~~~} B 
\]
between sets, or the intersection and union set 
\[
A \cap B, \; A \cup B 
\]
and the power set
\[
\mathrm{Pow}(A)=\text{the set of all subsets of } A. 
\]
Set theory\index{Set theory} has a disadvantage that will be significant for us 
later. If we replace the person Jan with the person Erik in the set $A$, the 
set 
\[
A'=\{\text{Erik's dog}, \text{ Heike's cat}\} 
\]
is created. If we now form the intersection with the set
\[
B=\{\text{Jan's pets}\}, 
\]
then $A\cap B$ and $A' \cap B$ are not equal, although $A$ and $A'$ are 
isomorphic\index{Isomorphism} as sets, because they have only changed into each 
other by exchanging one element. The intersection is therefore a delicate 
operation. It is called intensional\index{Intensionality}, as the intersection 
depends on the respective elements.

The question of the existence\index{Existence} of infinite sets played a major 
role in set theory\index{Set theory} from the beginning. 
Dedekind\index{Dedekind, Richard} believed that he could prove the 
existence\index{Existence} of an infinite set with the help of the totality of 
all thoughts possible for a human being. His proof is considered not strictly 
mathematical and therefore
not as correct. Dedekind\index{Dedekind, Richard} used infinite sets to 
construct infinite chains
\[
0, S(0), S(S(0)), \ldots 
\]
of elements in an infinite set with a starting element $0$ and a mapping $S$ 
of the set onto itself, which is called the successor function. Assuming that 
$S$ sufficiently separates the points of the set, these elements are pairwise 
distinct from each other and form a model\index{Model theory} of the natural 
numbers, usually denoted by $\mathbb{N}$.

Both Dedekind\index{Dedekind, Richard} and Cantor\index{Cantor, Georg} 
considered the collection of all natural numbers in such a model\index{Model 
theory} $\mathbb{N}$ as a new actual infinite object of mathematics. For 
Dedekind\index{Dedekind, Richard}, it was clear that there are many different 
models\index{Model theory} of the natural numbers. So he proved a theorem 
stating that all models\index{Model theory} are isomorphic to each 
other.\footnote{Dedekind\index{Dedekind, Richard} used second-order predicate 
logic, so non-standard models\index{Non-standard model} of $\mathbb{N}$ could 
not occur in his work.} This situation is a good example of mathematical 
objects 
being different, yet still being identified. Similar uniqueness theorems can be 
proven for the real numbers $\mathbb{R}$.

After the emergence of Russell's paradox\index{Paradox} around 1902, which 
arose from the consideration of the supposed set
\[
M=\{x \mid x \notin x\},
\]
the naive set theory\index{Set theory} of Cantor\index{Cantor, Georg} and 
Dedekind\index{Dedekind, Richard} fell into a crisis. From a letter from 
Hilbert\index{Hilbert, David} to Frege\index{Frege, Gottlob} in 1903, it is 
clear that Ernst Zermelo\index{Zermelo, Ernst} had discovered Russell's 
paradox\index{Paradox} before Russell\index{Russell, Bertrand} 
did.\footnote{See \cite{frege1976}.} He also conceived the axiom of 
choice in 1904. It states that from a family
\[
(M_i)_{i \in I}  
\]
of non-empty sets over an index set $I$, a tuple of elements $m_i$ can be 
selected such that $m_i \in M_i$ for all indices $i \in I$. Although this 
statement seems obvious, it constitutes an independent axiom because the set 
$I$ can be infinite. Zermelo\index{Zermelo, Ernst} succeeded in formulating a 
(conjecturally) contradiction-free axiom system for set theory\index{Set 
theory} in 1907. This list was modified by Abraham Fraenkel\index{Fraenkel, 
Abraham} in 1921 and completed by him in 1930 to the Zermelo-Fraenkel axiom 
system\index{Zermelo-Fraenkel axioms}. It is abbreviated as ZF, or ZFC if the 
axiom of choice is added.\footnote{The axiom of choice is equivalent to the 
well-ordering theorem and Zorn's lemma. See \cite[Ch. II]{manin} for the 
foundations of set theory.}

Set theory\index{Set theory} became the most popular of the possible 
foundations of mathematics with this axiom system. This is because sets seem 
intuitively graspable to many people, even if they are infinite. This 
perception is deceptive, as the 
consistency\index{Consistency}\index{Contradiction-freeness} 
of set theory\index{Set theory} is not provable\index{Provability} by its own 
means due to G\"odel's\index{Goed@G\"odel, Kurt} incompleteness 
theorems\index{Incompleteness theorem}. The set of real numbers $\mathbb{R}$ 
raises further deep questions. These include the 
undecidability\index{Undecidability} of the continuum 
hypothesis\index{Continuum hypothesis}. It states that every infinite subset of 
the real numbers is either bijective to $\mathbb{N}$ or bijective to 
$\mathbb{R}$. There are models\index{Model theory} of set theory in which the 
continuum hypothesis\index{Continuum hypothesis} holds and others in which it 
does not. Only by using additional axioms can it be decided. Thus, the view 
that infinite sets are well understood objects is at least questionable.

\section{Frege's\index{Frege, Gottlob} Concept of Number and
Logicism\index{Logicism}}

Frege\index{Frege, Gottlob} outlined his ideas on the concept of number in 
\enquote{The foundations of arithmetic} from 1884 and \enquote{Basic laws of 
arithmetic} from 1893.\footnote{See \cite{frege1884,frege1893}.} The first book 
was intended, among other things, to provide a justification for 
logicism\index{Logicism}. This idea, which Frege\index{Frege, Gottlob} was 
convinced of in his younger years, asserts that mathematics is not a 
synthetic\index{Synthetic judgement} a priori science in the sense 
of Immanuel Kant\index{Kant, Immanuel}, but analytic\index{Analytical  
judgement} a posteriori. This means that mathematics is based on more 
fundamental principles and derives its justification from the more fundamental 
logic. Since arithmetic is based on the principle of complete 
induction\index{Complete induction}, Frege\index{Frege, 
Gottlob} had to derive this principle from fundamental logical principles 
and introduce a precise concept of number. For this purpose, 
Frege\index{Frege, Gottlob} used a method called abstraction. In 
his considerations, he recognised that it is crucial to speak of concepts: 
\begin{quote}
This suggests to us ... that the number statement contains a statement about a 
concept.\footnote{Machine translation of \cite[\S 46]{frege1884}.}  
\end{quote}
To obtain number statements, he went beyond logic and used sets to describe 
concepts. He used the concept of conceptual scope\index{Conceptual scope} 
(German Begriffsumfang), which goes back to Leibniz\index{Leibniz, 
Gottfried Wilhelm}. The conceptual scope\index{Conceptual scope} of $F$ is the 
set of all objects that fall under the concept $F$. With these designations, 
Frege\index{Frege, Gottlob} defined the number of a concept $F$ as its 
conceptual scope\index{Conceptual scope}. 
The number $0$ can be defined, for example, by a concept that can never be 
fulfilled and whose conceptual scope\index{Conceptual scope} is the empty set. 
As a further illustration, Frege\index{Frege, Gottlob} gave the example of the 
large moons of Jupiter. This concept apparently has the number $4$. Since the 
same number can apply to several concepts, Frege\index{Frege, Gottlob} needed 
Hume's principle\index{Hume's principle}: 
\begin{quote}
The expression \enquote{the concept $F$ is equinumerous with the concept $G$} 
shall be equivalent to the expression \enquote{there is a relation $\varphi$, 
which uniquely assigns the objects falling under the concept $F$ to the objects 
falling under $G$}.\footnote{Machine translation of \cite[\S 72]{frege1884}.} 
\end{quote}
This principle states that two concepts $F$ and $G$ represent the same 
number exactly when there is a bijective mapping between the associated 
scope\index{Conceptual scope} sets, which we will also denote by $F$ and $G$ 
for simplicity. An obvious notation for the equinumerosity of the concepts $F$ 
and $G$ is thus
\[
F \cong G.
\]
Frege\index{Frege, Gottlob} noted in his investigations that this approach does 
not yet sufficiently define numbers. A famous example of his illustrates the 
difficulty of distinguishing numbers from objects like the person Julius 
Caesar\index{Caes@Caesar, Julius}. He sought a way out of this situation and 
had the idea of grouping all representatives of concepts with the same number 
by abstraction into an equivalence\index{Equivalence} class. In 
Frege's\index{Frege, Gottlob} own words: 
\begin{quote}
The number which belongs to the concept $F$ is the scope of the concept 
\enquote{equinumerous to the concept $F$}.\footnote{Machine translation of 
\cite[\S 72]{frege1884}.}  
\end{quote}
The number that belongs to the concept\index{Conceptual scope} $F$ is thus the 
set 
\[
\text{Number to the concept } F = \{G \mid G \cong F\}, 
\]
where a set of concepts $G$ is postulated from a concept $F$ that are 
equinumerous to $F$. For example, in today's language, the number $5$ is 
defined as the class of all sets consisting of $5$ elements. Although there are 
many sets with $5$ elements, there is only one such class. In the end, 
Frege\index{Frege, Gottlob} had constructed an equivalence\index{Equivalence} 
class of sets of equal cardinality with his method and thus the problem of 
defining numbers seemed to be solved.  

However, Frege's\index{Frege, Gottlob} approach was problematic. He used parts 
of naive set theory\index{Set theory}, which he apparently considered as 
an unproblematic part of logic. Furthermore, he used the axiom of unrestricted 
comprehension of set theory\index{Set theory} in a crucial way in the formation 
of a set of sets in the concept of number. This axiom allows the 
impredicative formation of sets, i.e., the formation of a 
set of new objects $G$ from an a priori not determined reservoir, which fulfil 
a given predicate with a free variable. The use of 
unrestricted comprehension inevitably led to Russell's paradox\index{Paradox}, 
as Frege\index{Frege, Gottlob} learned from Bertrand Russell\index{Russell, 
Bertrand} in a famous letter dated  June 16th 1902. This letter describes the 
problem that Russell\index{Russell, Bertrand} discovered, but also expresses 
his appreciation for Frege\index{Frege, Gottlob}:
\begin{quote}
Dear Colleague! \\ 
For a year and a half I have known your \enquote{Basic laws of arithmetic}, but 
only now have I been able to find the time for the thorough study that I intend 
to devote to your writings. I find myself in full agreement with you on all 
main points ... Only in one point have I encountered a difficulty. You claim 
(p.17) that the function can also form the indefinite element. I used to 
believe this, but now this view seems doubtful to me, because of the 
following contradiction: Let $w$ be the predicate, a predicate that 
cannot be predicated of itself. Can $w$ be predicated of itself? From any 
answer the opposite follows. Therefore, one must 
conclude that $w$ is not a predicate. Similarly, there is no class (as 
a whole) of those classes that do not belong to themselves as a whole ... with 
you I find the best of what I know from our time, and therefore I have 
allowed myself to express my deep respect to you.\footnote{Machine translation 
of \cite{frege1976}.}
\end{quote}
Russell's paradox\index{Paradox} can be represented in 
Frege's\index{Frege, Gottlob} thought world by the supposed set 
\[
M=\{x \mid x \notin x\},
\]
which uses the predicate $P(x)=x \notin x$, i.e., $x$ 
does not contain itself. If then $M$ contains itself, it does not contain 
itself and if it conversely does not contain itself, then it contains itself. 
This is obviously a antinomy\index{Antinomy}, i.e., a contradiction.

Frege\index{Frege, Gottlob} responded a few days later on June 22nd and 
indicated that he had recognised the dramatic consequences:
\begingroup
\addtolength\leftmargini{0.1in}
\begin{quote}
Your discovery of the contradiction has surprised me to the utmost and, almost 
I would say, dismayed, because thereby the ground on which I thought to build 
arithmetic is shaken. It seems that the transformation of the generality of an 
equality into a value-range equality (\S 9 of my Basic Laws) is not always 
permitted, that my Law V (\S 20, p. 36) is false and that my explanations in 
\S 31 are not sufficient to secure a meaning for my symbol combinations in all 
cases ... It is all the more serious as with the disappearance of my Law V not 
only the basis of my arithmetic, but the only possible basis of arithmetic 
seems to sink.\footnote{Machine translation of \cite{frege1976}.}
\end{quote}
\endgroup
Frege\index{Frege, Gottlob} was unable to solve the problems caused by 
Russell's\index{Russell, Bertrand} letter and was very unhappy about it. 
The Russellian antinomy\index{Antinomy} revealed a serious gap in the naive set 
theory\index{Set theory} founded by Cantor\index{Cantor, Georg} and 
Dedekind\index{Dedekind, Richard} from 1870 onwards. This was fixed in 
1930 through the final establishment of the Zermelo-Fraenkel axiom 
system\index{Zermelo-Fraenkel axioms} and, in particular, by later theories 
that contain classes that are not sets.

Bertrand Russell\index{Russell, Bertrand} tried another way out of the 
antinomies\index{Antinomy} of set theory\index{Set theory}, by hierarchising 
sets into an ascending sequence of infinitely many types, so that the 
self-reference on which the Russellian antinomy\index{Antinomy} is based, was 
excluded. Later he published these ideas, which founded type 
theory\index{Type theory}, together with Alfred North 
Whitehead\index{Whitehead, Alfred North} in the three-part book 
\enquote{Principia mathematica} and tried to underpin the idea of 
logicism\index{Logicism}.

Frege's\index{Frege, Gottlob} logicism\index{Logicism} was criticised in many 
ways. The deficits of naive set theory\index{Set theory} and its use as part of 
logic are already serious points of criticism. Henri 
Poincar\'e\index{Poincar\'e, Henri} was of the opinion that mathematical  
complete induction\index{Complete induction} is synthetic\index{Synthetic 
judgement} a priori and therefore cannot be derived from other principles in 
principle: 
\begin{quote} This law, which is as inaccessible to analytical proof 
as it is to experience, gives the actual type of synthetic judgement a 
priori\index{Apr@A priori judgement}.\footnote{See 
\cite[Part One, Ch. 1, Sec. VI]{poincare1902} for the French text. We provide a 
machine translation of the German translation by Elisabeth 
K\"ussner-Lindemann\index{Kuess@K\"ussner-Lindemann, Elisabeth} from the 
digitally available German edition of 1904 by Teubner Publ. House. See also 
\cite[p. 151]{skolem}.}
\end{quote} 
Poincar\'e's\index{Poincar\'e, Henri} statement can be interpreted to mean that 
infinite constructs like the natural numbers $\mathbb{N}$ and their 
arithmetic properties cannot be justified from purely logical principles. 
Instead, the existence\index{Existence} of the infinite set $\mathbb{N}$ and 
the 
validity of the proof principle of complete induction\index{Complete induction} 
must be demanded as an axiom. Since here the 
consistency\index{Consistency}\index{Contradiction-freeness} of the 
entire arithmetic is on trial, Poincar\'e's\index{Poincar\'e, Henri} objections 
form a serious critique of logicism\index{Logicism}.

Can logicism\index{Logicism} be saved? According to the view of a 
group of neologicists\index{Neologicism}\footnote{This group includes Crispin 
Wright\index{Wright, Crispin}, Charles Parsons\index{Parsons, Charles}, and Bob 
Hale\index{Hale, Bob}. For Hume's principle see \cite[Ch. 9]{linnebo}.} 
this was already achieved by Frege\index{Frege, Gottlob} himself by using  
Hume's principle\index{Hume's principle}, which equates the equinumerosity 
of sets with their isomorphism\index{Isomorphism}. They observed that 
Frege\index{Frege, Gottlob} essentially no longer used conceptual scopes 
after he had proven Hume's principle\index{Hume's principle}. If this 
principle is postulated as an axiom, then Russell's paradox\index{Paradox} 
and Law V can be avoided and Frege's\index{Frege, Gottlob} remaining proofs of 
the properties of the natural numbers are acceptable. In particular, 
Frege\index{Frege, Gottlob} had correctly defined the successor function $S$ 
and all axioms of Dedekind-Peano arithmetic\index{Dedekind-Peano arithmetic}. 
However, Hume's principle\index{Hume's principle} is not a genuine logical, but 
a set-theoretical structural principle\index{Structuralism}. This rescue 
attempt by the neologicists\index{Neologicism} thus does not fully correspond 
to the original intention of logicism\index{Logicism}.

Another way of approaching logicism\index{Logicism} is Russell's\index{Russell, 
Bertrand} type theory\index{Type theory}. In its current form, it can be 
understood as the syntax\index{Syntax} of all mathematics in the form of a 
higher-level logic. This allows the natural numbers $\mathbb{N}$ to be defined 
in a way that is closer to Dedekind's\index{Dedekind, Richard} definition and 
to logicism\index{Logicism}, because the problem of infinity does not occur in 
the calculus of type theory\index{Type theory} itself, but only in 
set-theoretical\index{Set theory} or categorical\index{Category theory} 
models. But even within type theory\index{Type theory} the question of the 
consistency\index{Consistency}\index{Contradiction-freeness} of elementary 
arithmetic remains. Thus, Poincar\'e's\index{Poincar\'e, Henri} objection 
remains in some way and logicism\index{Logicism} must yield to 
Kant's\index{Kant, Immanuel} view that mathematics is an a priori 
science.\footnote{On logicism and type theory see 
\cite{lambekscott1986,russell}.} 

\chapter{Scientific Languages}

The search for a universal language of mankind is an old dream. 
The Babylonian language confusion, which makes it difficult for humanity to 
communicate, has always posed problems. Over the 
centuries, there have been different approaches to inventing new 
artificial interlinguistic planned languages for easier 
communication in everyday life and formal languages\index{Formal language} as  
support for science. 

Leibniz\index{Leibniz, Gottfried Wilhelm} pursued such ideas more 
systematically 
than many others in his time. In connection with his concept 
of a universal science, called scientia generalis\index{Scientia generalis}, he 
hoped for formal symbolic calculi that enable correct proofs 
and calculations\index{Computability} and eliminate irrationality, 
inaccuracies and injustice in all sciences and in human
life. This contribution by Leibniz\index{Leibniz, Gottfried Wilhelm}, 
which went significantly beyond mathematics and logic, was an important step 
in the history of the philosophical concept of truth\index{Truth}. He 
also coined the motto \enquote{theoria cum praxi}. It means that 
science and applications should be closely linked. 

Umberto Eco\index{Eco, Umberto} has recorded the history of the most important 
approaches for interlinguistic planned languages and universal scientific 
languages in a book.\footnote{See \cite{eco1993}.} In it, he 
found that this search was a history of frequent failure, which 
has produced remarkable partial successes, such as the development of 
mathematics and logic, computers and artificial 
intelligence\index{Artificial intelligence}. It is a fascinating 
and still unsolved task to explore the fundamental limits of this idea. 

\section{Leibniz\index{Leibniz, Gottfried Wilhelm} and the Scientia 
Generalis\index{Scientia generalis}} 

Gottfried Wilhelm Leibniz\index{Leibniz, Gottfried Wilhelm} left behind 
an enormous oeuvre, which has not yet been fully edited to this day. His 
partly unpublished works only had an appropriate effect at the beginning of the 
20th century after the publication of a small part by Louis 
Couturat\index{Couturat, Louis}.\footnote{See \cite{couturat1901}.} 

The philosophical system of Leibniz\index{Leibniz, Gottfried Wilhelm}, 
especially that of his \enquote{Monadology},\footnote{Accessible at 
www.projekt-gutenberg.org.} is an alphabet of thoughts reminiscent of Giordano 
Bruno's\index{Bruno, Giordano} concept of monads. His approach consists in 
reducing thinking to primitive concepts and using these building blocks, as 
well 
as with the help of a logical calculus, to reduce the truth\index{Truth} of 
composite statements to simpler questions and ultimately to answer them. Many 
works by Leibniz\index{Leibniz, Gottfried Wilhelm} can therefore -- in addition 
to his works on theology and metaphysics\index{Metaphysics} -- be seen as the 
beginning of epistemology\index{Epistemology} with the help of symbolic 
methods. He believed that human thinking -- in contrast to his opinion of God's 
intuitive and comprehensive perception -- due to the limits of human reason, 
relied on symbolic knowledge. With this mutual positioning of theology and 
mathematics, Leibniz\index{Leibniz, Gottfried Wilhelm} was in the 
tradition\footnote{See \cite{eco1993}.} of Christian theology. 

Leibniz\index{Leibniz, Gottfried Wilhelm} anticipated the ideas of logicians 
like Gottlob Frege\index{Frege, Gottlob} in his thinking and built on 
insights that were in circulation long before his time. Although 
algorithms\index{Algorithm} have been used since antiquity, such as the 
Euclidean algorithm\index{Euclid} in antiquity or the calculations of Easter in 
the Middle Ages, it is only the Majorcan scholar Ram\'on Llull\index{Llull, 
Ram\'on} (Latin Raimundus Lullus) in the 13th century who is attributed with 
the concept of deductive and algorithmic\index{Algorithm} thinking, which 
is the basis of proofs and calculating machines. Leibniz\index{Leibniz, 
Gottfried Wilhelm}, in addition to his numerous other interests, not only dealt 
intensively with logic, but also with other parts of mathematics. In many ways, 
Leibniz's\index{Leibniz, Gottfried Wilhelm} dream of the scientia 
generalis\index{Scientia generalis} was also a dream of a generalised 
mathematics. 

Leibniz\index{Leibniz, Gottfried Wilhelm} refers in his 
\enquote{Dissertatio de arte combinatoria}\footnote{See \cite[Ser. VI, 
Vol. 1, No. 8]{leibniz}.} from 1666 explicitly to the book \enquote{Ars 
magna}\footnote{See the shorter version \enquote{Ars brevis} \cite{lullus2001},
the \enquote{Ars magna} and the work \enquote{Logica nova} \cite{lullus1985}.} 
by Llull\index{Llull, Ram\'on} from 1290. Leibniz was also familiar with the 
ideas of Ren\'e Descartes\index{Descartes, Ren\'e}, Thomas 
Hobbes\index{Hobbes, Thomas} and John Wilkins\index{Wilkins, John}. 

Hobbes\index{Hobbes, Thomas} had recognised in his work \enquote{De 
corpore}\footnote{See \cite{hobbes1655}.} that thinking in a generalised sense 
can be equated with calculating. Descartes\index{Descartes, Ren\'e} had 
already dreamt of a universal philosophical language before 
Leibniz\index{Leibniz, Gottfried Wilhelm}, as he wrote in a letter to the 
number theorist and monk Marin Mersenne\index{Mersenne, Marin} on November 20th 
1629:
\begin{quote}
One should arrange all thoughts methodically, just as the natural 
sequence of numbers is methodically arranged. Just as one can learn in a single 
day in any foreign language to name and write all numbers to infinity, the 
numbers, which certainly form an endless series of 
combinations, one must find the possibility to construct all 
words that are necessary to express everything that 
comes to the human mind and can come ... The invention of such a language 
depends on true philosophy.\footnote{Machine translation of the German version 
in \cite{blanke1996,couturat1903}. Louis Couturat\index{Couturat, Louis} was a 
significant Leibniz researcher.} 
\end{quote}
John Wilkins\index{Wilkins, John} published a book\footnote{See 
\cite{wilkins1668}.} in 1668 about a universal philosophical language, which 
should be superior to natural language. Most of the other universal 
languages\index{Lingua universalis} conceived at the time were primarily 
designed for recording knowledge and communication and anticipated today's 
interlinguistic planned languages.\footnote{Here Johann Joachim 
Becher\index{Becher, Johann Joachim}, George Dalgarno\index{Dalgarno, George}, 
Athanasius Kircher\index{Kircher, Athanasius} and Philippe 
Labb\'e\index{Labb\'e, Philippe} are to be mentioned, see \cite{blanke1996, 
eco1993}. Interlinguistic planned languages were developed by Johann Martin 
Schreyer\index{Schreyer, Johann Martin} (Volap\"uk), Giuseppe 
Peano\index{Peano, Giuseppe} (Latino sine flexione), Louis 
Couturat\index{Couturat, Louis} (Ido) and Ludwik Zamenhof\index{Zamenhof, 
Ludwik} (Esperanto).} 

Leibniz\index{Leibniz, Gottfried Wilhelm} thought more deeply than many 
others about a comprehensive universal science, which he 
called scientia generalis\index{Scientia generalis}. It should have an 
underlying universal scientific language, which we will call 
lingua universalis\index{Lingua universalis} in this text. This 
was for Leibniz\index{Leibniz, Gottfried Wilhelm} a tool to carry out 
thinking in a symbolic way and to perform conclusions and proofs with a 
computational method by term substitutions within the calculus. This art of 
judgement was called ars judicandi at that time. The proximity of proof and 
calculation, which will become important for us, was thus already present in 
the ideas of Leibniz\index{Leibniz, Gottfried Wilhelm}, 
Descartes\index{Descartes, Ren\'e}, Hobbes\index{Hobbes, Thomas} and 
Wilkins\index{Wilkins, John}. Furthermore, it should be possible to develop new 
creative thoughts and concepts in the form of an art of invention (ars 
inveniendi) by using this language.\footnote{See \cite{jost2019} for 
Leibniz\index{Leibniz, Gottfried Wilhelm}'s scientific ideas.}

Presumably due to time constraints, Leibniz\index{Leibniz, Gottfried Wilhelm} 
never fully developed these ideas himself in a satisfactory way. However, he 
regularly described them in his correspondence, for example in a letter to 
Nicolas R\'emond\index{Rem@R\'emond, Nicolas} dated January 10th 1714:
\begin{quote}
I should venture to add that if I had been less distracted, or if I were 
younger 
or had talented young men to help me, I should still hope to create a kind of 
universal symbolistic [sp\'ecieuse g\'en\'erale] in which all truths of reason 
would be reduced to a kind of calculus. At the same time this could be a kind 
of universal language or writing, though infinitely different from all such 
languages which have thus far been proposed, for the characters and the words 
themselves would give directions to reason, and the errors -- except those of 
fact -- would be only mistakes in calculation. It would be very difficult to 
form or invent this language or characteristic but very easy to learn it 
without any dictionaries. When we lack sufficient data to arrive at certainty 
in our truths, it would also serve to estimate degrees of probability and to 
see what is needed to provide this certainty.\footnote{See \cite[Ser. 
I]{leibniz} for the French original. English translation in \cite{leibniz1989}.}
\end{quote}
The calculus dreamed of by Leibniz\index{Leibniz, Gottfried Wilhelm} was 
realised within mathematical logic in the 19th and 20th centuries by 
Frege\index{Frege, Gottlob} and others. However, in its intended generality, it 
has not been constructed to this day and was not applied by 
Leibniz\index{Leibniz, Gottfried Wilhelm} himself within his 
metaphysics\index{Metaphysics}. One reason for this may be that his 
time would probably not have accepted a calculus of clearly derivable 
truth\index{Truth} anyway.

\section{Leibniz\index{Leibniz, Gottfried Wilhelm} and Mathematics}

Although Leibniz\index{Leibniz, Gottfried Wilhelm} was a self-taught 
mathematician, he had significant ideas. For many years he commissioned the 
construction of a calculating machine with a stepped drum for the four basic 
arithmetic operations and it can be speculated whether he dreamed of a powerful 
universal machine to put his symbolic calculus into practice and thus answer 
all philosophical questions rationally. 

Leibniz\index{Leibniz, Gottfried Wilhelm} introduced computing based on the 
binary number system with the digits $0$ and $1$, which would become the basis 
of digital computers. Equally successful was his invention of the differentials 
$\mathrm{dx}$ and the differential quotients
\[
\mathrm{\frac{df}{dx}} 
\]
in the form of a calculus. The differential $\mathrm{dx}$ symbolises an 
infinitesimally small quantity. The calculus\footnote{See the essays in 
\cite[Ser. VII, Vol. 4A]{leibniz}.} of Leibniz\index{Leibniz, Gottfried 
Wilhelm} has been maintained to this day and included the calculation rules 
in Fig.~\ref{fig:differentials}. 

\begin{figure}[ht!]
\begin{align*}
\mathrm{d(f+g)} &=\mathrm{df +dg} & &\text{(Additivity)} \cr 
\mathrm{d(af)} &=\mathrm{a df} & &\text{(Linearity)}\cr
\mathrm{d(fg)} &=\mathrm{f dg+g df} & &\text{(Product rule)}\cr 
\mathrm{d\left(\frac{f}{g}\right)} &=\mathrm{\frac{g df -f dg}{g^2}} & & 
\text{(Quotient rule)}\cr
\mathrm{d(f \circ g)} &= \mathrm{\frac{df}{dg} dg} & & \text{(Chain rule)} 
\end{align*}
\caption{\label{fig:differentials}Rules of calculus using the language of 
differentials.}
\end{figure}

In the development of infinitesimal calculus, Leibniz\index{Leibniz, Gottfried 
Wilhelm} was in a sometimes bitter competition with Isaac Newton\index{Newton, 
Isaac} and his theory of fluxions. Although he was inferior to 
Newton\index{Newton, Isaac} and the Royal Society on some occasions, his status 
in science today is on par with Newton's\index{Newton, Isaac}.

Calculus presupposes a concept of real numbers. Such a concept was certainly 
held by Leibniz\index{Leibniz, Gottfried Wilhelm} and others in his time. 
However, it was not until the 19th century that Cantor\index{Cantor, Georg} and 
Dedekind\index{Dedekind, Richard} precisely defined the set $\mathbb{R}$ of 
real numbers and proved theorems about it. Modern calculus, with its precise 
concept of limits, was only established in the 19th century by Karl 
Weierstra{\ss}\index{Weierstra{\ss}, Karl} and others. Only much later could a 
mathematically flawless definition of infinitesimally small quantities be given 
with the non-standard analysis\index{Non-standard model} of Abraham 
Robinson\index{Robinson, Abraham}, which revived the calculus of 
Leibniz\index{Leibniz, Gottfried Wilhelm}.

Leibniz\index{Leibniz, Gottfried Wilhelm} formulated some remarkable statements 
about real numbers. He already knew the difference between rational, 
irrational, algebraic and transcendental numbers and proved results about the 
transcendence of certain integrals of differential forms.\footnote{See the 
letter to Christiaan Huygens\index{Huygens, Christiaan} from April 1691 
\cite[Ser. III, Vol. 5, No. 17]{leibniz}.}

\section{The Rise of Mathematical Logic}

Aristotle\index{Aristotle} is considered the founder of mathematical logic 
through his study of syllogisms. These consist of a chain of statements, also 
called a mode (Latin modus):
\begin{quote}
Major premise: All B are C \\
Minor premise: All A are B \\
Conclusion: All A are C.
\end{quote}
In this example, the modus barbara is given. Aristotelian logic included a 
whole range of other such modes and had a dominant position for centuries in 
the scholastic tradition\index{Scholasticism} until the Middle Ages. 
Aristotle\index{Aristotle} already knew the law of excluded middle, 
i.e., the principle underlying proof by 
contradiction.\footnote{The law states that for any statement $A$, either $A$ 
or 
$\neg A$ is true. For historical remarks on this, see \cite[Ch. 
II]{lambekscott1986}. Only indirect proofs are known for the Thue-Siegel-Roth 
theorem and K\"onig's lemma.}

From today's perspective, the Aristotelian\index{Aristotle} modes are only 
fragments of modern logic. Aristotelian logic was revived in the 20th century 
by Jan {\L}ukasiewicz\index{Luka@{\L}ukasiewicz, Jan} and compared with 
then-modern 
theories. After ancient culture, logic was practised as a fundamental cultural 
technique in Europe, Asia and the Islamic world. These discoveries anticipated 
many modern developments, but often remained isolated. In early modern Europe, 
it was primarily Leibniz\index{Leibniz, Gottfried Wilhelm} who had the most 
significant insights in logic. He provided tools of mathematical logic in the 
form of a calculus, which anticipated the algebraic logic of the 19th century, 
including a version of quantifiers and thus predicate 
logic. Leibniz\index{Leibniz, Gottfried Wilhelm} 
anticipated modal logic in its beginnings, which is related 
to his theory of possible worlds\index{Theory of possible worlds}. However, all 
this remained hidden to most people of his time because Leibniz\index{Leibniz, 
Gottfried Wilhelm} did not publish his thoughts on logic.\footnote{See 
\cite[Ch. 8]{jost2019} and the text \enquote{Generales inquisitiones de analysi 
notionum et veritatum} from 1686 in \cite[Ser. VI, Vol. 4A, No. 165]{leibniz}.}

Mathematical logic experienced after Leibniz\index{Leibniz, Gottfried Wilhelm}
a great upswing in the 19th century, particularly at universities in 
Europe and the United States. The two English logicians 
George Boole\index{Boole, George} and Augustus de Morgan\index{Demor@De Morgan, 
Augustus}, as well as somewhat later Gottlob Frege\index{Frege, Gottlob} and 
Giuseppe Peano\index{Peano, Giuseppe}, particularly shaped the basic concepts 
of algebraic propositional logic and the beginnings of predicate 
logic.\footnote{See \cite{boole1847,demorgan1847,frege1879,peano}. Also 
contributing were John Venn\index{Venn, John}, William Stanley 
Jevons\index{Jevons, William Stanley}, Charles S. Peirce\index{Peirce, 
Charles}, Christine Ladd-Franklin\index{Ladd-Franklin, Christine} and Ernst 
Schr\"oder's\index{Schr\"oder, Ernst} book \enquote{Algebra of logic} 
\cite{schroeder1890}.} Logic as a mathematical theory only fully unfolded in 
the 
20th century. 

Frege\index{Frege, Gottlob} recognised from examples that natural language is 
not a good basis for the foundations of mathematics and therefore developed 
mathematical logic as 
\begin{quote}
The science of the most general laws of truth.\footnote{Machine translation of 
\cite[p. 139]{frege1983}.}
\end{quote}

He shaped mathematical logic in his \enquote{Begriffsschrift}.\footnote{See 
\cite{frege1879}.} Frege's\index{Frege, Gottlob} variant of a Leibnizian 
language included logical symbols such as the judgement stroke\index{Judgement} 
\[
\vdash A
\]
as well as the universal quantifier and the existential quantifier
\[
\forall x, \; \exists x.    
\]
With this, he introduced the first deductive system\index{Deductive system} 
with a formal language\index{Formal language} in the history of logic, which we 
would today call second-order predicate logic or 
simply second-order logic.\footnote{See \cite[Ch. I/II]{manin} for such 
concept formations.} Some of the common logical notations today, however, are 
more due to Peano\index{Peano, Giuseppe} than to Frege\index{Frege, Gottlob}. 

Leibniz\index{Leibniz, Gottfried Wilhelm} had a remarkable influence on Gottlob 
Frege\index{Frege, Gottlob}, Giuseppe Peano\index{Peano, Giuseppe} and 
Kurt G\"odel\index{Goed@G\"odel, Kurt}. At the beginning of the 
\enquote{Begriffsschrift}, Frege\index{Frege, Gottlob} explicitly referred to 
Leibniz\index{Leibniz, Gottfried Wilhelm} as a role model:
\begin{quote}
Leibniz too recognised the advantages of an appropriate notation, perhaps 
overestimated. His idea of a general characteristic, a calculus philosophicus 
or calculus ratiocinator was too gigantic for the attempt to realise it to get 
beyond mere preparations.\footnote{Machine translation of \cite{frege1879}.}
\end{quote}

The beginning of the classical logic of the 19th century was the calculus of 
propositional logic, to mathematically handle logical statements. In it, there 
are the possible truth values\index{Truth} 
\[
\top \text{ (true)}, \bot \text{ (false)}
\]
and statements are connected to new statements through the connections 
\[
A \wedge B \text{ (and)}, A \vee B \text{ (or)}, 
\]
\[
\neg A \text{ (negation)},A \xRightarrow{~~} 
B \text{ (conditional or implication)}.
\]

\begin{figure}[ht!]
\begin{tabular}{|c|c|c|c|}\hline
$A$ & $B$ & $A \Rightarrow B$ & $\neg A \vee B$ \\\hline
true & true & true & true \\\hline
true & false & false & false \\\hline
false & true & true & true \\\hline
false & false & true & true \\\hline
\end{tabular} 
\caption{\label{fig:truth_table}The truth tablet shows the identical truth 
values of $A \Rightarrow B$ and $\neg A \vee B$ depending on the truth values 
of $A$ and $B$.} 
\end{figure}

While the conjunction $\wedge$ is familiar to us in everyday life, the 
disjunction $\vee$ is often used in an exclusive sense. In mathematical logic, 
however, the statement $A \vee B$ means that at least one of the two statements 
$A$ or $B$ is fulfilled, i.e., they can both be fulfilled at the same time. In 
classical propositional logic, calculation rules such as 
\[
\neg A \vee \neg B=\neg (A \wedge B) \text{ and } 
(A \xRightarrow{~~} B) = \neg A \vee B
\]
apply, so that some symbols are in principle redundant. Proofs for such 
calculation rules can be conducted with truth tables (see 
Fig.~\ref{fig:truth_table}). 

The calculation rules of propositional logic for the 
unary or binary operations $\wedge, \vee, \Rightarrow$ and $\neg$ result in an 
algebraic structure, which in honour of George Boole\index{Boole, George} is 
also called Boolean algebra\index{Boolean algebra}. Such a structure $\Omega$ 
is in general not a ring in the usual mathematical sense, even though the 
operations $\vee$ and $\wedge$ have similar properties to addition and 
multiplication. In fact, $\vee$ and $\wedge$ allow a lattice 
structure\footnote{A lattice has a partial order structure $\le$ and two 
operations $\wedge$ (or $\cap$) and $\vee$ (or $\cup$), analogous to 
intersection (infimum) and union (supremum). The partial order is connected 
with the two operations. For example, $A \le B$ holds exactly when $A=A \wedge 
B$. In addition to the usual associative and commutative laws, the two 
absorption laws $A \vee (A \wedge B)=A$ and $A \wedge (A \vee B)=A$ also apply. 
Both operations are only semigroups, therefore a lattice is not a ring in 
general. Using the symmetric difference instead of $\vee$ yields a ring 
though, in fact a $\mathbb{Z}/2\mathbb{Z}$-algebra.} 
on $\Omega$ by the union or intersection of subobjects. The objects 
$A \Rightarrow B$ and $\neg A$ are called exponential object\index{Exponential 
object} $B^A$ and negation of $A$. On a Boolean algebra\index{Boolean algebra} 
there also exists a partial order $\le$, which can be defined by the relation 
\[
A \le B \text{ exactly when } A=A \wedge B 
\] 
through the operation $\wedge$. In $\Omega$ the equation 
\[
\neg A \vee B=(A \xRightarrow{~~} B) 
\]
holds and $\neg$ is an involution:
\[
\neg \neg A=A.
\]
Heyting algebras\index{Heyting algebra} are generalisations of Boolean 
algebras\index{Boolean algebra}, in which these two rules are abandoned. In 
them only\footnote{See \cite[Ch. I, Sec. 8]{maclane_moerdijk1992} on Boolean 
algebras and Heyting algebras.} 
\[
(\neg A \vee B) \le (A \xRightarrow{~~} B) \text{ and } 
A \le \neg \neg A 
\]
apply. Heyting algebras\index{Heyting algebra} play a major role in 
intuitionism\index{Intuitionism}.

\section{The Question of Truth\index{Truth}} 

Aristotle\index{Aristotle} also dealt with the concept of truth\index{Truth}. 
In his \enquote{Metaphysics}\index{Metaphysics}\footnote{See 
\cite{aristoteles}.} he made fundamental considerations on which two 
significant concepts of truth\index{Truth} are based. The first of these is the 
adequation theory\index{Adequation theory of truth} of truth. It sees the 
concept of truth\index{Truth} in the agreement of thinking with reality. Thomas 
Aquinas\index{Aquinas, Thomas} expressed this as follows:  
\begin{quote}
Veritas consistit in adaequatione intellectus et rei.\footnote{See 
\cite[Art. 1]{aquin}. Thomas Aquinas\index{Aquinas, Thomas} systematically 
linked Aristotelian\index{Aristotle} philosophy with Christian theology.}
\end{quote} 
By thinking, we mean the contents of consciousness and thoughts in 
our brain. We regard reality as the material physical 
reality together with its natural laws, but also as other phenomena, 
which may not be reducible to physics, such as questions of 
existence\index{Existence} and metaphysics\index{Metaphysics}. 

The correspondence theory\index{Correspondence theory of truth} of 
truth\index{Truth} -- related to the adequation\index{Adequation 
theory of truth} theory of truth -- replaces thoughts with the corresponding 
linguistic propositions in our thinking, i.e., the formulations of our 
thoughts. Aristotle\index{Aristotle} gave an explicit example of the 
correspondence theory\index{Correspondence theory of truth}: 
\begin{quote}
Not because our opinion, that you are white, is true, are you white, but 
because you are white, we tell the truth, by asserting this.\footnote{Machine 
translation of \cite[1051b]{aristoteles}.}
\end{quote}
The correspondence theory\index{Correspondence theory of truth} in its  
formulation thus means the correspondence between statements and facts of 
reality. Both theories are to be questioned, because they connect two a priori 
incomparable categories such as firstly thinking and the contents of 
consciousness in our brain -- or associated propositions -- and secondly 
reality in an inexplicable way. Problems lie in the relationship of the 
truth bearer\index{Truth bearer}, which could either be a thought or a 
linguistic statement, with the truth maker\index{Truth maker}, i.e., 
a fact of reality, because it is not obvious how facts could be determined 
without already having a definition of truth\index{Truth}. For these reasons, 
there have always been controversial discussions about the 
correspondence theory\index{Correspondence theory of truth}. Thomas 
Hobbes\index{Hobbes, Thomas} wrote in his \enquote{Leviathan} around 1651: 
\begin{quote}
For true and false are attributes of speech, not of things.\footnote{See 
\cite[Sec. 4.11]{leviathan}.}  
\end{quote}
When such concepts of truth\index{Truth} are applied to mathematics, 
serious problems arise, because mathematical concepts usually have 
no correspondence in reality. Moreover, the concept of reality   
itself is problematic, because it is not clear what exactly is meant by it. 
When people talk about reality, they are usually guided by aspects of matter 
and there is a basic trust in the concreteness of this physical reality. Upon 
deeper reflection, however, we encounter doubts and will come to the conclusion 
that even the physical reality has a very abstract nature, 
which we cannot comprehend as a whole and therefore can only describe with 
mathematical language. An example are electromagnetic fields, which 
spread in complete vacuum according to mathematical laws. Such 
fields are only tangible through measurements or their effect on other physical 
objects. Phenomena of quantum mechanics such as the dualism between 
particles and waves as well as the Heisenberg uncertainty relation elude 
even more our intuition. The abstractness of reality and its 
characteristics is the subject of the philosophical question of 
realism\index{Realism} and allows for contradictory positions.\footnote{See the 
article on realism and antirealism in \cite[Ch. IV]{schrenk} as well as the 
papers \cite{hawking} and \cite{penrose2006} on phenomena of physics.}

These problems with the adequation\index{Adequation theory of truth} 
and the correspondence theory\index{Correspondence theory of truth} of truth
had various consequences. On the one hand, this 
has led to the search for further theories of truth\index{Truth} that are free 
from connections with any kind of reality. An example of this is the coherence 
theory\index{Coherence theory of truth} of truth\index{Truth}. It states that 
for the truth\index{Truth} of statements of a new theory, the 
coherence\index{Coherence theory of truth} within the scope of the already 
existing or underlying theoretical structure is primarily important. A conditio 
sine qua non for this is again the 
consistency\index{Consistency}\index{Contradiction-freeness} -- i.e., the 
contradiction-freeness -- of the theory under consideration itself. The 
coherence theory\index{Coherence theory of truth} is also controversial in 
philosophy due to the existence\index{Existence} of ambiguous coherent systems, 
as Bertrand Russell\index{Russell, Bertrand} already noted.\footnote{See 
\cite{halbach1996,tarski1935,tarski1969,reclam,schrenk} on theories of truth.} 
We will return to whether this view is at least fruitful within mathematics. We 
do not provide an exhaustive description of all existing theories of 
truth\index{Truth}.

\section{Leibniz\index{Leibniz, Gottfried Wilhelm} and Truth\index{Truth}}

Leibniz\index{Leibniz, Gottfried Wilhelm} had his own conception of the 
correspondence theory\index{Correspondence theory of truth}. The concept of 
truth\index{Truth} for him was a property of thoughts, not things. In the 
\enquote{Dialogus}, he accordingly wrote in 1677:
\begin{quote}
Veritatem esse cogitationum non rerum.\footnote{See \cite[Ser. VI, Vol. 4A, 
No. 8]{leibniz}.} 
\end{quote}
The Leibnizian\index{Leibniz, Gottfried Wilhelm} theory of truth\index{Truth} 
uses the distinction between factual truths and rational truths. The 
truth\index{Truth} of propositions is examined by decomposing 
composite statements (notio composita) into simpler, indivisible statements 
using the symbolic calculus of the lingua universalis\index{Lingua 
universalis}. Leibniz\index{Leibniz, Gottfried Wilhelm} wrote about this in his 
\enquote{Monadology}:
\begin{quote}
There are also two kinds of truths: those of reasoning, which are necessary 
and whose opposite is impossible, and those of fact, which are contingent and 
whose opposite is possible. If a truth is necessary, one can find the reason 
for it by analysis, by resolving it into simpler ideas and truths, until one 
reaches the initial ones.\footnote{Machine translation of \S 33 of the 
\enquote{Monadology}, available at www.projekt-gutenberg.org, as well as the 
further essays in \cite[Ser. VI, Vol. 4A]{leibniz}.}   
\end{quote}
Factual truths apply in our world and are gained through experience. Rational 
truths apply in all possible worlds\index{Theory of possible worlds}. 
Leibniz\index{Leibniz, Gottfried Wilhelm} had invented a theory of possible 
worlds\index{Theory of possible worlds} to locate other realisations of 
truth\index{Truth}. It served him to represent the reality in which we live as 
the one preferred by God as the best of all possible worlds\index{Theory of 
possible worlds}.

Leibniz\index{Leibniz, Gottfried Wilhelm} incidentally spoke of God as a being 
that does not need the symbolic language because it can directly grasp the 
truth\index{Truth}. Moreover, in his opinion, God could experience the world 
and its history, i.e., the concept of space-time, as a whole all at once. If 
God is not considered as a being, but is understood as a universal 
principle, synonymous with the world, then Leibniz's\index{Leibniz, Gottfried 
Wilhelm} perspective leads to the realisation that, for us ordinary people, 
reality could only be a semantic\index{Semantics} internal representation, 
because we do not fully comprehend it. While Leibniz's\index{Leibniz, Gottfried 
Wilhelm} image of God naturally eludes all proof, his theory of possible 
worlds\index{Theory of possible worlds} was in a way realised 
in the 20th century through modal logic.

\section{Truth\index{Truth} and Kantian\index{Kant, Immanuel} 
Judgements\index{Judgement}}

In his considerations on the concept of truth\index{Truth}, 
Immanuel Kant\index{Kant, Immanuel} significantly developed 
Leibniz's\index{Leibniz, Gottfried Wilhelm} theory of truth\index{Truth}, 
without explicitly mentioning Leibniz\index{Leibniz, Gottfried 
Wilhelm}.\footnote{See \cite{boehm,martin-loef1994,willaschek2023}.} 
Kant\index{Kant, Immanuel} used the term 
judgement\index{Judgement} for statements whose truth\index{Truth} is 
to be investigated. He made an epistemic\index{Epistemology} 
distinction between a priori judgements\index{Apr@A priori judgement} -- 
necessary and universally valid judgements\index{Judgement} independent of 
experience -- and a posteriori 
judgements\index{Apos@A posteriori judgement} -- experiential judgements 
through empirical knowledge --, which in a way goes back to the 
Aristotelian\index{Aristotle} dualism between proteron (condition) and hysteron 
(conditioned).

Mathematical axioms are good examples of a priori\index{Apr@A priori 
judgement} judgements. A related example underlies any attempt at a 
proof of God\index{Proof of God}, because the existence\index{Existence} 
of divine beings cannot be derived from experience or observation of the world, 
even if some people claim this. Serious proofs of God\index{Proof of God}, 
like the proofs of Anselm of Canterbury\index{Canterbury, Anselm of} and Kurt 
G\"odel\index{Goed@G\"odel, Kurt}, are based on logical arguments and use an 
axiom that connects divine beings with maximal positive properties. 

Immanuel Kant\index{Kant, Immanuel} further distinguished between 
analytical\index{Analytical judgement} and synthetic\index{Synthetic 
judgement} judgements,\footnote{See \cite[Ch. 1]{linnebo}.} which express a 
reason for the validity of truth\index{Truth}. An analytical 
judgement\index{Analytical judgement} is an explanatory judgement, i.e., it can 
be directly inferred from the definition of the object or concept under 
consideration. In contrast, a synthetic judgement\index{Synthetic judgement} is 
an extension judgement, i.e., it can only be explained by means of a further 
justification. An example of an analytical judgement\index{Analytical 
judgement} is: 
\begin{quote}
It is dark at night. 
\end{quote}
This is because it is precisely the definition of night that it is dark at this 
time, unless there is a light source present. 

The pairs of properties analytical-synthetic\index{Analytical judgement}
\index{Synthetic judgement} and a priori-a posteriori\index{Apr@A priori 
judgement}\index{Apos@A posteriori judgement} are initially in principle  
independent of each other. In this system of thought, Leibniz's rational 
judgements become analytical judgements\index{Analytical judgement} a priori 
and the factual judgements become synthetic judgements\index{Synthetic 
judgement} a posteriori. Analytical judgements\index{Analytical judgement} a 
posteriori are excluded, as analytical judgements\index{Analytical judgement} 
are always a priori,\footnote{See \cite[p. 15-16]{linnebo}.} so that a third
category of synthetic\index{Synthetic judgement} judgements a priori 
remains, which according to Kant\index{Kant, Immanuel} is an existing category. 

Kant attributed the property pairing of synthetic\index{Synthetic judgement} a 
priori\index{Apr@A priori judgement} to mathematics. Mathematics is certainly a 
priori\index{Apr@A priori judgement}, as its judgements are not empirical 
knowledge. Kant\index{Kant, Immanuel} claimed that mathematics is 
synthetic\index{Synthetic judgement}, as it cannot be derived from anything 
else. As an example, he used the equation
\[
7+5=12. 
\]
It is the question whether all possible such equations, i.e., the entire 
arithmetic, directly follow from the definition of numbers or 
whether they can be derived from other principles. Such a principle 
could be a logical principle or a mathematical axiom. 
Leibniz\index{Leibniz, Gottfried Wilhelm} would have spoken of 
a sufficient reason in this context. 

From the question of whether mathematics is synthetic\index{Synthetic 
judgement} or not, logicism\index{Logicism} emerged, which aimed to reduce 
mathematics completely to logic, so that it 
would consequently be analytic\index{Analytical judgement}. Especially 
Gottlob Frege\index{Frege, Gottlob} and Bertrand Russell\index{Russell, 
Bertrand} propagated logicism\index{Logicism} and thus rejected 
the view of Kant\index{Kant, Immanuel}. Willard van Orman 
Quine\index{Quine, Willard van Orman} and others later had doubts 
that the distinction between analytic \index{Analytical judgement} 
and synthetic\index{Synthetic judgement} judgements is clear-cut 
and concluded that logicism\index{Logicism} is not sufficiently 
justified.\footnote{Quine\index{Quine, Willard van Orman} doubted in 
his essay \cite{quine1951} not only Kant's\index{Kant, Immanuel} 
distinction, but also the logical empiricism of the Vienna Circle.} 

In mathematics, the Kantian\index{Kant, Immanuel} distinction between 
analytic \index{Analytical judgement} and synthetic\index{Synthetic 
judgement} judgements\index{Judgement} is not common. But it is quite 
relevant, as synthetic\index{Synthetic judgement} judgements\index{Judgement}  
encounter us through the postulation of new axioms and the introduction 
of primitive concepts that cannot be derived from the remaining assumptions.  

\section{Truth\index{Truth}, Sense and Meaning According to Frege\index{Frege, 
Gottlob}}

We have already touched on the fact that Frege\index{Frege, Gottlob} 
dealt with the concept of equality\index{Equality}\index{Identity} and in this 
context examined the concepts of sense and meaning of linguistic statements. He 
also thought about how the truth\index{Truth} of statements can be proven and 
how the statements and their potential truth\index{Truth} are related to the 
thoughts and the content of our consciousness. The concepts 
truth\index{Truth}, sense and meaning all belong to the field of 
semantics\index{Semantics}, which is located in the language philosophy and 
linguistics ultimately going back to Frege\index{Frege, Gottlob}. With 
his investigations, in which he carefully distinguished all occurring concepts 
from each other, he became a significant philosopher of modern times after his 
works on logic.\footnote{See \cite{vossenkuhl}.} Frege\index{Frege, Gottlob} 
was 
of the opinion that the equation\index{Equality}\index{Identity}
\[
A=A 
\]
is an a priori\index{Apr@A priori judgement} true analytic\index{Analytical  
judgement} statement (in the sense of Kant\index{Kant, Immanuel}). However, the 
statement 
\[
A=B 
\]
is an extension of a completely different form, possibly 
synthetic\index{Synthetic judgement} and not automatically a 
priori\index{Apr@A priori judgement}. Frege\index{Frege, Gottlob} explained the 
difference between the concepts of sense and meaning using this example, as 
we have already explained using the triangle example.

Frege\index{Frege, Gottlob} particularly intensively dealt with the related 
concept of truth\index{Truth}. As early as 1892, he formulated the following 
thought:
\begin{quote}
One could even say: \enquote{The thought that $5$ is a prime number is true.} 
But if one looks more closely, one notices that actually no more is said than 
in the simple sentence \enquote{$5$ is a prime number}.\footnote{Machine 
translation of \cite[p. 31]{frege1962}.}
\end{quote}
Because of such statements, Frege\index{Frege, Gottlob} is often accused of a 
redundancy-theoretical attitude, which considers truth\index{Truth} irrelevant 
for the meaning of statements. Basically, he considered truth\index{Truth} 
undefinable at that time, as he noted in a note from 1906:
\begin{quote}
What is true, I consider unexplainable.\footnote{Machine translation of 
\cite[p. 189]{frege1983}.}
\end{quote}
However, he continued his reflections and tried to precisely explore aspects of 
truth\index{Truth}. He wrote in an essay in 1918:
\begin{quote}
So with every property of a thing, a property of a thought is connected, namely 
the truth.\footnote{Machine translation of \cite[p. 39]{frege2003}.}
\end{quote}
The bearer of truth\index{Truth bearer} lies, according to 
Frege's\index{Frege, Gottlob} opinion, in the statements of thoughts and not in 
the understanding of humans or in reality. The moment the truth\index{Truth} of 
a thought is asserted or shown, a judgement\index{Judgement} is made. 
Frege\index{Frege, Gottlob} wrote about thinking and the world of thoughts:
\begin{quote}
Thoughts are neither things of the outside world nor ideas. A third realm must 
be recognised.\footnote{Machine translation of \cite[p. 50]{frege2003}.}
\end{quote}
Frege's\index{Frege, Gottlob} realm of thoughts as a separate place of 
objective facts reminds us of the Platonic world of ideas\index{Platonic 
idealism}. At an earlier point, Frege\index{Frege, Gottlob} had written on this 
topic:
\begin{quote}
But if the subjective has no place, how is it possible that the objective 
number $4$ is nowhere? Now I assert that there is no contradiction in this. It 
is indeed the same for everyone who deals with it; but this has nothing to do 
with spatiality. Not every objective object has a place.\footnote{Machine 
translation of \cite[\S 61]{frege1884}.}
\end{quote}
So was Frege\index{Frege, Gottlob} a Platonist\index{Platonic idealism} 
or not? Similarly, one can ask why Frege\index{Frege, Gottlob} mentioned 
Kant\index{Kant, Immanuel}, Leibniz\index{Leibniz, Gottfried Wilhelm}, 
Mill\index{Mill, John Stuart} and others in his book \enquote{Foundations of 
arithmetic} and elsewhere, but not German idealism\index{Idealism}, represented 
by Fichte\index{Fichte, Johann Gottlieb}, Schelling\index{Schelling, Friedrich 
Wilhelm Joseph} and Hegel\index{Hegel, Georg Wilhelm Friedrich}, who  
had advanced conceptual philosophical thinking in their works.

However, Frege\index{Frege, Gottlob} was very sceptical about concepts like 
existence\index{Existence}, as his \enquote{Dialogue with P\"unjer on 
existence} proves. Thus, he was certainly rather sceptical towards both 
Platonism\index{Platonic idealism} and idealism\index{Idealism}. He is now 
considered the founder of analytical philosophy and had an influential student 
in Rudolf Carnap\index{Carnap, Rudolf}, who was a member of the Vienna 
Circle\index{Vienna circle} led by Moritz Schlick\index{Schlick, 
Moritz}.\footnote{See \cite[p. 60]{frege1983} and \cite{edmonds2021}.} The 
members advocated logical empiricism and vehemently rejected 
metaphysics\index{Metaphysics}. 

We will get to know the ideas of Hermann Weyl\index{Weyl, Hermann}, who
pondered the symbolic construction of reality in his philosophical works and 
had no fear of contact with idealism\index{Idealism} and 
metaphysics\index{Metaphysics}. He maintained a correspondence with Edmund 
Husserl\index{Husserl, Edmund} and was an admirer of Johann Gottlieb 
Fichte\index{Fichte, Johann Gottlieb}. In today's time, analytical philosophy 
and metaphysics\index{Metaphysics} are not easily distinguishable and their 
differences lie more in the approach than in the philosophical topics. Already 
Leibniz\index{Leibniz, Gottfried Wilhelm} had recognised that scientific 
metaphysics\index{Metaphysics} necessarily must be based on a deductive 
calculus\index{Deductive system}. 

\section{Tarski's\index{Tarski, Alfred} Theory of Truth\index{Truth}}

When investigating the truth\index{Truth} of statements in non-scientific 
contexts, there are fundamental difficulties associated with the concept of 
truth\index{Truth}. This is due to self-referential sentences in natural 
languages such as the paradox\index{Paradox} of Epimenides\index{Epimenides} 
\begin{quote}
A Cretan says: All Cretans lie\footnote{Also called the liar's paradox. 
It probably goes back to Eubulides of Miletus\index{Eubulides of Miletus}, see 
\cite[Part II]{brendel2015}.} 
\end{quote}
or sentences like the paradox\index{Paradox} attributed to Bertrand 
Russell\index{Russell, Bertrand} 
\begin{quote}
The barber who shaves all men who do not shave themselves.
\end{quote} 
These well-known paradoxes\index{Paradox} have been widely 
questioned\footnote{See \cite[Part I]{brendel2015}.} and it is better to use a 
self-referential statement of the form 
\begin{quote}
A: The statement A is false
\end{quote}
in which the paradox\index{Paradox} comes into sharper focus. In formal 
languages\index{Formal language}, such problems can be avoided, as 
Alfred Tarski\index{Tarski, Alfred} has shown in his language-analytical 
theory of truth\index{Truth}. 

His theory can be applied to formal object languages\index{Formal language} 
such as those underlying mathematics. For this purpose, an 
interpretation\index{Interpretation} of the statements in a 
metalanguage\index{Metalanguage} $M$ is used, which usually contains $L$ as 
a sublanguage, extending beyond a given formal language\index{Formal language} 
$L$. By separating the two languages $L$ and $M$, self-reference is avoided 
and a concept of truth\index{Truth} in the form of a truth 
predicate\index{Truth predicate} $T(x)$ in $M$ is defined.

Tarski\index{Tarski, Alfred} has given a syntactic\index{Syntax} adequacy 
condition to be fulfilled in biconditional form that reminds of 
Aristotle\index{Aristotle}:
\begin{quote}
The sentence \enquote{Snow is white} is true if and only if 
snow is white.\footnote{See \cite{tarski1935,tarski1969}.}
\end{quote}
Here, the left part of the sentence within the quotation marks is an atomic 
statement $p$ of the formal object language\index{Formal language} $L$. On the 
left side, the truth predicate\index{Truth predicate} $T$ of the
metalanguage\index{Metalanguage} $M$ applied to $p$ and is equivalent to 
the right-hand side, which consists of the translation $\tilde p$ of $p$ into 
the metalanguage\index{Metalanguage} $M$. So, a bit more formally, in $M$ 
it holds
\begin{quote}
$T(p)$ exactly when $\tilde p$,
\end{quote}
where a standard name for $p$ must be substituted as the argument of 
$T(x)$.\footnote{See \cite[Part III]{brendel2015} and \cite{tarski1969}.}  

The adequacy condition is not sufficient as a definition of 
truth\index{Truth}, but only a necessary prerequisite. The actual definition of 
truth\index{Truth} is made via the satisfiability in a suitable model, which is 
also called semantics\index{Semantics}. In it, the 
metalanguage\index{Metalanguage} $M$ is used as a deductive 
system\index{Deductive system}. The most natural example of mathematical 
semantics\index{Semantics} is the formal object language\index{Formal language} 
$L_\mathrm{ar}$ of Dedekind-Peano arithmetic\index{Dedekind-Peano arithmetic} 
and its interpretation\index{Interpretation} in the 
set-theoretical\index{Set theory} standard model\index{Model theory} 
$\mathbb{N}$ or the various non-standard models\index{Non-standard model}. In 
these models, the metalanguage\index{Metalanguage} $M$ is usually given by the  
Zermelo-Fraenkel set theory\index{Zermelo-Fraenkel axioms} with first-order 
predicate logic. From this it follows that the concept 
of truth\index{Truth} corresponds to the concept of 
provability\index{Provability} in a deductive system\index{Deductive system} 
corresponding to the model\index{Model theory} which has a richer formal 
metalanguage\index{Metalanguage}\index{Formal language} $M$ than $L$. 

Already Leibniz\index{Leibniz, Gottfried Wilhelm} and Frege\index{Frege, 
Gottlob} recognised that methods of mathematical logic can be used in 
inferences to prove the truth\index{Truth}. Later, an important reason was 
found why the proof must be carried out in a strictly richer 
metalanguage\index{Metalanguage}. In Dedekind-Peano 
arithmetic\index{Dedekind-Peano arithmetic}, self-referential constructs of 
the form 
\begin{quote}
A: The statement A is false
\end{quote}
can be cleverly generated with the trick of G\"odel 
numbering\index{Goed@G\"odel, Kurt}. This results in the (first) 
incompleteness theorem\index{Incompleteness theorem} by Kurt 
G\"odel\index{Goed@G\"odel, Kurt} and the related theorem on the 
undefinability of truth\index{Truth} by Alfred 
Tarski\index{Tarski, Alfred}. These theorems demonstrate the 
existence\index{Existence} of undecidable\index{Decidability} statements in 
formal object languages\index{Object language}\index{Formal language} $L$, 
which in $L$ are neither provable\index{Provability} nor refutable, as well as 
the undefinability of the truth predicate\index{Truth predicate} in $L$ 
itself.\footnote{See \cite[Vol. I, 1931]{goedel} and \cite{tarski1969}.}

\section{What Became of the Lingua Universalis\index{Lingua universalis}?} 

The development of logic and formal scientific 
languages\index{Formal language} has brought us within mathematically oriented 
sciences close to a realisation of a 
Leibnizian\index{Leibniz, Gottfried Wilhelm} lingua 
universalis\index{Lingua universalis}. In contrast, outside of mathematics, the 
attitude is widespread that this idea is idealistic, cannot encompass all human 
thinking and many truths\index{Truth} in the sciences and in everyday 
life fundamentally cannot be captured with a formal language\index{Formal 
language}. Even the optimistic Leibniz\index{Leibniz, Gottfried Wilhelm} 
believed in the limits of knowledge and thus of the lingua 
universalis\index{Lingua universalis}. Through his mill example\index{Mill 
example}, he expressed that phenomena such as the consciousness of humans and 
animals and other qualia\index{Qualia} are difficult to explain. In the life 
sciences and in philosophy, the phenomenon of consciousness is controversially 
discussed, so that the exact boundaries in this case remain unclear to this 
day. 

Digitisation has brought forth new 
algorithms\index{Algorithm} in recent years that can excellently imitate human 
thinking. In almost all areas, the formal verifiability of knowledge is 
necessary due to large amounts of data and made possible by intelligent 
digital assistant systems. However, most of the systems that have emerged 
do not yet have the judgement that Leibniz\index{Leibniz, 
Gottfried Wilhelm} had envisaged. It is an interesting and difficult 
question to explore the nature of universal scientific languages. They 
certainly extend -- despite the mentioned counter-positions -- far into areas 
beyond mathematics and computer science. It remains an exciting question to 
explore their limits and the connection with the foundations of our world. 

\chapter{Mathematical Thinking}

Statements like \enquote{I was always bad at maths} or \enquote{What is there 
still to research in mathematics?} are commonplace in the media and at private 
meetings, regularly causing unpleasant feelings for professionals. Many people 
are proud of their lack of knowledge in mathematics, even if they are ashamed 
of their other intellectual deficits. Even if we cannot expect mathematics 
to be important to all people, we should continue to develop the mathematical 
curriculum in schools and try to change this situation.

This is not hopeless, as the mathematical way of thinking is more similar 
to that in art. Godfrey Harold Hardy\index{Hardy, Godfrey Harold}, a famous 
British number theorist, wrote perhaps the best-known book about the inner 
view of mathematics. It is titled \enquote{A mathematician's 
apology}.\footnote{See \cite{hardy}.} In this book, the beauty of mathematics 
and its kinship with art play a special role. Hardy\index{Hardy, Godfrey 
Harold} emphasises therein that the usefulness of mathematics in applications 
is no measure of its sense and quality. He points out that areas like number 
theory possess an inner beauty and have no applications. On this point, 
Hardy\index{Hardy, Godfrey Harold} was mistaken, as number theory now has many 
applications, such as blockchains\index{Cryptography} and public-key 
cryptography\index{Cryptography}.

Outside of number theory, the way of thinking of mathematics can also be 
impressively explained using the world of transfinite sets, the universal, free 
constructions of algebra and the description of topological\index{Topology} 
space forms. 

\section{Natural Numbers}

Number theory, also called arithmetic, is a field of mathematics that deals 
with various types of numbers and their properties. This includes in particular 
the natural numbers, the integers and the algebraic numbers, which are the 
roots of polynomials with integer coefficients. Number theory and geometry 
existed already in antiquity. Euclid\index{Euclid} and 
Diophantus\index{Diophantus} were leaders in these fields. 

The multi-part work \enquote{Arithmetica} by 
Diophantus\index{Diophantus}\footnote{See \cite{diophant}.} contains a 
fascinating arithmetic problem which consists in 
the search for integer solutions of diophantine polynomial equations. The most 
famous equation of this kind is Fermat's equation\index{Fermat, Pierre de}
\[
x^n+y^n=z^n.
\]
According to a theorem by Andrew Wiles\index{Wiles, Andrew} from 1993, it only 
has strictly positive integer solutions when $n=1$ or $2$, as Pierre de 
Fermat\index{Fermat, Pierre de} had suspected. Solutions for the case $n=2$ can 
already be found on Babylonian clay tablets\footnote{See \cite{struik} for the 
history of ancient mathematics.} around 1800 BC, such as the famous Plimpton 
322 tablet. Slightly more general than solutions of a polynomial equation are 
algebraic varieties, which arise as simultaneous zero sets of 
several polynomials. These frequently occur in mathematical 
modelling.\footnote{Polynomials and algebraic varieties appear in mathematical 
models of almost all natural and social sciences. Two beautiful examples of 
this are quantum field theory and algebraic statistics.} 

Euclid\index{Euclid} is the originator of the Euclidean 
algorithm\index{Algorithm}, the prototype of recursive 
thinking\index{Recursion}. It is based on the fact that two natural numbers 
$a$ and $b$ (with $b \neq 0$) can be written in the form  
\[
a=qb+r \text{ with } 0 \le r < b,
\]
where $r$ is the remainder of $a$ modulo $b$ and $qb$ is the maximum multiple 
of $b$ that is still less than or equal to $a$. The Euclidean 
algorithm\index{Algorithm} replaces the pair $(a,b)$ with the pair $(b,r)$ in 
each step and stops as soon as the remainder $r$ becomes zero for the first 
time. The number $b$ in the last pair $(b,0)$ is the sought-after greatest 
common divisor of $a$ and $b$. The algorithm\index{Algorithm} also works if $b$ 
is larger than $a$, because then $r=a$ and $q=0$ and the pair $(a,b)$ is 
replaced by $(b,a)$ in the first step. Since $r<b$ in each step, the process 
terminates after a finite number of iterations. This algorithm\index{Algorithm} 
can be easily programmed and executed quickly in any common programming 
language.\footnote{In the Python programming language: 
\\
def gcd(a,b):\\
\noindent\hspace*{5mm} while b>0:\\
\noindent\hspace*{10mm} a,b=b,a\%b\\
\noindent\hspace*{5mm} return a}

The calculation of $q$ as the floor of $\frac{a}{b}$ can be omitted, because 
the remainder $r$ is reached when $b$ is subtracted from $a$ until a positive 
number less than $b$ is obtained. This variant is called the fast Euclidean 
algorithm\index{Algorithm}. 

We are interested in the properties of the entirety $\mathbb{N}$ of all natural 
numbers with $0$ as the smallest number among them. The infinite sequence 
\[
0,1,2,3,4,5,\ldots 
\]
is based on the principle of counting on. It is generated by the successor 
function $S$, which maps each number $n$ to the following number
\[
S(n)=n+1. 
\]
It holds 
\[
\begin{array}{rl}
1&=S(0) \\
2&=S(S(0)) \\
\vdots & 
\end{array}
\]
If we write 
\[
S(n)=n+1, 
\]
this suggests an addition mapping. The natural numbers indeed possess an 
addition
\[
m+n
\]
and a multiplication 
\[
m \cdot n.
\]
These two arithmetic operations can be precisely defined using the successor 
function $S$ and recursion\index{Recursion}. 

Very large numbers can be easily written in mathematics, for example by using 
powers of ten
\[
\begin{array}{rl}
10^2&=100 \\
10^3&=1000 \\
10^6&=1000000, \text{ one million} \\
\vdots & \\
10^{100}&=1\underbrace{000\ldots 000}_{100\text{ zeros}}, 
\; \text{a googol} \\ \vdots& 
\end{array}
\]
The infinite sequence of natural numbers has no counterpart in the real 
world. This can be easily understood, because the number of particles in 
the universe is very large, but finite. According to reasonably realistic 
estimates based on current models, this number in the visible universe is about 
the order of magnitude $10^{N}$ with $N$ between $80$ and $90$.  

\section{Prime Numbers}

The sequence of ever larger numbers is not mysterious in itself. Only when we 
consider a concept like divisibility does the world of prime numbers unfold 
from 
the natural numbers
\[
2,3,5,7,11,13,\ldots,101,103,\ldots 
\]
These are exactly the natural numbers that are divisible only by $1$ and 
themselves. There are infinitely many prime numbers and every natural 
number can be uniquely written as a product of such. A 
well-known proof for the infinity of the set of all prime numbers was already 
given by Euclid\index{Euclid}.\footnote{Assuming there are only finitely many 
prime numbers $p_1 < p_2 < \cdots < p_n$, consider the number $N=p_1 p_2 \cdots 
p_n+1$ and show that none of the prime factors of $N$ appear in the finite 
list. Contradiction.} Another wonderful proof considers the 
Fermat\index{Fermat, Pierre de} numbers of the form
\[
F_m=2^{2^m}+1 
\]
for $m=0,1,2,\ldots$. The first of these numbers are 
\[
\begin{array}{rl}
F_0&=3 \\
F_1&=5 \\
F_2&=17\\
F_3&=257\\
F_4&=65537.
\end{array}
\]
Fermat\index{Fermat, Pierre de} had speculated in a letter to the 
mathematician Bernard Frenicle de Bessy\index{Frenicle de Bessy, Bernard} 
in 1640 that all Fermat\index{Fermat, Pierre de} numbers are primes. 
However, the number $F_5=641 \cdot 6700417$ is not a prime number. We 
claim, nevertheless, that all these numbers are coprime to each other, 
i.e., they have no common divisor. This implies the infinity of 
prime numbers, because the prime factors in the numbers $F_n$ are all different 
from each other. For this claim, we first show the formula\footnote{First show 
$F_m-2=2^{2^m}-1=(2^{2^{m-1}}+1)(2^{2^{m-1}}-1)=F_{m-1}(F_{m-1}-2)$ and then 
apply complete induction.}  
\[
F_m=F_0 F_1 \cdots F_{m-1}+2. 
\]
If $F_k$ and $F_m$ with $k<m$ had a common divisor $p$, then 
$p$ would also divide the $2$. This is a contradiction, because all $F_m$ and 
thus also $p$ are odd.

Computers today can handle numbers of the order of several hundred decimal 
places and with great effort investigate the factorisation of such numbers 
into prime factors. The problem of finding the factorisation of a number is 
presumably much harder than checking whether a given number is prime, because 
there are prime number tests that have a limited runtime. An example of a prime 
number with $50$ digits is the number
\[
p=53542885039615245271174355315623704334284773568199. 
\]
The currently largest known prime numbers are the rarely occurring  
Mersenne\index{Mersenne, Marin} primes of the form $2^p-1$ with $p$ prime. 
At the time of writing this book, the number 
\[
2^{136.279.841}-1
\]
with $41.024.320$ digits was the largest proven
Mersenne\index{Mersenne, Marin} prime number.\footnote{See www.mersenne.org and 
\cite[Kap. 12]{mspiontkow2011}.} 

Since the time of Bernhard 
Riemann\index{Riemann, Bernhard}, a surprising amount is known about the 
distribution of prime numbers. Although he 
has a comparatively small oeuvre, his works were all the more 
influential. Among them is the formulation of the Riemann 
hypothesis\index{Riemann hypothesis}\footnote{Contained in \cite{riemann}.} 
in an unpublished work of less than seven pages. It states that 
the zeros of the complex-valued Riemann $\zeta$-function
\[
\zeta(s)=\sum_{n=1}^\infty n^{-s}
\]
are either the so-called trivial zeros at the negative even numbers 
$n=-2,-4,-6,\ldots$ or lie on the vertical line with the equation 
$\mathrm{Re}(s)=\frac{1}{2}$ in the complex plane. 
This hypothesis\index{Riemann hypothesis} has not yet been proven, 
although there is plenty of evidence for it. The assumption of further zeros on 
$\mathrm{Re}(s)=\frac{1}{2}$ was a great insight by 
Riemann\index{Riemann, Bernhard}. In the aforementioned short work, he proved 
a connection between the location of the zeros of $\zeta(s)$ and the 
distribution of prime numbers by applying the methods of Fourier theory to 
the $\zeta$-function. This resulted in an analytical formula for 
the number $\pi(x)$ of prime numbers below a limit $x$, which can be expressed 
with the help of the non-trivial zeros. Only in 1932 was a formula discovered 
in Riemann's\index{Riemann, Bernhard} posthumous notes with which he had 
calculated zeros of the $\zeta$-function.\footnote{The Riemann-Siegel formula 
is contained in \cite{siegel}.}

While such formulas are very precise and can in principle find every prime 
number, there are also results like the prime number theorem, which predicts 
the asymptotic distribution of prime numbers. It states that $\pi(x)$ is well 
approximated by the function $\frac{x}{\log(x)}$ and -- even somewhat  
better -- by the logarithmic integral $\mathrm{Li}(x)=\int_2^x 
\frac{dt}{\log(t)}$. This theorem was proven in 1896 by Salomon 
Hadamard\index{Hadamard, Salomon} and Charles-Jean de La Vall\'ee 
Poussin\index{La Vall\'ee Poussin, Charles-Jean de}. Assuming 
the Riemann hypothesis\index{Riemann hypothesis}, it can be shown that the 
error is bounded by 
\[
| \pi(x) - \mathrm{Li}(x) | < \frac{\sqrt{x} \log(x)}{8 \pi}. 
\]
The zeros and special values of the Riemann $\zeta$-function and its 
generalisations in the form of $L$-functions carry information about 
numerous arithmetic objects.\footnote{The Langlands programme\index{Langlands, 
Robert} addresses such questions.} 

\section{The Structure of the Number System}

From the natural numbers $\mathbb{N}$, the entire number system from school can 
be built. This includes the transition from $\mathbb{N}$ to the integers 
$\mathbb{Z}$ and from there to the rational numbers $\mathbb{Q}$, the 
fractions. Finally, we want to construct the real numbers $\mathbb{R}$ from 
$\mathbb{Q}$. The transition from $\mathbb{N}$ to $\mathbb{Z}$ is called the 
Grothendieck group. For this, we imagine any integer, especially the negative 
numbers, as a difference of natural numbers 
\[
a-b, \quad \text{ with } a,b \in \mathbb{N}. 
\]
This representation is not unique, so we equate pairs $(a,b)$ with 
the same difference, such as 
\[
(1,2)=(2,3).  
\]
Negative numbers $-b$ correspond in this view to the pairs $(0,b)$. 
In general, we identify 
\[
(a,b)=(c,d) \text{ exactly when } a+d=b+c.
\]
The equation on the right side only contains natural numbers due to this 
rearrangement. Consequently, the integers $\mathbb{Z}$ are pairs of natural 
numbers that are identified in a precise manner. They form a commutative ring, 
i.e., they carry an addition
\[
(a,b)+(c,d)=(a+c,b+d)
\]
and a multiplication
\[
(a,b) \cdot (c,d)= (ac+bd,ad+bc), 
\]
both of which are connected by the distributive law. 

\medskip 
When transitioning from $\mathbb{Z}$ to $\mathbb{Q}$, fractions 
\[
\frac{a}{b} 
\]
are formed from integers in a similar way, where $a$ and $b$ are integers and 
$b\neq 0$. Here too, we work more precisely with pairs $(a,b)$ with $b \neq 0$ 
and identify
\[
\frac{a}{b}= \frac{c}{d} \text{, if } ad=bc.
\]
The usual rules of fraction calculation apply:
\begin{align*}
\frac{a}{b}+\frac{c}{d}&=\frac{ad+bc}{bd} \cr
\frac{a}{b} \cdot \frac{c}{d}&=\frac{ac}{bd}.
\end{align*}
This construction presented here can be carried out for every commutative ring 
$R$ free of zero divisors instead of $\mathbb{Z}$ and is called the field of 
quotients of $R$. Like the Grothendieck group, it first appeared in an 
unpublished manuscript\footnote{Printed in \cite{ms2023}.} by Richard 
Dedekind\index{Dedekind, Richard}. Of course, these calculation rules for the 
transition from $\mathbb{Z}$ to $\mathbb{Q}$ were in use long before 
Dedekind\index{Dedekind, Richard} and have long been part of the school 
curriculum.

Fields of quotients are special cases of a more general construction, in which 
all objects in a multiplicatively closed subset $S$ of a ring $R$ are made 
invertible. This is called the localisation of a ring after $S$, denoted 
$S^{-1}R$. In the field of quotients of a commutative ring free of zero 
divisors, localisation is done after all ring elements that are not equal to 
$0$. The Grothendieck group and more generally the localisations are examples 
of quotient constructions, where elements are identified so that new 
mathematical objects can be created.

From the rational numbers $\mathbb{Q}$, the real numbers $\mathbb{R}$ can be 
constructed using Dedekind cuts, which are certain subsets of $\mathbb{Q}$, or 
with Cauchy sequences. The approach with Cauchy sequences is called completion 
and is defined as a quotient after the null sequences with respect to the 
norm $|-|$, i.e., the absolute value. With norms not fulfilling the 
Archimedean axiom, this method leads to the more exotic $p$-adic numbers 
$\mathbb{Q}_p$, which are a completion of the rational numbers $\mathbb{Q}$ 
with respect to the $p$-adic norm $|-|_p$ for a prime number $p$. The $p$-adic 
norm of $p$ is given by $1/p$, so that high powers of $p$, unlike in the 
absolute value case, have a small norm.\footnote{See \cite[Kap. 
13]{mspiontkow2011}.} 

\section{Transfinite Numbers}

Georg Cantor\index{Cantor, Georg} defined a fascinating world of infinite 
numbers using set theory\index{Set theory}, which are called 
ordinal\index{Ordinal number} and cardinal numbers\index{Cardinal number}. 
Ordinal numbers\index{Ordinal number} are well-ordered  
sets, i.e., totally ordered\footnote{A set with an order relation $\le$ is 
totally ordered if for any two elements $a,b$, either $a \le b$ or $b \le a$ 
holds.} sets, such that every non-empty subset has a smallest element. 
Cardinal numbers\index{Cardinal number} are certain
ordinal numbers\index{Ordinal number}, the smallest among all 
ordinal numbers\index{Ordinal number} that possess the same cardinality. 
We want to explain this using examples.

All natural numbers are simultaneously ordinal numbers\index{Ordinal number} 
and cardinal numbers\index{Cardinal number}. John von Neumann\index{Neumann, 
John von} represented them as sets of the respective preceding numbers:
\[
0=\emptyset, 1=\{0\}, 2=\{0,1\}, \ldots, n+1=\{0,1,\ldots,n\},\ldots 
\]
The smallest infinite ordinal number\index{Ordinal number} is $\omega$. It 
corresponds to the normal ordering 
\[
0<1<2<3<\cdots 
\]
of the natural numbers. The associated cardinal number\index{Cardinal number} 
of $\mathbb{N}$ is denoted by $\aleph_0$. One of the most remarkable 
insights of Cantor\index{Cantor, Georg} was that counting can continue beyond 
the ordinal number\index{Ordinal number} $\omega$:
\[
\omega<\omega+1<\omega+2<\cdots<2\omega<\cdots<\omega^2<\cdots<\omega^\omega
<\cdots<\varepsilon_0<\cdots 
\]
In addition, Cantor\index{Cantor, Georg} introduced an addition and 
multiplication of ordinal numbers\index{Ordinal number}, both of which are 
generally not commutative. The addition $\alpha+\beta$ is defined by 
\enquote{placing $\alpha$ and $\beta$ side by side}. For example, 
\[
\omega+1=\omega \cup \{\ast\} \neq \omega,
\]
where $\ast$ is greater than all $n \in \omega$. On the other hand, it is 
apparently true that 
\[
1+\omega=\omega. 
\]
The multiplication $\alpha \cdot \beta$ is realised by replacing each element 
of $\beta$ with a copy of the entire well-ordered set $\alpha$. The 
exponentiation is even more complicated. 

In this way, an exotic ordinal number arithmetic is created, which neither 
fulfils commutativity nor usual cancellation rules. Also, fascinating ordinal 
numbers\index{Ordinal number} can be defined, such as the number 
$\varepsilon_0$, which is given as an unfathomably large power
\[
\varepsilon_0=\omega 
\uparrow \omega=\underbrace{\omega^{\omega^{\iddots^\omega}}}_{\omega 
\text{ times}}. 
\]
This number fulfils the unusual fixed point equation 
\[
\omega^{\varepsilon_0}=\varepsilon_0. 
\]
Due to the well-ordering of $\varepsilon_0$, however, every strictly descending 
sequence of ordinal numbers\index{Ordinal number}
\[
\varepsilon_0 > \alpha_1 > \alpha_2 > \cdots 
\]
is necessarily finite. Ordinal numbers\index{Ordinal number} are the basis for 
the incredibly effective proof technique of transfinite 
induction\index{Transfinite induction}, in which -- as with complete 
induction\index{Complete induction} -- the series of ordinal 
numbers\index{Ordinal number} is climbed.

All ordinal numbers\index{Ordinal number} $\alpha < \varepsilon_0$ can be 
associated with finite trees that have a distinguished vertex  -- the root. 
Trees are graphs that do not contain multiple edges and no loops, i.e., closed 
paths that return to the starting point. A tree with a root and attached 
subtrees $B_1,\ldots, B_m$ corresponds to the ordinal number\index{Ordinal 
number} in Cantor's\index{Cantor, Georg} normal form
\[
\alpha=\omega^{\beta_1}+  \cdots + \omega^{\beta_m},
\]
where $\beta_1 \geq \cdots \geq \beta_m$ are ordinal numbers\index{Ordinal 
number} that correspond to the subtrees. In Fig.~\ref{fig:ordinals}, we 
have indicated some examples.\footnote{See \cite{dershowitz} for definitions 
and the well-ordering on trees.}

\begin{figure}[ht!]
\begin{tikzpicture}
\draw[-] (-5,-1) node [below,yshift=-1mm] {$0$};
  
\draw[-] (-4,0) node {}
  -- (-4,-1) node [below,yshift=-1mm] {$1$};

\draw[-] (-2.7,0) node {}
  -- (-3,-1) node [below,yshift=-1mm] {$2$}
  -- (-3.3,0) node {};
  
\draw[-] (-1.7,0) node {}
  -- (-2,-1) node [below,yshift=-1.9mm] {$n$} 
  -- (-2.3,0) node [right] {$...$};
  
\draw[-] (-1,1) node {}
  -- (-1,0) node {}
  -- (-1,-1) node [below,yshift=-1.9mm] {$\omega$};

\draw[-] (-0.3,1) node {}
  -- (-0.3,0) node {}
  -- (0,-1) node [below,yshift=-1.1mm] {$\omega+1$}
  -- (0.3,0) node {};

\draw[-] (0.7,1) node {}
  -- (0.7,0) node {}
  -- (1,-1) node [below,yshift=-1.1mm] {$\omega \cdot 2$}
  -- (1.3,0) node {}
  -- (1.3,1) node {};
  
\draw[-] (1.7,1) node {}
  -- (2,0) node {}
  -- (2,-1) node [below,yshift=-0.5mm] {$\omega^2$}
  -- (2,0) node {}
  -- (2.3,1) node {};

\draw[-] (2.7,1) node {}
  -- (3,0) node {}
  -- (3,-1) node [below,yshift=-1.1mm] {$\omega^n$}
  -- (3,0) node {}
  -- (3.3,1) node [left] {$...$};

\draw[-] (4,2) node {}
  -- (4,1) node {}
  -- (4,0) node {}
  -- (4,-1) node [below,yshift=-1.1mm] {$\omega^\omega$};
  
\draw[-] (5,1) node [above] {$\vdots$}
  -- (5,0) node {}
  -- (5,-1) node [below,yshift=-0.8mm] {$\omega \uparrow n$};
  
\draw[fill=black] (-5,-1) circle (.5ex);

\draw[fill=black] (-4,-1) circle (.5ex);
\draw[fill=black] (-4,0) circle (.5ex);

\draw[fill=black] (-3,-1) circle (.5ex);
\draw[fill=black] (-3.3,0) circle (.5ex);
\draw[fill=black] (-2.7,0) circle (.5ex);

\draw[fill=black] (-2,-1) circle (.5ex);
\draw[fill=black] (-2.3,0) circle (.5ex);
\draw[fill=black] (-1.7,0) circle (.5ex);

\draw[fill=black] (-1,-1) circle (.5ex);
\draw[fill=black] (-1,0) circle (.5ex);
\draw[fill=black] (-1,1) circle (.5ex);

\draw[fill=black] (0,-1) circle (.5ex);
\draw[fill=black] (-0.3,0) circle (.5ex);
\draw[fill=black] (0.3,0) circle (.5ex);
\draw[fill=black] (-0.3,1) circle (.5ex);

\draw[fill=black] (1,-1) circle (.5ex);
\draw[fill=black] (0.7,0) circle (.5ex);
\draw[fill=black] (1.3,0) circle (.5ex);
\draw[fill=black] (1.3,1) circle (.5ex);
\draw[fill=black] (0.7,1) circle (.5ex);

\draw[fill=black] (2,-1) circle (.5ex);
\draw[fill=black] (2,0) circle (.5ex);
\draw[fill=black] (2.3,1) circle (.5ex);
\draw[fill=black] (1.7,1) circle (.5ex);

\draw[fill=black] (3,-1) circle (.5ex);
\draw[fill=black] (3,0) circle (.5ex);
\draw[fill=black] (3.3,1) circle (.5ex);
\draw[fill=black] (2.7,1) circle (.5ex);

\draw[fill=black] (4,-1) circle (.5ex);
\draw[fill=black] (4,0) circle (.5ex);
\draw[fill=black] (4,1) circle (.5ex);
\draw[fill=black] (4,2) circle (.5ex);

\draw[fill=black] (5,-1) circle (.5ex);
\draw[fill=black] (5,0) circle (.5ex);
\draw[fill=black] (5,1) circle (.5ex);
\draw[fill=black] (5,2) circle (.5ex);

\end{tikzpicture}
\caption{\label{fig:ordinals}The graphical correspondence between ordinal 
numbers\index{Ordinal number} $\alpha < \varepsilon_0$ and root trees.}
\end{figure}

John von Neumann\index{Neumann, John von} defined around 1928 the cumulative 
hierarchy\index{Cumulative hierarchy} of sets
\begin{align*}
V_0=&\emptyset \\
V_{\alpha+1}=&\mathrm{Pow}(V_\alpha) \\
V_\lambda=&\bigcup_{\alpha<\lambda} V_\alpha, 
\end{align*}
which is based on the empty set. Here, $\alpha$ is any ordinal 
number\index{Ordinal number} and $\lambda$ is a limit ordinal number, i.e., it 
is not of the form $\lambda=\lambda'+1$ for an ordinal 
number\index{Ordinal number} $\lambda'$. The smallest limit ordinal number 
is obviously $\omega$. For very large ordinal numbers\index{Ordinal number} 
$\kappa$, which are called strongly inaccessible cardinal 
numbers\index{Cardinal number}, the sets $V_\kappa$ of the cumulative 
hierarchy\index{Cumulative hierarchy} are also called 
Grothendieck universes\index{Universe} and are models\index{Model theory} of 
the Zermelo-Fraenkel axioms\index{Zermelo-Fraenkel axioms}, as can be shown by 
transfinite induction. The postulation of strongly 
inaccessible cardinal numbers\index{Cardinal number} represents one of many 
possible additional axioms, which imply the 
consistency\index{Consistency}\index{Contradiction-freeness} of the 
Zermelo-Fraenkel axioms\index{Zermelo-Fraenkel axioms}.\footnote{An 
uncountable cardinal number $\kappa$ is strongly inaccessible if it cannot be 
represented by cardinal arithmetic from smaller cardinal numbers. In 
particular, for every cardinal number $\alpha < \kappa$ it still holds that 
$2^\alpha<\kappa$. The Morse-Kelley axioms also imply the consistency of the 
Zermelo-Fraenkel axioms.}  

In his paradise, Cantor\index{Cantor, Georg} formulated the continuum 
hypothesis\index{Continuum hypothesis}. It states that there is no set whose 
cardinality is strictly between the cardinality $\aleph_0$ of the natural 
numbers and the cardinality $2^{\aleph_0}$ of the real numbers. This conjecture 
cannot be proven because there are models\index{Model theory} of set 
theory\index{Set theory}, which satisfy the continuum 
hypothesis\index{Continuum 
hypothesis}, as Kurt G\"odel\index{Goed@G\"odel, Kurt} proved with the help of 
the cumulative hierarchy\index{Cumulative hierarchy}, and other 
models\index{Model theory} that do not satisfy it. The latter was shown by Paul 
Cohen\index{Cohen, Paul} with the help of the forcing method\index{Forcing}. In 
this method, new sets are 
constructed from a given and sufficiently small model\index{Model theory} of 
set 
theory\index{Set theory}, by the refined trick of adjunction, which also occurs 
in algebra.\footnote{See \cite{cohen1963, cohen1964} and \cite[Ch. 
7]{hoffmann2018}. The forcing method is inspired by G\"odel's completeness 
theorem\index{Goed@G\"odel, Kurt} and the L\"owenheim-Skolem theorem. A modern 
proof of Cohen's theorem comes from Lawvere\index{Lawvere, William} and 
Tierney\index{Tierney, Myles}. It can be found in \cite[Ch. 
VI]{maclane_moerdijk1992}.}

Interestingly, Skolem\index{Skolem, Thoralf} had anticipated the idea of 
the forcing method\index{Forcing} in set theory\index{Set theory} 
long before. He wrote in a paper from 1922: 
\begin{quote} 
It would certainly be of much greater interest if one could prove 
that a new subset $Z_0$ could be adjoined without 
contradictions arising; but this will probably be very 
difficult.\footnote{Machine translation of \cite[p. 149]{skolem}.}
\end{quote}

\section{Free Objects and Adjunction}

The method of adjunction is very common in mathematics and consists of adding 
new independent generating elements to a given structure, so that certain 
properties are preserved, but the new elements acquire the desired properties 
in a controlled manner. 

The adjunction method consists of two steps. In the first step, a generic new 
element is added, which is independent of the elements of the initial 
situation. In this way, a free structure is created. In the second step, 
conditions in the form of relations are imposed on the new elements, so that 
the calculation rules of the old elements are preserved and the newly 
constructed object fulfils the desired rules.

A beautiful example of a free object are the natural numbers $\mathbb{N}$. They 
form a free monoid\footnote{A monoid is a set $M$ together with an associative 
operation $\ast \colon M \times M \longrightarrow M$ and a neutral element $e 
\in M$. It is free, when there is a subset $E \subseteq M$, such that every 
element of $M$ can be uniquely written as a combination of finitely many 
elements of $E$.} with addition as the operation and a single generating 
element $1$, because every natural number is a sum of values of the number $1$:
\[
5=1+1+1+1+1. 
\]
The number $0$ also falls under this definition, as it arises when 
the number $1$ is summed up zero times. 

Free monoids can have more than one generating element. The free monoid over 
two elements consists -- written multiplicatively -- of finite words in the 
letters $a$ and $b$. Such words can be graphically visualised as nodes in a 
binary tree (see Fig.~\ref{fig:words}). 

\begin{figure}[ht!]
\begin{tikzpicture}[level distance=1cm,
  level 1/.style={sibling distance=2cm},
  level 2/.style={sibling distance=1cm}]
  \node {Empty word}
    child {node {a}
      child {node {$a^2$}}
      child {node {$ab$}}
    }
    child {node {$b$}
    child {node {$ba$}}
      child {node {$b^2$}}
    };

\end{tikzpicture}
\caption{\label{fig:words}Binary tree for words in two letters $a,b$.}
\end{figure}

Similarly, there are free groups\footnote{A (multiplicative) group $G$ is 
free, when there is a subset $E \subseteq G$, such that every element of 
$G$ can be uniquely written as a combination of finitely many elements of $E$ 
and their inverses.} over one or more generating elements. The free group with 
a single generating element is the set of integers $\mathbb{Z}$. Free groups 
with two or more generating elements are no longer commutative and are 
therefore preferably noted as multiplicative groups, as the symbol $+$ is 
usually reserved for commutative groups.

The free group $F_2$ over two generating elements 
$a$ and $b$ consists of arbitrary finite words in the two letters 
and their inverses $a^{-1}$ and $b^{-1}$, such as the word
\[
w=aba^{-1}b^{-1}. 
\]
Free objects possess a characteristic universal property. This 
states for $F_2$ that for every group $G$ with two generating 
elements, whether commutative or not, there is a mapping 
\[
F_2 \xlongrightarrow{~~~~~~} G  
\]
which maps the two generating elements of $F_2$ onto the generating 
elements of $G$. This mapping is surjective, i.e., it reaches 
every element of $G$, but it is generally not injective, unless 
$G$ is itself a free group. A similar universal property applies 
for free monoids and for free commutative groups $\mathbb{Z}^n$ of 
rank $n$. 

With this knowledge, we now want to look at a concrete example of the 
adjunction method. Given a ring $R$, there is the possibility to add a new 
element $T$ to add to the object $R$, and to form the free object
\[
R[T],
\]
which is called the polynomial ring over $R$. The elements in $R[T]$ 
are polynomials, i.e., finite sums 
\[
f=a_dT^d+ \cdots +a_2T^2 +a_1T+a_0, 
\]
where the coefficients $a_i$ are elements of $R$ and the powers $T^i$ 
are new expressions. All powers form infinitely many new elements and 
the finite sums of such expressions form an infinitely large set. The 
polynomial ring $R[T]$ is therefore a much larger object than $R$ itself. The 
symbol $T$ is usually called an unknown or variable. The addition and 
multiplication in this new object $R[T]$ is explained in an obvious way.

The polynomial ring $R[T]$ is a free structure, because no conditions were 
imposed on the variable $T$ and the powers $T^i$ form a copy of the natural 
numbers, because the exponent $i$ can exactly run through the natural numbers. 

In the second step of the adjunction, a relation $f=0$ is imposed. This 
constructs a new object, denoted by 
\[
S=R[T]/(f).
\]
In $S$, two polynomials from $R[T]$ are identified when their difference is a 
multiple of $f$. This procedure thus enforces the equation $f=0$.

As an example, we take the polynomial $f=T^2-2$ and form the object
\[
S=R[T]/(T^2-2). 
\]
What corresponds to the element $T$ in this construct? Since we have enforced 
the equation 
\[
T^2-2=0
\]
with the relation, it holds for the image of the variable $T$ in $S$ that 
$T^2=2$ is, i.e., $T$ is a square root of $2$:
\[
T=\sqrt{2}. 
\]
Through two steps, we have constructed a new ring $S$ from the ring $R$ in 
which the element $2$ has a square root, even though it did not necessarily 
have a square root in $R$ before. We write this as  
\[
S=R[\sqrt{2}]. 
\]
The ring $S$ is different from $R$ exactly when $\sqrt{2}$ did not exist in 
$R$. In the case $R=\mathbb{Z}$, this is obviously the case.

Another example is given by the polynomial $f=T^2$. In this case,
\[
S=R[T]/(T^2) 
\]
is the ring of infinitesimal numbers over the ring $R$, which was used by 
Leibniz\index{Leibniz, Gottfried Wilhelm} in his infinitesimal calculus. The 
element $T$ in this ring can be interpreted as an infinitesimal object.

\section{Geometry and the Concept of Space}

Geometric shapes model many situations in applications of mathematics. We first 
consider as simple geometric structures as possible. The simplest 
$0$-dimensional geometric object is obviously the point, which we denote by 
$\Delta_0$ for our following purposes. In dimension $1$ there is the unit 
interval $[0,1]$, which we denote by $\Delta_1$. The circle $S^1$ arises from 
the interval $\Delta_1$, when the starting and ending point are identified. In 
this respect, $\Delta_1$ is the prototypical geometric object of dimension one. 
In two dimensions there is the triangle $\Delta_2$, the unit square $\Box_2$ 
and many other polygons. By triangulation, all these objects can be decomposed 
into sums of the simplest triangles $\Delta_2$. For example, the unit square 
$\Box_2$ is made into the sum of two triangles by inserting a diagonal. 

In higher dimensions it becomes more complicated. Already in antiquity, the 
$5$ Platonic and the $18$ Archimedean bodies were found. The Platonic bodies 
are, besides the cube $\Box_3$ and the tetrahedron $\Delta_3$, the octahedron, 
the dodecahedron and the icosahedron. In addition to Archimedean 
solids, there are also truncated bodies, such as the truncated icosahedron or 
football body, which consists of 5- and 6-sided polygons and has been used as 
a model for the football since 1970. From dimension 4, not all of these 
geometric shapes generalise and only the three series $\Delta_n$ (the 
$n$-simplex\index{Simplex} or hypersimplex), $\Box_n$ (the hypercube) and 
$\Diamond_n$ (the cross polytope, the generalisation of the octahedron) remain. 
Of these three series, $\Delta_n$ is the most common and both $\Box_n$ and 
$\Diamond_n$ can be decomposed into simplices\index{Simplex} $\Delta_n$. It is 
practical to coordinate $\Delta_n$ by $n+1$ coordinates $t_0,\ldots,t_n$ that 
satisfy the equations 
\[
0 \le t_i \le 1 \text{ and } t_0+ \cdots + t_n=1. 
\]
The boundary $\partial \Delta_n$ of a single simplex\index{Simplex} $\Delta_n$ 
is given by the union of all parts of the form $\Delta_{n-1}$ that arise by 
setting a coordinate $t_i$ to zero. 

All the spaces we have considered so far are vivid elementary examples of 
topological spaces\index{Topology}. A general topological space\index{Topology} 
in mathematics is usually described as a set of points in which certain subsets 
are distinguished as open sets.\footnote{A topological space is a set $X$, in 
which a collection of open subsets $U$ is distinguished. Here, $X$ and the 
empty set $\emptyset$ are open and (arbitrary) unions and finite intersections 
of open sets are open.} A useful class of topological spaces\index{Topology} 
are the metric spaces, on which a concept of distance exists between any two 
points. The three-dimensional Euclidean space surrounding us $\mathbb{R}^3$ 
defines for two points $a=(a_1,a_2,a_3)$ and $b=(b_1,b_2,b_3)$ the Euclidean 
distance\footnote{Due to the general theory of relativity, space is curved, so 
this statement is only approximately correct.} 
\[
d(a,b)=\sqrt{(a_1-b_1)^2+(a_2-b_2)^2+(a_3-b_3)^2}. 
\]
However, in mathematics there are numerous examples of topological 
spaces\index{Topology} that are quite different from Euclidean space. For 
example, there are the non-Euclidean geometries, discovered in the 18th 
century by J\'anos Bolyai\index{Bolyai, J\'anos}, Nikolai 
Lobachevsky\index{Lobachevsky, Nikolai} and Carl Friedrich 
Gau{\ss}\index{Gau{\ss}, Carl Friedrich}.\footnote{See \cite{struik} for the 
history of non-Euclidean geometry.} In these geometries, Euclid's\index{Euclid} 
parallel postulate does not hold. With the help of differential geometry, which 
Bernhard Riemann\index{Riemann, Bernhard} subsequently developed, such spaces 
could be realised in particular as hyperbolic manifolds with constant negative 
curvature, in which the straight lines are given as geodesics, i.e., shortest 
connections between points. Manifolds are -- figuratively speaking -- geometric 
spaces with good smoothness properties that locally look like subsets of 
$\mathbb{R}^n$. In particular, a tangent space can be defined at each point. 

The structure-preserving mappings between topological spaces\index{Topology} 
are the continuous mappings. They are defined by the fact that they preserve 
the concept of distance in a certain way. This is equivalent to the fact that 
for any two points $a,b \in X$ that are close enough together, the image points 
$f(a), f(b)$ are not far apart.\footnote{Mathematically, this is defined so 
that the pre-images of open sets are again open. In a metric space, this 
property is defined by the well-known $\varepsilon$-$\delta$-definition.} Two 
spaces have an indistinguishable topological structure -- i.e., they are 
homeomorphic\index{Homeomorphism} -- when there is a bijective continuous 
mapping $f$ between them, such that both $f$ itself and the inverse mapping of 
$f$ are continuous. 

\section{Symmetries}

Topological spaces\index{Topology} often possess symmetries. For example, the 
three-dimensional Euclidean space $\mathbb{R}^3$ that surrounds us has a 
symmetry group consisting of translations and rotations that preserve 
lengths and angles. Anyone who jumps from a diving board in 
a swimming pool utilises these possibilities of movement. The Platonic 
and Archimedean solids and many ornaments in artworks and buildings since 
antiquity possess finite symmetries. Symmetry groups in physics 
are closely related to conservation quantities according to a theorem by Emmy 
Noether\index{Noether, Emmy}.\footnote{Contained in \cite{noether}.} 

A symmetry group $G$ operates on a topological space\index{Topology} $X$ by 
each group element $g$ forming a homeomorphism\index{Homeomorphism} $g \colon X 
\longrightarrow X$. The sets $\{gx \mid g\in G \}$ for a fixed point $x\in X$ 
are referred to as orbits and there is a continuous mapping $\pi \colon X 
\longrightarrow X/G$ from $X$ to the orbit space $X/G$. An example is the real 
number line $X=\mathbb{R}$ and the additive group $G=\mathbb{Z}$, which 
operates 
as a symmetry group through translation. The orbit of each point $x$ is given 
by 
all translates $x+n$ with $n \in \mathbb{Z}$. Through the quotient mapping 
$\pi(x)=\exp(2\pi i x)$, the orbit space $\mathbb{R}/\mathbb{Z}$ is identified 
with the unit circle $S^1$.

There are conditions for so-called good group operations on 
manifolds $X$, so that $X/G$ still remains a manifold 
and $\pi$ is an unbranched covering. The latter means that all 
orbits are isomorphic to $G$ and for the mapping $\pi$ for each $x \in X$ 
there exists a sufficiently small open set $U$ in $X$, so that $\pi$ 
forms a homeomorphism\index{Homeomorphism} from $U$ to $\pi(U)$. 
For less good operations, the quotients $X/G$ form a class of spaces under 
certain conditions on the operation of $G$, which are 
called orbifolds or stacks. Individual points in $X/G$ possess 
automorphisms as additional information and can become singular. 
A simple example of this is the operation of the group $G=\mu_n$ of the 
complex $n$-th roots of unity on $X=\mathbb{C}$. Then $X/G$ is again 
homeomorphic to  $\mathbb{C}$, it holds $\pi(z)=z^n$ and the zero point 
possesses the entire group $\mu_n$ as an automorphism group. Outside of the 
zero point, the operation is good. 

\section{Simplicial Spaces\index{Simplicial space}} 

One way to go from a geometric object to a simpler structure is the method of 
surveying, i.e., triangulation. Here, support points are chosen, connection 
curves between them, surface pieces between the connection lines, and so on up 
to higher dimensions. In Fig.~\ref{fig:triangulation}, the triangulation of 
a sphere $S^2$ is indicated. 

\begin{figure}[ht!]
\tdplotsetmaincoords{60}{111}
\begin{tikzpicture}[line join=round,tdplot_main_coords,declare function={a=5;}]
\coordinate (O) at (0,0,0);
\draw[tdplot_screen_coords] (0,0,0) circle (2);
\path (xyz spherical cs:radius=2,latitude=90,longitude=00)
node[fill,circle,inner sep=2pt]{}; 
\foreach \X in {0,1,2} 
{\draw plot[variable=\t,domain=90:90-109.5,smooth] 
(xyz spherical cs:radius=2,latitude=\t,longitude=\X*120)
node[fill,circle,inner sep=2pt]{};}
\draw plot[variable=\t,domain=0:360,smooth] 
(xyz spherical 
cs:radius=2,latitude={90-109.5-14*abs(sin(1.5*\t))},longitude=\t);
\end{tikzpicture}
\hskip2cm
\begin{tikzpicture}[line join=round,tdplot_main_coords,declare 
function={a=3.7;}]
\begin{scope}[canvas is xy plane at z=0,transform shape]
 \path foreach \X [count=\Y] in {A,B,C}
 {(\Y*120:{a/(2*cos(30))}) coordinate(\X)};
\end{scope}
\path (0,0,{a*cos(30)}) coordinate (D);
\draw foreach \X/\Y [remember=\X as \Z (initially D)] in {A/B,B/C,C/D,D/A}
 {(\X) -- (\Z) -- (\Y)};
\end{tikzpicture}
\caption{\label{fig:triangulation}Sphere triangulated with four points together 
with a tetrahedron which resembles the simplicial abstraction of this 
triangulation.} 
\end{figure}

A triangulation of a topological space\index{Topology} is an attempt to break 
it down into the simplest components, which are given by continuous mappings $f 
\colon \Delta_n \longrightarrow X$. The resulting combinatorial object, which 
is glued together from simplices\index{Simplex}, is called a simplicial 
space\index{Simplicial space}. Popular geometric bodies, such as spheres $S^n$, 
can be easily assembled in a combinatorial way from two or more copies of 
simplices\index{Simplex} $\Delta_n$ that are glued along the edges 
$\partial \Delta_{n}$.

Graphs are very simple simplicial spaces\index{Simplicial space}, consisting 
only of $0$- and $1$-dimensional simplices\index{Simplex}, i.e.,  of nodes 
(vertices) and edges. They form important structures in mathematics and are as 
fundamental as numbers. Graphs often describe combinatorial relationships. 
Dependencies between people, organisations, and objects in the world can also 
be described using graphs. The simplest graphs are the trees. They contain no 
multiple edges and no loops, i.e., closed paths, that return to the starting 
point.

Simplicial spaces\index{Simplicial space} are special examples 
of geometric realisations of simplicial 
sets\index{Simplicial set}, the totality of which is denoted by 
$\mathbf{sSet}$. These consist of an abstract set of 
simplices\index{Simplex} in each degree $n$ and from gluing data between 
the simplices\index{Simplex}. Simplicial sets\index{Simplicial set} 
can be complicated and contain many simplices\index{Simplex}. Thus 
all continuous mappings $f \colon \Delta_n \longrightarrow X$ 
for a topological space\index{Topology} $X$ can be considered and we 
obtain a simplicial set\index{Simplicial set}, which is denoted 
by $\mathrm{Sing}_\bullet(X)$ and plays an important role. 
Simplicial sets\index{Simplicial set} can be seen as 
combinatorial-algebraic abstraction of topological spaces\index{Topology}. 
Their geometric realisations possess good 
properties.\footnote{Simplicial sets are contained in the class {\bf 
CGHaus} of compactly generated Hausdorff spaces, which include CW-complexes 
and locally compact Hausdorff spaces. See \cite[Ch. I]{goerssjardine} on 
simplicial sets.}

\section{Paths, Fundamental Group\index{Fundamental group} 
and Homotopies\index{Homotopy theory}}

The actual shape of topological spaces\index{Topology} and 
their generalisations is not decisive in many contexts, 
but rather the equivalence class\index{Equivalence}, which arises when 
continuous deformations, called homotopies\index{Homotopy theory}, 
are allowed. Their equivalence classes\index{Equivalence} 
are called homotopy types\index{Homotopy theory}.

\begin{figure}[ht!]
\begin{tikzpicture}
\draw [ultra thick,gray]  plot [smooth, tension=4] coordinates
{ (0,0) (1,1.5) (2,0.5) (3,0)};
\draw [ultra thick,gray] plot [smooth, tension=2] coordinates 
{ (0,0) (1,-0.3) (2,-1.5) (3,0)};
\draw[fill=black] (0,0) circle (.5ex);
\node[] at (0,-0.3) {$a$};
\draw[fill=black] (3,0) circle (.5ex);
\node[] at (3,0.3) {$b$};
\node[] at (0.5,1.5) {$p$};
\node[] at (2,-1.7) {$q$};
\end{tikzpicture}
\caption{\label{fig:path}Paths $p,q$ from one point $a$ to another point $b$ 
can be quite different and can have self-intersections.}
\end{figure} 

What is a homotopy\index{Homotopy theory}? We will first explain paths in 
topological spaces\index{Topology} and then expand the concept. A path $p$ in 
a topological space\index{Topology} $X$ is given by a continuous mapping
\[
p \colon [0,1] \xlongrightarrow{~~~~~~} X, 
\]
where $[0,1]$ is the real unit interval. Here, $p(0)=a$ is the starting point 
and $p(1)=b$ is the end point of the path (see Fig.~\ref{fig:path}).

Paths $p$ and $q$ can naturally be composed to a path $q \circ p$ 
when the end point of $p$ coincides with the starting point of $q$ 
(see Fig.~\ref{fig:path_composition}). Expressed in mathematical language, 
\[
(q \circ p)(t)= 
\begin{cases} 
p(2t) & t < \frac{1}{2} \cr
b=p(1)=q(0) & t=\frac{1}{2} \cr
q(2t-1) & t > \frac{1}{2}.
\end{cases}
\]  
            
\begin{figure}[ht!]
\begin{tikzpicture}
\draw [ultra thick,gray]  plot [smooth, tension=4] coordinates 
{ (0,0) (1,1.5) (2,0.5) (3,0)};
\draw [ultra thick,gray] plot [smooth, tension=2] coordinates 
{ (3,0) (1,-0.3) (2,-1.5) (6,0)};
\draw[fill=black] (0,0) circle (.5ex);
\node[] at (0,-0.3) {$a$};
\draw[fill=black] (3,0) circle (.5ex);
\node[] at (3,0.3) {$b$};
\node[] at (0.5,1.5) {$p$};
\node[] at (3,-1.2) {$q$};
\draw[fill=black] (6,0) circle (.5ex);
\node[] at (6,0.3) {$c$};
\end{tikzpicture}
\caption{\label{fig:path_composition}The composition $q \circ p$ of two paths 
$p,q$ is a concatenation.}
\end{figure} 

The inverse path $p^{-1}$ to a path $p$ from $a$ to $b$ is a path from $b$ to 
$a$ and is given by reversing the direction of path $p$. It is given by the 
function $p^{-1}(t)=p(1-t)$, where $t$ is the coordinate in the interval 
$[0,1]$.

\medskip
Two paths $p,q$ with the same starting and ending points are 
homotopic\index{Homotopy theory} to each other, denoted $p \simeq q$, when 
there is a homotopy\index{Homotopy theory} between them. Such is given by 
a continuous mapping 
\[
h \colon [0,1]\times [0,1] \xlongrightarrow{~~~~~~} X, 
\]
so that $p(t)=h(0,t)$ and $q(t)=h(1,t)$ define the two paths (see 
Fig.~\ref{fig:path_homotopy}). For the defined operations, the following 
calculation rules apply up to homotopy\index{Homotopy theory}:
\begin{align*}
p^{-1} \circ p & \simeq 1_a \\
p \circ p^{-1} & \simeq 1_b \\
(p \circ q) \circ r & \simeq p \circ (q \circ r). 
\end{align*}
Here, $1_a$ and $1_b$ are the constant paths at $a$ and $b$. The 
homotopy\index{Homotopy theory} of paths in these formulas is based on a 
reparametrisation of the paths as mappings from the unit interval to $X$. The 
formulas become incorrect if the symbol $\simeq$ is replaced by 
equality\index{Equality}\index{Identity}, because already the first formula is 
only valid up to homotopy\index{Homotopy theory}.

\begin{figure}[ht!]
\begin{tikzpicture}
  \coordinate (O) at (0,0,0);
  \coordinate (A) at (3,0,0);

  \draw[color=gray] (O)--(A);
  \draw[color=gray] (O) to [bend left=5] (A);
  \draw[color=gray] (O) to [bend right=5] (A);
  \draw[color=gray] (O) to [bend left=10] (A);
  \draw[color=gray] (O) to [bend right=10] (A);
  \draw[color=gray] (O) to [bend left=15] (A);
  \draw[color=gray] (O) to [bend right=15] (A);
  \draw[color=gray] (O) to [bend left=20] (A);
  \draw[color=gray] (O) to [bend right=20] (A);
  \draw[color=gray] (O) to [bend left=25] (A);
  \draw[color=gray] (O) to [bend right=25] (A);
  \draw[color=gray] (O) to [bend left=30] (A);
  \draw[color=gray] (O) to [bend right=30] (A);
  \draw[color=gray] (O) to [bend left=35] (A);
  \draw[color=gray] (O) to [bend right=35] (A);
  \draw[color=gray] (O) to [bend left=40] (A);
  \draw[color=gray] (O) to [bend right=40] (A);
  \draw[color=gray] (O) to [bend left=45] (A);
  \draw[color=gray] (O) to [bend right=45] (A);
  \draw[color=gray] (O) to [bend left=50] (A);
  \draw[color=gray] (O) to [bend right=50] (A);
  \draw[color=gray] (O) to [bend left=55] (A);
  \draw[color=gray] (O) to [bend right=55] (A);
  \draw[color=gray] (O) to [bend left=60] (A);
  \draw[color=gray] (O) to [bend right=60] (A);
  \draw[color=gray] (O) to [bend left=65] (A);
  \draw[color=gray] (O) to [bend right=65] (A);
  \draw[color=gray] (O) to [bend left=70] (A);
  \draw[ultra thick,color=darkgray] (O) to [bend right=70] (A);
  \draw[ultra thick,color=darkgray] (O) to [bend left=70] (A);
\draw[fill=black] (0,0) circle (.5ex);
\node[] at (-0.2,0) {$a$};
\draw[fill=black] (3,0) circle (.5ex);
\node[] at (3.2,0) {$b$};
\node[] at (1.5,1) {$p$};
\node[] at (1.5,-1) {$q$};
\end{tikzpicture}
\caption{\label{fig:path_homotopy}A homotopy between two paths $p$ and $q$ from 
$a$ to $b$ is a continuous deformation.}
\end{figure} 

The fundamental group\index{Fundamental group} $\pi_1(X,*)$ is the group of 
homotopy classes\index{Homotopy theory} of all paths in $X$ with fixed start 
and end point $*$, the base point. It is usually not commutative, as the 
example of the \enquote{lying figure eight} (infinity symbol) shows, whose 
fundamental group\index{Fundamental group} is the free group $F_2$ over two 
generating elements, which correspond to the two obvious paths in 
Fig.~\ref{fig:lying_eight}. 

\begin{figure}[ht!]
\begin{tikzpicture}
\draw[] (0,0) circle (0.5cm);
\draw[] (1,0) circle (0.5cm);
\draw[fill=black] (0.5,0) circle (1mm);
\end{tikzpicture}
\caption{\label{fig:lying_eight}The \enquote{lying figure eight} given by two 
circles touching at a  
base point.}
\end{figure}

The fundamental group\index{Fundamental group} has some disadvantages. It 
obscures individual properties of paths in $X$ and the dependence on the base 
point is cumbersome through the homotopies\index{Homotopy theory}. Moreover, 
there are interesting topological spaces\index{Topology} that hardly allow real 
paths and for which our definition therefore provides nothing. For example, the 
usual definition is usually worthless in algebraic geometry. In this case, it 
is better to define the fundamental group\index{Fundamental group} using 
so-called \'etale coverings of $X$, which are particularly unbranched. The 
standard definition is only suitable for topological spaces\index{Topology} 
that have enough real paths. 

Often it is better to consider paths between all points, without 
necessarily limiting oneself to homotopy classes\index{Homotopy theory} and a 
base point. In the general case, we therefore work best 
with a higher categorical structure\index{Category theory} 
$\Pi_\infty(X)$, which is called the fundamental infinity groupoid or more 
generally infinity groupoid\index{Infinity groupoid}. In it, all 
calculation rules only apply up to certain equivalences\index{Equivalence}, 
i.e., homotopies\index{Homotopy theory} and higher generalisations thereof. 

What is a general homotopy\index{Homotopy theory} and a 
homotopy\index{Homotopy theory} equivalence\index{Equivalence}? For this, we 
consider two continuous mappings 
\[
f,g \colon X \xlongrightarrow{~~~~~~} Y
\]
and call both homotopic\index{Homotopy theory}, if there is a continuous 
mapping
\[
h \colon [0,1] \times X \xlongrightarrow{~~~~~~} Y
\]
with $h(0,-)=f$ and $h(1,-)=g$. Two spaces $X,Y$ are then called 
homotopy equivalent\index{Homotopy theory}, if there are continuous mappings 
\[
F \colon X \xlongrightarrow{~~~~~~} Y \text{ and } G \colon Y 
\xlongrightarrow{~~~~~~} X,
\]
so that $G \circ F$ and $F \circ G$ are each homotopic\index{Homotopy theory} 
to the identity\index{Identity}\index{Equality} on $X$ and $Y$. There is a 
somewhat more general concept of weak homotopy 
equivalence\index{Equivalence}\index{Homotopy theory} 
between two spaces\index{Topology}. They preserve essential 
characteristics and invariants\index{Invariant} of topological 
spaces\index{Topology}. 

\section{Topological Invariants\index{Invariant}} 

Properties of equivalent\index{Equivalence} objects in mathematics, which are 
common to all representatives, are called invariants\index{Invariant}. 
In mathematics and physics, there are many invariants\index{Invariant} that 
carry crucial information about the objects under consideration. Invariants of 
topological spaces\index{Topology} play a special role, as they are preserved 
under homeomorphism\index{Homeomorphism} or under more general operations 
like homotopy equivalence\index{Homotopy theory}\index{Equivalence}. 

The mathematician Emmy Noether\index{Noether, Emmy} has made a decisive 
contribution to ensuring that invariants\index{Invariant} of algebraic, 
topological and physical nature were assigned their appropriate significance. 
She defined homological\index{Homology group} invariants\index{Invariant} of 
abstract algebraic chain complexes and thus initiated the algebraisation of 
topology\index{Topology}.\footnote{See \cite{ms2027} and the collected works 
\cite{noether}.} 

How does this work? Each topological space\index{Topology} can be assigned a 
natural simplicial set\index{Simplicial set} and thus a natural singular chain 
complex.\footnote{A chain complex is a sequence $\cdots \longrightarrow C_{n+1} 
\longrightarrow C_{n}\longrightarrow C_{n-1} \longrightarrow \cdots $ of 
abelian groups $C_i$ (or more generally of modules over a ring $R$) together 
with $\mathbb{Z}$-linear mappings $d_i \colon C_i \longrightarrow C_{i-1}$, 
which satisfy $d_{i-1} \circ d_i=0$.} For this, we consider the set of all 
continuous mappings of the simplices\index{Simplex} $\Delta_n$ to $X$:
\[
\mathrm{Sing}_n(X)=\{f \colon \Delta_n \longrightarrow X \text{  continuous}\}.
\]
This defines a canonical simplicial set\index{Simplicial set}, 
which is denoted by 
\[
\mathrm{Sing}_\bullet(X).
\]
It also includes boundary mappings 
\[
\partial_n \colon \mathrm{Sing}_n(X) \xlongrightarrow{~~~~~~} 
\mathrm{Sing}_{n-1}(X),
\]
which are given by restriction to the boundaries $\partial \Delta_n$ of the 
simplices\index{Simplex}. The elements of $\mathrm{Sing}_n(X)$ are called 
$n$-simplices\index{Simplex}, just like the objects $\Delta_n$ themselves. The 
$n$-simplices\index{Simplex} degenerate when $n$ is larger than the dimension 
of $X$. They are omitted for many purposes, but not always, because the 
simplicial set\index{Simplicial set} $\mathrm{Sing}_\bullet(X)$ contains 
valuable homotopy-theoretical\index{Homotopy theory} information about $X$, 
which is then lost.   

To every simplicial set\index{Simplicial set} $S_\bullet$ -- and therefore 
to every topological space\index{Topology} $X$ through the simplicial 
set\index{Simplicial set} $\mathrm{Sing}_\bullet(X)$ -- homology 
groups\index{Homology group} can be assigned, as Emmy 
Noether\index{Noether, Emmy} had introduced in general. In this 
view, we consider the chain complex of abelian groups underlying the 
simplicial set\index{Simplicial set} $S_\bullet$  
\[
\mathbb{Z}S_\bullet \colon \cdots \longrightarrow \mathbb{Z}S_n \longrightarrow 
\mathbb{Z}S_{n-1} \longrightarrow \cdots \longrightarrow \mathbb{Z}S_1 
\longrightarrow \mathbb{Z}S_0,
\]
where $\mathbb{Z}S_i$ is the free abelian group generated by the 
simplices\index{Simplex} in $S_i$. The linear mappings between these free 
groups are generated by the given boundary mappings $\partial_i \colon S_i \to 
S_{i-1}$ in $S_\bullet$. Then we form at each point the homology 
groups\index{Homology group} $H_n(S_\bullet)$ as the quotients
\[
H_n(S_\bullet)=\frac{\mathrm{Ker}(\mathbb{Z}S_n \longrightarrow 
\mathbb{Z}S_{n-1})}{\mathrm{Im}(\mathbb{Z}S_{n+1} \longrightarrow 
\mathbb{Z}S_{n})}. 
\]
The $n$-th Betti number $b_n$ is then defined as the rank of $H_n(S_\bullet)$. 
The zeroth Betti number $b_0$ gives the number of connected components of 
$S_\bullet$ and the first Betti number $b_1$ the number of loops. The $n$-th 
Betti number indicates whether a higher-dimensional \enquote{hole} exists. 

We want to give some simple applications of this formal algebraic 
definition. First, let's consider a triangle, i.e., the boundary of the 
$2$-simplex\index{Simplex} (see Fig.~\ref{fig:boundary}). 

\begin{figure}[ht!]
\begin{tikzpicture}
\draw[->] (0,0) node {}
  -- (2,0) node {}
  -- (1,1.5) node {}
  -- cycle;
\draw[fill=black] (0,0) circle (.5ex);
\draw[fill=black] (2,0) circle (.5ex);
\draw[fill=black] (1,1.5) circle (.5ex);
\end{tikzpicture}
\caption{\label{fig:boundary}An equilateral triangle is the boundary $\partial 
\Delta_2$ of a $2$-simplex $\Delta_2$.}
\end{figure}

The triangle consists of three copies of $\Delta_1$ as edges and three vertices 
$\Delta_0$. We give all edges a counter-clockwise orientation. This defines 
the end and the starting point of each edge. The chain complex $S_\bullet$  
which calculates the homology groups\index{Homology group} is of the form
\[
\mathbb{Z}^3 \xlongrightarrow{~~ \partial ~~} \mathbb{Z}^3, 
\]
where $\partial$ calculates \enquote{end point minus starting point} from each 
edge and is represented in a suitable basis by the matrix
\[
\left(\begin{matrix}
-1 & 0 & 1 \cr 
1 & -1 & 0 \cr
0 & 1  & -1 
\end{matrix}\right).
\]
The entries $+1$ stand for the endpoints and $-1$ 
for the starting points. This matrix has determinant $0$ and rank $2$, because 
the first two columns are linearly independent and the sum of the columns or 
rows is $0$. From this it follows that both the kernel of $\partial$ and the 
cokernel of $\partial$ are isomorphic to $\mathbb{Z}$. Therefore, the Betti 
numbers are $b_0=b_1=1$.

A slightly different example is the already mentioned \enquote{lying figure 
eight}, where two loops $\Delta_1$ and only one point $\Delta=0$ occur. In 
this figure, the chain complex is given by 
\[
\mathbb{Z}^2 \xlongrightarrow{~~ \partial ~~} \mathbb{Z}, 
\]
where $\partial$ is given by the zero mapping, because the end and 
starting point of the $1$-simplices\index{Simplex} coincide. Thus 
$b_1=2$ and $b_0=1$ hold. More generally, in graphs, the Betti number 
$b_1$ is equal to the number of loops.

The higher-dimensional \enquote{holes} are given by the higher 
Betti numbers $b_n$. For spheres, for example, it holds 
\[
b_i(S^n)=\begin{cases} 1 & \text{for } i=0 \text{ and } i=n \cr  
0 & \text{otherwise}. \end{cases} 
\]
The homotopy groups\index{Homotopy group} $\pi_n(X)$ of a topological 
space\index{Topology} $X$ are defined as 
homotopy classes\index{Homotopy theory} of continuous mappings 
\[
f \colon S^n \xlongrightarrow{~~~~~~} X.
\]
A suitable composition can be defined on the 
homotopy groups\index{Homotopy group}, so that they carry a 
group structure and for $n\ge 2$ are even commutative. The most important 
of these groups is the fundamental group\index{Fundamental group} 
$\pi_1(X,*)$, where $*$ is a base point in $X$. A space $X$ is called 
a homotopy $n$-type\index{Homotopy theory}, if the 
homotopy groups\index{Homotopy group}  $\pi_i(X)$ for $i > n$ are zero.
 
A famous example of a homotopy group\index{Homotopy group} is 
provided by the $2$-sphere $S^2$, for which there is the Hopf fibration
\[
S^3 \xlongrightarrow{~~~~~~} S^2 
\]
which contains exactly the crucial information about the generator of the third 
homotopy group\index{Homotopy group} $\pi_3(S^2)\cong \mathbb{Z}$.
The sphere $S^2$ is therefore not a homotopy $2$-type\index{Homotopy theory}. 

\chapter{Mathematics in our Culture}

Mathematics is a multicultural science with a rich past. We do not know exactly 
when and where it first originated. Beginnings that go beyond simple counting 
methods can be found in the ancient Orient long before Greek antiquity, among 
the Maya, in China and India, and in many other places in the world. It was 
certainly a cultural technique from the outset, which had concrete applications 
in practical areas of the societies of the time. Only in Greek antiquity did a 
mathematical science slowly develop with the fundamental areas of arithmetic, 
geometry and logic.\footnote{See \cite{struik} for the history of mathematics. 
The book \cite{rademachertoeplitz} refers to the historical significance of 
arithmetic and geometry.} 

Most people are only aware of a few applications of mathematics. In 
recent times, blockchain algorithms\index{Cryptography}, 
cryptographic\index{Cryptography} protocols and machine 
learning\index{Artificial intelligence} have gained a high level of recognition 
due to digitalisation. Numerous other examples can be found in 
technology and natural sciences as well as in the social sciences. 

Contrary to widespread beliefs, however, mathematics is not limited to
quantitative and algorithmic methods. Its actual task is to provide concepts 
and buildings for other sciences. The interlocking with physics is remarkable, 
as we have already indicated in the introduction. In the future, further 
challenges await mathematics. One of them comes from the exploration 
of complex systems, where multidisciplinary phenomena need to be understood and 
managed.  

\section{Cryptography\index{Cryptography} or the Art of 
Encryption} 

Encryption methods\index{Cryptography} are used to exchange secrets 
securely or to generate tamper-proof certificates and 
authentications.\footnote{See \cite{wigderson2019} on 
these topics.} These protocols often calculate in residue classes modulo 
a very large number $N$, i.e., with the numbers 
\[
\mathbb{Z}/N\mathbb{Z}=\{0,1,2,3,\ldots,N-1\} 
\]
and these residue classes are added and multiplied in the usual way, 
with multiples of $N$ possibly being subtracted if the result is strictly 
greater than $N-1$. This is familiar to most people when $N=7$ 
is, because the seven days of the week 
\[
\text{Mon, Tue, Wed, Thu, Fri, Sat, Sun} 
\]
form residue classes modulo $N=7$. 

In cryptography\index{Cryptography} or when entering passwords into computer 
systems, one-way functions are used, which are often also called hash 
functions. With these, it is possible to verify the confidentiality or the 
verification of the truth\index{Truth} of statements -- such as the 
authentication of a person -- without jeopardising security and anonymity. With 
mathematical methods, many hash functions can be generated and confidentiality 
protocols can be technologically realised and made secure. This is where the 
actual value of number theory in cryptography\index{Cryptography} lies. 

Hash functions are functions with the property that from a function value 
$f(x)$ the argument $x$ can only be calculated with great difficulty. 
They are called one-way functions when the corresponding function $f$ 
is injective, i.e., when different arguments have different function values. 
Hash functions are comparable to classical or genetic fingerprints, as they 
possess similar properties. The powers
\[
g^m \hskip-3mm \mod N
\]
form the basis of number-theoretic hash functions and are used in RSA 
encryption, ElGamal encryption, in Diffie-Hellman key exchange and in Shamir's 
three-pass protocol. Here we can either consider the base $g$ or the exponent 
$m$ as a secret. 

The calculation of $m$ given the power $g^m$ of $g$ is  
called the discrete logarithm problem. Neither trying out all 
possibilities nor index calculations or the sophisticated 
algorithms\index{Algorithm} of Pohlig-Hellman, Pollard and Shanks
significantly alleviate this difficulty. The security of 
ElGamal encryption and the Diffie-Hellman key exchange rely on this. 

In many protocols, $N=p$ is a large prime number. In practice, often only the 
elements in $\mathbb{Z}/N\mathbb{Z}$ are used that are coprime to $N$. These 
form a multiplicative subgroup $U_N$ of units in $\mathbb{Z}/N\mathbb{Z}$. In 
the case of RSA encryption, which we will now go into in a little more detail, 
$N=pq$ is the product of two different large prime numbers with more than a 
hundred digits.\footnote{See \cite[Kap. 5]{mspiontkow2011} for details on some 
protocols. In the ElGamal protocol, a suitable subgroup $G \subseteq U_p$ with 
$q$ elements is often used, where $q$ is a Sophie-Germain prime, i.e., 
$p=2q+1$ is also a prime number. The mathematician Sophie 
Germain\index{Germain, 
Sophie} dealt with Fermat's conjecture for such prime numbers as exponents.}

\section{The Security of RSA Encryption}

In the case of RSA encryption, the secret is given by a number 
$g \in \mathbb{Z}/N\mathbb{Z}$. The security of the procedure is guaranteed if 
it is difficult for large $N$ to determine the base $g$ from the knowledge of 
$g^m$, even if $m$ is known. 

In the case $N=p$ prime, there is a simple decryption method for this problem, 
called the backdoor. For this, we search with the help of the Euclidean 
algorithm\index{Algorithm} for a $m'$ with 
\[
mm' \equiv 1 \hskip-3mm \mod p-1
\]
and Fermat's little theorem\footnote{Fermat's little theorem states that 
modulo a prime number $p$, always $g^{p} \equiv g$ for all $g$. From this 
follows $g^{p-1} \equiv 1$ for all $g$ coprime to $p$.} guarantees that
\[
(g^m)^{m'} \equiv g^{mm'} \equiv g \hskip-3mm \mod N. 
\]
The decryption is therefore also carried out with a simple 
exponential function. This is obviously not a particularly secure method if 
$m$ is known, because the number $m'$ is easy to calculate from $m$ if 
$N=p$ is a prime number. 

This behaves quite differently when $N=pq$ is the product of two large 
different prime numbers. The decryption $g$ of $g^m$ 
given exponent $m$ is also obtained in this case by 
exponentiation with a suitable $m'$, if the factorisation of $N=pq$ 
is known. This time, $m'$ is not so easy to determine, but requires 
Euler's theorem, which generalises Fermat's little theorem. This 
theorem states that 
\[
g^{\varphi(N)} \equiv 1 \hskip-3mm \mod N, 
\]
if $g$ is coprime to $N$, where $\varphi(N)=(p-1)(q-1)$ is the Euler 
$\varphi$-function\footnote{The value of the Euler $\varphi$-function 
$\varphi(N)$ at the point $N$ gives the number of numbers coprime to $N$ 
between $1$ and $N-1$, i.e., the order of $U_N$.} of $N$. The 
decryption (or backdoor) is then given by exponentiation with $m'$, 
where 
\[
mm' \equiv 1 \hskip-3mm \mod \varphi(N).
\]
RSA encryption is therefore secure because it is just as difficult to factorise 
a large number $N$ of the form $N=pq$ as it is to calculate 
$\varphi(N)=(p-1)(q-1)$. A malicious person can, even if they know $N$, only 
calculate the number $m'$ and thus $g$ from the knowledge of $g^m$ and $m$ with 
unreasonable time and computational effort. 

Factoring large natural numbers poses a problem, which is presumably more 
difficult than verifying whether a number is prime. In 2002, it was shown that 
prime number tests in polynomial time\index{Polynomial time} are possible and 
for this the expression \enquote{Primes is in P} was coined. For the 
factorisation of numbers, there are algorithms\index{Algorithm}, such as the 
quadratic sieve or the number field sieve, which have subexponential 
complexity\index{Complexity theory}. However, no deterministic 
polynomial time\index{Polynomial time} algorithm\index{Algorithm} is known. 
Peter Shor\index{Shor, Peter} proved in 1994 that a quantum 
computer\index{Quantum computing} would achieve this in polynomial 
time\index{Polynomial time} in a probabilistic way,\footnote{See 
\cite{shor,wigderson2019}.} i.e., the computation must be repeated until the 
output becomes reliable. Because of this fact, cryptographic\index{Cryptography} 
postquantum protocols are being researched.

\section{A Protocol for Exchanging Secrets}

In Shamir's three-pass protocol, two people, usually referred to as Alice and 
Bob in cryptography\index{Cryptography}, want to exchange a secret without a 
malicious person, often called Eve (from English evil), being able to uncover 
it. The idea can be easily explained without mathematics. The concept is that 
the secrets are sent via a regular, lockable suitcase, to which Alice and Bob 
can attach locks with their own keys. Alice puts the secret in the suitcase and 
locks it with her own key, which no one else possesses, not even Bob and 
certainly not a malicious Eve. Alice hands the suitcase to Bob, who in turn 
locks the suitcase a second time and returns it to Alice. Now Alice has the 
opportunity to unlock her own lock without the suitcase being able to be 
opened, because Bob's lock is still attached. Finally, she gives the suitcase 
back to Bob, who can now unlock it and share the secret with Alice. During this 
process, Eve never had the opportunity to find an unlocked suitcase.

If Alice and Bob's locks are realised with natural numbers, Alice chooses 
private numbers $m$ and $m'$ with 
\[
mm' \equiv 1 \hskip-3mm \mod p-1 
\]
and Bob his private numbers $n$ and $n'$ with 
\[
nn' \equiv 1 \hskip-3mm \mod p-1.
\]
Here, $p$ is a large prime number which was previously agreed 
upon. Then the protocol is realised by Alice first sending the number $g^m$ to 
Bob, then Bob forming the power $(g^m)^{n}=g^{mn}$ and sending it back to 
Alice, who in turn forms the $m'$-th power:
\[
(g^{mn})^{m'}=g^{mm'n} \equiv g^{n} \hskip-3mm \mod p.
\]
Finally, Bob exponentiates this number again with the $n'$-th power and obtains
\[
(g^n)^{n'}=g^{nn'} \equiv g \hskip-3mm \mod p.  
\]
Alice has thus successfully sent the secret $g$ to Bob. Along the way, the 
numbers $g^m$, $g^{mn}$ and $g^{n}$ appear, from which $g$ cannot be obtained 
with realistic effort without essentially trying out all possible powers.

Similar methods can be used to realise digital signatures\index{Cryptography}, 
i.e., signatures and other authenticity certificates, and enable 
zero-knowledge proofs\index{Cryptography}, with which the truth\index{Truth} 
of statements can be checked without seeing the contents.

\section{Blockchains\index{Cryptography}}

The fascinating new technology of blockchains\index{Cryptography} is also based 
on hash functions. Commercial hash functions such as the secure-hash algorithm 
SHA 256 are used, which differ from the hash functions we have considered.

Blockchains\index{Cryptography} were invented by Satoshi 
Nakamoto\index{Nakamoto, Satoshi}. This name is a pseudonym and it is not known 
who is behind it. In the only work under this name, this person wrote:
\begin{quote}
\hskip2mm A purely peer-to-peer version of electronic cash would allow online 
payments to be sent directly from one party to another without going through a 
financial institution. Digital signatures provide part of the solution, but the 
main benefits are lost if a third party is still required to prevent 
double-spending. We propose a solution to the double-spending problem using a 
peer-to-peer network. The network timestamps transactions by hashing them 
into an ongoing chain of hash-based proof-of-work, forming a record that 
cannot be changed without redoing the proof-of-work.\footnote{The unpublished 
work is available at www.bitcoin.org.} 
\end{quote}
A blockchain\index{Cryptography} is a list of blocks that represent 
confidential data and are gradually extended. The processing of 
blockchains\index{Cryptography} takes place in communities of computer 
networks, whose nodes possess private 
identities\index{Identity}\index{Equality}. There is no central organisation, 
such as a bank, that monitors the proceedings. The appending of new blocks 
takes place in a mathematical consensus procedure in the form of a 
proof-of-work. The basic idea of blockchains\index{Cryptography} is that each 
block is connected with the past blocks, i.e., the data of previous processes 
and the identities\index{Identity}\index{Equality} of the participants are 
preserved at each step in the form of values of hash functions and can be 
individually stored as evidence (see Fig.~\ref{fig:blockchain}). The 
verification of new transactions in the proof-of-work is extremely 
computationally intensive. This requires a mathematical performance, which is 
called mining. 

\begin{figure}[ht!]
\begin{tikzpicture}[
  scale=0.8,
  node distance=9mm,typetag/.style={rectangle,draw=black!50,font=\scriptsize, 
anchor=west}
]
  \node (a) {transaction 1};

  \node [below=of a.west, typetag, xshift=2mm] (a1) { owner 1 public key };
  \node [below=of a1.west, typetag] (a2) { hash };
  \node (a3) [below=of a2.west, typetag] { owner 0 signature };
  \node (a4) [below=of a3] 
  { owner 1  private key}; 
  
  \node [draw=black!50, fit={(a) (a1) (a2) (a3)}] {};
  \node [draw=black!50, fit =(a4)] {};

  \node (b) at (5cm, 0) {transaction 2};

  \node (b1) [below=of b.west, typetag, xshift=2mm] { owner 2 public key };
  \node (b2) [below=of b1.west, typetag] { hash };
  \node (b3) [below=of b2.west, typetag] { owner 1 signature };
  \node (b4) [below=of b3] { owner 2 private key }; 

  \node [draw=black!50, fit={(b) (b1) (b2) (b3)}] {};
  \node [draw=black!50, fit =(b4)] {};
  
  \node (c) at (10cm, 0) {transaction 3};

  \node (c1) [below=of c.west, typetag, xshift=2mm] { owner 3 public key };
  \node (c2) [below=of c1.west, typetag] { hash };
  \node (c3) [below=of c2.west, typetag] { owner 2 signature };
  \node (c4) [below=of c3] { owner 3 private key }; 

  \node [draw=black!50, fit={(c) (c1) (c2) (c3)}] {};
  \node [draw=black!50, fit =(c4)] {};
  
  \draw [->] (a1) -- (a2) {}; 
  \draw [->] (a2) -- (a3) {};  
  \draw [->] (b1) -- (b2) {}; 
  \draw [->] (b2) -- (b3) {};  
  \draw [->] (c1) -- (c2) {}; 
  \draw [->] (c2) -- (c3) {}; 
  \draw [dashed,->] (a1) -- (b3) node[midway,sloped] { verify };
  \draw [dashed,->] (a4) -- (b3) node[midway,sloped] { sign };
  
  \draw [->] (a) -- (b2) {};
  \draw [->] (b) -- (c2) {};
  \draw [dashed,->] (b1) -- (c3) node[midway,sloped] { verify };
  \draw [dashed,->] (b4) -- (c3) node[midway,sloped] { sign };

\end{tikzpicture}
\caption{\label{fig:blockchain}Flowchart of blockchains\index{Cryptography} 
according to Satoshi Nakamoto\index{Nakamoto, Satoshi}. It shows three 
transactions in which signatures and verifications occur. }
\end{figure}

An important application of blockchains\index{Cryptography} are 
cryptocurrencies\index{Cryptography} such as bitcoins. The 
total volume of all bitcoins is capped and the exchange rate is quite 
volatile. Often the question is asked, what is the equivalent value of bitcoins 
and other cryptocurrencies\index{Cryptography}, as they ultimately only exist 
as data. The answer is like cash. Perhaps with a historical coin, the issue 
value still roughly corresponded to the metal value. At the latest, banknotes 
have broken with this tradition and the value of a carefully produced and 
difficult to counterfeit banknote is that it can be used to purchase other 
things without any problems. The value of money thus corresponds to the trust 
associated with it. The same is true for cryptocurrencies\index{Cryptography}. 
Whether and in what form they will establish themselves in the future is still 
unclear, especially, since the high energy consumption is extremely harmful to 
the environment.

\section{Factorisation and Quantum Computing\index{Quantum computing}}

Classic digital computers calculate with bits, which are elements of the set 
$\{0,1\}$. This idea ultimately goes back to Leibniz\index{Leibniz, Gottfried 
Wilhelm}. They calculate functions that map $m$-bit vectors to $n$-bit 
vectors:
\[
f \colon \{0,1\}^m \xlongrightarrow{~~~~~~} \{0,1\}^n. 
\]
Such functions can be practically arbitrary. However, all conceivable 
calculations\index{Computability} are realised using a few standard gates that 
mimic logical operations. The AND, OR and NOT gates form a generating set 
from which any other function $f$ can be generated. These gates correspond to 
the logical operations $\wedge$, $\vee$ and $\neg$. The NAND gate, i.e., the 
negated AND gate, already generates these alone.

Quantum computers\index{Quantum computing} work with qubits instead of bits and 
utilise the laws of quantum mechanics. Physical states in quantum mechanics are 
vectors in a Hilbert space of finite or infinite dimension. The temporal 
dynamics of quantum mechanical processes are described by the Schr\"odinger 
equation\index{Schr\"odinger, Erwin}. From the theory of quantum mechanics, it 
follows that the physically permissible operations $f$ that can occur between 
such Hilbert spaces are either projection mappings or unitary operators. 
Projections correspond to measurements on a quantum mechanical system.

A qubit is an element in the two-dimensional Hilbert space $H=\mathbb{C}^2$ 
with the standard scalar product. The basis of $H$ is often denoted in physics 
as
\[
\rvert 0 \rangle, \rvert 1 \rangle.
\]
A vector in $H$ is given by a linear combination
\[
\alpha \rvert 0 \rangle + \beta  \rvert 1 \rangle,
\]
where $\alpha$ and $\beta$ are complex numbers. When $n$ such qubits are 
considered together, the appropriate Hilbert space that describes the overall 
situation is the $n$-fold tensor product
\[
H^{\otimes n}=\underbrace{H \otimes H \otimes \cdots \otimes H}_{n \text{ 
times}}. 
\]
This complex vector space has the dimension $2^n$ and every vector in $H$ is of 
the form 
\[
u=\lambda_{00 \ldots 0} \rvert 00\ldots 0 \rangle + \cdots 
+ \lambda_{11 \ldots 1} \rvert 11\ldots 1 \rangle.
\]
This large dimension is the actual reason why quantum 
computers\index{Quantum computing} can be so powerful. A state $u$ is called 
entangled if it cannot be written as a tensor product $v \otimes w$. The 
simplest case of this kind occurs at $n=2$, because 
\[
u=\rvert 01\rangle + \rvert 10 \rangle  
\]
is not an unentangled tensor product of the form 
\[
\left( \alpha \rvert 0 \rangle + \beta  \rvert 1 \rangle \right) \otimes 
\left( \gamma \rvert 0 \rangle + \delta  \rvert 1 \rangle \right).
\]
This follows from the fact that the coefficients $\lambda_{ij}$ of the 
unentangled vectors for $n=2$ satisfy the cone equation  
\[
\lambda_{00} \cdot \lambda_{11} - \lambda_{01} \cdot \lambda_{10}=0
\]
but $u$ does not. There are suitable physical experiments that produce 
entangled particles. For example, entangled pairs of photons can be generated 
via polarisation, which even at a great distance from each other retain this 
property at times. This leads to interesting transmission channels and 
cryptographic protocols\index{Cryptography} in quantum information theory. 

Quantum gates are the necessary gates used in the context of 
quantum computers\index{Quantum computing}. Unlike with 
digital computers, gates are necessary that induce unitary operations. 
Examples of such gates are the Hadamard gate H and the 
CNOT gate, i.e., the controlled NOT gate:
\[
\text{H}=\frac{1}{\sqrt{2}} \left( \begin{matrix} 1 & 1 \cr 1& -1 \end{matrix} 
\right), \quad  \text{CNOT}=\left( \begin{matrix} 1 & 0 & 0 & 0 \cr 
0 & 1 & 0 & 0 \cr 0 & 0 & 0 & 1 \cr 0 & 0 & 1 & 0  \end{matrix} \right).
\]
In 1994, Peter Shor\index{Shor, Peter} found the Shor 
algorithm\index{Algorithm} for quantum computers\index{Quantum computing}, 
which is capable of factorising natural numbers in a probabilistic way in 
polynomial time\index{Polynomial time}, provided a quantum 
computer\index{Quantum computing} can calculate largely error-free with a 
sufficiently large number of qubits. Shor\index{Shor, Peter} solved this task 
by implementing the so-called discrete Fourier transformation as a quantum 
algorithm and reading out the information, i.e., the coefficients in the base 
representation in $H^{\otimes n}$, in the final state of the system through a 
clever mix of quantum mechanical measurement and postprocessing using 
elementary number theory.\footnote{See \cite{shor,wigderson2019}.} 

\section{Mathematical Physics}

In his habilitation thesis, Bernhard Riemann\index{Riemann, Bernhard} 
developed the mathematical theory of manifolds together with their metric 
and differential geometric structure. This theory has particularly contributed 
to the discovery of Einstein's field equations\index{Einstein, Albert} in 
general relativity theory. 

In mathematical physics, the concept of space-time is usually modelled on 
Riemannian\index{Riemann, Bernhard} manifolds, where in cosmology singularities 
must also be allowed, which explain phenomena such as black holes. It is common 
to describe a manifold $M$ and its metric-differential geometric structure 
locally using coordinates. However, global properties play a more decisive 
role. On $M$ and its tensor fields often operates a group of gauge 
transformations as a symmetry group, which preserves physical quantities of 
the given theory. The physics to be described with such a differential 
geometric model does not depend on the specifically chosen coordinates and the 
gauge transformations. This property is called covariance\index{Covariance} and 
is a good example that the mathematical description of physical objects is only 
relevant up to isomorphism\index{Isomorphism} or 
equivalence\index{Equivalence}.\footnote{See \cite{hawking,penrose2006} 
for covariance and space-time.}

The idea of covariance\index{Covariance} goes in parts back to 
Leibniz\index{Leibniz, Gottfried Wilhelm}, who in his correspondence with 
Samuel Clarke\index{Clarke, Samuel} from the years 1715-1716 discussed the 
structure of physical space. Leibniz\index{Leibniz, Gottfried Wilhelm} wrote in 
a letter dated February 25th 1716:
\begin{quote}
As for my own opinion, I have said more than once that I hold space 
to be something purely relative, as time is -- that I hold it to be an order of 
coexistences, as time is an order of successions ... I have many demonstrations 
to confute the fancy of those who take space to be a substance or at least an 
absolute being.\footnote{Printed in \cite[p. 14]{clarke2000}.}
\end{quote}
In the same letter, Leibniz\index{Leibniz, Gottfried Wilhelm} considered the 
possibility of swapping the cardinal directions East and West in space and then 
wrote further in loc. cit.:
\begin{quote}
But if space is nothing else but this order or relation, and is nothing 
at all without bodies but the possibility of placing them, then those two 
states, the one such as it is now, the other supposed to be the quite 
contrary way, would not at all differ from one another. 
\end{quote}
Leibniz\index{Leibniz, Gottfried Wilhelm} thus linked his concept of space with 
the concept of equality\index{Equality}\index{Identity} and 
isomorphism\index{Isomorphism}. Albert Einstein\index{Einstein, Albert} 
formulated the principle of covariance\index{Covariance} as follows:
\begin{quote}
The general laws of nature are to be expressed by equations that are valid 
for all coordinate systems, i.e., they are covariant with respect to arbitrary 
substitutions (generally covariant).\footnote{Author's translation of \cite[p. 
776]{einstein1916}.}
\end{quote}
Einstein's\index{Einstein, Albert} field equations describe the 
general theory of relativity. The decisive matrix-valued equation 
among them describes the relationship of the curvature tensors with the 
mass ratios in space and reads:\footnote{See \cite[Ch. 19]{penrose2006}.}
\[
R_{\mu \nu}-\frac{1}{2} g_{\mu \nu} R + \Lambda g_{\mu \nu}= 
\frac{8 \pi G}{c^4} T_{\mu \nu}.
\]
It is noteworthy that the ideas of Riemann\index{Riemann, Bernhard} together 
with Einstein's\index{Einstein, Albert} theory of relativity found an 
application in GPS navigation only after more than 150 years. The location 
determination using the radio signals of several satellites would be 
significantly less precise if the special and general theory of relativity were 
not taken into account at the same time.

Some remarks from Riemann's\index{Riemann, Bernhard} habilitation thesis are 
also worth mentioning. In the last section he wrote: 
\begin{quote}
The question about the validity of the assumptions of geometry in 
the infinitesimal is linked to the question of the inner reason for the 
mass ratios in space. With this question, which may still be counted as part of 
the doctrine of space, the above remark applies that in a discrete manifold the 
principle of mass ratios is already contained in the concept of this manifold, 
but in a continuous one it must come from elsewhere. It must therefore either 
be the case that the real underlying space forms a discrete manifold, or the 
reason for the mass ratios must be sought outside, in binding forces acting on 
it.\footnote{Machine translation of \cite{riemann}.} 
\end{quote}
These remarkable sentences show deep insights and question the continuous 
space-time models that prevail to this day. Roger Penrose\index{Penrose, Roger} 
has developed spin networks as a new idea for discrete space-time 
models from around 1971. Spin networks in dimension $d$ are based on a 
decorated graph whose edges represent $(d-1)$-dimensional wall structures (see 
Fig.~\ref{fig:spin_network}). The decoration of vertices and edges encodes in a 
discrete way representation-theoretic data. When a spin network evolves in 
discrete time-steps, the graph becomes a two-dimensional simplicial complex 
which is called spin foam. Another approach is non-commutative geometry, where 
$C^\ast$-algebras take the role of local functions and no underlying 
topological space\index{Topology} $X$ is available. Both concepts could prove 
useful for the unification of quantum field theory with the general theory of 
relativity into quantum gravity.\footnote{See \cite{connes1995,hawking} and 
\cite[Ch. 32]{penrose2006}.} 

\begin{figure}[ht!]
\begin{tikzpicture}[line join = round, line cap = round]
\pgfmathsetmacro{\factor}{-0.5};
\coordinate (A) at (2,0,-2*\factor);
\coordinate (B) at (-3,0,-2*\factor);
\coordinate (C) at (0,2,2*\factor);
\coordinate (D) at (0,-2,2*\factor);
\node [black] at (-0.15,-0.9) {\textbullet}; 
\node [black] at (0.8,0.22) {\textbullet}; 
\node [black] at (-0.21,0.6) {\textbullet}; 
\node [black] at (-0.7,0.27) {\textbullet}; 
\draw[-,dashed] (-0.1,0.2) -- (0.5,-2.5,2) node[left] {$e_1$}; 
\draw[-,dashed] (-0.1,0.2) -- (4,1.5,3) node[right] {$e_2$};  
\draw[-,dashed] (-0.1,0.2) -- (0,3.5,2) node[above] {$e_3$}; 
\draw[-,dashed] (-0.1,0.2) -- (-1,3,6) node[left] {$e_4$};
\draw[-] (A)--(D)--(B)--cycle;
\draw[-] (A) --(D)--(C)--cycle;
\draw[-] (B)--(D)--(C)--cycle;
\end{tikzpicture}
\caption{\label{fig:spin_network}A typical piece of a $3$-dimensional spin 
network is given by a $4$-valent decorated graph. The four edges at each vertex 
and their wall structures form a tetrahedron. Further decorations are omitted 
in the picture.}
\end{figure}

Similar to Einstein's\index{Einstein, Albert} field equation in general 
relativity, physical laws can often be described with one or more mathematical 
equations. In quantum mechanics, the Schr\"odinger 
equation\index{Schr\"odinger, Erwin}  
\[
 i \hbar \frac{\partial}{\partial t} \Psi(t) = H \, \Psi(t)
\]
plays an important role for the wave function $\Psi(t)$. Another example is the 
heat conduction equation 
\[
\frac{\partial}{\partial t} u(x,t) -a \Delta u(x,t)=f(x,t). 
\]
It describes the propagation of temperature $u$ in three-dimensional 
space (with the $x$-coordinate) depending on time $t$ given a heat source 
$f(x,t)$ and conductivity index $a>0$. Another example is the system of three 
Euler partial differential equations, which describe the flow of a fluid and 
their extension to the Navier-Stokes equations. The Clay Foundation in Boston 
named seven Millennium problems\index{Millennium problems} in 2000, the 
solutions of which pose challenges in mathematics and six of which have not 
been solved to this day. They are a reminiscence of the Hilbert 
problems\index{Hilbert, David} from 1900. One of these problems asks about the 
existence\index{Existence} and smoothness of solutions to the Navier-Stokes 
equations.\footnote{See www.claymath.org.}

\section{Artificial Neural Networks\index{Neural network} and Deep 
Learning\index{Artificial intelligence}}

Artificial intelligence\index{Artificial intelligence} is an umbrella term for 
many types of algorithmically\index{Algorithm} controlled systems that simulate 
human thinking or can support and replace humans. The range is broad and 
includes suitable algorithms\index{Algorithm} as well as machines, such as 
robots. After the beginnings of artificial intelligence\index{Artificial 
intelligence} in the times of Alan Turing\index{Turing, Alan} and Norbert
Wiener\index{Wiener, Norbert}, the discovery of Hebb's rule\index{Hebb's rule} 
by Donald Hebb \index{Hebb, Donald} in 1949 and the invention of the 
perceptron\index{Perceptron} by Frank Rosenblatt\index{Rosenblatt, Frank} in 
the late 1950s, numerous attempts over several decades have been made. During 
this time, statistical learning algorithms\index{Algorithm} such as the PAC 
learning algorithm\index{Algorithm} by Leslie Valiant\index{Valiant, Leslie}, 
optimisation algorithms\index{Algorithm} of various kinds, genetic and 
evolutionary algorithms\index{Algorithm} and many others have been designed. 
Since around 2006, enormous breakthroughs have been achieved in a field called 
machine learning\index{Artificial intelligence}. A 
specific form of this includes deep learning\index{Artificial intelligence}, 
which involves training artificial neural networks\index{Neural 
network}.\footnote{See \cite{hebb1949,rosenblatt1958} for Hebb's rule and the 
perceptron, \cite{goodfellow2016} and \cite{kutyniok2023} for deep learning, 
\cite{valiant2013} for learning algorithms and \cite{wigderson2019} for the 
theoretical computer science behind all this.}

\begin{figure}[ht!]
\begin{tikzpicture}[x=1.3cm, y=1.3cm, >=stealth]

\tikzset{%
  every neuron/.style={
    circle,
    draw,
    minimum size=0.4cm
  },
  neuron missing/.style={
    draw=none, 
    scale=4,
    text height=0.333cm,
    execute at begin node=\color{black}$\vdots$
  },
}

\foreach \m/\l [count=\y] in {1,2,3,missing,4} {
\ifnum \y=4 \node [every neuron/.try, neuron \m/.try] (input-\m) at 
(0,2.5-\y) {};
\else \node [every neuron/.try, neuron \m/.try,fill=lightgray] (input-\m) at 
(0,2.5-\y) {};
\fi  }
  
\foreach \m [count=\y] in {1,missing,2}
  \node [every neuron/.try, neuron \m/.try] (hidden-\m) at 
(2,2-\y*1.25) {};

\foreach \m [count=\y] in {1,missing,2}  
  \node [every neuron/.try, neuron \m/.try] (output-\m) at (4,1.5-\y) {}; 

\foreach \l [count=\i] in {1,2,3,100}
  \draw [<-] (input-\i) -- ++(-1,0)
    node [above, midway] {$E_{\l}$};

\foreach \l [count=\i] in {1,N} {
  \ifnum \i=1  \node [above] at (hidden-\i.north) {$V_1$};
  \else \node [below] at (hidden-\i.south) {$V_\l$};
  \fi}
  
\foreach \l [count=\i] in {1,4} 
    \draw [->] (output-\i) --  (6,-0.5)  
    node [every neuron/.try,fill=lightgray] at (6.15,-0.5) {}; 

\foreach \i in {1,...,4}
  \foreach \j in {1,...,2}
    \draw [->] (input-\i) -- (hidden-\j);
       
\draw [->] (hidden-1) -- (output-1) node [above=0.2] {$V'_1$};
\draw [->] (hidden-2) -- (output-2) node [below=0.2] {$V'_N$};
\draw [->] (hidden-1) -- (output-2);
\draw [->] (hidden-2) -- (output-1);
\draw [->] (6.3,-0.5) -- (7.2,-0.5) node [above,midway] {$A$} ;
\end{tikzpicture}
\caption{\label{fig:neural_net}Artificial neural network\index{Neural network} 
with $100$ input nodes, two hidden layers, and one output node.}
\end{figure}

Graph-like functional networks also occur in nature. Already in the 19th 
century, parts of the brain could be dissected so precisely that connected 
structures of neurons (nerve cells) and glial cells were discovered, with the 
neurons exchanging via synapses and being able to exhibit electrical 
activation. Around 1943, a simplified mathematical model of this system in the 
form of an artificial neural network\index{Neural network} was designed by 
Warren McCulloch\index{McCulloch, Warren} and Walter Pitts\index{Pitts, 
Walter}.\footnote{See \cite{cullochpitts}.} 

A simple form of an artificial neural network\index{Neural network} is the 
multilayer perceptron\index{Perceptron}. It consists of a finite 
directed graph, which includes input nodes $E_\bullet$ and output nodes 
$A_\bullet$ as well as nodes $V_\bullet$ from several hidden layers with 
their connecting edges. The number of hidden layers is referred to as the depth 
of the network and the number $N$ its width. Fig.~\ref{fig:neural_net} suggests 
a blueprint for a network, with which from a digital square image with $10 
\times 10$ pixels with grey values between $0$ (white) and $1$ (black) the 
letter Y is to be recognised with the help of a single output node.

\begin{figure}[ht!]
\begin{tikzpicture}
\draw[step=0.5cm,color=gray] (0,0) grid (5,5);
\node[fill=black,inner sep=7pt] at (2.75,+0.75) {};
\node[fill=black,inner sep=7pt,minimum size=0.1] at (2.25,+1.25) {};
\node[fill=gray,inner sep=7pt,minimum size=0.1] at (2.25,+1.75) {};
\node[fill=black,inner sep=7pt,minimum size=0.1] at (2.25,+2.25) {};
\node[fill=black,inner sep=7pt,minimum size=0.1] at (2.25,+2.75) {};
\node[fill=black,inner sep=7pt,minimum size=0.1] at (2.25,+0.25) {};
\node[fill=gray,inner sep=7pt,minimum size=0.1] at (0.75,4.75) {}; 
\node[fill=black,inner sep=7pt,minimum size=0.1] at (3.25,4.75) {}; 
\node[fill=black,inner sep=7pt,minimum size=0.1] at (2.75,3.25) {}; 
\node[fill=gray,inner sep=7pt,minimum size=0.1] at (+1.25,3.75) {};
\node[fill=black,inner sep=7pt,minimum size=0.1] at (+0.75,3.75) {};
\node[fill=black,inner sep=7pt,minimum size=0.1] at (+0.75,4.25) {};
\node[fill=black,inner sep=7pt,minimum size=0.1] at (+1.75,3.25) {};
\node[fill=gray,inner sep=7pt,minimum size=0.1] at (+3.25,3.75) {};
\node[fill=black,inner sep=7pt,minimum size=0.1] at (+3.25,4.25) {};
\foreach \x in {0.25,2.25,3.75,4.25} 
\node[fill=lightgray,inner sep=7pt,minimum size=0.1] at (\x,4.75) {};
\foreach \x in {1.75,4.75} 
\node[fill=lightgray,inner sep=7pt,minimum size=0.1] at (\x,4.25) {};
\foreach \x in {2.25,4.25} 
\node[fill=lightgray,inner sep=7pt,minimum size=0.1] at (\x,3.75) {};
\foreach \x in {0.75,1.25,2.25,4.25} 
\node[fill=lightgray,inner sep=7pt,minimum size=0.1] at (\x,3.25) {};
\foreach \x in {0.25,1.25,2.75,3.75} 
\node[fill=lightgray,inner sep=7pt,minimum size=0.1] at (\x,2.75) {};
\foreach \x in {0.75,2.75,3.75,4.75} 
\node[fill=gray,inner sep=7pt,minimum size=0.1] at (\x,2.25) {};
\foreach \x in {1.75,3.25,4.25} 
\node[fill=lightgray,inner sep=7pt,minimum size=0.1] at (\x,1.75) {};
\foreach \x in {0.25,1.25,3.75} 
\node[fill=lightgray,inner sep=7pt,minimum size=0.1] at (\x,1.25) {};
\foreach \x in {0.75,1.75,3.25,4.75} 
\node[fill=lightgray,inner sep=7pt,minimum size=0.1] at (\x,0.75) {};
\foreach \x in {0.75,2.75,3.25,4.25} 
\node[fill=lightgray,inner sep=7pt,minimum size=0.1] at (\x,0.25) {};
\end{tikzpicture}
\caption{\label{fig:training_image}Array with $10 \times 10$ grey pixels which 
resembles a training image for the letter $Y$.}
\end{figure}

How is it calculated? The depicted neural network\index{Neural network} is 
intended to compute a real function $f\colon \mathbb{R}^{100} \longrightarrow 
\mathbb{R}$, whose arguments are the grey values $e_i$ of the $100$ pixels of a 
digital image (See Fig.~\ref{fig:training_image}). If the letter Y is 
recognised, then the network should approximate the function value 
$f(e_1,\ldots,e_{100}) \approx 1$ and otherwise yield $f(e_1,\ldots,e_{100}) 
\approx 0$.

Important for the proper functioning of a neural network\index{Neural network} 
are the transitions between the individual layers, also called feedforward. If 
we focus on the transition between the input nodes $E_\bullet$ and a node $V_i$ 
of the first hidden layer, this is realised for each $i$ between $1$ and $N$ by 
formulas
\[
v_i=\sigma(W_i+b_i), 
\]
where $W_i$ are suitable functions of the input variables $e_j$, $b_i$ are 
rational constants -- called the bias or negative threshold -- and $\sigma$ is 
a non-linear activation function, for example of the form 
\[
\sigma(x)=\max(x,0)=
\begin{cases}
0 & \text{ for } x \le 0 \cr 
x & \text{ for } x \ge 0.
\end{cases} 
\]
In our example, $W_i$ can be a simple linear function 
\[
W_i(e_1,\ldots,e_{100})=w_{i,1} e_1 + \cdots +w_{i,100} e_{100}, 
\]
where the coefficients $w_{i,j}$ are positive or negative rational numbers 
called weights. There are other choices for $W_i$ and $\sigma$, which may be 
advantageous in some situations. Analogous formulas apply for the transitions 
in the hidden layers and in the final step to the output nodes. 

How does the training of such a neural network\index{Neural network} work? At 
the beginning of the training, the weights $w_{i,j}$ and the 
(negative) threshold values $b_i$ are initially set. The network is then 
trained by input data, by calculating for each input image a 
cost function $K$ depending on the weights and threshold values. In the 
simplest case, $K$ measures the sum of the squares of the deviations 
from $f(e_1,\ldots,e_{100})$ from the desired value $0$ or $1$, which -- 
in our example -- is represented by the actual letter on the image. By 
successive adjustment of the weights $w_{i,j}$ and the 
threshold values $b_i$, the cost function $K$ is slowly minimised. This 
method is called backpropagation or backtracking. Mathematically 
speaking, this process consists of a gradient descent, similar to the way 
of descending from a mountain as quickly as possible against the direction of 
the steepest slope. The goal is a local minimum of the cost function. To avoid 
nonsensical or unrealistic local minima of the cost function in this method, 
which do not come close to the solution of the problem, the gradient descent is 
usually combined with a stochastic component, which is randomly combined from 
several choices. Related to this is the method of genetic 
algorithms\index{Algorithm}, with which unwanted local minima are avoided by 
evolutionary jumps in the parameters.

By inputting large amounts of data, the neural network\index{Neural network} 
can be adjusted over a long period of time until the weights and threshold 
values are altered so that it can recognise desired structural patterns from 
any new data with good accuracy. This makes the network\index{Neural network} a 
good predictor in appropriate new situations. Essentially, the method of deep 
learning\index{Artificial intelligence} involves reducing a high-dimensional 
problem to a few essential parameters, with the correct reduction found by 
training through approximation. 

In addition to the simple case of the multilayer perceptron\index{Perceptron}, 
there are more powerful recurrent neural networks\index{Neural network}, where
certain directed edges transmit information backwards, thus enabling 
recursive\index{Recursion} calculations of sequential functions $f(n)$, which 
depend on a discrete parameter $n$. Other variants use sophisticated filters in 
the hidden layers. Such networks are called convolutional neural 
networks\index{Neural network}. They are particularly effective in image 
processing. An even more radical generalisation of neural 
networks\index{Neural network} uses architectures which include 
higher-dimensional simplices\index{Simplex} and which take this  
topological\index{Topology} structure into account in information processing. 
It is conjectured that such architectures, which run under the label 
\enquote{topological deep learning}, are more efficient and reliable in the 
detection of topological\index{Topology} structures in data. Together with the 
enormous computing power of today's processors, such sophisticated refinements 
of the multilayer perceptron\index{Perceptron} contribute to the current and 
upcoming success of artificial intelligence\index{Artificial intelligence}.

\section{Pitfalls of Artificial Intelligence\index{Artificial intelligence}}

A priori, there are many systematic errors in neural 
networks\index{Neural network} and other learning algorithms\index{Algorithm}. 
Traditional results of analysis usually cannot accurately show the quality of 
the approximation, as the prerequisites of usual mathematical theorems, such as 
convexity, linearity or smoothness, are hardly ever met in the applications. 
Therefore, it is surprising that this method is so successful in many cases. 
The depth of the network, more than its width, plays a crucial role. At the 
moment, a new field of mathematical analysis of deep learning\index{Artificial 
intelligence} is emerging. On the one hand, it provides estimates of various 
types of systematic approximation errors and, on the other hand, it 
investigates fundamental problems and sources of error.\footnote{See 
\cite{kutyniok2023} for the mathematics of deep learning.}

It is noteworthy that deep learning\index{Artificial intelligence}, unlike 
other learning algorithms, is not significantly affected in its function by the 
phenomena of overfitting and overparameterisation. In the case of overfitting, 
the algorithm\index{Algorithm} adapts optimally to the training data due to the 
numerous adjustable parameters, so that it can make perfect decisions in the 
given cases, but may not generalise well in many other situations and thus 
possibly produce incorrect results. The related overparameterisation is caused 
by the fact that many more adjustable parameters, such as the weights in neural 
networks\index{Neural network}, are present in the algorithm than the dimension 
of the essential parameters of the structures being examined.

Learning algorithms\index{Algorithm} may produce solutions that are incorrect 
or have an undesirable bias, for example if the training data have a bias. A 
fundamental problem is that there are undecidable 
problems\index{Undecidability} for deep learning\index{Artificial 
intelligence}, 
just as there are for ordinary Turing machines\index{Turing 
machine}.\footnote{See \cite{ben-david2019}.} A danger comes from incorrect 
applications that arise from ignorance of the pitfalls of the methods. 
Questioning such algorithms\index{Algorithm} is therefore justified. For these 
reasons, it is necessary to continue developing the theory quickly. The 
exploration of learning is a more general goal than the method of deep 
learning\index{Artificial intelligence} and includes other algorithms already 
discovered\index{Algorithm}. The mathematical correctness of all these methods 
will one day presumably be verifiable and the use therefore safer and ethically 
less questionable than today.

Is artificial intelligence\index{Artificial intelligence} really a form 
of intelligence? This is a difficult question. Artificial 
intelligence\index{Artificial intelligence} in the form of deep 
learning\index{Artificial intelligence} is a method that can draw conclusions 
from complex problems with large amounts of data. This involves a reduction of 
complexity and an understanding in the form of predictors is achieved. Human 
intelligence is also capable of navigating complex situations and making good 
decisions. However, our intelligence seems to possess more properties than just 
the ability to reduce complexities in particularly complicated situations. The 
discovery of insights in the sciences, creative artistic work, or even 
understanding a subtle joke, are achievements of human intelligence that do not 
seem to be explainable. The last word has not yet been spoken. So there is room 
for further research.

\section{Topological Data Analysis\index{Topological data analysis} in 
Bioinformatics}

The life sciences are increasingly applying methods from mathematics and 
computer science. The amount of data that can be generated from gene sequences 
of organisms, from molecular data, from image data in high-resolution imaging 
techniques, from configurations of neural networks\index{Neural network} in 
model organisms or from other methods is enormous and can only be managed with 
sophisticated methods. We want to address the problem of data 
analysis\index{Topological data analysis}, which can be approached with 
mathematical methods from algebraic topology.

Mathematically, the aim is to recognise a geometric pattern from a large set of 
discrete data points in a high-dimensional space. Now, topological 
spaces\index{Topology} and their invariants\index{Invariant} are precisely the 
objects that reflect typical geometric patterns. In applications, large 
collections of data points with a given distance concept, called metric, are 
relevant, which can have different manifestations. For genetic data, this could 
be the Hamming distance between vectors with combinatorial data of genetic 
sequences, which by definition is the number of all 
non-matching components. In neural networks\index{Neural network}, the 
reciprocal of the strength of the synaptic connections provides a metric.

There is a novel method to gain topological information from data point clouds. 
It is called persistent homology or topological data analysis\index{Topological 
data analysis}. For this, the metric structure on the space of data is used to 
define a chain complex. There are two popular variants, firstly the 
Vietoris-Rips complex and secondly the \v{C}ech complex.

We first consider the Vietoris-Rips complex. The points $P$ of the data cloud 
form the $0$-simplices\index{Simplex}. Around each of these points $P$ we 
place closed metric spheres $B_{P,\varepsilon}$ with radius $\varepsilon>0$. 
Two different $0$-simplices\index{Simplex} $P$ and $Q$ are connected by an 
edge if their distance is less than or equal to $\varepsilon$, i.e., if $Q$ is 
contained in $B_{P,\varepsilon}$. Given pairwise different 
$0$-simplices\index{Simplex} $P_0, \ldots, P_k$, we define a 
$k$-simplex\index{Simplex} if all points have a distance less than or equal to 
$\varepsilon$ from each other, i.e., if each point is in the closed 
$\varepsilon$-sphere of each other. In the case of the \v{C}ech complex, we 
define a $k$-simplex\index{Simplex} for $k+1$ pairwise different points 
$P_0,\ldots,P_k$, if the intersection
\[
B_{P_0,\varepsilon} \cap \cdots \cap B_{P_k,\varepsilon} \neq \emptyset
\]
is not the empty set. With both methods, we get a simplicial 
set\index{Simplicial set} $V(\varepsilon)_\bullet$ for each $\varepsilon$. The 
associated chain complex $\mathbb{Z}V(\varepsilon)_\bullet$ of abelian groups 
is called the Vietoris-Rips complex or \v{C}ech complex.

Persistent homology arises when we calculate the homology\index{Homology group} 
of the chain complexes $\mathbb{Z} V(\varepsilon)_\bullet$ and investigate 
which homology groups\index{Homology group} and Betti numbers these chain 
complexes have when $\varepsilon$ varies. These Betti numbers 
$b_n=b_n(\varepsilon)$ can be calculated using methods of linear 
algebra and can be visually represented in the form of barcodes 
depending on $\varepsilon$. The persistent homology is then the result of the 
barcodes that occur over longer distances (see Fig.~\ref{fig:barcode}). 
Mathematically speaking, the homology\index{Homology group} of filtered chain 
complexes, which are filtered over the parameter $\varepsilon$, is calculated 
in an efficient manner. 

\begin{figure}[ht!]
\begin{tikzpicture}[x=0.7cm,y=0.35cm]
 
\draw[dashed,->]  (-11,6)--(4,6) node[anchor=west] {$\varepsilon$};
\draw[-]  (-10,-1)--(-7,-1) node[anchor=west] {$b_0$};
\draw[-]  (-10,-2)--(-7,-2) node[anchor=west] {$b_0$};
\draw[-]  (-10,-3)--(-7,-3) node[anchor=west] {$b_0$};
\draw[-]  (-10,-4)--(-7,-4) node[anchor=west] {$b_0$};
\draw[-]  (-10,-5)--(-7,-5) node[anchor=west] {$b_0$};
\draw[-]  (-10,-6)--(-7,-6) node[anchor=west] {$b_0$};
\draw[-]  (-10,-7)--(-7,-7) node[anchor=west] {$b_0$};
\draw[-]  (-10,-8)--(-4,-8) node[anchor=west] {$b_0$};
\draw[-]  (-10,-9)--(2,-9) node[anchor=west] {$b_0$};
\draw[-]  (-4,-10)--(2,-10) node[anchor=west] {$b_1$};
    
\draw (-11,1.5) node {.};
\draw (-11,3) node {.};
\draw (-10.5,4) node {.};
\draw (-10,1) node {.};
\draw (-10,5) node {.};
\draw (-9.5,1.5) node {.};
\draw (-9,2.5) node {.};
\draw (-9,4) node {.};
    
\draw (-8,1.5) node {.};
\draw (-8,3) node {.};
\draw (-7.5,4) node {.};
\draw (-7,1) node {.};
\draw (-7,5) node {.};
\draw (-6.5,1.5) node {.};
\draw (-6,2.5) node {.};
\draw (-6,4) node {.};
 
\draw[-]  (-8,1.5)--(-8,3);
\draw[-]  (-8,3)--(-7,5);

\draw[-]  (-7,1)--(-6.5,1.5);
\draw[-]  (-6.5,1.5)--(-6,2.5);
\draw[-]  (-6,2.5)--(-6,4);

\draw (-5,1.5) node {.};
\draw (-5,3) node {.};
\draw (-4.5,4) node {.};
\draw (-4,1) node {.};
\draw (-4,5) node {.};
\draw (-3.5,1.5) node {.};
\draw (-3,2.5) node {.};
\draw (-3,4) node {.};

\draw[-]  (-5,1.5)--(-5,3);
\draw[-]  (-4,1)--(-5,1.5);

\draw[-]  (-5,1.5)--(-4,1);
\draw[-]  (-4,5)--(-3,4);
\draw[-]  (-5,3)--(-4,5);

\draw[-]  (-4,1)--(-3.5,1.5);
\draw[-]  (-3.5,1.5)--(-3,2.5);
\draw[-]  (-3,2.5)--(-3,4);
    
\draw (-2,1.5) node {.};
\draw (-2,3) node {.};
\draw (-1.5,4) node {.};
\draw (-1,1) node {.};
\draw (-1,5) node {.};
\draw (-0.5,1.5) node {.};
\draw (-0,2.5) node {.};
\draw (0,4) node {.};

\draw[-]  (-2,1.5)--(-2,3);
\draw[-]  (-1,1)--(-2,1.5);
\draw[-]  (-1,5)--(0,4);
\draw[-]  (-2,3)--(-1,5);

\draw[fill=darkgray!10] (0,4) -- (-0.5,1.5) -- (0,2.5) -- cycle; 

\draw[-]  (-1,1)--(-0.5,1.5);
\draw[-]  (-0.5,1.5)--(0,2.5);
\draw[-]  (0,2.5)--(0,4);
    
\draw (1,1.5) node {.};
\draw (1,3) node {.};
\draw (1.5,4) node {.};
\draw (2,1) node {.};
\draw (2,5) node {.};
\draw (2.5,1.5) node {.};
\draw (3,2.5) node {.};
\draw (3,4) node {.};

\draw[-]  (1,1.5)--(1,3);
\draw[-]  (1,1.5)--(2,1);
\draw[-]  (2,5)--(3,4);
\draw[-]  (1,3)--(2,1);
\draw[-]  (1,3)--(2,5);

\draw[fill=darkgray!10] (3,2.5) -- (2.5,1.5) -- (2,1) -- cycle; 
\draw[fill=darkgray!10] (1,3) -- (1,1.5) -- (2,1) -- cycle; 
\draw[fill=darkgray!10] (3,4) -- (2.5,1.5) -- (3,2.5) -- cycle; 

\draw[-]  (2,1)--(2.5,1.5);
\draw[-]  (2.5,1.5)--(3,2.5);
\draw[-]  (3,2.5)--(3,4);
    
\end{tikzpicture}
\caption{\label{fig:barcode}Barcodes of persistent simplices\index{Simplex} 
with increasing $\varepsilon$ for a circular point cloud. From the third image 
on, the Betti numbers $b_0$ and $b_1$ stabilise to $1$.}
\end{figure}

Topological data analysis has shown remarkable applications, particularly in 
the life sciences. However, it remains to be seen how successful 
it will be in competition with artificial intelligence\index{Artificial 
intelligence}, in particular with \enquote{topological deep learning}.

\section{Opportunities for Mathematics}

In the past, the sciences mostly worked in a monodisciplinary manner. Modern 
questions are often interdisciplinary and have complex systems as their 
mathematical basis. A good example of this is climate research and the Earth's 
ecosystem as a whole. This requires studying the physical and chemical 
foundations of meteorology and geosciences, as well as the influence of humans, 
which are interdependent. In the various layers of the Earth, physical 
interfaces spatially meet. Even if the dynamics occurring in each layer are 
understood, the spatial boundary transition must be modelled. This requires 
spatial multiscale methods and coarsening tricks as well as a high 
computational effort. In climate research, temporal multiscale methods 
are necessary to handle the dynamics on the time axis. While the entire climate 
operates on a large time scale, there are weather effects and human behaviour, 
both of which occur on a narrower time scale. For the mathematical explanation 
of the connection of these two scales, Klaus Hasselmann\index{Hasselmann, 
Klaus} was awarded the Nobel Prize in 2021. His model demonstrated that the 
climate generally has a great stability, but in recent decades human-made 
effects outside this normal range of fluctuation have been added.

Another future task is the development of mathematical consensus protocols. Due 
to their immense computational effort and energy consumption, 
blockchains\index{Cryptography} are no longer contemporary. Novel and more 
efficient data networks in the style of distributed ledgers are required, which 
are not monitored by any central authority. They are a special case of 
distributed systems, where many participants exchange data and execute 
algorithms without a central instance having a control function. Financial 
transactions in the areas of mobility, energy and health are good examples of 
this. The requirements for transparency, security and data protection are 
enormous. Complex and distributed systems are just two significant examples 
among others. Mathematics will play a crucial role in all future challenges, as 
it can design universally applicable models.\footnote{See \cite{jost2017}.} 

\chapter{Computability\index{Computability} and 
Decidability\index{Decidability}}

The idea of recursion\index{Recursion} and complete 
induction\index{Complete induction} is very old. But it was not until the 19th 
century that Richard Dedekind\index{Dedekind, Richard} with his 
books\footnote{Contained in \cite{dedekind}, see \cite{ms2023}.} 
\enquote{Continuity and irrational numbers} and \enquote{What are and what 
should the numbers be?} established the fundamental properties of natural 
numbers on a solid foundation and defined the construction of the number system 
with integers, rational and real numbers from natural numbers. One of his most 
important achievements was the discovery and proof of the recursion 
theorem\index{Recursion}. 

Building on this, a theory of computability\index{Computability} only 
crystallised around 1936. At this time, almost simultaneously, all theoretical 
computability models\index{Computability} were developed that we know today. 
They turned out to be equivalent and describe in an abstract fashion the 
strength of the later developed digital computers. In particular, the functions 
of $\mathbb{N}^n$ to $\mathbb{N}$ that can be computed with a Turing 
machine\index{Turing machine} are exactly all partially defined 
recursive\index{Recursion} functions that can be constructed from elementary 
functions using primitive recursion\index{Recursion} and the $\mu$-operator. 
From these contexts, the original Church-Turing thesis\index{Church-Turing 
thesis} was derived, which informally states that every effectively 
computable\index{Computability} function is recursive\index{Recursion} 
according to this definition. The question remains open whether there are 
stronger computability models\index{Computability} that can be physically 
realised. To date, such expectations, referred to as 
hypercomputing\index{Hypercomputing}, have not been realised, despite 
promising approaches for computer architectures in the field of 
neuromorphic computing.\footnote{See \cite{siegelmann1999}.} 

There are mathematical problems that cannot be solved algorithmically at all. 
These include in particular the halting problem\index{Halting problem} for 
Turing machines\index{Turing machine}, the word problem\index{Word problem} and 
the problem of solvability of integer diophantine equations. Such 
undecidable\index{Undecidability} problems show us the limits of 
computability\index{Computability}.  

\section{The Method of Recursion\index{Recursion}} 

Dedekind's recursion theorem\index{Recursion} states that a mapping $f$ from 
the natural numbers to another set can be uniquely constructed by specifying 
$f(0)$ and providing a rule for how the value of $f(n)$ can be obtained from 
the values $f(m)$ for $m \le n-1$ for each $n \ge 1$. 

We have already seen the method of recursion\index{Recursion} in the form 
of the Euclidean algorithm\index{Algorithm}, where the crucial step is to 
reduce the problem to smaller numbers. Another example is the calculation of 
the factorial function, which is defined by the formula
\[
f(n)=n!=1 \cdot 2 \cdots (n-1) \cdot n. 
\]
For this function, $f(0)=1$ and 
\[
f(n)=n \cdot f(n-1) \text{ for } n \ge 1. 
\]
The calculation of the function $f$ at the argument $n$ thus uses the product 
of $n$ with the value of the function at the argument $n-1$. Such a call of a
value of a function at a smaller argument is called recursion\index{Recursion}. 
In most programming languages, recursive\index{Recursion} programs can be 
written so that the program code for the factorial function contains a similar 
line like 
\[
\text{return } n * f(n-1), 
\]
i.e., the function $f$ calls itself. 

\medskip
To better illustrate the principle of recursion\index{Recursion}, we want to 
explain the illustrative and at first glance unmathematical example of the 
Towers of Hanoi. There are three rods and $n$ holed discs given, which are 
stacked on the first rod. The task is to move the discs of the tower one after 
the other so that they are stacked on the third rod at the end. The middle rod 
can be used as a help stack and a larger disc must never lie on a smaller one. 
There is a recursive\index{Recursion} solution strategy $L(n)$ for $n$ discs. 
It consists of placing the top $n-1$ discs on the help stack using the solution 
strategy $L(n-1)$, then placing the largest disc on the third rod and then 
applying the solution strategy $L(n-1)$ again to transport the discs from the 
help stack to the third rod. Fig.~\ref{fig:hanoi} illustrates the case $n=4$.

\newcommand{\pyramid}[1]{%
    \begin{tikzpicture}[x=0.25cm,y=0.25cm]
        \foreach \b [count=\n] in {#1} {
            \if \b5 \draw [black,fill=gray]  
                (-\b,{(\n-1)*0.4cm}) rectangle (\b,\n*0.4cm);
            \else \draw [black,fill=lightgray]
                (-\b,{(\n-1)*0.4cm}) rectangle (\b,\n*0.4cm);
            \fi
        }
    \end{tikzpicture}
}

\begin{figure}[ht!]
\hskip-5cm
\pyramid{5,4,3,2}
\vskip0.4cm 
\hskip-2.7cm
\pyramid{5}
\pyramid{4,3,2}
\vskip0.4cm 
\hskip2.8cm \pyramid{4,3,2}
\pyramid{5}
\vskip0.4cm 
\hskip5.1cm\pyramid{5,4,3,2}
\caption{\label{fig:hanoi}For the recursive solution in the case $n=4$, one 
first (recursively) places the upper three disks on the help stack. Then one 
puts the lowest disc on the third rod. After that, one places (recursively) the 
three smaller discs from the help stack there too.}
\end{figure}

In the solution for $n$ discs, the solution for $n-1$ discs is applied twice 
and the largest disc is moved exactly once. Therefore, for a complete solution, 
\[
1+2(1+2(1+...)...)=1+2+4+8+ \cdots + 2^{n-1}=2^n-1
\]
discs must be moved. It can be shown that it does not work with fewer moves. 

Dedekind\index{Dedekind, Richard} showed with the recursion 
theorem\index{Recursion} that the natural numbers are unique up to 
isomorphism\index{Isomorphism}. However, he used quantifiers over all subsets 
of the natural numbers, so from a modern point of view, he based his proofs on 
second-order predicate logic. Thoralf Skolem\index{Skolem, Thoralf} later 
showed that in first-order predicate logic, non-standard 
models\index{Non-standard model} of the natural numbers exist that have 
different properties than we usually expect. Skolem's\index{Skolem, Thoralf} 
exclusive use of first-order logic became a prevailing standard after some 
time. 

Dedekind's\index{Dedekind, Richard} approach can be seen as the beginning of 
structuralism\index{Structuralism}, a mathematical-philosophical view. He used 
the concept of mapping in a modern sense and thus anticipated the emergence of 
category theory\index{Category theory} in the 20th century, which emphasised 
the structural properties of mathematical objects more than the objects 
themselves, which are interchangeable.\footnote{See \cite{ms2023,ms2027}.}

Recursive\index{Recursion} calculations and proofs with complete 
induction\index{Complete induction} were already known before 
Dedekind\index{Dedekind, Richard}, for example by Franziskus 
Maurolicus\index{Maurolicus, Franziskus}, Blaise Pascal\index{Pascal, Blaise} 
and Jakob Bernoulli\index{Bernoulli, Jakob}. In England, Lady Ada 
Lovelace\index{Lovelace, Ada} developed a program for 
recursive\index{Recursion} calculations of Bernoulli numbers\index{Bernoulli, 
Jakob} on the unfinished mechanical calculator planned by Charles 
Babbage\index{Babbage, Charles} since 1822, which he called the Analytical 
Engine. Herbert Gra{\ss}mann\index{Gra{\ss}mann, Herbert}, Charles S. 
Peirce\index{Peirce, Charles} and Giuseppe Peano\index{Peano, Giuseppe} also 
considered the axioms of natural numbers, without using the recursion 
theorem\index{Recursion} in the way Dedekind\index{Dedekind, Richard} did. 
Peano\index{Peano, Giuseppe}, like Russell\index{Russell, Bertrand}, 
extensively developed the language of logic from Frege\index{Frege, Gottlob} 
and refined the axioms of natural numbers based on the preliminary work of 
Dedekind\index{Dedekind, Richard}.

Between 1900 and 1930, recursion theory\index{Recursion} was slowly developed 
further and a theory of primitive recursive functions\index{Recursion} emerged, 
particularly in the works of David Hilbert\index{Hilbert, David}, 
Thoralf Skolem\index{Skolem, Thoralf} and R\'ozsa P\'eter\index{Peter@P\'eter, 
R\'ozsa}. Primitive recursive functions\index{Recursion} form the smallest 
class of everywhere defined functions
\[
f \colon \mathbb{N}^n \xlongrightarrow{~~~~~~} \mathbb{N}
\]
in several variables, which contains constant functions and projection mappings 
and is closed under the recursion scheme\index{Recursion}
\begin{align*}
f(0,y)&=g(y) \cr
f(x+1,y)&=h(f(x,y),x,y), 
\end{align*}
where $g$ and $h$ are also primitive recursive functions\index{Recursion} and 
the variable $x$ denotes one of the arguments of $f$ and $y$ the rest. In 
addition to the factorial function, many other elementary functions are 
primitive recursive\index{Recursion}. Surprisingly, the function $f$ with
\[
f(n) = p_n \text{ (the $n$-th prime number)}
\]
is also primitive recursive\index{Recursion}.

The example of recursion theory\index{Recursion} and the beginnings of the 
proof-theoretic work of Hilbert\index{Hilbert, David} and 
Ackermann\index{Ackermann, Wilhelm} clearly show how closely proof 
theory\index{Proof theory} and computability theory\index{Computability} are 
related. Wilhelm Ackermann\index{Ackermann, Wilhelm} and R\'ozsa 
P\'eter\index{Peter@P\'eter, R\'ozsa} dealt with functions that are 
computable\index{Computability}, but not primitive recursive\index{Recursion}. 
This includes the Ackermann-P\'eter function $A(m,n)$, which in 
P\'eter's\index{Peter@P\'eter, R\'ozsa} representation from 1935 is given by
\begin{align*} 
A(0,n)&=n+1 \cr 
A(m+1,0)&=A(m,1) \cr 
A(m+1,n+1)&=A(m,A(m+1,n)).
\end{align*}
This function is not primitive recursive\index{Recursion}, because 
its growth is stronger than that of any primitive recursive\index{Recursion} 
function. The discovery of rapidly growing functions can be traced back to 
Godfrey Harold Hardy\index{Hardy, Godfrey Harold}. 

Primitive recursive\index{Recursion} functions can be classified with the help 
of the Grzegorczyk hierarchy\index{Grzegorczyk, Andrzej}.\footnote{See 
\cite{grzegorczyk1964,rose}.} For more general computable \index{Computability} 
functions, an analogous transfinite hierarchy had first shown to be difficult 
by 
results of Stephen C. Kleene\index{Kleene, Stephen}, John 
Myhill\index{Myhill, John}, Norman Arthur Routledge\index{Routledge, Norman 
Arthur}, and Solomon Feferman\index{Feferman, Solomon}. However, after an idea 
of Georg Kreisel\index{Kreisel, Georg}, this works for the subclass of 
so-called  ordinal-recursive  resp. $\epsilon_0$-recursive\index{Recursion} 
functions. The multiple recursive functions of R\'ozsa 
P\'eter\index{Peter@P\'eter, R\'ozsa} correspond in this respect to the ordinal 
number\index{Ordinal number} $\omega^\omega$ and the  
Ackermann-P\'eter function $A(m,n)$ to the ordinal number\index{Ordinal 
number} $\omega^2$. In principle, all functions which are 
provably computable\index{Computability} in Dedekind-Peano 
arithmetic\index{Dedekind-Peano arithmetic} can be classified via 
recursion\index{Recursion} over ordinal numbers\index{Ordinal number} 
$<\epsilon_0$ inside the L\"ob-Schwichtenberg-Wainer 
hierarchy\index{Loeb@L\"ob-Schwichtenberg-Wainer hierarchy}. 
G\"odel\index{Goed@G\"odel, Kurt} used a variant of this characterisation in 
his famous system $T$ to study the 
consistency\index{Consistency}\index{Contradiction-freeness} 
of Dedekind-Peano arithmetic\index{Dedekind-Peano arithmetic}.\footnote{See 
\cite{rose,schwichtenberg} for aspects of recursive functions and their 
hierarchies.}

\section{The Theory of Computability\index{Computability}} 

The first machine-like models of computability\index{Computability} were 
invented around 1910 by Axel Thue\index{Thue, Axel} and are now called Thue 
systems\index{Thue system}. His standpoint was that every form of computation 
or proof arises as a sequence of term substitutions. Thue's\index{Thue, Axel} 
ideas had a tremendous influence in linguistics, especially through the works 
of Noam Chomsky\index{Chomsky, Noam} and Richard Montague\index{Montague, 
Richard}.\footnote{See \cite{thue}.}

After Kurt G\"odel\index{Goed@G\"odel, Kurt} and Jacques 
Herbrand\index{Herbrand, Jacques} had dealt with a generalisation of primitive 
recursive functions\index{Recursion} in a brief correspondence in 1931, which 
was later improved, it was not until the annus mirabilis 1936 that the concept 
of computability\index{Computability} was clarified with the involvement of 
several people. The most famous among them is the concept of the Turing 
machine\index{Turing machine} and thus the term Turing 
computability\index{Computability} by Alan Turing\index{Turing, 
Alan}.\footnote{See \cite{turing1936}.} 

\begin{figure}[ht!]
\begin{tikzpicture}[start chain=1 going right,start chain=2 going below,node 
distance=-0.15mm]
    \node [on chain=2] {tape};
    \node [on chain=1] at (-1.16,-.4) {\ldots};  
    \foreach \x in {1,2,...,11} {
        \x, \node [draw,on chain=1,fill=gray!20] {};
    } 
    \node [name=r,on chain=1] {\ldots}; 
    \node [name=k, draw, on chain=2] at (-0.335,-.65) {};    
    \node at (1.6,-.85) {read/write head};
    \node [on chain=2] {};
    \node [draw,on chain=2,fill=gray!20] {program};
    \chainin (k) [join]; 
\end{tikzpicture}
\caption{\label{fig:turing}A blueprint of a Turing machine consists of a 
steering unit (program), a read/write head, and the tape on which labels from 
an alphabet are written.}
\end{figure}

Fig.~\ref{fig:turing} outlines the functioning of a Turing machine\index{Turing 
machine}. The unit labelled as program contains a deterministic function that 
describes how to proceed from the states and the labels of the tape with an 
alphabet, i.e., how the read/write head moves and how the tape is newly 
labelled.

Turing\index{Turing, Alan} noted that universal Turing 
machines\index{Turing machine} exist that can simulate any other Turing 
machine\index{Turing machine}. This is similar to a compiler that makes an 
executable file from program code. The theoretical concept of Turing machines 
is therefore very similar to a programmable computer. However, the history of 
computer development was lengthy. After Leibniz's\index{Leibniz, Gottfried 
Wilhelm} calculating machine and Babbage's\index{Babbage, Charles} rather 
unsuccessful attempt to develop his programmable Analytical Engine, the Z3 
computer by Konrad Zuse\index{Zuse, Konrad} from 1941 and the American computer 
Mark I from 1944 were the first functioning mechanical freely programmable 
computers that were universal Turing machines\index{Turing machine}. The Eniac 
computer from 1946 at the University of Pennsylvania was the first electronic 
computer of this kind. Today, computers in their internal structure are mostly 
realised by the von Neumann\index{Neumann, John von} architecture, which was 
conceived by John von Neumann\index{Neumann, John von} and his team during 
their work on the Eniac.

Emil Post\index{Post, Emil}, who was familiar with Thue's\index{Thue, Axel} 
work, invented a mathematical theory of automata that resembles Thue 
systems\index{Thue system} and Turing machines\index{Turing machine}. Alonzo 
Church\index{Church, Alonzo} founded the $\lambda$-calculus -- a third 
alternative computability theory and today the basis of functional programming 
languages -- and invented the so-called simple theory of types\index{Type 
theory}. 

Finally, Stephen C. Kleene\index{Kleene, Stephen} defined recursive 
functions\index{Recursion}, also known as $\mu$-recursive functions, 
as an extension of the primitive recursive functions\index{Recursion} and 
demonstrated the equivalence of these four definitions.\footnote{See 
\cite{church1940,turing1936,post1936} and the book 
\cite{kleene1952}.} For this, he used the $\mu$-operator, which for each 
partially defined recursive\index{Recursion} function 
\[
f(t,x_1,\ldots,x_n) \colon \mathbb{N}^{n+1} \xlongrightarrow{~~~~~~} 
\mathbb{N}
\]
forms a new -- possibly only partially defined -- function $\mu f$ with
\[
\mu f(x_1,\ldots,x_n) \colon \mathbb{N}^{n} \xlongrightarrow{~~~~~~} 
\mathbb{N},
\]
where the value 
\[
\mu f(x_1,\ldots,x_n)=k,  
\]
if there is a smallest natural number $k$ such that
\[
f(k,x_1,\ldots,x_n)=0 
\]
and $f(t,x_1,\ldots,x_n)$ is defined for all $0 \le t \le k$. If such a $k$ 
does not exist, $\mu f(x_1,\ldots,x_n)$ is undefined. The $\mu$-operator can 
be implemented in programming languages using while-loops by counting up. The 
class of $\mu$-recursive functions is generated from the primitive recursive 
functions\index{Recursion} together with the $\mu$-operator and forms a class 
that includes partially defined functions. For example, the 
Ackermann-P\'eter function can be expressed through the $\mu$-operator.

The $\mu$-operator is by definition a kind of search operator that seeks a 
root of $f$ in a variable, if one exists. If this is not the case, the 
resulting function $\mu f$ is undefined. For example, integer solutions of 
diophantine equations such as the Catalan equation
\[
x^2-y^3=1
\]
can be calculated by nested application of the $\mu$-operator -- or 
equivalently using while-loops. In this example, it is best to search on all 
lines of the form $x+y=n$, where $n$ runs through all natural 
numbers.\footnote{In the Python programming language, we search for the 
solution $(3,2)$ as follows: 
\\
k,n=1,2\\
while k**2-(n-k)**3!=1:\\
\noindent\hspace*{5mm}  if k<n-1: \\
\noindent\hspace*{10mm} k=k+1\\
\noindent\hspace*{5mm}  else: \\
\noindent\hspace*{10mm} k,n=1,n+1\\
print(k,n)}
Preda Mih\u{a}ilescu\index{Mih\u{a}ilescu, Preda} proved the Catalan 
conjecture\index{Catalan's conjecture} in 2002, i.e., this equation only has 
the solution $(x,y)=(3,2)$ in natural numbers $x,y \ge 1$. 

The Church-Turing thesis\index{Church-Turing thesis}, as originally 
formulated by Alonzo Church\index{Church, Alonzo} and others, states that all 
effectively computable\index{Computability} functions are always 
recursive\index{Recursion}. This thesis has been confirmed many times, but it 
is very informal in terms of the concept of computability\index{Computability} 
and its physical realisation and is therefore considered 
unprovable\index{Provability}. A refutation of the Church-Turing 
thesis\index{Church-Turing thesis} would only be possible through stronger 
physically realisable computability models\index{Computability}, which do not 
exist to this day.

Artificial neural networks\index{Neural network} with rational (or 
computable\index{Computability}) weights, as well as the still largely 
fictional quantum computers\index{Quantum computing}, are fascinating 
computability models\index{Computability}. Both fulfil the Church-Turing 
thesis\index{Church-Turing thesis}, as they can in principle be simulated by 
Turing machines\index{Turing machine}. Non-deterministic Turing 
machines\index{Turing machine} represent another computability 
model\index{Computability}, which is predominantly of a theoretical nature. 
Mathematically, this is achieved by replacing the deterministic function in the 
program of a Turing machine\index{Turing machine} with a relation, i.e., a 
multi-valued function.

Algorithms\index{Algorithm} are classified in various ways according to 
objectives or properties.\footnote{See \cite{wigderson2019}. Concrete 
goals are realised by optimisation, sorting, and search algorithms. Well-known 
sorting algorithms are Quicksort and Bubblesort. Another class of algorithms 
are cryptographic algorithms and protocols, such as the symmetric 
protocols AES (advanced encryption standard) and DES (digital encryption 
standard) or the asymmetric public-key protocols like RSA cryptography, ECC 
(elliptic curve cryptography) and ElGamal cryptography.} An interesting 
class are the stochastic algorithms\index{Algorithm}, also called probabilistic 
algorithms\index{Algorithm}. Effective stochastic algorithms\index{Algorithm} 
are the Monte Carlo algorithms\index{Algorithm}. An example of this is Monte 
Carlo integration, where high-dimensional integrals are evaluated by randomly 
selecting a sufficient number of random support points. Stochastic 
algorithms\index{Algorithm} are often more efficient than deterministic 
algorithms\index{Algorithm}. In the case of prime number tests, the 
Agrawal-Kayal-Saxena prime number test is a deterministic 
algorithm\index{Algorithm} that provably runs in polynomial 
time\index{Polynomial time}. However, probabilistic prime number tests, such as 
the Miller-Rabin test or the Solovay-Strassen test, are often much faster in 
practice. This is no coincidence, as the Miller-Rabin test has a deterministic 
variant that runs in polynomial time\index{Polynomial time} under the 
assumption of a generalisation of the Riemann hypothesis\index{Riemann 
hypothesis}, because the set of necessary samples can be limited.

A special class of algorithms\index{Algorithm} are the quantum mechanical 
algorithms\index{Algorithm}. These include the Shor algorithm\index{Algorithm}, 
which was found by Peter Shor\index{Shor, Peter} in 1994 for the factorisation 
of natural numbers, as well as the Deutsch-Jozsa algorithm\index{Algorithm} 
and Grover's search algorithm\index{Algorithm}. They could be effectively 
implemented on freely programmable universal quantum computers\index{Quantum 
computing} as soon as they can process enough qubits. Quantum annealers and 
other commercial quantum computers\index{Quantum computing} that are not 
universal cannot usually execute such algorithms.

\section{Complexity Theory\index{Complexity theory}}

Even if an algorithm\index{Algorithm} exists for a given problem, it matters 
greatly how quickly the algorithm\index{Algorithm} works and how much storage 
space it requires. To increase the performance of algorithms\index{Algorithm}, 
parallelisation is used, which arranges program parts so that they can be 
executed simultaneously. However, not every algorithm\index{Algorithm} is 
suitable for this.

The effort required by an algorithm\index{Algorithm} is considered in 
complexity theory\index{Complexity theory}, which is a separate field of 
research with many open questions.\footnote{See \cite{wigderson2019} for 
algorithms and complexity theory.} Computable\index{Computability} functions 
and algorithms\index{Algorithm} that calculate these can be divided into 
complexity classes\index{Complexity theory} depending on the growth or degree 
of nesting of the functions and the computational and storage effort generated. 
An important classification of functions is provided by the Grzegorczyk 
hierarchy\index{Grzegorczyk, Andrzej}. One of the lowest levels of this are the 
so-called elementarily computable\index{Computability} functions, which are 
sufficient for surprisingly many situations.\footnote{See \cite[Ch. 5]{rose}.}

A completely different classification is given by the runtime of 
algorithms\index{Algorithm}. The complexity class\index{Complexity theory} 
$\mathbf{P}$ is the class of polynomial time 
algorithms\index{Algorithm}\index{Polynomial time}. These are defined by the 
requirement that the runtime (and thus the memory requirement) is polynomial 
in the effort given by the coding length of the input data.\footnote{There are 
other characterisations of recursive functions in $\mathbf{P}$ through 
predicative recursion, see \cite{bellantonicook1992} and the literature 
therein.} It should be noted that the input of a natural number $n$ into a 
computable\index{Computability} function (or a Turing machine\index{Turing 
machine}) means an effort of $O(\log(n))$ input data, because the number of 
digits of $n$ is proportional to $\log(n)$.\footnote{The effort is $O(h(n))$, 
if there is a function $h$ and a positive constant $C$ such that for large 
values $n$ the effort grows like $C \cdot h(n)$.} Algorithms\index{Algorithm} 
in the class $\mathbf{P}$ thus have the effort $O(\log^k(n))$ for a natural 
number $k$.

Another famous complexity class\index{Complexity theory} is the larger class 
$\mathbf{NP}$ of computable\index{Computability} functions that can be computed 
by non-deterministic Turing machines\index{Turing machine} with 
polynomial\index{Polynomial time} runtime. A modern definition of this  
complexity class\index{Complexity theory} states that a problem is in 
class $\mathbf{NP}$ if there is a deterministic polynomial 
time\index{Polynomial time} algorithm $V$ (verifier) which can verify the 
correctness of a solution of the problem. One of the Millennium 
problems\index{Millennium problems} of the Clay Foundation asks whether
\[
\mathbf{P}=\mathbf{NP}  
\]
holds or not. Remarkable about this question is that there are problems that 
are in $\mathbf{NP}$ and are $\mathbf{NP}$-complete. If such a problem is in 
$\mathbf{P}$, then $\mathbf{P}=\mathbf{NP}$ already holds and vice versa. The 
Millennium problem\index{Millennium problems} therefore only needs to be tested 
on one of them. These include the SAT problem, the traveling salesman problem, 
and the clique problem.

The SAT problem, also known as the satisfiability problem of propositional 
logic, consists in deciding for each formula of propositional logic whether 
there are truth values\index{Truth} (true or false) for all variables occurring 
in it, so that the whole formula becomes true after substituting these values. 
Obviously, any solution to this problem is verifiable in polynomial 
time\index{Polynomial time}. On the other hand, the number of truth values to 
be substituted for $n$ variables is of the order of $2^n$ and thus exponential 
in $n$. Stephen Cook\index{Cook, Steven} and Leonid Levin\index{Levin, Leonid} 
independently showed that the SAT problem is 
$\mathbf{NP}$-complete.\footnote{See \cite[Ch. 3]{wigderson2019} for 
$\mathbf{NP}$-completeness.}

Quantum computers\index{Quantum computing} are potentially much more powerful 
than classical digital computers, as the Shor algorithm\index{Algorithm} shows. 
Their polynomial time\index{Polynomial time} complexity class\index{Complexity 
theory} is denoted by $\mathbf{BQP}$. It lies in the class 
$\mathbf{PSPACE}$ of decidable problems with polynomial effort. In recent 
years, quantum computers\index{Quantum computing} with more than $50$ qubits 
have been built, which can execute specially adapted 
algorithms\index{Algorithm} faster than classical computers. Presumably, 
algorithmic problems always require significantly greater effort with classical 
computers than with quantum computers, even if this is not provable at this 
point. This popular assumption refines the Church-Turing 
thesis\index{Church-Turing thesis} and is called the supremacy of quantum 
computers\index{Quantum computing}. Currently, the technology is still far from 
such expectations.

\section{Undecidable Problems}
 
Surprisingly, there are mathematical problems that cannot be solved 
algorithmically and are called undecidable\index{Undecidability}. They are 
based on the existence\index{Existence} of undecidable\index{Undecidability} 
subsets of the natural numbers. A set $S \subseteq \mathbb{N}$ is 
decidable\index{Decidability} if $S$ and its complement $\mathbb{N} \setminus 
S$ are both recursively enumerable, i.e., they are either empty or the set of 
values of a computable\index{Computability} function 
$f \colon \mathbb{N} \longrightarrow \mathbb{N}$. In particular (and 
equivalently), the characteristic function of $S$ is 
computable\index{Computability}.

Hilbert\index{Hilbert, David} and Ackermann\index{Ackermann, Wilhelm} 
formulated the decision problem\index{Decision 
problem}\index{Entscheidungsproblem} in 1928. It asks 
whether for every statement in a formal language\index{Formal language} over 
first-order predicate logic, it can be decided whether 
it is provable\index{Provability} or not. Alan Turing\index{Turing, Alan} and 
Alonzo Church\index{Church, Alonzo}\footnote{See \cite{church1936,turing1936}.} 
independently showed in 1936 that the decision problem is algorithmically 
unsolvable and thus undecidable\index{Undecidability}, by showing that the 
subset $S$ of provable\index{Provability} statements is 
undecidable\index{Undecidability}. The proof is nowadays carried out using the  
connection with the halting problem\index{Halting problem} for Turing 
machines\index{Turing machine}. This result relativised the ideas of 
Leibniz\index{Leibniz, Gottfried Wilhelm} on the decidability of all scientific 
questions. 

The halting problem\index{Halting problem} asks whether there is an 
algorithm\index{Algorithm} that decides for a given 
G\"odel numbering\index{Goed@G\"odel, Kurt} $T_n$ of Turing 
machines\index{Turing machine}, which compute partially defined functions $f_n 
\colon \mathbb{N} \longrightarrow \mathbb{N}$, whether the function $f_n$ is 
defined at the argument $n$. Today's proofs for the 
undecidability\index{Undecidability} of the halting problem\index{Halting 
problem} use a variant of Cantor's diagonal argument\index{Cantor, Georg}. In 
these, assuming that the halting problem\index{Halting problem} is 
decidable\index{Decidability}, the total computable\index{Computability} 
function $h$ with 
\[
h(k)=
\begin{cases}
f_k(k)+1 & \text{if } f_k(k) \text{ is defined} \cr
0  & \text{otherwise}
\end{cases}
\]
is considered. There is then a number\index{Goed@G\"odel, Kurt} $n$, such that 
$h$ is computed by the Turing machine\index{Turing machine} $T_n$ corresponding 
to the function $f_n$. This yields a contradiction because essentially 
\[
h(n)=f_n(n) \text{ and } h(n)=f_n(n)+1
\]
would simultaneously hold.

There are numerous other significant problems in mathematics for which no 
general solution algorithm exists. These include Hilbert's 10th 
problem\index{Hilbert problems} and the word problem\index{Word problem} for 
semigroups and groups. Somewhat more generally, a deep theorem of Henry 
Gordon Rice\index{Rice, Henry Gordon} states that any non-trivial 
semantic\index{Semantics} property of Turing machines\index{Turing machine} is 
undecidable\index{Undecidability}. By this we mean a non-trivial property of 
the output of the Turing machine\index{Turing machine}, for example that the 
input $n$ yields the $n$-th decimal digit of $\pi$ as output. The proof of 
Rice's theorem\index{Rice, Henry Gordon} can be pursued with the help of the 
halting problem\index{Halting problem}, which also provides a beautiful 
semantic\index{Semantics} property.\footnote{See \cite{bridges1994,davis1965}.}

The 10th Hilbert problem\index{Hilbert problems}\index{Hilbert, David} 
asks for an algorithm\index{Algorithm} that, for every 
integer diophantine equation 
\[
F(x_1,...,x_n)=0
\]
decides whether the equation has an integer solution or not. 
After preliminary work by Martin Davis\index{Davis, Martin}, Julia 
Robinson\index{Robinson, Julia} and Hilary Putnam\index{Putnam, Hilary}, 
Yuri Matiyasevich\index{Matiyasevich, Yuri} showed in 1970 that the 10th 
problem\index{Hilbert problems} is undecidable\index{Undecidability}, 
i.e., there is no such algorithm\index{Algorithm}. The 
obvious search algorithm for finding solutions is not effective, 
because there is no a priori estimate for the size of possible 
solutions as a termination condition.

The negative solution of the 10th Hilbert problem\index{Hilbert problems} 
follows from the existence\index{Existence} of 
undecidable\index{Undecidability} 
sets and a theorem by Davis\index{Davis, Martin}, 
Matiyasevich\index{Matiyasevich, Yuri}, Robinson\index{Robinson, Julia} and 
Putnam\index{Putnam, Hilary}, which states that every recursively 
enumerable\index{Recursion} set is diophantine and therefore given by the 
projection of a zero set of a diophantine equation. 

Davis\index{Davis, Martin}, Putnam\index{Putnam, Hilary} and 
Robinson\index{Robinson, Julia} first proved that certain 
exponential equations, such as the zero set 
\[
x^y=z 
\]
in the $(x,y,z)$-space, are diophantine, although this does not seem so 
at first glance. Julia Robinson\index{Robinson, Julia} had the idea to use 
the classical Pell's equation
\[
x^2-dy^2=\pm 1 
\]
in the case of $d=a^2-1$, whose infinitely many solutions have suitable 
growth behaviour. The proof by Matiyasevich\index{Matiyasevich, Yuri} 
proceeded somewhat differently, as he used an equation of the form 
\[
x=F_{2y} 
\]
in the $(x,y)$-space, where $F_m$ is the sequence of Fibonacci numbers 
\[
0,1,1,2,3,5,8,13,21,\ldots  
\]
defined by $F_0=0$, $F_1=1$ and $F_{m+2}=F_{m}+F_{m+1}$. Modern proofs usually 
use Pell's equation as Julia Robinson\index{Robinson, Julia} did.\footnote{See 
\cite[Ch. V]{manin}.}

The word problem\index{Word problem} for finitely presented 
groups or semigroups
\[
G=\{a_1,\ldots,a_s \mid r_1,\ldots,r_t\}, 
\]
which are given by finitely many generating elements $a_1,\ldots,a_s$ and 
relations\footnote{A relation $r$ is an equation $u=v$, where $u, v$ are words 
in the generating elements. In the case of groups, $v=1$ without loss of 
generality.} $r_1,\ldots,r_t$, consists in deciding whether an arbitrary word 
$w$ in $G$ is equal to a given word $w_0$: 
\[
w=w_0.
\]
So an algorithm\index{Algorithm} is sought that finds finitely many 
relations, so that the word $w$ goes over into the word $w_0$ after applying 
these relations as term substitutions.

The word problem\index{Word problem} was first formulated by Axel 
Thue\index{Thue, Axel} in his already mentioned works on 
computability\index{Computability} and he was aware of the 
difficulty of this problem. Max Dehn\index{Dehn, Max}, who had already solved 
the 3rd Hilbert problem\index{Hilbert problems} in 1903, formulated the word 
problem\index{Word problem} independently from Thue\index{Thue, Axel}. The 
presentation of semigroups and groups is closely related to the 
computability model\index{Computability} of Thue systems\index{Thue system} 
because the relations $r_j$ can be understood as term substitutions. 

The proof for the undecidability\index{Undecidability} of the 
word problem\index{Word problem} for semigroups was independently provided by 
Emil Post\index{Post, Emil} and Andrey Markov\index{Markov, Andrey} 
in 1947. Indeed, the connection with computation models\index{Computability} 
can be used to demonstrate undecidability\index{Undecidability}. Pyotr S. 
Novikov\index{Novikov, Pyotr} and William Boone\index{Boone, William} solved 
the word problem\index{Word problem} for groups a few years later. 

There are important classes of semigroups and groups for which the 
word problem\index{Word problem} is solvable. On the other hand, there are 
specific finitely presented groups for which the problem has no solution. A 
handy example was found by Gregory Tseytin\index{Tseytin, Gregory}. For 
the semigroup 
\[ 
G=\langle a,b,c,d,e \mid ac=ca, ad=da, bc=cb,bd=db,ce=eca, 
\]
\[
de=edb, cdca=cdcae, caaa=aaa, daaa=aaa \rangle 
\]
the problem is undecidable\index{Undecidability}, whether an arbitrary word 
matches $w=aaa$.\footnote{See \cite{davis1965,ms2024,post1947,thue}.}

\section{Artificial Intelligence\index{Artificial intelligence} 
and Hypercomputing\index{Hypercomputing}}

We have dealt with the Platonic world of ideas\index{Platonic idealism} and its 
connection with abstract concepts and truth\index{Truth}. All this thinking 
takes place in our heads. It is an interesting and still unsolved question 
whether human intelligence is superior to a computer or not. 
G\"odel\index{Goed@G\"odel, Kurt} was convinced that this is the case. His 
argument used ideas around his own 
incompleteness theorem\index{Incompleteness theorem} and is often referred to 
in slight variations as the G\"odel-Lucas-Penrose\index{Penrose, 
Roger}\index{Lucas, John}\index{Goed@G\"odel, Kurt} argument. John 
Searle's\index{Searle, John} thought experiments go in the same direction. 

The answer to this question was not at all obvious in Alan 
Turing's\index{Turing, Alan} time. He proposed a related test in an 
article,\footnote{See \cite{turing1950}.} which is now called the Turing 
test\index{Turing test}. It essentially involves whether a human can be 
distinguished from a computer through a question-and-answer game. The 
slowness of the brain and its error-proneness make it clearly inferior to 
today's computers in many respects. On the other hand, the test neglects some 
crucial aspects of human intelligence, such as our consciousness, the presence 
of which may distinguish us from a machine. The Turing test\index{Turing, 
Alan}\index{Turing test} in its historical version is no longer the right 
question for various reasons. 

In recent years, the field of artificial 
intelligence\index{Artificial intelligence} has come to the fore. 
A prophet in this field was the Leibniz\index{Leibniz, Gottfried 
Wilhelm}-inspired mathematician Norbert Wiener\index{Wiener, Norbert}, who 
developed crucial ideas for it as early as the 1950s and 
at the same time directed a critical, forward-looking view of the 
future.\footnote{See \cite{wiener}.} After many decades of 
investigating promising approaches, there have only been decisive breakthroughs 
in implementation in recent years. The currently dominant technology is 
primarily based on the training of various forms of recurrent artificial neural
networks\index{Neural network} using large amounts of data as input. In 
many areas, such as image recognition, this method has significantly surpassed 
any other approach. It is known that recurrent 
artificial neural networks\index{Neural network} can simulate any 
Turing machines\index{Turing machine}. 

Recurrent neural networks\index{Neural network} in their simplest form are 
based -- similar to the multilayer perceptron\index{Perceptron} -- on systems 
of equations 
\begin{align*}
v_\bullet(n) & = \sigma_V \left( W_E \, e_\bullet(n) + W_V v_\bullet(n-1) + b_V 
\right) \cr 
a_\bullet(n) & = \sigma_A \left( W_A \, v_\bullet(n) + b_A \right),  
\end{align*}
where $n=0,1,2,\ldots$ is a discrete time parameter, $\sigma_V$ and $\sigma_A$ 
are activation functions, $b_V$ and $b_A$ are biases, and  $W_E,W_V$ and $W_A$ 
are transition matrices with weights that operate on the input vectors 
$e_\bullet(n)$ and the hidden variables $v_\bullet(n)$ to finally yield the 
output $a_\bullet(n)$ (see Fig.~\ref{fig:recurrent_net}). 

\begin{figure}[ht!]
\begin{tikzpicture}[x=1.3cm, y=1.3cm, >=stealth]

\tikzset{every loop/.style={min distance=10mm,in=120,out=60,looseness=10}}

\tikzset{%
  every neuron/.style={
    circle,
    draw,
    minimum size=0.4cm
  },
  neuron missing/.style={
    draw=none, 
    scale=4,
    text height=0.333cm,
    execute at begin node=\color{black}$\vdots$
  },
}

\foreach \m/\l [count=\y] in {1,2,3} {
\ifnum \y=3 \node [every neuron/.try, neuron \m/.try,fill=lightgray] (input-\m) 
at 
(0,1.5-\y) {};
\else \node [every neuron/.try, neuron \m/.try,fill=lightgray] (input-\m) at 
(0,1.5-\y) {};
\fi  }
  
\foreach \m [count=\y] in {1,2,3}
  \node [every neuron/.try, neuron \m/.try] (hidden-\m) at (2,2-\y*1.25) {};

\foreach \m [count=\y] in {1,2,3}  
  \node [every neuron/.try, neuron \m/.try] (output-\m) at (4,1.5-\y) {}; 

\foreach \l [count=\i] in {1,2,3}
  \draw [<-] (input-\i) -- ++(-1,0)
    node [above, midway] {$E_{\l}$};
  
\foreach \l [count=\i] in {1,2,3} 
    \draw [->] (output-\i) --  (6,-0.5)  
    node [every neuron/.try,fill=lightgray] at (6.15,-0.5) {}; 

\foreach \i in {1,...,3}
  \foreach \j in {1,...,3}
    \draw [->] (input-\i) -- (hidden-\j);
       
\draw [->] (hidden-1) -- (output-1) node [above=0.2] {};
\draw [->] (hidden-2) -- (output-2) node [below=0.2] {};
\draw [->] (hidden-2) -- (output-3) node [below=0.2] {};
\draw [->] (hidden-3) -- (output-3) node [below=0.2] {};

\draw [->,dashed] (output-2) -- (hidden-3) node [below=0.2] {};
\draw [->,dashed] (output-1) -- (hidden-3) node [below=0.2] {};
\draw [->,dashed] (output-3) -- (hidden-1) node [below=0.2] {};
\draw [->,dashed] (output-2) -- (output-1) node [below=0.2] {};
\draw [->,dashed] (output-2) -- (output-3) node [below=0.2] {};
\draw [->,dashed] (hidden-1) -- (hidden-2) node [above=0.2] {};
\draw [->,dashed] (hidden-3) -- (hidden-2) node [below=0.2] {};

\path (hidden-1) edge [->,anchor=center,loop above,dashed] node {} (hidden-1);

\draw [->] (hidden-1) -- (output-2);
\draw [->] (hidden-2) -- (output-1);
\draw [->] (6.3,-0.5) -- (7.2,-0.5) node [above,midway] {$A$} ;
\end{tikzpicture}
\caption{\label{fig:recurrent_net}A recurrent neural network\index{Neural 
network} also has feedback edges (dashed) that point backwards or sideways.} 
\end{figure}

There is a fascinating new research direction which combines neural 
networks\index{Neural network} and category theory\index{Category theory}. The 
idea is that states of neural networks\index{Neural network} represent  
morphisms in some (higher) category\index{Category theory} which depend on 
parameters (i.e., the weights). The training process corresponds to changes of 
weights and thus results in a flow of morphisms until one reaches an optimal 
choice.

\medskip 
It is an open question whether there are other physically realisable 
computability models\index{Computability}, besides quantum 
computers\index{Quantum computing} or the different variants of artificial 
neural networks\index{Neural network}, that yield strictly more 
computing power than classical computers or can even compute functions that are 
not Turing computable\index{Computability}\index{Turing machine}. Such 
expectations are referred to as hypercomputing\index{Hypercomputing}. Possible 
architectures for this are generalisations of recurrent neural 
networks\index{Neural network} which contain analogue physical systems as a 
network structure in the hidden layers, possibly with quantum mechanical 
components. also known as reservoir computing. As with many other approaches, 
in addition to the question of hypercomputing\index{Hypercomputing}, a 
significant reduction in energy consumption is also sought. 

\chapter{Deductive Systems\index{Deductive system} and 
Incompleteness\index{Incompleteness theorem}}

In antiquity, scientific discourse developed rapidly. Arithmetic and geometric 
theorems were precisely derived under the assumption of axioms. 
Euclid's\index{Euclid} influential book \enquote{Elements} is an expression 
of this culture. From today's perspective, every mathematical proof is based on 
a syntactic\index{Syntax} calculus, which we refer to as a deductive 
system\index{Deductive system}. This concept includes a formal 
language\index{Formal language} and logical inference rules. After 
further developments of Aristotelian\index{Aristotle} logic by 
Llull\index{Llull, Ram\'on}, Leibniz\index{Leibniz, Gottfried Wilhelm}, 
Bolzano\index{Bolzano, Bernard} and others, it was only Frege\index{Frege, 
Gottlob} in his \enquote{Begriffsschrift} who introduced a deductive 
system\index{Deductive system} of today's kind. At the same time, recursion 
theory and the axiomatics of arithmetic were developed by 
Dedekind\index{Dedekind, Richard}. Shortly afterwards, Peano\index{Peano, 
Giuseppe} laid the foundations of today's notation in logic.\footnote{See 
\cite{bolzano1837,frege1879,ms2023,peano,wardhaugh}.} 

After the turn of the 20th century, further perspectives emerged. Axel 
Thue\index{Thue, Axel} designed his Thue systems\index{Thue system}, which 
represent elementary operations on trees. Hilbert\index{Hilbert, David} worked 
together with Wilhelm Ackermann\index{Ackermann, Wilhelm} and Paul 
Bernays\index{Bernays, Paul} on the mathematisation of proofs in axiomatic 
theories and developed the Hilbert calculus\index{Hilbert calculus}. Movements 
such as Hilbert's\index{Hilbert, David} formalism\index{Formalism} of axiomatic 
theories, intuitionism\index{Intuitionism} and 
constructivism\index{Constructivism} emerged.\footnote{Significant 
examples of axiomatic theories appeared already at the end of the 19th century, 
like Euclidean geometry, whose axioms were closely examined by Moritz 
Pasch\index{Pasch, Moritz} and David Hilbert\index{Hilbert, David}, 
Zermelo-Fraenkel set theory and Dedekind-Peano arithmetic. See 
\cite{linnebo} for formalism.} These ideas led to the emergence of the field of 
proof theory\index{Proof theory}, in which one can talk about mathematical 
proofs in the form of a metamathematics. A highlight of these developments are 
the two incompleteness theorems\index{Incompleteness theorem} by 
G\"odel\index{Goed@G\"odel, Kurt}. They form a fundamental obstacle for 
Hilbert's dream\index{Hilbert, David} of a proof of the
consistency\index{Consistency}\index{Contradiction-freeness} of mathematical 
theories, which only appears partially surmountable through the use of 
transfinite numbers or other axioms. 

\section{Formal Languages\index{Formal language} and Deductive 
Systems\index{Deductive system}}

A deductive system\index{Deductive system} is a calculus in which proofs can be 
conducted. The basis of every deductive system\index{Deductive system} is a 
formal language\index{Formal language} along with some axioms and inference 
rules. A proof in a deductive system\index{Deductive system} uses the axioms 
and inference rules in a suitable order until the desired result is 
achieved.\footnote{See \cite{boole1847,manin,takeuti2013}.}

Let's first turn to formal languages\index{Formal language}. The classical, 
two-valued propositional logic is one of the simplest formal 
languages\index{Formal language} with the special characters
\[
\neg, \wedge, \vee, \xRightarrow{~~}, 
\]
which we have already learned. Through the relationships 
\[
A \vee B = \neg(\neg A \wedge \neg B), \; (A \xRightarrow{~~} B) = \neg(A 
\wedge 
\neg B)
\]
it is sufficient to only consider the logical symbols $\wedge$ and $\neg$. 
Propositional logic has some variants, most of which are unknown to most 
people. There are multi-valued logics, modal logic and intuitionistic 
logic\index{Intuitionism}.

First-order predicate logic extends the 
propositional logic and consists of a reservoir of 
symbols for constants, free and bound variables, function and relation symbols 
(predicates ), as well as the additional symbols  
\[
\forall \text{ and } \exists, 
\]
which we have already learned. While function symbols can take all permissible 
values within the language, predicates provide the truth 
values true or false. A vivid example of a unary 
predicate is a property of the form  
\[
P(x) \colon x \text{ is red.} 
\]
The equality sign\index{Equality}\index{Identity} $=$ is usually part of every 
formal language\index{Formal language} as a binary relation or as a binary 
predicate . In a formal language\index{Formal language} 
$L$, there are usually other non-logical symbols. Depending on the language, 
these can be constants like $0$ and $1$, unary function symbols like $S$ 
(successor mapping), minus sign, inverse symbol or binary function symbols like 
$+, \cdot, \circ$ or relation symbols like $<$ and $\in$.  

From the constants, the variables and the function or relation symbols of a 
formal language\index{Formal language} $L$, terms can be formed. All 
constants and free variables are terms and the values of functions 
$f(t_1,\ldots,t_n)$, into which terms $t_1,\ldots,t_n$ have been inserted, are 
again terms. 

From terms, formulas can be formed, with atomic formulas arising by inserting 
terms into predicates $P(t_1,\ldots,t_n)$. At this point, the 
equality sign\index{Equality}\index{Identity} $=$ can occur. General formulas 
arise from atomic formulas by using the logical symbols, i.e., if $P$ and $Q$ 
are formulas, then $\neg P$, $P \wedge Q$, $P \vee Q$ and $P \Rightarrow Q$ are 
formulas. If $P(x)$ is a formula containing a free variable $x$, then 
\[
\exists x \, P(x) \text{ and } \forall x \, P(x)
\]
are again formulas. After all, a sentence is a formula which is closed, i.e., 
without free variables. Such expressions can be viewed as statements (or 
propositions, theorems) which have a truth value.

Calculi include, in addition to a formal language\index{Formal language}, 
additional axioms and inference rules. In all deductive systems\index{Deductive 
system}, substitution is one of the inference rules. The notation 
$P[\frac{t}{x}]$ is common when the expression $t$ is substituted for $x$ in 
the formula $P(x)$. When performing substitutions in sequence, the 
equality\index{Identity}\index{Equality} 
$P[\frac{t}{x}][\frac{s}{t}]=P[\frac{s}{x}]$ is usually required. The calculi 
of Frege\index{Frege, Gottlob}, Russell\index{Russell, Bertrand} and the 
Hilbert calculus\index{Hilbert calculus} possess numerous logical axioms and 
have only the modus ponens
\[
\frac{P \quad P \Longrightarrow Q}{Q} 
\]
as a new inference rule. The Hilbert calculus\index{Hilbert calculus} and 
preceding calculi of Frege\index{Frege, Gottlob} and others correspond little 
to the usual methods of mathematical reasoning and the logical axioms are 
unnatural in some versions. Therefore, from 1926, Jan 
{\L}ukasiewicz\index{Luka@{\L}ukasiewicz, Jan}, Stanis{\l}aw 
Ja\'skowski\index{Ja\'skowski, Stanis{\l}aw} and Gerhard Gentzen\index{Gentzen, 
Gerhard} introduced new calculi, of which we will use the calculus of natural 
deduction\index{Calculus of natural deduction}.\footnote{See \cite{gentzen}.} 
This calculus uses judgements\index{Judgement} of the form
\[
\Gamma \vdash A, 
\]
where $\Gamma$ denotes a finite number of hypotheses and $A$ is a single 
formula. The symbol $\vdash$ goes back to the Fregean\index{Frege, Gottlob} 
judgement stroke\index{Judgement} and expresses the derivability of $A$ 
based on the assumption $\Gamma$ in the calculus. It is the metalinguistic 
version of implication and must be distinguished from the object 
language\index{Object language} symbol $\Rightarrow$. If $\Gamma$ is empty, 
then $\vdash A$ denotes a formula $A$ that is derivable without hypotheses, 
i.e., an axiom or a theorem.

\medskip
In propositional logic, there are two tautological inference rules
\[
\frac{}{A \vdash A}, \;  \frac{\Gamma \vdash A}{\Gamma' \vdash A} 
\text{ for } \Gamma' \supseteq \Gamma 
\]
and two symmetry rules
\[
\frac{\Gamma \vdash (A \wedge B)}{\Gamma \vdash (B \wedge A)}, \; 
\frac{\Gamma \vdash (A \vee B)}{\Gamma \vdash (B \vee A)}.
\]
In addition, there are four introduction rules each for 
$\wedge, \vee, \neg, \Rightarrow$
\[
\frac{\Gamma \vdash A \quad \Gamma \vdash B}{\Gamma \vdash A \wedge B}, 
\; \frac{\Gamma \vdash A}{\Gamma \vdash A \vee B}, \; \frac{\Gamma, A \vdash B 
\quad \Gamma, A \vdash \neg B}{\Gamma \vdash \neg A},
\; \frac{\Gamma, A \vdash B}{\Gamma \vdash (A \Rightarrow B)}   
\]
and the corresponding four elimination rules
\[ 
\frac{\Gamma \vdash A \wedge B}{\Gamma \vdash A},\; 
\frac{\Gamma \vdash A \vee B \quad \Gamma, A \vdash C \quad \Gamma, B \vdash 
C}{\Gamma \vdash C},\;\frac{\Gamma \vdash \neg \neg A}{\Gamma \vdash A}
\]
and  
\[
\frac{\Gamma \vdash A \quad  \Gamma \vdash A \Rightarrow B}{\Gamma \vdash 
B} 
\text{ (modus ponens)}.
\]
One of the most important discoveries of Gentzen\index{Gentzen, Gerhard} is his 
main theorem of cut elimination with the provable inference rule
\[
\frac{\Gamma \vdash A \quad \Delta, A \vdash B}{\Gamma, \Delta \vdash B}, 
\]
which allows the formula $A$ to be eliminated. This rule in turn implies the 
two symmetry rules. Gentzen\index{Gentzen, Gerhard} used the cut elimination 
rule to show the consistency\index{Consistency}\index{Contradiction-freeness}, 
i.e., the contradiction-freeness, of Dedekind-Peano 
arithmetic\index{Dedekind-Peano arithmetic}.

The ability to find and elegantly conduct correct mathematical proofs is a 
great art. George P\'olya\index{Polya@P\'olya, George} wrote a wonderful book 
about it in 1945 titled \enquote{How to solve it}.\footnote{See 
\cite{polya1945}.} An unspectacular proof of a simple theorem like 
\[
\vdash (A \wedge B) \wedge C \xRightarrow{~~} A \wedge (B \wedge C), 
\]
which -- due to the rules applicable to $\Rightarrow$ -- is equivalent to 
\[
(A \wedge B) \wedge C \vdash A \wedge (B \wedge C),           
\] 
consists of the derivation tree shown in Fig.~\ref{fig:proof_tree}. The proof 
uses the abbreviation $\Gamma=(A \wedge B) \wedge C$ and -- starting from the 
tautology $\Gamma \vdash \Gamma$ -- uses conclusion rules at the solid lines. 

\begin{figure}[ht!] 
\[
\cfrac{
\cfrac{
\cfrac{\Gamma \vdash (A \wedge B) \wedge C}{\Gamma \vdash A \wedge B}
}{
\Gamma  \vdash A
} 
\quad 
\cfrac{
\cfrac{
\cfrac{\Gamma \vdash (A \wedge B) \wedge C}{ \cfrac{\Gamma \vdash A \wedge B 
}{\Gamma \vdash B \wedge A} }
}
{\Gamma \vdash B} 
\quad 
\cfrac{ \cfrac{\Gamma \vdash (A \wedge B) \wedge C}{\Gamma \vdash C \wedge (A 
\wedge B) } }{\Gamma \vdash C}
}
{\Gamma \vdash B \wedge C}
}
{
\Gamma \vdash A \wedge (B \wedge C)
}
\]
\caption{\label{fig:proof_tree}Derivation tree of $(A \wedge B) \wedge C 
\vdash A \wedge (B \wedge C)$.}
\end{figure}

In predicate logic, additional introduction rules for 
$\exists$, $\forall$ and $=$ are needed. These are 
\[ 
\frac{\Gamma \vdash P[\frac{t}{x}]}{\Gamma \vdash \exists x \, P(x)}, \; 
\frac{\Gamma \vdash P[\frac{t}{x}]}{\Gamma \vdash \forall x \, P(x)} \; (t 
\text{ 
free}), \; \frac{}{\Gamma \vdash P=P}.
\]
The corresponding elimination rules are 
\[
\frac{\Gamma \vdash \exists x \, P(x) \quad 
P[\frac{t}{x}] \vdash Q}{\Gamma \vdash Q} \; (t \text{ free}), \frac{\Gamma 
\vdash  \forall x \, P(x)}{\Gamma \vdash P[\frac{t}{x}]},  \; \frac{\Gamma 
\vdash s=t \quad \Gamma \vdash P[\frac{t}{x}]}{\Gamma \vdash P[\frac{s}{x}]}.
\]
The adjective free here means that the variable $t$ appears independently of 
the other occurring variables and constants and can take any values.

\section{The Dedekind-Peano Arithmetic\index{Dedekind-Peano arithmetic}}

We now want to deal with the formal language\index{Formal language} of 
Dedekind-Peano arithmetic\index{Dedekind-Peano arithmetic}. This language 
is used to describe the foundations of arithmetic, i.e., the theory of natural 
numbers $\mathbb{N}$ and their additive and multiplicative properties. 

The language of Dedekind-Peano arithmetic\index{Dedekind-Peano arithmetic}
in addition to the logical symbols $\wedge$, $\vee$, $\neg$, $\Rightarrow$, 
$\forall$, $\exists$, consists of the symbols $0$ (zero), 
$S$ (successor function), $+$ (addition), $\cdot$ (multiplication) and $=$ 
(equality). A typical formula in it is given by the $4$-squares theorem
\[
\forall n \, \exists a \, \exists b \, \exists c \, \exists d \quad 
n=a^2+b^2+c^2+d^2, 
\]
which states that every natural number $n$ can be represented as the sum of 
four squares. In the calculus of this language, the usual conventions 
such as \enquote{dot before line}, bracket symbols and other abbreviations are 
allowed. In Dedekind-Peano arithmetic\index{Dedekind-Peano arithmetic} there 
are the axioms listed in Fig.~\ref{fig:peano_axioms}. 

\begin{figure}[ht!]
\begin{align}
\forall n \quad &  S(n) \neq 0 \\
\forall m \forall n \quad  & S(m)=S(n) \xLeftrightarrow{~~~} m=n \\
\forall n \quad & n+0=n=0+n \\
\forall m \forall n \quad  & m+S(n)=S(m+n)=m+S(n) \\
\forall n \quad & n \cdot 0=0=0 \cdot n \\
\forall m \forall n \quad  & m \cdot S(n)=m \cdot n+m \\ 
& P(0) \wedge \forall n \left(P(n) \xRightarrow{~~} P(S(n))\right) 
\xRightarrow{~~} \forall n \; P(n). 
\end{align}
\caption{\label{fig:peano_axioms}The axioms of Dedekind-Peano 
arithmetic\index{Dedekind-Peano arithmetic}.}
\end{figure}

The successor function $S$ is noted as 
\[
S(n)=n+1 
\]
when addition has already been defined by recursion\index{Recursion}. So it 
holds:
\[
1=S(0), \; 2= S(S(0)), 3=S(S(S(0))), \; \ldots 
\]
The most important axiom is the last axiom (7) of complete 
induction\index{Complete induction}. For every property $P$ of natural numbers 
which can be formulated in the formal language\index{Formal language} of 
Dedekind-Peano arithmetic\index{Dedekind-Peano arithmetic} it expresses the 
following statement:

\begin{quote}
If $P(0)$ holds and if $P(n)$ always implies $P(S(n))$, then $P$ is fulfilled 
for every natural number.
\end{quote}

This statement can be interpreted as a rule of inference 
\[
\frac{P(0) \quad \forall n \left(P(n) \Longrightarrow 
P(S(n))\right)}{\forall n \; P(n)}. 
\]
Since there are a priori infinitely many such properties $P$, there are 
infinitely many axioms in first order predicate logic. If one quantifies over 
all $P$, and thus changes to second order logic, then there are only finitely 
many axioms.  

\medskip
In Dedekind-Peano arithmetic\index{Dedekind-Peano arithmetic} the order 
relation $<$ is also of great importance. This relation can be added to the 
language. 

\section{G\"odel's\index{Goed@G\"odel, Kurt} First  
Incompleteness Theorem\index{Incompleteness theorem}} 

Kurt G\"odel\index{Goed@G\"odel, Kurt} spent a lot of his life intensely 
studying Leibniz\index{Leibniz, Gottfried Wilhelm} and the idea of the scientia 
generalis\index{Scientia generalis}. He became so engrossed in the reading that 
he occasionally believed that the editions of Leibniz's works were deliberately 
withholding certain content when publishing. 

G\"odel\index{Goed@G\"odel, Kurt} revolutionised mathematical logic around 
1930 already in his dissertation supervised by Hans Hahn\index{Hahn, Hans} 
with the proof of the completeness theorem\index{Completeness theorem}. He 
succeeded shortly afterwards, in his work \enquote{On formally 
undecidable\index{Undecidability} propositions of Principia mathematica} from 
1931, to form the self-referential G\"odel sentence\index{Goed@G\"odel, Kurt} 
$Q$ within the Dedekind-Peano arithmetic\index{Dedekind-Peano arithmetic} 
which -- colloquially speaking -- is of the form 
\begin{quote}
I am not provable. 
\end{quote}
He wrote in the introduction to this work, which summarises the entire proof, 
that he was inspired by the paradox\index{Paradox} of 
Epimenides\index{Epimenides} for his proof.\footnote{See \cite[Vol. I, 
1931]{goedel}. In fact, G\"odel\index{Goed@G\"odel, Kurt} did not use 
Dedekind-Peano arithmetic, but the formal system in \enquote{Principia 
mathematica} \cite{russell}. For his proof, fragments of the Dedekind-Peano 
arithmetic are sufficient if they meet the Hilbert-Bernays-L\"ob criteria.} 

In it, G\"odel\index{Goed@G\"odel, Kurt} used the method of  
G\"odelisation\index{Goed@G\"odel, Kurt}, to number sequences of symbols in 
formulas and proofs of Dedekind-Peano 
arithmetic\index{Dedekind-Peano arithmetic} with the help of elementary number 
theory. There are different ways to define G\"odel numbers\index{Goed@G\"odel, 
Kurt} for sequences of symbols. One variant uses the sequence of prime numbers
\[
p_1=2, \; p_2=3, \; p_3=5, \ldots 
\]
and assigns a natural number to each symbol, for example 
\[
A \to 1, B \to 2, C \to 3, \ldots, H \to 8, \ldots 
\]
These numbers are used as exponents of the prime numbers. Thus, a sequence of 
symbols like BACH is given the G\"odel number\index{Goed@G\"odel, Kurt}
\[
2^2 \cdot 3^1 \cdot 5^3 \cdot 7^8= 8647201500.
\]
Similar numberings of expressions in formal languages\index{Formal language} 
already occur with Leibniz\index{Leibniz, Gottfried Wilhelm}. The G\"odel 
number\index{Goed@G\"odel, Kurt} of a formula $P$ is denoted by $\ulcorner P 
\urcorner$. The unique prime factor decomposition of natural numbers makes it 
possible to reconstruct the sequence of symbols BACH from the G\"odel 
number\index{Goed@G\"odel, Kurt} $8647201500$. 

Since proofs of statements of arithmetic are finite sequences of formulas, 
G\"odel\index{Goed@G\"odel, Kurt} had to show that the 
proof operations in the deductive system\index{Deductive system} of 
Dedekind-Peano arithmetic\index{Dedekind-Peano arithmetic} correspond to 
primitive recursive\index{Recursion} functions of the G\"odel 
numbers\index{Goed@G\"odel, Kurt}. This property allowed him to form the  
G\"odel sentence\index{Goed@G\"odel, Kurt} $Q$ that corresponds to the 
colloquial statement above in a syntactic\index{Syntax} way. 

How does G\"odel's\index{Goed@G\"odel, Kurt} proof work in detail? The proof 
of the theorem uses class signs (German \enquote{Klassenzeichen}), i.e., 
formulas of the form $\varphi(x)$ which depend on a free variable. These 
can also be assigned G\"odel numbers\index{Goed@G\"odel, Kurt} $\ulcorner 
\varphi(x) \urcorner$. It is possible to number the class signs with an index 
$k$ so that $\varphi_k(x)$ receives the G\"odel number\index{Goed@G\"odel, 
Kurt} $k$. Not all natural numbers appear as an index.
G\"odel\index{Goed@G\"odel, Kurt} showed that the mapping $\mathrm{diag}$, 
which sends the number $k$ to the G\"odel number\index{Goed@G\"odel, Kurt} 
$\ulcorner \varphi_k(k) \urcorner$ of $\varphi_k(k)$, is primitive 
recursive\index{Recursion}. This part of the proof is reminiscent of Cantor's 
diagonal argument\index{Cantor, Georg}. With these preconditions, he proved 
that the two-place provability relation\index{Provability} 
\[
\mathrm{Proof}(x,y) = \begin{cases} x \text{ encodes a 
proof } \cr \text{of the formula } \varphi_y(y) \end{cases} 
\]
is primitive recursive\index{Recursion} and can be formulated 
syntactically\index{Syntax} in the Dedekind-Peano 
arithmetic\index{Dedekind-Peano arithmetic}. Thus, using $\mathrm{diag}$, 
$\neg$ and $\mathrm{Proof}$, the class sign
\[
\varphi(y) = \forall x \, \neg \mathrm{Proof}(x,y)
= \begin{cases} \text{ For all } x \text{ it is not the case that } x  \cr 
\text{ encodes a proof of the formula } \varphi_y(y) 
\end{cases} 
\]
can be defined. Therefore, there is an $n$ with $\varphi=\varphi_n$ and the 
sought-after G\"odel sentence\index{Goed@G\"odel, Kurt} can be written as 
\[
Q=\varphi_n(n). 
\]
It now easily follows that $Q$ is provable if and only if $Q$ is not 
provable and vice versa. This results in a contradiction as in the 
paradox\index{Paradox} of Epimenides\index{Epimenides}.\footnote{For a modern 
presentation see \cite{hoffmann2017} and \cite[Ch. 4]{hoffmann2018}.}  

With this, G\"odel\index{Goed@G\"odel, Kurt} had proven the first 
incompleteness theorem\index{Incompleteness theorem}. It asserts that in any 
consistent\index{Consistency}\index{Contradiction-freeness} deductive 
system\index{Deductive system} which can formulate the elementary features of 
Dedekind-Peano arithmetic\index{Dedekind-Peano arithmetic}, there are 
sentences $Q$ which are neither provable nor refutable. Barkley 
Rosser\index{Rosser, Barkley} improved G\"odel's\index{Goed@G\"odel, Kurt} 
assumptions on consistency\index{Consistency}\index{Contradiction-freeness} in 
1936 and arrived at the assumptions stated here. A very short proof was found 
by Raymond Smullyan\index{Smullyan, Raymond}. He constructed a minimal formal  
language in which the self-referential G\"odel sentence\index{Goed@G\"odel, 
Kurt} $Q$ can be formulated.\footnote{See \cite{rosser1936,smullyan2013} and 
\cite[Ch. II]{manin}.}

In modern proofs, G\"odel's\index{Goed@G\"odel, Kurt} trick is 
usually explained with the diagonal lemma\index{Goed@G\"odel, Kurt} which 
was implicitly contained in G\"odel's\index{Goed@G\"odel, Kurt} proof. It was 
first explained by G\"odel\index{Goed@G\"odel, Kurt} in a lecture at the 
Institute for Advanced Study in Princeton in 1934 and he said that it was based 
on an idea by Rudolf Carnap\index{Carnap, Rudolf}.\footnote{See 
\cite{carnap1934}.} The diagonal lemma\index{Goed@G\"odel, Kurt} states that 
every class sign $P(x)$ has a sentence $Q$ as a fixed point, for which 
\[
Q \xLeftrightarrow{~~~} P(\ulcorner Q \urcorner) 
\] 
is provable. As class sign $P(x)$, with the help of the unary 
provability\index{Provability} predicate 
$\mathrm{Proof}(y)=\exists x \, \mathrm{Proof}(x,y)$, the predicate
\[
P=\neg \mathrm{Proof}
\]
of non-provability\index{Provability} can be chosen, i.e., $P(n)$ is true if and 
only if $n$ is not the G\"odel number\index{Goed@G\"odel, Kurt} of a provable 
formula $Q$. From the diagonal lemma now results the famous formula
\[
Q \xLeftrightarrow{~~~} \neg \mathrm{Proof}(\ulcorner Q \urcorner) 
\]
from which -- as in the original proof by G\"odel\index{Goed@G\"odel, Kurt} -- 
it follows that $Q$ is neither provable nor refutable.

G\"odel's\index{Goed@G\"odel, Kurt} first incompleteness 
theorem\index{Incompleteness theorem} has exciting implications beyond 
mathematics. Per Martin-L\"of\index{Martin-L\"of, Per} has connected the 
G\"odelian\index{Goed@G\"odel, Kurt} concept of 
incompleteness\index{Incompleteness theorem} with Kantian synthetic 
judgements\index{Synthetic judgement}.\footnote{See \cite{martin-loef1994}.} 
The incompleteness\index{Incompleteness theorem} corresponds to the information 
hidden in a synthetic\index{Synthetic judgement} judgement which cannot be 
obtained analytically\index{Analytical judgement}. The diagonal lemma also has
a variant in computer science in the fix point theorem of 
Kleene-Rogers\index{Kleene, Stephen}\index{Rogers, Hartley}. It 
implies, for example, that in any programming language there is a program which 
has its own code as output.

\section{G\"odel's\index{Goed@G\"odel, Kurt} Second  
Incompleteness Theorem\index{Incompleteness theorem}} 

The consistency\index{Consistency}, i.e., the 
contradiction-freeness\index{Contradiction-freeness}, of any mathematical 
theory that can describe Dedekind-Peano arithmetic\index{Dedekind-Peano 
arithmetic} can be expressed by the syntactic\index{Syntax} arithmetic formula 
\[
\mathrm{Con} = \neg \mathrm{Proof}(\ulcorner 0=1\urcorner).
\]
The second incompleteness theorem\index{Incompleteness theorem} of 
G\"odel\index{Goed@G\"odel, Kurt} states that the proof of the formula 
$\mathrm{Con}$ is not possible within the same theory. This theorem has a 
surprising history. A proof for it was not included in the submitted manuscript 
and was only announced at the end of G\"odel's\index{Goed@G\"odel, Kurt} 1931 
publication in a somewhat cryptic remark as a corollary for a second part that 
never appeared. Possibly, G\"odel\index{Goed@G\"odel, Kurt} had this passage 
inserted shortly before printing because he feared a competing publication by 
John von Neumann\index{Neumann, John von}, who wrote him a letter outlining 
this result on November 20th 1930. In his responses to von 
Neumann\index{Neumann, John von}, G\"odel\index{Goed@G\"odel, Kurt} presented 
the argument himself and a proof was later published by Hilbert\index{Hilbert, 
David} and Bernays\index{Bernays, Paul}.\footnote{John von Neumann's letters to 
G\"odel\index{Goed@G\"odel, Kurt} can be found in \cite[Vol. V]{goedel}. 
G\"odel's\index{Goed@G\"odel, Kurt} reply letters are partly only preserved in 
the form of handwritten drafts. See 
\cite{hilbert1939,vonplato2020,schuette2020} for further aspects.} The proof of 
the incompleteness theorem\index{Incompleteness theorem} only requires the 
contradiction-freeness\index{Contradiction-freeness}\index{Consistency} 
$\mathrm{Con}$. Therefore, with additional arguments, the implication 
\[
\mathrm{Con} \xRightarrow{~~} Q 
\]
can be proven. Since the G\"odel sentence\index{Goed@G\"odel, Kurt} $Q$ 
is unprovable within the theory, this also applies to $\mathrm{Con}$. The 
consistency\index{Consistency}\index{Contradiction-freeness} of the more 
general Zermelo-Fraenkel set theory\index{Zermelo-Fraenkel axioms} also 
cannot be derived as a formula, because it would -- by virtue of the standard 
model\index{Model theory} $\mathbb{N}$ -- imply the 
consistency\index{Consistency}\index{Contradiction-freeness} of the 
subsystem of arithmetic.

\medskip 
This was a spectacular result, as it dampened the hopes of 
Hilbert's\index{Hilbert, David} programme to establish the 
proof theory\index{Proof theory} and thus metamathematics on exact 
foundations. Neither Hilbert\index{Hilbert, David} nor 
G\"odel\index{Goed@G\"odel, Kurt} saw the second 
incompleteness theorem\index{Incompleteness theorem} as devastating 
for Hilbert's\index{Hilbert, David} programme though. Instead, both 
suggested independently shortly afterwards to modify the finite standpoint 
and continue Hilbert's\index{Hilbert, David} programme. 
Hilbert\index{Hilbert, David} wrote in 1934: 
\begin{quote}
With regard to this goal, I would like to emphasise that the temporarily arisen 
opinion, that certain recent results by G\"odel imply the infeasibility of my 
proof theory, has proven to be erroneous. That result indeed only shows that 
for the more advanced proofs of consistency one must exploit the finite 
standpoint in a sharper way than is necessary when considering elementary 
formalisms.\footnote{Machine translation of \cite{hilbert1934}. The lectures 
from 1930 and 1931 appeared in \cite[App. D]{ewaldsieg}, and the correspondence 
between G\"odel\index{Goed@G\"odel, Kurt} and Bernays\index{Bernays, Paul} in 
\cite[Vol. IV]{goedel}. More about G\"odel's\index{Goed@G\"odel, Kurt} 
theorems can be found in \cite{cheng2021,sieg,tapp}.}
\end{quote}
G\"odel expressed this in 1931 as follows:
\begin{quote}
It should be expressly noted that theorem XI (...) is in no contradiction to 
Hilbert's formalistic standpoint. For this only presupposes the 
existence of a proof of consistency conducted with finite means and it would be 
conceivable that there are finite proofs that cannot be represented in $P$ 
[G\"odel's formal system] (...).\footnote{Machine translation of \cite[Vol. I, 
1931]{goedel}.}
\end{quote}
Even many years later, around 1961-1962, 
G\"odel\index{Goed@G\"odel, Kurt}, according to Gerald E. Sacks\index{Sacks, 
Gerald}, occasionally remarked that Hilbert's\index{Hilbert, David} 
programme was still open. 

The above-mentioned requirements for $\mathrm{Con}$ as 
a syntactic\index{Syntax} arithmetic formula, which also apply in all 
non-standard models\index{Non-standard model}, are very strong. Therefore, 
attempts have been made to express 
consistency\index{Consistency}\index{Contradiction-freeness} 
either through a different formula with a deviating proof predicate or through 
an infinite schema of formulas. As a second possibility, 
Gentzen\index{Gentzen, Gerhard} and Ackermann\index{Ackermann, Wilhelm} used 
transfinite induction\index{Transfinite induction} up to the ordinal 
number\index{Ordinal number} $\varepsilon_0$. In the even more general case of 
set theory\index{Set theory}, 
consistency\index{Consistency}\index{Contradiction-freeness} can be 
demonstrated in a third way using Grothendieck universes\index{Universe} which 
are based on additional axioms about inaccessible cardinal 
numbers\index{Cardinal number}. The price for consistency 
proofs\index{Consistency}\index{Contradiction-freeness} of the last two types 
are new axioms which again 
generate incompleteness\index{Incompleteness theorem}. In the first case, the 
concept of consistency\index{Consistency}\index{Contradiction-freeness} is 
different.\footnote{See \cite{artemov2019}, \cite[Sec. 5]{cheng2021}, 
\cite{gentzen} and \cite{stepien2017}.}

Vladimir Voevodsky\index{Voevodsky, Vladimir} expressed doubts about the 
consistency\index{Consistency}\index{Contradiction-freeness} of mathematics 
in a remarkable lecture in 2010. He saw type theory\index{Type theory} and 
intuitionism\index{Intuitionism} as a possible way out. 
This happened around the same time as Edward Nelson\index{Nelson, Edward} from 
Princeton -- a committed opponent of Platonism\index{Platonic idealism} 
and semantics\index{Semantics} -- unsuccessfully tried to refute the 
consistency\index{Consistency}\index{Contradiction-freeness} of 
arithmetic\index{Dedekind-Peano arithmetic}.
Hilbert's\index{Hilbert, David} programme is, all things considered, 
indeed still partially open and new insights into it are conceivable.

\section{The G\"odel\index{Goed@G\"odel, Kurt}-Lucas\index{Lucas, 
John}-Penrose\index{Penrose, Roger} Argument}

Most people are probably convinced that human 
qualities and qualia\index{Qualia} such as empathy, consciousness, intuition,
understanding, feelings and other factors distinguish our thinking from 
a Turing machine\index{Turing machine}. However, such opinions are 
far from a scientifically tenable thesis. The 
opposite idea is that human thinking and behaviour as 
a whole is equivalent to a universal Turing machine\index{Turing machine}. 
This position is called the mechanistic thesis\index{Mechanistic thesis}.

Such questions were what Alan Turing\index{Turing, Alan} wanted to decide. He 
formulated the Turing test\index{Turing test}, a hypothetical experiment, 
in which targeted questions are used to determine whether an 
interlocutor is a human or a computer. John Searle's\index{Searle, John}  
thought experiment of the Chinese Room has a similar 
approach.\footnote{See \cite{turing1950,searle1992}.}

In G\"odel's\index{Goed@G\"odel, Kurt} estate, an unpublished 
speech manuscript for the Gibbs Lecture of the American Mathematical Society at 
Christmas 1951 was found, from which it became clear that he believed in the 
superiority of human thinking. He outlined an argument that used his 
incompleteness theorem\index{Incompleteness theorem}:
\begin{quote}
Evidently no well-defined system of correct axioms can comprise all objective 
mathematics, since the proposition which states the consistency of the system 
is true, but not demonstrable in the system. However, as to subjective 
mathematics, it is not precluded that there should exist a finite rule 
producing all its evident axioms ... The assertion, however, that they are all 
true could at most be known with empirical certainty, on the basis of a 
sufficient number of instances or by other inductive inferences. If 
it were so, this would mean that the human mind (in the realm of pure 
mathematics) is equivalent to a finite machine that, however, is unable to 
understand completely its own functioning ... So the following disjunctive 
conclusion is inevitable: Either mathematics is incompletable in this sense, 
that its evident axioms can never be comprised in a finite rule, that is to 
say, the human mind (even within the realm of pure mathematics) infinitely 
surpasses the powers of any finite machine, or else there exist absolutely 
unsolvable diophantine problems of the type specified.\footnote{See 
\cite[Vol. III, 1951]{goedel}.}
\end{quote}

G\"odel\index{Goed@G\"odel, Kurt} derived here an alternative 
between the opposite of the mechanistic thesis\index{Mechanistic thesis}
and the existence\index{Existence} of certain arithmetic propositions. 
Moreover, he considered the possibility that the human brain could be 
equivalent to a Turing machine\index{Turing machine} that does not fully 
understand its own functioning. At the end of his speech, 
G\"odel\index{Goed@G\"odel, Kurt} tried to substantiate the view of Platonic 
idealism\index{Platonic idealism} with similar arguments. 

The exact proof argument is not given precisely in 
G\"odel's\index{Goed@G\"odel, Kurt} text. The idea is to consider the deductive 
system\index{Deductive system} $S$ of Dedekind-Peano 
arithmetic\index{Dedekind-Peano arithmetic}. 
G\"odel's\index{Goed@G\"odel, Kurt} first incompleteness 
theorem\index{Incompleteness theorem} provides a method to generate a more 
complete deductive system\index{Deductive system} $S'$ by adjoining the 
G\"odel sentence\index{Goed@G\"odel, Kurt} $Q_S$ of $S$. Alternatively, 
the statement $\mathrm{Con}_S$ of the  
consistency\index{Consistency}\index{Contradiction-freeness} of $S$ 
can be adjoined and a new G\"odel sentence\index{Goed@G\"odel, Kurt} $Q_{S'}$ 
in $S'$ results. The extension $S'$ is called consistency 
extension\index{Consistency}\index{Contradiction-freeness}. Such 
consistency extensions\index{Consistency}\index{Contradiction-freeness} can 
be repeated through infinite transfinite\index{Transfinite induction} 
iterations over ordinal numbers\index{Ordinal number}.\footnote{See 
\cite{feferman1962,turing1939}.} Among all the true theorems in the standard 
model\index{Model theory} in the Dedekind-Peano arithmetic\index{Dedekind-Peano 
arithmetic} there are still those after such iterations that are not provable. 
There are now two options. Either humans can prove all true theorems of the 
Dedekind-Peano arithmetic\index{Dedekind-Peano arithmetic}. In this case, 
G\"odel's first incompleteness theorem\index{Incompleteness theorem} 
implies that the human brain is infinitely superior to any Turing 
machine\index{Turing machine} and the mechanistic thesis\index{Mechanistic 
thesis} is false. Otherwise, there is a true arithmetic theorem that is 
absolutely unsolvable, as G\"odel\index{Goed@G\"odel, Kurt} puts it.

John Randolph Lucas\index{Lucas, John} tried in 1961 to develop a similar 
argument in his own way. Roger Penrose\index{Penrose, Roger} also published 
variants of G\"odel's\index{Goed@G\"odel, Kurt} and Lucas'\index{Lucas, John} 
arguments in two books from 1989 onwards. Lucas\index{Lucas, John} and 
Penrose\index{Penrose, Roger} went beyond G\"odel\index{Goed@G\"odel, Kurt} and 
claimed that the mechanistic thesis\index{Mechanistic thesis} is always false. 
We do not want to fully reproduce the exact arguments here, as there are gaps 
in the proof that were uncovered by the two logicians Martin Davis\index{Davis, 
Martin} and Solomon Feferman\index{Feferman, Solomon} and could never be 
eliminated despite many debates. This was merely criticism of the proof itself, 
as the two aforementioned experts, like G\"odel\index{Goed@G\"odel, Kurt}, 
Lucas\index{Lucas, John} and Penrose\index{Penrose, Roger}, were convinced of 
the superiority of human thought.\footnote{See 
\cite{lucas1961,penrose1989,penrose1994,davis1993, feferman1996}.} 

The thought processes of Penrose\index{Penrose, Roger} and 
Lucas\index{Lucas, John} are occasionally used in popular science literature in 
an unscientific way as arguments to refute the mechanistic 
thesis\index{Mechanistic thesis} and to draw all possible conclusions from it. 
Without proof, however, such claims remain mere speculation.

In this context, Roger Penrose\index{Penrose, Roger} conjectured that the human 
brain needs quantum mechanical mechanisms to exceed the limit of 
Turing computability\index{Computability}. At the same time, this could 
explain the phenomena of consciousness and free will. He argued that the 
collapse of the wave function in quantum mechanics plays a crucial role in 
this. However, quantum mechanical effects in the functioning of the brain have 
never been proven to this day. The free will of humans is at least questionable 
due to neuroscientific experiments, as human decisions are often preceded by 
signals in the brain by fractions of a second. 

\section{Intuitionism\index{Intuitionism}}

Proofs of classical logic occasionally use proofs by contradiction and thus the 
law of excluded middle, called tertium non datur in Latin. This method generates 
mathematical objects which cannot be described constructively anymore. Leopold 
Kronecker\index{Kronecker, Leopold} was one of the first who pointed this out. 
In his papers and opinions, he rejected the upcoming abstract mathematics of 
Cantor\index{Cantor, Georg}, Dedekind\index{Dedekind, Richard}, and 
Hilbert\index{Hilbert, David} in which uncountable infinite sets like the real 
numbers were looked at even if the single elements are in general not 
accessible. 

Intuitionism\index{Intuitionism} is a movement in mathematics coined by Hermann 
Weyl\index{Weyl, Hermann}, Luitzen Egbertus Jan Brouwer\index{Brouwer, Luitzen 
Egbertus Jan} and Arend Heyting\index{Heyting, Arend} which takes such doubts 
into account. It has a related manifestation in 
constructivism\index{Constructivism} and is an antirealistic\index{Realism} 
position in philosophical terms, at least in its extreme variants. Both views 
result from the desire for mathematics to be carried out through a 
comprehensible thought process of humans. In it, the existence\index{Existence} 
of mathematical objects should be verifiable and the truth\index{Truth} of 
statements should be exclusively obtainable by constructive 
proofs\index{Provability}. The use of the law of excluded middle is rejected and 
the axiom of choice is considered a tautology.\footnote{See the appendix in 
\cite{beeson1985} for the history of these ideas.} This attitude is historically 
explainable from the antinomies of set theory at the beginning of the 20th 
century and led to the development of intuitionistic\index{Intuitionism} logic. 
 
The law of excluded middle states that  
\[
\vdash A \vee \neg A   
\]
for every statement $A$. This is equivalent to the equivalent judgements 

\begin{minipage}{\textwidth}
\begin{itemize}
\item $\neg \neg A \vdash A$
\item $\vdash ((A \Rightarrow B)\Rightarrow A) \Rightarrow A$ \; 
(Peirce's\index{Peirce, Charles} law). 
\end{itemize}
\end{minipage}

in the calculus of natural deduction\index{Calculus of natural deduction}, in 
which $A$ and $B$ are arbitrary. Both inference rules therefore do not apply in 
intuitionistic\index{Intuitionism} logic. The reversal 
\[
A \vdash \neg \neg A
\]
applies without additional assumptions and follows from the calculus with the 
introduction rule for $\neg$ (with $\neg A$ instead of $A$ and $\Gamma=B=A$).

Brouwer\index{Brouwer, Luitzen Egbertus Jan} was a charismatic person with a 
varied life. He already dealt with intuitionism\index{Intuitionism} in his 
dissertation. After a creative phase in which he was able to prove deep results 
of topology\index{Topology}, such as the invariance of dimension and his fixed 
point theorem, he turned to the foundations of mathematics around 1913. The 
famous foundational dispute between 1921 and 1928 led to a bitter dispute with 
Hilbert\index{Hilbert, David}, who had little appreciation for 
intuitionism\index{Intuitionism}. This was followed by the exclusion of 
Brouwer\index{Brouwer, Luitzen Egbertus Jan} from the editorial board of the 
Mathematical Annals. Arend Heyting\index{Heyting, Arend} developed 
intuitionistic\index{Intuitionism} propositional logic and 
Heyting arithmetic\index{Heyting arithmetic} as the 
intuitionistic\index{Intuitionism} counterpart of Dedekind-Peano 
arithmetic\index{Dedekind-Peano arithmetic}. 
Hermann Weyl\index{Weyl, Hermann} adapted the ideas of 
intuitionism\index{Intuitionism} in his own way and also acted against 
Hilbert\index{Hilbert, David} during the foundational dispute.\footnote{See 
\cite{bridges1984,bridges1994,feferman1978,schuette2020} and 
\cite[Appendix]{beeson1985} on constructivism and intuitionism.} 

The ideas behind intuitionism\index{Intuitionism} and some metamathematical 
positions of Brouwer\index{Brouwer, Luitzen Egbertus Jan} were controversial 
from the start. His concern was that humans as creative subjects construct 
the objects of mathematics in their mind and work out the 
proofs\index{Provability} of statements. In particular, he thought a lot about 
the continuum and wanted to understand it better with the help of so-called 
choice sequences. Many ideas of Brouwer\index{Brouwer, Luitzen Egbertus Jan} 
were influential and have been further developed later. Especially Stephen C. 
Kleene\index{Kleene, Stephen}, George Kreisel\index{Kreisel, Georg}, John R. 
Myhill\index{Myhill, John} and Richard E. Vesley\index{Vesley, Richard} have 
further elaborated intuitionistic\index{Intuitionism} mathematics in the 1960s. 
In this way, less dogmatic variants were created und 
intuitionistic\index{Intuitionism} logic got the appropriate importance.

All logical connections and quantifiers are considered as conditions to be 
verified in intuitionism\index{Intuitionism}. For example, a statement of the 
form $A \wedge B$ is provable\index{Provability} if there is a proof of $A$ and 
a proof of $B$. All other logical operations can be interpreted in a similar 
form. This concept is called Brouwer-Heyting-Kolmogorov interpretation
\index{Brouwer-Heyting-Kolmogorov interpretation} since it was discovered by 
Brouwer\index{Brouwer, Luitzen Egbertus Jan}, Heyting\index{Heyting, Arend}, 
and Kolmogorov\index{Kolmogorov, Andrey}. If proofs are described using the 
$\lambda$-calculus, the Brouwer-Heyting-Kolmogorov interpretation
\index{Brouwer-Heyting-Kolmogorov interpretation} provides a connection 
between proofs and $\lambda$-expressions, and is related to the Curry-Howard 
correspondence\index{Curry-Howard correspondence}.

G\"odel\index{Goed@G\"odel, Kurt} noted that while 
intuitionism\index{Intuitionism} appears as a restriction of classical 
mathematics, it can also be seen as an extension. To this end, he constructed 
-- with the help of double negation -- an interpretation\index{Interpretation} 
of first-order classical logic in intuitionistic logic\index{Intuitionism} 
using a variant of the recursive assignment $\varphi \mapsto 
\varphi^\mathrm{n}$ in Fig.~\ref{fig:double_negation}.

\begin{figure}[ht!]
\begin{align*}
\varphi^\mathrm{n} &= \neg \neg \varphi, \text{ if } \varphi 
\text{ is atomic}\cr 
(\varphi \wedge \psi)^\mathrm{n} & = \varphi^\mathrm{n} \wedge 
\psi^\mathrm{n} \cr 
(\varphi \vee \psi)^\mathrm{n} & = \neg(\neg \varphi^\mathrm{n} \wedge 
\neg \psi^\mathrm{n}) \cr 
(\varphi \Rightarrow \psi)^\mathrm{n} & = (\varphi^\mathrm{n} \Rightarrow 
\psi^\mathrm{n}) \cr 
(\neg \varphi)^\mathrm{n} & = \neg \varphi^\mathrm{n} \cr
(\forall x \, \varphi)^\mathrm{n} & = \forall x \, \varphi^\mathrm{n} \cr
(\exists x \, \varphi)^\mathrm{n} & = \neg \forall x \, \neg  
\varphi^\mathrm{n}. 
\end{align*}
\caption{\label{fig:double_negation}The rules of double negation for all 
logical operations.}
\end{figure}

G\"odel\index{Goed@G\"odel, Kurt} applied this 
interpretation\index{Interpretation} to arithmetic and showed that 
every contradiction derivable in 
Dedekind-Peano arithmetic\index{Dedekind-Peano arithmetic} 
is also derivable in Heyting arithmetic\index{Heyting arithmetic}. 
Consequently, Heyting arithmetic\index{Heyting arithmetic} 
and Dedekind-Peano arithmetic\index{Dedekind-Peano arithmetic} are 
equiconsistent\index{Consistency}\index{Contradiction-freeness}.\footnote{ See 
\cite[Vol. I, 1933e]{goedel}.} Moreover, G\"odel's 
incompleteness theorems\index{Incompleteness theorem}
\index{Goed@G\"odel, Kurt} with their consequences also apply in an 
intuitionistic\index{Intuitionism} version. This makes it clear that 
intuitionism\index{Intuitionism} can only alleviate concerns 
about classical logic to a limited extent. 

\section{Constructivism\index{Constructivism} and the Question of Existence}

An influential variant of intuitionism\index{Intuitionism} is 
constructivism\index{Constructivism}. The existence\index{Existence} 
of mathematical objects is a crucial concept in both currents. 
Behind this is the idea that mathematical objects, whose 
existence\index{Existence} is asserted in a theorem, must always be concretely 
constructed. This implies that proofs by contradiction and the axiom of choice 
in its classical formulation are not permitted. It is common to regard the 
axiom of choice as a tautology by interpreting the specification of any number 
of non-empty sets as the stronger assumption of specifying just as many 
explicit elements. 

The connection of constructivistic\index{Constructivism} and 
intuitionistic\index{Intuitionism} ways of thought in their original form 
seems to contradict the existence\index{Existence} of abstract mathematical 
objects, similar to Benacerraf's\index{Benacerraf, Paul} 
dilemma\index{Benacerraf's dilemma}. On the one hand, these movements are of 
antirealist\index{Realism} resp. nominalistic\index{Nominalism} nature, as 
we already remarked, since they make the existence\index{Existence} depend on 
mathematically acting humans. On the other hand, they treat the 
existence\index{Existence} of some of the constructed objects as preferred 
compared to classical approaches, which is rather a Platonic\index{Platonic 
idealism} viewpoint. The concept of mathematical 
construction\index{Constructivism} free of ontological\index{Ontology} 
attributes does not seem to be problematic though. 

What exactly construction means is interpreted differently in different 
variants. An unproblematic requirement for constructability is that concrete 
objects should be named for mathematical existence claims\index{Existence} to 
generate evidence. A good example of this is the $4$-squares theorem, which 
provides a representation 
\[
n=a^2+b^2+c^2+d^2
\]
for every natural number $n$. For a given $n$, the explicit specification of an 
algorithm\index{Algorithm} that calculates the quadruple $(a,b,c,d)$ depending 
on $n$ is such evidence. The proof of this theorem can be conducted 
constructively\index{Constructivism} in this sense.\footnote{See 
\cite[Kap. 9]{mspiontkow2011} for a constructive proof.}

Hilbert\index{Hilbert, David} stirred up the movements of 
intuitionism\index{Intuitionism} and constructivism\index{Constructivism} 
through his attitude towards mathematics, which was particularly manifested in 
formalism\index{Formalism}. In a letter dated December 29th 1899 to 
Frege\index{Frege, Gottlob}, he wrote:
\begin{quote}
If the arbitrarily set axioms do not contradict each other with all their 
consequences, then they are true, then the things defined by the axioms exist. 
This is for me the criterion of truth and 
existence\index{Existence}.\footnote{Machine translation of \cite{frege1976}.}
\end{quote}
Hilbert\index{Hilbert, David} experienced strong counter-reactions due to his 
clear views. Especially Brouwer\index{Brouwer, Luitzen Egbertus Jan} and Oskar 
Becker\index{Becker, Oskar} as well as occasionally Weyl\index{Weyl, Hermann} 
fought against Hilbert\index{Hilbert, David}. In retrospect, it is not 
necessary to assume that Hilbert\index{Hilbert, David} meant the word 
existence\index{Existence} in a Platonic\index{Platonic idealism} or similar 
sense. He was probably more concerned with being able to work with such objects 
without contradiction, even if they were not 
constructively\index{Constructivism} generated. This is a legitimate 
standpoint, as Hilbert\index{Hilbert, David} had clearly recognised that the 
consistency\index{Consistency}\index{Contradiction-freeness} of mathematics is 
the real problem and not the supposedly missing existence\index{Existence} of 
its objects. This is substantiated with the fact that 
Hilbert\index{Hilbert, David} later had a tendency to constructive methods in 
mathematics in his finitistic proof theory\index{Proof theory}. 

It seems that the existence of mathematical objects, in particular of infinite 
sets, is a steady source of misunderstandings. Many ideological fights in the 
philosophy of mathematics suffer from the claim of set-theoretical 
semantics\index{Semantics} in a too strong ontological\index{Ontology} sense, 
although reasonable realistic\index{Realism} positions would suffice. In fact, 
as Hilbert\index{Hilbert, David} had realised, the 
consistency\index{Consistency}\index{Contradiction-freeness} of mathematics is 
the real problem. 

What are the consequences of the constructivist\index{Constructivism} 
viewpoint? Let us first consider the real numbers. Any $a \in \mathbb{R}$ is 
constructively\index{Constructivism} given by a Cauchy sequence $(a_n)_{n \in 
\mathbb{N}}$ that converges to $a$ and fulfils a suitable error estimate. At 
least from a classical point of view, the trichotomy 
\[
a<0 \text{ or } a=0 \text{ or } a>0 
\]
certainly applies. However, for a given Cauchy sequence 
$(a_n)_{n \in \mathbb{N}}$, it is not constructively\index{Constructivism} 
decidable which of the three cases applies. This phenomenon is well-known in 
numerical mathematics. The step function 
\[
f(x)=\begin{cases} 0 & \text{ for } x<0 \cr 1 & \text{ for } x\ge 
0 \end{cases}  
\]
is -- in other words -- not meaningfully definable in a 
constructive\index{Constructivism} sense. This implies, for example, that there 
is no characteristic function $\chi_\mathbb{Q}$ of the subset of rational 
numbers in the real numbers. Brouwer\index{Brouwer, Luitzen Egbertus Jan} had 
already seen this and he proved that in his constructive\index{Constructivism} 
formalism every total real function is continuous and uniformly continuous on 
compact intervals. Not least because of such consequences, the 
intuitionistic\index{Intuitionism} and the constructive\index{Constructivism} 
view have not prevailed widely to this day.

There are mathematical propositions for which no constructive proof is known. 
An example of this is K\"onig's lemma, which was found by D\'enes 
K\"onig\index{Koen@K\"onig, D\'enes} in 1936.\footnote{See \cite[Ch. 
IV]{beeson1985}.} It states that a connected graph, in which only
infinitely many edges spring from each node, is exactly infinite, 
when there is an infinite path in it that avoids itself. In 
the special case of a finitely branched, infinite tree, we are looking for an 
infinite path. The proof of this lemma in classical logic is not 
very difficult. It determines at each step by means of complete 
induction\index{Complete induction} a new edge, over which infinitely 
many other nodes can be reached. The law of excluded middle is used to 
show that such exists.\footnote{The proof implements a 
weak form of the axiom of choice. A constructive proof outside of 
special cases is not known.}

In the time after Brouwer\index{Brouwer, Luitzen Egbertus Jan}, the schools 
of Andrey A. Markov\index{Markov, Andrey}, Nicolai A. Shanin\index{Shanin, 
Nicolai}, and Errett Bishop\index{Bishop, Errett} further developed 
constructive\index{Constructivism} mathematics. This resulted in new 
approaches that can handle large parts of mathematics. The search for those 
theorems of mathematics that have constructive proofs provides extremely 
valuable insights.\footnote{See \cite{bridges1984,bridges1994,feferman1978}.} 
Even though intuitionism\index{Intuitionism} and 
constructivism\index{Constructivism} still play a subordinate role, they 
are nevertheless significant as they form the foundations of dependent type 
theory\index{Type theory} and underlie the internal syntax\index{Syntax} of 
categories\index{Category theory} as we will see.

\chapter{Category Theory\index{Category theory}}

Category theory\index{Category theory} is -- alongside set theory\index{Set 
theory} and type theory\index{Type theory} -- one of the three possible 
foundations of mathematics. It incorporates the structural mathematical 
thinking in an abstract way and is well-suited for the justification of a 
mathematical semantics\index{Semantics}. Richard Dedekind\index{Dedekind, 
Richard} has in a certain sense founded structuralism\index{Structuralism} in 
mathematics. He  realised that mathematical objects, such as numbers, can have 
many isomorphic\index{Isomorphism} set realisations. It is only their 
operations and structural properties that are unique. 
Dedekind\index{Dedekind, Richard} has axiomatically characterised the natural 
numbers in his book \enquote{What are and what should the numbers be?} as 
chains and proved their uniqueness up to isomorphism\index{Isomorphism} with 
the recursion theorem\index{Recursion}.\footnote{Contained in 
\cite{dedekind,ms2023}.}

Emmy Noether\index{Noether, Emmy} has significantly contributed to establishing 
these ideas of Dedekind\index{Dedekind, Richard} in modern 
algebraand in topology\index{Topology}, where they were 
particularly fruitful.\footnote{On Emmy Noether\index{Noether, Emmy} and her 
influence see \cite{ms2027}.} In the form of category theory\index{Category 
theory}, structuralism\index{Structuralism} was further developed by the 
topologists Saunders MacLane\index{MacLane, Saunders} and Samuel 
Eilenberg\index{Eilenberg, Samuel} from about 1945 onwards and then found its 
way into the works of Alexander Grothendieck\index{Grothendieck, Alexander}, 
Daniel Quillen\index{Quillen, Daniel} and others. Especially the development of 
model and path categories\index{Model category}\index{Path category} has shown 
the importance of category theory\index{Category theory} in homotopy 
theory\index{Homotopy theory} and in related algebraic areas, because they 
conceptualise the role of (weak) equivalences\index{Equivalence} in a 
conceptual way. An important role in current research is played by higher 
categories\index{Category theory} and infinity categories\index{Infinity 
category} as further developments. From model categories\index{Model category} 
infinity categories\index{Infinity category} can be obtained through 
localisation.\footnote{See \cite{cisinski2015,dwyerkan,groth2015,lurie2009,
quillen}.}

\section{Categories\index{Category theory}}

A category\index{Category theory} $\mathcal{C}$ consists of a collection 
$\mathrm{Ob}(\mathcal{C})$ of objects that share common structural features of 
a given mathematical concept. Between any two objects $A$ and $B$ of 
$\mathcal{C}$ there is a collection $\mathrm{Hom}_\mathcal{C}(A,B)$ of 
structure-preserving morphisms, which are represented in the form of arrows  
\[
A \xlongrightarrow{~~ f ~~} B.  
\]
For each object $A$ in $\mathcal{C}$ there is a 
natural morphism $\mathrm{id}_A \colon A \longrightarrow A$, which 
is called identity\index{Identity}\index{Equality}. 

The collections of objects and morphisms do not necessarily have to be 
sets in a given set theory\index{Set theory}, even though this is often 
required. In such a case, the category\index{Category theory} $\mathcal{C}$ is 
called small. We will try as far as possible to consider category 
theory\index{Category theory} as an alternative basis for mathematics and to 
avoid sets when studying categories\index{Category theory}.

Publications in category theory\index{Category theory} are usually full of 
diagrams of arrows. The corresponding morphisms are chained by composition (see 
Fig.~\ref{fig:composition}).

\begin{figure}[ht!]
\begin{tikzpicture}
\node (a) at (0,0) {$A$};
\node (b) at (3,0) {$B$};
\node (c) at (6,0) {$C$};
\draw[->] (a) to node[scale=.7] [above] {$f$} (b);
\draw[->] (b) to node[scale=.7] [above] {$g$} (c);
\draw[->] (a) to [bend right] node[scale=.7] [above] {$g \circ f$} (c);
\end{tikzpicture}
\caption{\label{fig:composition}Composition $g \circ f$ of two morphisms $f,g$.}
\end{figure}

For this, the associative law -- similar to multiplication -- applies: 
\[
h \circ (g \circ f)= (h \circ g) \circ f.
\]
A morphism $f \colon A \longrightarrow B$ is called an isomorphism if there 
is an inverse morphism $f^{-1} \colon B \longrightarrow A$ that 
fulfils $f^{-1} \circ f=\mathrm{id}_A$ and $f \circ f^{-1}=\mathrm{id}_B$.

Arrows can be chained in different ways. A typical situation is a commutative 
diagram (see Fig.~\ref{fig:commutative_diagram}) for which 
\[
i \circ h=g \circ f 
\]
applies, i.e., it does not matter which path the square diagram is traversed. 
This results in a unique new arrow, which is drawn dashed as a diagonal.

\begin{figure}[ht!]
\begin{tikzpicture}
  \matrix (m) [matrix of math nodes,row sep=3em,column sep=4em,minimum 
width=2em]
  {
     A & B \\
      C & D \\};
  \path[-stealth]
    (m-1-1) edge node [left] {$h$} (m-2-1)
            edge node [above] {$f$} (m-1-2)
    (m-2-1.east|-m-2-2) edge node [below] {$i$}
            node [above] {} (m-2-2)
    (m-1-2) edge node [right] {$g$} (m-2-2)
    (m-1-1) edge [dashed] node [below] {} (m-2-2);
\end{tikzpicture}
\caption{\label{fig:commutative_diagram}A commutative diagram is a square of 
four objects and four morphisms as edges. The composition of composable edges 
coincide (with the dashed diagonal morphism).}
\end{figure}

There are many different categories\index{Category theory} in 
mathematics. Important examples are the categories\index{Category theory} 
$\mathbf{Set}$ and $\mathbf{Top}$. In $\mathbf{Set}$ the objects are sets
and the morphisms are the set-theoretic\index{Set theory} mappings. In 
$\mathbf{Top}$, the objects are the topological spaces\index{Topology} and 
the morphisms are the continuous mappings. A special category is assigned to 
each topological space\index{Topology} $X$, namely the category 
$\mathbf{Off}(X)$ of open sets\index{Category theory} in $X$. The objects 
therein are the open sets $U$ in $X$ and the morphisms are the inclusions 
between them. Another important category\index{Category theory} $\mathbf{Ch}$ 
consists of chain complexes of abelian groups. It contains objects such as the 
singular chain complexes $\mathbb{Z}\mathrm{Sing}_\bullet(X)$ of topological 
spaces\index{Topology}.

\medskip
It is possible, in a given category\index{Category theory} 
$\mathcal{C}$, to reverse all arrows. This results in the opposite 
category\index{Category theory} $\mathcal{C}^{\mathrm{op}}$. 

\section{Groupoids\index{Groupoid}}

A very simple category\index{Category theory} is assigned to each group $G$. It 
has only one object $\ast$ and for each $g \in G$ an arrow, which is also 
denoted by $g$ (see Fig.~\ref{fig:groupoid}).

\begin{figure}[ht!]
\begin{tikzpicture}
\node (A) at (0,0) {$\ast$};
\path (A) edge [anchor=center,loop above] node {$g$} (A);
\path (A) edge [anchor=center,loop right] node {$g^{-1}$} (A);
\end{tikzpicture}
\caption{\label{fig:groupoid}The category\index{Category theory} of a group 
has only one object $\ast$ but one morphism for each element $g$ of $G$.}
\end{figure}

This category\index{Category theory} is a groupoid\index{Groupoid}, i.e., 
a category\index{Category theory}, in which every morphism $f$ is 
an isomorphism\index{Isomorphism}. A groupoid\index{Groupoid} is a 
far-reaching generalisation of a group. 

\medskip 
An interesting groupoid\index{Groupoid}, which does not come from a 
group in this form, is the 
fundamental groupoid\index{Fundamental groupoid} $\Pi_1(X)$ for a 
topological space\index{Topology} $X$. It is defined as the 
category\index{Category theory}, whose objects are the points of $X$ and 
whose morphisms between $a$ and $b$ are the 
homotopy classes\index{Homotopy theory} of paths from $a$ to $b$. From 
this fundamental groupoid\index{Fundamental groupoid} $\Pi_1(X)$, both 
the set $\pi_0(X)$ of connected components of $X$ and the 
fundamental group\index{Fundamental group} $\pi_1(X,*)$ for every base point 
$*$ can be reconstructed. 

This groupoid\index{Groupoid} contains the homotopy classes\index{Homotopy 
theory} of all paths in $X$ and thus more information than the fundamental 
group\index{Fundamental group} $\pi_1(X,*)$. To go from homotopy 
classes\index{Homotopy theory} of paths to all paths, a further step of 
abstraction is necessary, which will lead us from categories\index{Category 
theory} to infinity categories\index{Infinity category}.

\section{Functors\index{Functor}}

Like all mathematical objects, categories\index{Category theory}
like $\mathbf{Set}$ or $\mathbf{Top}$ are not unique. There are variants of 
$\mathbf{Set}$ which depend on chosen axioms. For example, the 
sets in $\mathbf{Set}$ are often restricted to subsets of a fixed 
Grothendieck universe\index{Universe} $U$. Some significant 
theorems of mathematics assert the equivalence\index{Equivalence} of 
two categories\index{Category theory}. We therefore need suitable 
tools to compare categories\index{Category theory} with each other.
A covariant functor\index{Functor} between two categories is an arrow 
\[
F \colon \mathcal{A} \xlongrightarrow{~~~~~~} \mathcal{B}  
\]
which sends objects $A$ in $\mathcal{A}$ to objects $F(A)$ in $\mathcal{B}$ 
and morphisms $f \colon A_1 \longrightarrow A_2$ in $\mathcal{A}$ to 
morphisms 
\[
F(f) \colon F(A_1) \xlongrightarrow{~~~~~~} F(A_2). 
\]
A contravariant functor\index{Functor} reverses the arrows, i.e., it holds 
\[
F(f) \colon F(A_2) \xlongrightarrow{~~~~~~} F(A_1).
\]
An example is the covariant functor\index{Functor} 
$F \colon \mathbf{Top} \longrightarrow \mathbf{Set}$, which sends a 
topological space\index{Topology} to the underlying set and continuous 
maps to the underlying map, i.e., it ignores continuity. Other examples are the 
covariant functors of homology \index{Homology group} and homotopy 
groups\index{Homotopy group}, which represent functors\index{Functor} 
$\mathbf{Top} \longrightarrow \mathbf{Gr}$ in 
the category of groups\index{Category theory}. At the level 
of the underlying chain complexes of abelian groups, which Emmy 
Noether\index{Noether, Emmy} already considered, the functor\index{Functor} 
$\mathbf{Top} \longrightarrow \mathbf{Ch}$ can be considered, which assigns to 
each topological space\index{Topology} $X$ the singular 
chain complex $\mathbb{Z}\mathrm{Sing}_\bullet(X)$.

Functors\index{Functor} can be compared with each other. A natural 
transformation\index{Natural transformation} $T$ between two 
functors\index{Functor} $F$ and $G$ assigns to each object $A$ of $\mathcal{A}$ 
a morphism $T_A \colon F(A) \longrightarrow G(A)$, so that the diagram (see  
Fig.~\ref{fig:nat_transformation}) in the category\index{Category theory} 
$\mathcal{B}$ commutes.

\begin{figure}[ht!]
\begin{tikzpicture}
  \matrix (m) [matrix of math nodes,row sep=3em,column sep=4em,minimum 
width=2em]
  {
     F(A) & F(B) \\
      G(A) & G(B) \\};
  \path[-stealth]
    (m-1-1) edge node [left] {$T_A$} (m-2-1)
            edge node [above] {$F(f)$} (m-1-2)
    (m-2-1.east|-m-2-2) edge node [below] {$G(f)$}
            node [above] {} (m-2-2)
    (m-1-2) edge node [right] {$T_B$} (m-2-2);
\end{tikzpicture}
\caption{\label{fig:nat_transformation}A natural transformation\index{Natural 
transformation} $T$ between two functors\index{Functor} $F,G$ yields a 
commutative diagram for every morphism $f \colon A \rightarrow B$ since $T$ 
commutes with the application of functors\index{Functor}.}
\end{figure}

Even with categories\index{Category theory} there is the basic question of 
equality\index{Equality}\index{Identity} and there can be equivalent 
variants. Equivalences\index{Equivalence} of two 
categories\index{Category theory} $\mathcal{C}$ and $\mathcal{D}$ arise 
through two functors\index{Functor} 
\[
F \colon  \mathcal{C} \xlongrightarrow{~~~~~~} \mathcal{D}\text{ and } 
G \colon \mathcal{D} \xlongrightarrow{~~~~~~} \mathcal{C},
\]
so that there are natural transformations 
\[
\varepsilon \colon \mathrm{id}_\mathcal{D} 
\xlongrightarrow{~~ \cong ~~} F \circ G \text{ and } 
\eta \colon G \circ F \xlongrightarrow{~~ \cong ~~} \mathrm{id}_\mathcal{C} 
\]
which induce isomorphisms\index{Isomorphism} on all objects.

\section{Presheaves and Sheaves\index{Sheaf}}

With the help of functors, sheaves\index{Sheaf} and presheaves can be defined 
as new objects on topological spaces\index{Topology}. Such objects are 
important as structure sheaves\index{Sheaf} and as coefficients for cohomology 
theories. In the history of the development of mathematics, they were 
introduced around 1950 within complex analysis by Henri Cartan\index{Cartan, 
Henri} and Kiyoshi Oka\index{Oka, Kiyoshi} in the form of sheaves\index{Sheaf} 
of holomorphic functions on complex-analytical spaces.\footnote{On open 
subsets $U \subset X$ there are commutative algebras $\mathcal{O}_X(U)$ of 
holomorphic functions, which together form the structure sheaf $\mathcal{O}_X$. 
An analogous algebraic example are affine spectra $X=\mathrm{Spec}(A)$ for 
commutative rings $A$, which consist of all prime ideals in $A$ and are endowed 
with the Zariski topology. By gluing together such affine spectra, 
the concept of the scheme by Alexander 
Grothendieck\index{Grothendieck, Alexander} is obtained. Schemes also carry a 
structure sheaf $\mathcal{O}_X$. Instead of the Zariski topology, there are 
other topologies on a scheme such as the \'etale topology and other 
Grothendieck topologies, where the open sets are usually not 
subsets of $X$ but morphisms $U \longrightarrow X$ in a category. Variants of 
schemes in algebraic geometry are the rigid-analytic spaces by John 
Tate\index{Tate, John} and the perfectoid spaces by Peter 
Scholze\index{Scholze, Peter}, which open up new possibilities for research. 
Scholze\index{Scholze, Peter} received the Fields Medal in 2018 for his ideas.} 

A presheaf\index{Sheaf} $\mathcal{F}$ is a contravariant functor from 
the category\index{Category theory} of open 
sets $\mathbf{Off}(X)$ in a topological space\index{Topology} 
$X$ with values in the category $\mathbf{Set}$. This means that for 
two open sets $V \subset U$ there is a restriction mapping 
\[
\rho^U_V \colon \mathcal{F}(U) \xlongrightarrow{~~~~~~} \mathcal{F}(V)
\]
which must satisfy two functor properties\index{Functor}. Firstly, the trivial 
inclusion $U \subset U$ induces the identity\index{Identity}\index{Equality} 
$\mathrm{id}_{\mathcal{F}(U)}$ as restriction mapping $\rho^U_U$ and secondly, 
for three nested open sets $W \subset V \subset U$ the composition rule for 
restriction mappings applies (see Fig.~\ref{fig:presheaf}). The category of 
presheaves\index{Category theory}\index{Sheaf} is denoted by $\mathbf{Psh}(X)$. 

\begin{figure}[ht!]
\begin{tikzpicture}
\node (a) at (0,0) {$\mathcal{F}(U)$};
\node (b) at (3,0) {$\mathcal{F}(V)$};
\node (c) at (6,0) {$\mathcal{F}(W)$};
\draw[->] (a) to node[scale=.7] [above] {$\rho^U_V$} (b);
\draw[->] (b) to node[scale=.7] [above] {$\rho^V_W$} (c);
\draw[->] (a) to [bend right] node[scale=.7] [above] 
{$\rho^U_W=\rho^V_W \circ \rho^U_V$} (c);
\end{tikzpicture}
\caption{\label{fig:presheaf}The functor property\index{Functor} takes care 
that the composition of restriction mappings is again a restriction mapping.}
\end{figure}

Sheaves\index{Sheaf} of sets\index{Category theory} on a 
topological space\index{Topology} $X$ are special 
presheaves\index{Sheaf} that fulfil an additional property.\footnote{The sheaf 
property refers to coverings. If $\mathcal{F}$ is a presheaf, $U$ is an open 
set in $X$ and $f \in \mathcal{F}(U)$ is an element, then open coverings 
$U=\bigcup_i U_i$ are considered. If we denote by $f_i$ the restriction of 
$f$ to $U_i$, then obviously the restrictions of $f_i$ and $f_j$ 
to the subset $U_i \cap U_j$ coincide for all pairs of indices $i,j$. 
A presheaf $\mathcal{F}$ is a sheaf if, conversely, given 
functions $f_i \in \mathcal{F}(U_i)$ with coinciding restrictions
on $U_i \cap U_j$, a section $f \in \mathcal{F}(U)$ can be constructed,
which delivers the function $f_i$ on each $U_i$. The sheaf property is
often symbolised by the exactness of the following equaliser sequence:
\[
\mathcal{F}(U) \to \prod_i \mathcal{F}(U_i) \rightrightarrows 
\prod_{i,j} \mathcal{F}(U_i \cap U_j).
\]
}
They form a category\index{Category theory}, denoted by $\mathbf{Sh}(X)$. A 
sheaf\index{Sheaf} or a presheaf forms a family of sets over $X$, because for 
every point $x \in X$ the stalk $\mathcal{F}_x$ of the sheaf\index{Sheaf} can 
be defined as a set.\footnote{To define the stalk $\mathcal{F}_x$, the limit of 
all $\mathcal{F}(U)$ with $x \in U$ must be considered and it results 
in $\mathcal{F}_x$ as a set that depends on $x$. This limit object 
receives incoming arrows from all $\mathcal{F}(U)$ with $x \in U$ and is 
a universal object with this property, i.e., for every other object 
$\mathcal{G}$ with such incoming arrows there is an arrow from 
$\mathcal{F}_x$ to $\mathcal{G}$.}

The category\index{Category theory} of presheaves\index{Sheaf} 
$\mathbf{Psh}(X)$ is an example of a functor category\index{Functor}. It 
consists of the collection
\[
\hat{\mathcal{C}}=\mathbf{Set}^{\mathcal{C}^{\mathrm{op}}}
\] 
of all functors\index{Functor} $\mathcal{C}^{\mathrm{op}} \longrightarrow  
\mathbf{Set}$ on a given category\index{Category theory} 
$\mathcal{C}$. The category\index{Category theory} $\mathbf{Set}$ 
is the special case of the presheaves\index{Sheaf} $\mathbf{Psh}(*)$ on the 
point. For the category\index{Category theory} 
$\mathcal{C}=\Delta_\bullet$ of finite ordered sets of the form
\[
\Delta_n=\{0 < 1 < \cdots < n\}
\]
with monotone mappings as morphisms, the functor category\index{Functor} 
$\hat{\Delta}_\bullet$ is the category\index{Category theory} $\mathbf{sSet}$
of simplicial sets\index{Simplicial set}.\footnote{The reason is that 
there are injective boundary mappings $d_i \colon \Delta_{n-1} \longrightarrow 
\Delta_n$ on the one hand, which omit an element $0 \le i \le n$, and on the 
other hand surjective degeneration mappings $s_i \colon \Delta_{n} 
\longrightarrow \Delta_{n-1}$, which hit the element $i$ twice. Each 
contravariant functor from $\Delta_\bullet$ to $\mathbf{Set}$ then provides a 
simplicial set $S_\bullet$ in $\mathbf{sSet}$.}

Functor categories\index{Functor} often have better properties 
than $\mathcal{C}$ itself and every small category\index{Category theory} 
possesses a useful embedding 
\[
y \colon \mathcal{C} \xlongrightarrow{~~~~~~} \hat{\mathcal{C}}, \;  B 
\mapsto \mathrm{Hom}_\mathcal{C}(-,B),
\]
which is also called Yoneda embedding\footnote{If the collections 
$\mathrm{Hom}_\mathcal{C}(A,B)$ do not form sets, then $\hat{\mathcal{C}}$ can 
be replaced by the cocompletion of $\mathcal{C}$ which consists not only of 
presheaves.} after Nobuo Yoneda\index{Yoneda, Nobuo}. 

\section{The Category {\bf Set} as Elementary Topos\index{Elementary 
topos}\index{Topos}}

How can $\mathbf{Set}$ be defined as a category\index{Category theory} 
without using axiomatic set theory\index{Set theory}? William 
Lawvere\index{Lawvere, William} has in the 1960s described the category of 
sets\index{Set theory}\index{Category theory} axiomatically in a structural 
way which provides an alternative to Zermelo-Fraenkel set 
theory\index{Zermelo-Fraenkel axioms}\index{Set theory}. 
Lawvere\index{Lawvere, William} has named his approach the \enquote{Elementary 
theory of the category of sets} (ETCS). Later, other versions such as SEAR 
(\enquote{Sets, elements and relations}) by Jeremy Avigad\index{Avigad, Jeremy} 
emerged.\footnote{See \cite{lawvere1964}.} 

The notion of an elementary topos\index{Elementary topos} has been coined by  
William Lawvere\index{Lawvere, William} and Myles Tierney\index{Tierney, Myles} 
around 1970, because they saw that certain properties of $\mathbf{Set}$, in the 
form Lawvere\index{Lawvere, William} has asked  for ETCS, form an important 
class of categories\index{Category theory}. 

\medskip 
The precise axioms of an elementary topos\index{Elementary topos}\index{Topos} 
$\mathcal{E}$ are:  
\begin{itemize}
\item $\mathcal{E}$ is (locally) cartesian closed\index{Cartesian closed 
category}. 
\item $\mathcal{E}$ has (finite) limits and colimits.
\item $\mathcal{E}$ has a subobject classifier\index{Subobject classifier}. 
\end{itemize}

With these axioms, $\mathbf{Set}$ is an elementary topos\index{Elementary 
topos}\index{Topos} which has a natural number object $\mathbb{N}$, i.e., a 
zero object $0$ and a successor map  $S \colon \mathbb{N} \longrightarrow 
\mathbb{N}$.\footnote{See \cite[Part II]{lambekscott1986}.} 

There are, in addition to $\mathbf{Set}$, two further examples of 
elementary topoi\index{Topos}\index{Elementary topos} that we have already 
encountered, namely the functor categories\index{Functor} 
\[
\hat{\mathcal{C}}=\mathbf{Set}^{\mathcal{C}^{\mathrm{op}}}
\]
for a category\index{Category theory} $\mathcal{C}$ and the category 
$\mathbf{Sh}(X)$ of sheaves of sets\index{Category theory} on a topological 
space\index{Topology} $X$.

\section{Locally Cartesian Closed Categories\index{Category 
theory}\index{Cartesian closed category}}

The first group of properties of an elementary topos\index{Elementary 
topos}\index{Topos} states that $\mathcal{E}$ is a cartesian closed 
category\index{Cartesian closed category}. In such categories\index{Category 
theory}, the element symbol is not used. Instead, the elements of $M$ are given 
by arrows 
\[
1 \xlongrightarrow{~~ x ~~} M 
\]
from a canonical terminal object $1$ to $M$. Terminal means that for all 
objects $A$ in $\mathcal{E}$ a morphism $A \longrightarrow 1$ exists. 
With the help of this object it can be tested whether two morphisms $f,g \colon 
M \longrightarrow M'$ in $\mathrm{Hom}(M,M')$ are equal by checking whether 
their composition with all elements $1 \xlongrightarrow{x} M$ coincide. Hence, 
we require that sufficiently many elements are present to guarantee the 
function extensionality\index{Extensionality} of morphisms. 

\medskip 
For any two objects $A$, $B$ in $\mathcal{E}$ one asks that the binary 
cartesian product $A \times B$ exists (see Fig.~\ref{fig:cat_product}). 

\begin{figure}[ht!]
\begin{tikzpicture}
\matrix (m) [matrix of math nodes,row sep=3em,column sep=4em,minimum width=2em]
  {
     A \times B & B \\
      A & 1 \\};
  \path[-stealth]
    (m-1-1) edge node [left] {$\mathrm{pr}_A$} (m-2-1)
            edge node [above] {$\mathrm{pr}_B$} (m-1-2)
  (m-2-1.east|-m-2-2) edge node [below] {}
            node [above] {} (m-2-2)
    (m-1-2) edge node [right] {} (m-2-2);
\draw +(-.2,0) -- +(0,0)  -- +(0,.2);
\end{tikzpicture}
\caption{\label{fig:cat_product}The binary cartesian product $A \times 
B$ is characterised by two projection morphisms to $A$ and $B$ and the 
universal property for any object $T$ with given morphisms to $A$ and $B$.}
\end{figure}

The binary cartesian product fulfils the universal property that every morphism 
\[
T \xlongrightarrow{~~~ h ~~~} A \times B  
\]
is uniquely given by two morphisms $T \longrightarrow A$ and $T \longrightarrow 
B$ into the two components in the form $h=(f,g)$. More generally, a 
commutative diagram is called cartesian, or pullback diagram, if all 
properties of the binary cartesian product, including the universal property, 
are satisfied.  

\medskip 
Furthermore, there should be an internal exponential object\index{Exponential 
object} $X^B$ in every locally cartesian closed category\index{Category 
theory}\index{Cartesian closed category} $\mathcal{E}$, which represents the 
functions from $B$ to $X$. This means that all morphisms $B \longrightarrow X$ 
are represented by morphisms $1 \longrightarrow X^B$ into an object $X^B$ and 
more generally every morphism 
\[
A \times B \xlongrightarrow{~~~ f ~~~} X
\]
corresponds to a unique morphism
\[
A \xlongrightarrow{~~~ F ~~~} X^B
\]
where $f(a,b)=F(a)(b)$ holds. This implies that morphisms with multiple 
arguments can successively be replaced by simple morphisms. In the literature, 
the exponential object\index{Exponential object} is also called an internal 
hom-object. The property of being a cartesian closed category\index{Cartesian 
closed category} holds even locally in an elementary topos\index{Elementary 
topos}\index{Topos}. This means that all comma categories\index{Category 
theory} $\mathcal{E}/_X$ are cartesian closed\index{Cartesian closed category}. 
Their objects $A$ are defined as objects in $\mathcal{E}$ along with a morphism 
$A \longrightarrow X$.

\medskip 
In addition to the terminal object $1$, which corresponds to a singleton set in 
$\mathbf{Set}$, there is an initial object $0$ in locally cartesian closed 
categories\index{Category theory}\index{Cartesian closed category}, which 
corresponds to the empty set in $\mathbf{Set}$. Initial means that there are 
morphisms $0 \longrightarrow A$ for all $A$ in $\mathcal{E}$. The objects $0$ 
and $1$ do not represent numbers or truth values\index{Truth}, although a 
connection exists.

\section{Limits and Colimits}

The second group of properties of an 
elementary topos\index{Elementary topos}\index{Topos} requires that the 
category\index{Category theory} $\mathcal{E}$ has so-called finite limits 
and finite colimits, denoted by lim and colim. 

\medskip 
The simplest example of a limit is the binary cartesian product $A \times B$ in 
$\mathbf{Set}$. A more general limit is a fiber product $A \times_X B$ 
for two morphisms $f \colon A \longrightarrow X$ and $g \colon B 
\longrightarrow X$ in $\mathcal{E}$ (see Fig.~\ref{fig:fiber_product}). 

\begin{figure}[ht!]
\begin{tikzpicture}
\matrix (m) [matrix of math nodes,row sep=3em,column sep=4em,minimum width=2em]
  {
     A \times_X B & B \\
      A & X \\};
  \path[-stealth]
    (m-1-1) edge node [left] {$\mathrm{pr}_A$} (m-2-1)
            edge node [above] {$\mathrm{pr}_B$} (m-1-2)
  (m-2-1.east|-m-2-2) edge node [below] {$f$}
            node [above] {} (m-2-2)
    (m-1-2) edge node [right] {$g$} (m-2-2);
\draw +(-.2,0) -- +(0,0)  -- +(0,.2);
\end{tikzpicture}
\caption{\label{fig:fiber_product}The cartesian diagram for the fiber product 
is the binary cartesian product in the comma category\index{Category theory} 
$\mathcal{E}/_X$.}
\end{figure}

The object $A \times_X B$ is the binary cartesian product in the comma 
category\index{Category theory} $\mathcal{E}/_X$. The universal property for  
$A \times_X B$ for two morphisms $T \longrightarrow A$ and 
$T \longrightarrow B$, whose compositions with $f$ and $g$ coincide as 
morphisms $T \longrightarrow X$, has an arrow 
\[
T \xlongrightarrow{~~~~~~} A \times_X B 
\]
as a consequence, whose respective compositions with the two projections again 
coincide with the given morphisms. In the case of the terminal object 
$X=1$ the binary cartesian product results in
\[
A \times_X B = A \times B. 
\]
A special case of a colimit is a pushout object $A +_X B$ for two morphisms $f 
\colon X \longrightarrow A$ and $g \colon X \longrightarrow B$ (see 
Fig.~\ref{fig:pushout}). 

\begin{figure}[ht!]
\begin{tikzpicture}
\matrix (m) [matrix of math nodes,row sep=3em,column sep=4em,minimum width=2em]
  {
    X &  A   \\
      B & A +_X B \\};
  \path[-stealth]
    (m-1-1) edge node [left] {$g$} (m-2-1)
            edge node [above] {$f$} (m-1-2)
  (m-2-1.east|-m-2-2) edge node [below] {$i_B$}
            node [above] {} (m-2-2)
    (m-1-2) edge node [right] {$i_A$} (m-2-2);
\draw +(.2,0) -- +(0,0)  -- +(0,-.2);
\end{tikzpicture}
\caption{\label{fig:pushout}The characterising cocartesian diagram for 
$A +_X B$ is dual to the diagram for the fiber product. It shows two morphisms 
from $A$ and $B$ to $A +_X B$ satisfying the universal property.}
\end{figure}
 
In the case of the initial object $X=0$, the binary sum $A+B$, also called 
coproduct, results. In $\mathbf{Set}$, the coproduct corresponds to the 
disjoint union. More generally, a commutative diagram is called cocartesian, 
or pushout diagram, if all properties of the pushout object, including the 
universal property, are satisfied.  

\section{Subobject Classifiers\index{Subobject classifier}}

The third group of properties of an 
elementary topos\index{Elementary topos}\index{Topos} demands that the 
category\index{Category theory} $\mathcal{E}$ has a subobject 
classifier\index{Subobject classifier} $\Omega$ together with an arrow $1 
\xlongrightarrow{~~ t ~~} \Omega$ ($t$ as in true).\footnote{See \cite[Ch. I, 
Sec. 4]{maclane_moerdijk1992} for an axiomatic definition of elementary topoi 
and examples of subobject classifiers.} This guarantees for every subobject $B 
\subset A$ a classifying arrow $\chi$, so that in the depicted diagram $t$ 
becomes a universal embedding (see Fig.~\ref{fig:subobject}).

\begin{figure}[ht!]
\begin{tikzpicture}
  \matrix (m) [matrix of math nodes,row sep=3em,column sep=4em,minimum 
width=2em]
  {
     B & 1 \\
      A & \Omega \\};
  \path[-stealth]
    (m-1-1) edge node [left] {} (m-2-1)
            edge node [above] {} (m-1-2)
    (m-2-1.east|-m-2-2) edge node [below] {}
            node [above] {$\chi$} (m-2-2)
    (m-1-2) edge node [right] {$t$} (m-2-2);
\end{tikzpicture}
\caption{\label{fig:subobject}The characterising diagram for the subobject 
classifier\index{Subobject classifier} views a subobject $B$ of $A$ as the
locus where the characteristic truth function $\chi$ has value $1$.}
\end{figure}

The idea here is that $\Omega$ is a collection of truth values\index{Truth} and 
contains a distinguished element $1$ that corresponds to the value true. In the 
example of the elementary topos\index{Elementary topos}\index{Topos} 
$\mathbf{Set}$, $\Omega=\{0,1\}$ with the usual classical truth 
values\index{Truth} $0$ (false) and $1$ (true). The morphism $\chi$ in this 
example is the characteristic function of the subset $B$, i.e., for all 
elements $a$ of $A$ it holds
\[
\chi(a)=\begin{cases} 1 &\text{ if } a \text{ is an element of } B \cr 0 & 
\text { otherwise.} 
\end{cases} 
\]
An elementary topos\index{Elementary topos}\index{Topos} admits internal 
power objects $\mathrm{Pow}(A)$ for every $A$ in $\mathcal{E}$. The 
subobject classifier\index{Subobject classifier} $\Omega$ is thus of the form 
$\mathrm{Pow}(1)$. All power objects $\mathrm{Pow}(A)$ and all collections 
\[
\mathrm{Sub}_{\mathcal{E}}(A)=\mathrm{Hom}_{\mathcal{E}}(A,\Omega)
\]
of subobjects of an object $A$ in $\mathcal{E}$ carry the structure of 
Heyting algebras\index{Heyting algebra}, which in both cases are usually not 
Boolean algebras\index{Boolean algebra}.

In a Heyting algebra\index{Heyting algebra} there are operations $\wedge, \vee, 
\Rightarrow$ and $\neg$. In an elementary topos\index{Topos}\index{Elementary 
topos}, the operations $\wedge$ and $\vee$ on the above Heyting 
algebras\index{Heyting algebra} are induced by fiber products 
and coproducts. The operation $A \Rightarrow B$ between two objects $A,B$ is 
traced back to the exponential object\index{Exponential object} $B^A$. The 
negation $\neg A$ is defined by $A \Rightarrow 0$. The operation $\wedge$ is 
related to the universal property of $\Omega$ as shown in 
Fig.~\ref{fig:subobject_algebra}. Summarizing, we say that an elementary 
topos\index{Topos}\index{Elementary topos} is a Heyting 
category\index{Category theory}.\footnote{See \cite[Ch. IV, 
Sec. 8]{maclane_moerdijk1992} and \cite[Part II]{lambekscott1986}.}

\begin{figure}[ht!]
\begin{tikzpicture}
  \matrix (m) [matrix of math nodes,row sep=3em,column sep=4em,minimum 
width=2em]
  {
     1 & 1 \\
      \Omega \times \Omega & \Omega \\};
  \path[-stealth]
    (m-1-1) edge node [left] {$t \times t$} (m-2-1)
            edge node [above] {} (m-1-2)
    (m-2-1.east|-m-2-2) edge node [below] {}
            node [above] {$\wedge$} (m-2-2)
    (m-1-2) edge node [right] {$t$} (m-2-2);
\end{tikzpicture}
\caption{\label{fig:subobject_algebra}The subobject classifier\index{Subobject 
classifier} $\Omega$ is a Heyting algebra\index{Heyting algebra} and thus has 
an intersection operation $\wedge \colon \Omega \times \Omega \rightarrow 
\Omega$.}
\end{figure}

\section{$2$-Categories\index{Category theory} and 
Bicategories\index{Bicategory}}

Higher categories\index{Category theory} and 
infinity categories\index{Infinity category} are modifications of 
categories\index{Category theory} that can handle the 
equivalence\index{Equivalence} of mathematical objects better 
than ordinary categories\index{Category theory}.\footnote{See 
\cite{cisinski2015,groth2015,lurie2008,lurie2009}.} In this process, 
categorical calculation rules are weakened in different ways and the coherence 
of this approach is secured by higher morphisms (see 
Fig.~\ref{fig:higher_morphism}).

\begin{figure}[ht!]
\begin{tikzpicture}
\node (a) at (0,0) {$A$};
\node (b) at (3,0) {$B$};
\draw[->] (a) to [bend left] node[scale=.7] (f) [above] {$f$} (b);
\draw[->] (a) to [bend right] node[scale=.7] (g) [below] {$g$} (b);
\draw[-{Implies},double distance=1.5pt,shorten >=2pt,shorten <=2pt] (f) 
to node[scale=.7] [right] {$\alpha$} (g);
\end{tikzpicture}
\caption{\label{fig:higher_morphism}Higher $2$-morphism $\alpha$ between two 
$1$-morphisms $f,g$.}
\end{figure}

There are precise definitions for higher $n$-categories\index{Category theory} 
for small natural numbers $n=2,3,\ldots$, which in the limit case $n 
\to \infty$ are called infinity categories\index{Infinity category}. We 
start by explaining the case of a strict $2$-category\index{Category theory}. 
In them, there are additional $2$-morphisms between the ordinary 
$1$-morphisms. The composition of arrows can occur in different ways. Thus, in 
a $2$-category, in addition to the composition of two ordinary morphisms, 
there is a vertical and a horizontal possibility to link the $2$-morphisms 
(see Fig.~\ref{fig:vertical} and Fig.~\ref{fig:horizontal}). 

\begin{figure}[ht!]
\begin{tikzpicture}
\node (a) at (0,0) {$A$};
\node (b) at (3,0) {$B$};

\draw[->] (a) to [bend left] node[scale=.7] (f) [above] {$f$} (b);
\draw[->] (a) to node[scale=0.1] (g) [right] {$g$} (b);
\draw[->] (a) to [bend right] node[scale=.7] (h) [below] {$h$} (b);
\draw[-{Implies},double distance=1.5pt,shorten >=2pt,shorten <=2pt] 
(f) to node[scale=.7] [right] {$\alpha$} (g);
\draw[-{Implies},double distance=1.5pt,shorten >=2pt,shorten <=2pt] 
(g) to node[scale=.7] [right] {$\beta$} (h);
\end{tikzpicture}
\caption{\label{fig:vertical}The vertical composition concatenates two 
$2$-morphisms between $1$-morphisms $f,g$ resp. $g,h$ to a $2$-morphism between 
$f$ and $h$.}
\end{figure}

\begin{figure}[ht!]
\begin{tikzpicture}

\node (a) at (0,0) {$A$};
\node (b) at (3,0) {$B$};
\node (c) at (6,0) {$C$};

\draw[->] (a) to [bend left] node[scale=.7] (f) [above] {$f$} (b);
\draw[->] (a) to [bend right] node[scale=.7] (g) [below] {$g$} (b);
\draw[-{Implies},double distance=1.5pt,shorten >=2pt,shorten <=2pt] 
(f) to node[scale=.7] [right] {$\alpha$} (g);
\draw[->] (b) to [bend left] node[scale=.7] (f') [above] {$f'$} (c);
\draw[->] (b) to [bend right] node[scale=.7] (g') [below] {$g'$} (c);
\draw[-{Implies},double distance=1.5pt,shorten >=2pt,shorten <=2pt] 
(f') to node[scale=.7] [right] {$\beta$} (g');
\end{tikzpicture}
\caption{\label{fig:horizontal}The horizontal composition of $2$-morphisms 
between $1$-morphisms $f, g$ and $f',g'$ yields a $2$-morphisms between the 
$1$-morphisms $f \circ g$ and $f' \circ g'$.}
\end{figure}

In a strict $2$-category\index{Category theory}, it is required that both 
types of compositions are associative and commute with each other (see 
Fig.~\ref{fig:commute}). Thus, the diagram shown has no ambiguities. In 
addition, there should be $2$-morphisms $\mathrm{id}_f$ for each morphism $f$ 
that are identities\index{Identity}\index{Equality} with respect to both 
compositions. 

\begin{figure}[ht!]
\begin{tikzpicture}

\node (a) at (0,0) {$A$};
\node (b) at (3,0) {$B$};
\node (c) at (6,0) {$C$};

\draw[->] (a) to [bend left] node[scale=.7] (f) [above] {$f$} (b);
\draw[->] (a) to node[scale=0.1] (g) [right] {$g$} (b);
\draw[->] (a) to [bend right] node[scale=.7] (h) [below] {$h$} (b);
\draw[-{Implies},double distance=1.5pt,shorten >=2pt,shorten <=2pt] 
(f) to node[scale=.7] [right] {$\alpha$} (g);
\draw[-{Implies},double distance=1.5pt,shorten >=2pt,shorten <=2pt] 
(g) to node[scale=.7] [right] {$\beta$} (h);

\draw[->] (b) to [bend left] node[scale=.7] (f) [above] {$f'$} (c);
\draw[->] (b) to node[scale=0.1] (g) [right] {$g$} (c);
\draw[->] (b) to [bend right] node[scale=.7] (h) [below] {$h'$} (c);
\draw[-{Implies},double distance=1.5pt,shorten >=2pt,shorten <=2pt] 
(f) to node[scale=.7] [right] {$\gamma$} (g);
\draw[-{Implies},double distance=1.5pt,shorten >=2pt,shorten <=2pt] 
(g) to node[scale=.7] [right] {$\delta$} (h);
\end{tikzpicture}
\caption{\label{fig:commute}Horizontal and vertical compositions commute with 
each other.}
\end{figure}

The most important example of a strict $2$-category\index{Category theory} 
is $\mathbf{Cat}$, the category\index{Category theory} of all (small) 
categories with functors\index{Functor} as morphisms and natural 
transformations\index{Natural transformation} as $2$-morphisms. Another 
example is the category $\mathbf{Gr}$ of groups\index{Category theory}, where 
the morphisms are the group homomorphisms
\[
G \xlongrightarrow{~~ f ~~} H, \; \text{with } f(ab)=f(a)f(b) 
\]
and the $2$-morphisms 
\[
f \xRightarrow{~~~~} f'
\]
between $f$ and $f'$ are given by the conjugations with elements of the 
group:
\[
f'(g)=hf(g)h^{-1}, \text{ with } h\in H.  
\]
As a rule, strict $2$-categories\index{Category theory} are rare, 
although it can be shown that every $2$-category\index{Category theory} 
is equivalent to a strict one\index{Equivalence}. More often, 
bicategories\index{Bicategory}\index{Category theory} are considered. These are 
weakenings of $2$-categories, in which the composition of 
$1$-morphisms only needs to be associative modulo certain higher expressions.

\begin{figure}[ht!]
\tikzset{zshift/.style={xshift={-0.3*#1},yshift={-0.9*#1}}}
\begin{tikzcd}[row sep=2cm,column sep=2cm,inner sep=1ex]
& B \arrow{d}{g} \arrow{dr}{h \circ g} & \\
A \arrow{ur}{f} 
\arrow{r}{g \circ f} \arrow{rr}{\hskip-11mm (h \circ g) \circ f \Longrightarrow
h \circ (g \circ f)}& |[zshift=-1.5cm]| C \arrow{r}{h}& |[zshift=1cm]|D
\end{tikzcd}
\caption{\label{fig:associator}The associator is a $2$-morphism which replaces 
the shift of brackets from left to right in the associativity law for the 
composition of morphisms. It is naturally represented by a tetrahedron.}
\end{figure}

As a replacement for the associative law, there is an invertible $2$-morphism 
\[
a_{h,g,f} \colon (h \circ g) \circ f \xRightarrow{~~~~} h \circ (g \circ f)
\]
-- called the associator -- which depends naturally on $h$, $g$ and $f$. 
It describes the shift of brackets from left to right and is naturally 
represented by a $3$-dimensional tetrahedron, the sides of which correspond to 
the possible compositions (see Fig.~\ref{fig:associator}). Associators for four 
objects $i,h,g,f$ satisfy the pentagon rule, which states that the two ways of 
going around the pentagon yield the same result (see 
Fig.~\ref{fig:pentagon_rule}).  

\begin{figure}[ht!]
\begin{tikzpicture}[rotate=90,scale=3]

\draw (-0*360/5:1) node (v1) {$(i \circ h) \circ (g \circ f)$};
\draw (-1*360/5:1) node (v3) {$i \circ (h \circ (g \circ f))$};
\draw (-2*360/5:1) node (v5) {$i\circ ((h \circ g) \circ f)$};
\draw (-3*360/5:1) node (v7) {$(i\circ (h \circ g)) \circ f$};
\draw (-4*360/5:1) node (v9) {$((i \circ h) \circ g) \circ f)$};

\draw[-{Implies},double distance=1.5pt] (v1) -- (v3) ;
\draw[-{Implies},double distance=1.5pt] (v5) -- (v3);
\draw[-{Implies},double distance=1.5pt] (v7) -- (v5);
\draw[-{Implies},double distance=1.5pt] (v9) -- (v7);
\draw[-{Implies},double distance=1.5pt] (v9) -- (v1);
\end{tikzpicture}
\caption{\label{fig:pentagon_rule}Four composable 
morphisms can be bracketed in different ways which are related by associators 
(shifts of brackets) in the form of a pentagon.} 
\end{figure}

Additionally, in a bicategory\index{Bicategory}, the 
existence\index{Existence} of invertible $2$-morphisms 
\[
l_f \colon \mathrm{id}_B \circ f \xRightarrow{~~~~} f, \quad r_f \colon f 
\circ \mathrm{id}_A \xRightarrow{~~~~} f
\] 
for every morphism $f \colon A \longrightarrow B$ is required, which are 
called unitors. 

\section{Higher Categories\index{Category theory} and 
Infinity Categories\index{Infinity category}}

The situation becomes even more interesting with 
$3$-categories\index{Category theory}, where the pentagon rule is not 
a equality\index{Identity}\index{Equality}, but is replaced by a 
$3$-morphism, which is called the associahedron. In general, there 
are $n$-categories\index{Category theory} for all $n$. In these,
higher $k$-morphisms are often described in a geometric way that reminds one 
of simplicial structures\index{Simplicial space}\index{Simplicial set}. The 
$k$-dimensional morphisms correspond to $k$-dimensional 
simplices\index{Simplex}, which describe the relationships between 
$(k-1)$-dimensional objects. For example, the composition $g \circ f$ of two 
morphisms $f$ and $g$ is described by the choice of a $2$-dimensional 
simplex\index{Simplex} $\Delta_2$ (see Fig.~\ref{fig:composition2}). 

\begin{figure}[ht!]
\begin{tikzpicture}
\node (a) at (0,-1) {$A$};
\node (b) at (2,0) {$B$};
\node (c) at (4,-1) {$C$};
\draw[->] (a) to node[scale=.7] [above] {$f$} (b);
\draw[->] (b) to node[scale=.7] [above] {$g$} (c);
\draw[->] (a) to node[scale=.7] (x) [below] {$g \circ f$} (c);
\end{tikzpicture}
\caption{\label{fig:composition2}Composition $g \circ f$ of concatenable 
morphisms $f,g$.}
\end{figure}

Infinity categories\index{Infinity category} are even more 
general than $n$-categories\index{Category theory}. They were found by Michael 
Boardman\index{Boardman, Michael} and Rainer Vogt\index{Vogt, Rainer} around 
1973, who discovered that higher $k$-morphisms for all $k \ge 1$ can be very 
useful in homotopy theory\index{Homotopy theory}.\footnote{See 
\cite{boardman}.} The calculation rules and relationships between higher 
morphisms are, as in the cases of $2$- and $3$-categories\index{Category 
theory} or bicategories\index{Bicategory}, to be determined. Thus, the 
composition of higher morphisms need not be unique nor strictly associative. 
But the basic principle applies that the choices for non-unique higher 
morphisms are in a suitable sense equivalent\index{Equivalence} to each other, 
if they are not unique. Such coherence conditions among $k$-morphisms are 
given by $(k+1)$-morphisms. In the simplest case, the 
equivalence\index{Equivalence} of two morphisms $f,g \colon A 
\longrightarrow B$ is described by a $2$-simplex\index{Simplex} (see 
Fig.~\ref{fig:equal_morphisms}). 

\begin{figure}[ht!]
\begin{tikzpicture}
\node (a) at (0,-1) {$A$};
\node (b) at (2,0) {$A$};
\node (c) at (4,-1) {$B$};
\draw[->] (a) to node[scale=.7] [above] {$\mathrm{id}_A \;\;$} (b);
\draw[->] (b) to node[scale=.7] [above] {$f$} (c);
\draw[->] (a) to node[scale=.7] (x) [below] {$g$} (c);
\end{tikzpicture}
\caption{\label{fig:equal_morphisms}Equivalence\index{Equivalence} of two 
morphisms $f,g$ is more general than judgemental 
equality\index{Equality}\index{Identity} $f=g$.}
\end{figure}

Infinity categories\index{Infinity category} are in a certain sense the limit 
case of a hierarchy of $(n,r)$-categories\index{Category theory}. The letter 
$n$ stands for the same counting as in $n$-categories\index{Category theory} 
and the letter $r$ with $0 \le r \le n+1$ for another index. In an 
$(n,r)$-category\index{Category theory}, $k$-morphisms become trivial if 
$k>n$, and they become invertible equivalences\index{Equivalence} if $k>r$. 
Also in infinity categories\index{Infinity category}\index{Category theory} 
this classification can be made and $(\infty,r)$-categories\index{Infinity 
category} result. Most often in this context, 
$(\infty,1)$-categories\index{Infinity category} are meant when we 
speak of infinity categories\index{Infinity category}. Overall, the scheme in 
Fig.~\ref{fig:cat_scheme} emerges, in which the concept of the set (or class) 
and the partially ordered set (poset) generalise.

\begin{figure}
\begin{tabular}{c|ccccc} 
& $n=0$ & $n=1$ & $n=2$ & $\cdots$ & $n=\infty$ \\\hline 
$r=0$ &  set   & groupoid  & $2$-groupoid     & 
$\cdots$ &  $(\infty,0)$-category \\ 
$r=1$ &  poset & category & $(2,1)$-category & 
$\cdots$ & $(\infty,1)$-category \\
$r=2$ &  poset & $2$-poset  & $2$-category     & 
$\cdots$ & $(\infty,2)$-category \\ 
$r=3$ & poset  & $2$-poset  & $3$-poset        & 
$\cdots$ & $(\infty,3)$-category \\ 
$\vdots$ & $\vdots$ & $\vdots$  & $\vdots$  & $\ddots$ & $\vdots$ 
\end{tabular}
\caption{\label{fig:cat_scheme}Two-dimensional infinite hierarchy of higher 
categories.}
\end{figure}    

Infinity categories\index{Infinity category} (or higher categories) are mostly 
studied in the form of set-theoretic models\index{Model theory} which are 
called quasi-categories, (weak) Kan complexes\index{Kan, 
Daniel}, complete Segal spaces\index{Segal, Graeme}, Segal categories, or 
simplicial categories.\footnote{See \cite{groth2015}.} However, this 
contradicts the general belief that category theory\index{Category theory} 
ought to be considered as an independent foundation of mathematics. So far, no 
clear axiomatic system has prevailed for infinity categories\index{Infinity 
category} since many properties of set-theoretic models\index{Model theory} 
depend on the chosen model. However, it seems desirable to define intrinsic 
axioms without resorting to models. Recently this has become a important 
research question by the work of Emily Riehl\index{Riehl, Emily}, Dominic 
Verity\index{Verity, Dominic} and the group around Denis-Charles 
Cisinski\index{Cisinski, Denis-Charles}. A comprehensive answer cannot be given 
at this point.\footnote{See \cite{cisinski2025,riehlverity}.} Attempts have 
also been made to define $(n,r)$-topoi or elementary 
$(\infty,r)$-topoi\index{Elementary topos}\index{Topos} in a model-independent 
way\index{Model theory}. Such categories\index{Category theory} should fulfil 
similar axioms as ordinary elementary topoi\index{Topos}\index{Elementary 
topos}.\footnote{See \cite{shulman2018,shulman2019}.}  

\section{Homotopical\index{Homotopy theory} Categories\index{Category theory}}

The above mentioned elementary topoi\index{Topos}\index{Elementary topos} 
$\mathbf{Set}$, $\hat{\mathcal{C}}=\mathbf{Set}^{\mathcal{C}^{\mathrm{op}}}$, 
and $\mathbf{Sh}(X)$ are so-called Grothendieck topoi\index{Grothendieck 
topos}\index{Topos}, which Alexander Grothendieck\index{Grothendieck, 
Alexander} introduced to define new topologies in algebraic geometry.
Grothendieck topoi\index{Grothendieck topos}\index{Topos} are special cases of  
elementary topoi\index{Topos}\index{Elementary topos} in which the notion of an 
open set is replaced by a categorical\index{Category theory} 
notion.\footnote{For Grothendieck topoi see \cite{illusie2004}.} The usual 
category\index{Category theory} $\mathbf{Top}$ is lacking properties of an 
elementary topos\index{Elementary topos}\index{Topos}, like the existence of 
exponential objects\index{Exponential object}, and therefore it is advantagous 
to pass to the category\index{Category theory} $\mathbf{sSet}$ of simplicial 
sets\index{Simplicial set}.\footnote{Exponential objects exist in the category 
$\mathbf{CGHaus}$ of compactly generated Hausdorff spaces.} The traditional 
definition of topological spaces\index{Topology} is not suitable for many 
situations in algebraic and arithmetic geometry. Therefore,  
Alexander Grothendieck\index{Grothendieck, Alexander} and others have invented 
smarter definitions of topologies\index{Topology}, which result in 
further categorical\index{Category theory} variants of $\mathbf{Top}$. 

Such categories\index{Category theory} are called 
homotopical\index{Homotopy theory} categories and should resemble 
the characteristics of an infinity structure\index{Category theory} which 
$\mathbf{Top}$ carries. By this we mean that, in addition to the usual objects 
and $1$-morphisms, there are $2$-morphisms, hence homotopies\index{Homotopy 
theory} $h$ between continuous mappings 
\[
f,g \colon X \xlongrightarrow{~~~~~~} Y, 
\]
which are defined by a continuous mapping 
\[
h \colon [0,1] \times X \xlongrightarrow{~~~~~~} Y.  
\]
Furthermore, there are homotopies\index{Homotopy theory} between 
homotopies\index{Homotopy theory} etc., resulting in typical 
characteristics of an $(\infty,1)$-category\index{Category theory}. 

Formally, the correct setup in which these observations make sense is given by 
model categories\index{Model category} of Daniel Quillen\index{Quillen, 
Daniel}. How did this idea come about? Many recent developments in category 
theory\index{Category theory} in mathematics can be traced back to a nearly 
$600$-page manuscript \enquote{A la poursuite des champs}\footnote{See 
\cite{grothendieck}.} by Alexander Grothendieck\index{Grothendieck, Alexander}, 
which began with a letter to Daniel Quillen\index{Quillen, Daniel} dated 
February 19th 1983. This manuscript greatly influenced the further development 
of homotopy theory\index{Homotopy theory} and (higher) category 
theory\index{Category theory}. 

Model categories\index{Model category} were introduced by 
Daniel Quillen\index{Quillen, Daniel} and his school and include, in addition 
to weak equivalences\index{Equivalence}, two more classes of morphisms, called 
fibrations and cofibrations. Closely related are the path categories\index{Path 
category} of Kenneth Brown\index{Brown, Kenneth} and the \enquote{tribes} of 
Andr\'e Joyal\index{Joyal, Andr\'e}. The latter concepts arise from fibred 
objects in model categories\index{Model category}, i.e., objects $A$ for which 
the morphism $A \longrightarrow \ast$ to the terminal object is a fibration. 
The structure-preserving equivalences\index{Equivalence} of model 
categories\index{Model category} -- and thus of homotopical 
categories\index{Category theory} -- are given by so-called Quillen 
equivalences\index{Equivalence}. Homotopical categories\index{Category theory} 
are closely related 
to $(\infty,1)$-categories\index{Infinity category} since with the help of 
Dwyer-Kan localisation, homotopical categories\index{Category theory} can be 
turned into $(\infty,1)$-categories\index{Infinity category} and higher 
elementary topoi\index{Elementary topos}\index{Topos}.\footnote{For model 
categories, path categories, and tribes see \cite{brown1973,joyal2017,quillen}. 
The Dwyer-Kan localisation is also called simplicial localisation. It localises 
with respect to weak equivalences. See \cite{dwyerkan}. For 
$(\infty,1)$-categories and higher elementary topoi 
see \cite{cisinski2015,cisinski2025}, \cite{lurie2009}, \cite{riehlverity}, 
\cite{shulman2018}, and \cite{shulman2019}.}

\medskip
We want to explain path categories\index{Path category} in more detail. In 
them, two types of morphisms are distinguished, called fibrations and weak 
equivalences\index{Equivalence}, and in which 
path objects exist. Their properties are modelled after the topological case. 
There are thus some axioms that these types of morphisms 
must fulfil. For example, isomorphisms\index{Isomorphism} are always 
weak equivalences\index{Equivalence} and the 2-out-of-3 
rule applies: if two of the three morphisms $f,g,g \circ f$ are weak 
equivalences\index{Equivalence}, so is the third. The composition of 
two fibrations and all isomorphisms\index{Isomorphism} are fibrations. In 
homotopy theory\index{Homotopy theory} a morphism $p \colon E \longrightarrow 
B$ of topological spaces\index{Topology} is called a  Hurewicz fibration if it 
satisfies a certain lifting property for all topological\index{Topology} spaces 
$X$, resp. only for CW-complexes in the case of Serre fibrations (see 
Fig.~\ref{fig:cat_lifting}).

\begin{figure}[ht!]
\begin{tikzpicture}
  \matrix (m) [matrix of math nodes,row sep=3em,column sep=4em,minimum 
width=2em]
  {
    X \times \{0\} & E \\
      X \times I & B \\};
  \path[-stealth]
    (m-1-1) edge node [left] {$\mathrm{id} \times \{0\}$} (m-2-1)
            edge node [above] {$f$} (m-1-2)
    (m-2-1.east|-m-2-2) edge node [below] {$F$}
            node [above] {} (m-2-2)
    (m-1-2) edge node [right] {$p$} (m-2-2)
    (m-2-1) edge [dashed] node [above] {$\tilde F$} (m-1-2);
\end{tikzpicture}
\caption{\label{fig:cat_lifting}Lifting property of a Hurewicz or a Serre 
fibration.}
\end{figure}

Path categories\index{Path category} possess path objects. This includes for 
each object $X$ a path object $\mathrm{Path}(X)$ together with a composition of 
mappings as indicated in the diagram in Fig.~\ref{fig:path_object}. Here we 
require that the morphism  
\[
r \colon X \xlongrightarrow{~~~~~~} \mathrm{Path}(X)
\]
is a weak equivalence\index{Equivalence}, and the morphism 
\[
p \colon \mathrm{Path}(X) \xlongrightarrow{~~~~~~} X \times X
\]
is a fibration. 

\begin{figure}[ht!]
\begin{tikzpicture}
  \matrix (m) [matrix of math nodes,row sep=3em,column sep=4em,minimum 
width=2em]
  {
      & \mathrm{Path}(X) \\
      X & X \times X \\};
  \path[-stealth]
    (m-2-1.east|-m-2-2) edge node [below] {$(\mathrm{id},\mathrm{id})$}
            node [above] {} (m-2-2)
    (m-1-2) edge node [right] {$p$} (m-2-2)
    (m-2-1) edge node [dashed, above] {$r$} (m-1-2);
\end{tikzpicture}
\caption{\label{fig:path_object}A path object for $X$ is a fibration 
$\mathrm{Path}(X)$ over $X \times X$ together with a morphism $X \rightarrow 
\mathrm{Path}(X)$ which is a weak equivalence\index{Equivalence}.}
\end{figure}

With the help of path objects, an abstract concept of 
homotopy\index{Homotopy theory} can be defined, in which two mappings $f,g 
\colon Y \longrightarrow X$ are homotopic\index{Homotopy theory} if there is a 
lift $h$ (see Fig.~\ref{fig:homotopy}).

\begin{figure}[ht!]
\begin{tikzpicture}
  \matrix (m) [matrix of math nodes,row sep=3em,column sep=4em,minimum 
width=2em]
  {
      & \mathrm{Path}(X) \\
      Y & X \times X \\};
  \path[-stealth]
    (m-2-1.east|-m-2-2) edge node [below] {$(f,g)$}
            node [above] {} (m-2-2)
    (m-1-2) edge node [right] {$p$} (m-2-2)
    (m-2-1) edge [dashed] node [above] {$h$} (m-1-2);
\end{tikzpicture}
\caption{\label{fig:homotopy}A homotopy\index{Homotopy theory} between 
morphisms $f,g \colon Y \rightarrow X$ can be viewed as a lifting $h$ of the 
map $(f,g) \colon Y \rightarrow X \times X$ to the path object 
$\mathrm{Path}(X)$.}
\end{figure}

Path objects exist in the case of the category\index{Category 
theory} of topological spaces\index{Topology}, because the space 
$\mathrm{Path}_a(X)$ of paths with a fixed starting point $a$ is contractible 
to the point path by reparametrisation and therefore the natural mapping 
(constant path)
\[
r \colon X \xlongrightarrow{~~~~~~} \mathrm{Path}(X)
\]
is a homotopy equivalence\index{Equivalence}\index{Homotopy theory} for every 
topological space\index{Topology} $X$. The projection mapping to the starting 
and ending point
\[
p \colon \mathrm{Path}(X) \xlongrightarrow{~~~~~~} X \times X
\]
is a fibration, in which all fibres correspond to the path spaces with a fixed 
starting and ending point.

\section{Fundamental Infinity Groupoid\index{Infinity groupoid}} 

The fundamental infinity groupoid\index{Infinity groupoid} $\Pi_\infty(X)$ of a 
topological\index{Topology} space $X$ with enough real paths consists of all 
paths in $X$ without considering their homotopy classes\index{Homotopy theory}. 
It is a generalisation of the fundamental groupoid\index{Fundamental groupoid} 
$\Pi_1(X)$. As a model\index{Model theory} for it the simplicial 
set\index{Simplicial set} $\mathrm{Sing}_\bullet(X)$ is used, which is created 
by continuous mappings $\Delta_n \longrightarrow X$ for all $n$. 
Homotopies\index{Homotopy theory} between paths are formed by continuous 
mappings 
\[
h \colon [0,1] \times [0,1] \xlongrightarrow{~~~~~~} X,
\]
i.e., by continuous mappings of a square $\Box^2$ to $X$. However, a square can 
be divided into two triangles by its diagonal. This trick ensures that 
homotopies\index{Homotopy theory} and their higher generalisations can be 
traced back to continuous mappings $\Delta_n \longrightarrow X$.

The object $\Pi_\infty(X)$ is an example of an 
$(\infty,0)$-category\index{Category theory}, whose objects consist of the 
points of $X$ and the morphisms consist of the paths. The $2$-morphisms are 
data that encode the homotopies\index{Homotopy theory} between the paths, which 
are used for the calculation rules. These are the rules for the inverse path 
$p^{-1}$, for the composition of paths and for the associative law, all of 
which only apply up to homotopy\index{Homotopy theory}. The $3$- and 
higher-dimensional objects are generated by the homotopies\index{Homotopy 
theory} between the homotopies\index{Homotopy theory} and so on.\footnote{See 
\cite{lurie2008} for an easily accessible introduction.}

The model\index{Model theory} $\mathrm{Sing}_\bullet(X)$ is a Kan 
complex\index{Kan, Daniel}. Such  structures are named after Daniel 
Kan\index{Kan, Daniel} and form abstractions of the combinatorics inside the 
models\index{Model theory} $\mathrm{Sing}_\bullet(X)$. They are defined by the 
existence\index{Existence} of fillings in the form of an 
$n$-simplex\index{Simplex} $\Delta_n$ for those subsets\index{Simplicial set} 
$\Lambda_n^i$ of $\partial \Delta_n$ where the $i$-th side was removed from 
$\partial \Delta_n$ for  $0 \le i \le n$ (see Fig.~\ref{fig:kan_filling}). This 
condition ensures that all morphisms of dimension $n \ge 1$ are invertible. 
Therefore, Kan complexes\index{Kan, Daniel} are models\index{Model theory} of 
$(\infty,0)$-categories\index{Infinity category} or infinity 
groupoids\index{Infinity groupoid}. 

\begin{figure}[ht!]
\begin{tikzpicture}
\node (a) at (0,-2) {$\Delta_n$};
\node (b) at (0,0) {$\Lambda_n^i$};
\node (c) at (3,0) {$X$};
\draw[->] (b) to node[scale=.7] [above] {} (a);
\draw[->] (b) to node[scale=.7] [above] {$x$} (c);
\draw[->,dashed] (a) to node[scale=.7,dashed] (x) [below] {} (c);
\end{tikzpicture}
\caption{\label{fig:kan_filling}A filling in a Kan complex is given by the 
extension of a map $x \colon \Lambda_n^i \rightarrow X$ to a map  
$\Delta_n \rightarrow X$.}
\end{figure}

In weak Kan complexes\index{Kan, Daniel}, it is required that only the inner 
subsets\index{Simplicial set} $\Lambda_n^i$ for $1 \le i \le n-1$ have 
fillings. As a result, in weak Kan complexes\index{Kan, Daniel} all morphisms 
of dimension $n \ge 2$ are invertible and weak Kan complexes\index{Kan, Daniel} 
are models\index{Model theory} of $(\infty,1)$-categories\index{Infinity 
category}.

\medskip There exists a structure of a Quillen\index{Quillen, Daniel} model 
category\index{Model category} on $\mathbf{sSet}$ such that the Kan 
fibrations\index{Kan, Daniel} between Kan complexes\index{Kan, Daniel}, which 
are also defined by filling property (see Fig.~\ref{fig:kan_fibration}), 
correspond precisely to 
the Serre-fibrations in the Quillen\index{Quillen, Daniel} model 
category\index{Model category} structure on $\mathbf{Top}$ under a suitable 
Quillen equivalence\index{Equivalence} of model structures. 

\begin{figure}[ht!]
\begin{tikzpicture}
\node (a) at (0,-2) {$\Delta_n$};
\node (b) at (0,0) {$\Lambda_n^i$};
\node (c) at (3,0) {$X$};
\node (d) at (3,-2) {$Y$}; 
\draw[->] (b) to node[scale=.7] [above] {} (a);
\draw[->] (b) to node[scale=.7] [above] {$x$} (c);
\draw[->,dashed] (a) to node[scale=.7,dashed] [below] {} (c);
\draw[->] (a) to node[scale=.7] [below] {$y$} (d);
\draw[->] (c) to node[scale=.7] [right] {$f$} (d);
\end{tikzpicture}
\caption{\label{fig:kan_fibration}A fibration of Kan complexes is defined by a 
lifting property with respect to the inclusion $\Lambda_n^i \subseteq 
\Delta_n$.}
\end{figure}

\section{Category Theory\index{Category theory} and the Concept of Space} 

The development of a most general concept of space\index{Topology} and 
number is a goal within and outside of mathematics. One of the 
prevailing trends is the categorification\index{Category theory} 
of mathematical structures. In particular, $(\infty,0)$- and  
$(\infty,1)$-categories\index{Infinity category} are the crucial 
structures with which we can generalise the concept of space and number. 

As we have seen, the fundamental group $\pi_1(X,\ast)$ is usefully extended to 
the infinity groupoid\index{Infinity groupoid} $\Pi_\infty(X)$ -- an 
$(\infty,0)$-category\index{Infinity category}. Alexander 
Grothendieck\index{Grothendieck, Alexander} conjectured 
that every infinity groupoid\index{Groupoid}\index{Infinity groupoid} 
is equivalent to one of the form $\Pi_\infty(X)$. Assuming this 
homotopy hypothesis\index{Homotopy hypothesis} by 
Grothendieck\index{Grothendieck, Alexander}, it follows that the 
homotopy types\index{Homotopy theory} of topological 
spaces\index{Topology} are given by the totality of all 
infinity groupoids\index{Infinity groupoid}, thus 
by an $(\infty,1)$-category\index{Infinity category}, 
which can be constructed from model categories\index{Model category} or 
path categories.\index{Path category}\footnote{Such $(\infty,1)$-categories 
result from a Dwyer-Kan localisation after the weak equivalences in model 
categories or the related path categories. Here, certain -- previously 
non-invertible -- morphisms are inverted and they become equivalences. This 
generalises the Gabriel-Zisman localisation of categories and originates in the 
construction of fractions, i.e., the localisation of a ring. See 
\cite{dwyerkan} and \cite[Ch. 7]{cisinski2015}.} This observation shows that 
infinity groupoids\index{Infinity groupoid} or 
$(\infty,0)$-categories\index{Category theory} and 
homotopy theory\index{Homotopy theory} are closely related.

We have already seen that the usual definition of $\Pi_\infty(X)$
only makes sense for certain topological spaces\index{Topology} $X$. For 
various interesting situations, a different approach is better to arrive at 
an infinity groupoid\index{Infinity groupoid} $\Pi_\infty(X)$. In algebraic 
geometry, the underlying topological spaces\index{Topology} 
almost never have real paths. However, there is an alternative definition of 
a corresponding simplicial object\index{Simplicial set} $\Pi_\infty(X)$, 
if $X$ is an algebraic scheme,\footnote{Schemes are generalisations of 
algebraic varieties that locally look like prime ideal spectra 
$\mathrm{Spec}(R)=\{ \mathfrak{p} \mid \mathfrak{p} \text{ is a prime ideal in 
} R \}$.} an algebraic space\footnote{Algebraic spaces are 
generalisations of schemes and are special cases of stacks. They 
naturally arise as quotients of schemes $U$ by 
identification via an equivalence relation $R \subseteq U \times U$, 
where the projection maps $R \longrightarrow U$ are each \'etale, i.e., 
in particular unbranched.} or an algebraic stack.\footnote{Stacks are 
general space concepts that exist in algebraic, topological and 
differentiable form. Grothendieck\index{Grothendieck, Alexander}
and Artin\index{Artin, Michael} developed algebraic stacks, which generalise 
algebraic schemes. They are more suitable as moduli spaces for algebraic 
objects and for forming quotients than schemes. A stack $\mathcal{X}$ is a 
category that is fibred over a base category $\mathcal{B}$ in the form of a 
functor $\pi \colon \mathcal{X} \longrightarrow \mathcal{B}$, where the fibres 
are groupoids. In addition, so-called descent data must be satisfied, which 
state that the family $\pi$ forms a $2$-sheaf (more precisely a 
$(2,1)$-sheaf) of groupoids over $\mathcal{B}$ in a suitable Grothendieck 
topology. In algebraic geometry, $\mathcal{B}$ is usually given by the category 
$\mathcal{B}=(\mathbf{Sch}/S)_{\mathrm{fppf}}$ of schemes over a fixed base 
scheme $S$ and equipped with the fppf-topology. A scheme $X$ over $S$ can be 
considered as a stack, where $\mathcal{X}=\mathbf{Sch}/X$ and under $\pi$ a 
scheme $U/X$ is mapped to $U/S$. The fibre over $U/S$ consists as a category of 
the objects $\mathrm{Hom}_S(U,X)$ with $\mathrm{id}: U \to U$ as the only 
morphism. All fibres are thus groupoids. In topology, $\mathcal{B}$ is 
typically the category $\mathbf{Top}$ and in differential geometry a category 
of manifolds.} It can be constructed using the \'etale homotopy theory and 
suitable \'etale -- especially unbranched -- coverings of $X$.\footnote{See 
\cite{artinmazur1969,carchedi2020,friedlander1982}.}

Gauge theories of mathematical physics often use differentiable stacks, on 
which gauge transformations operate in the form of groupoids\index{Groupoid} 
and preserve physical properties, such as the Lagrange functions in 
Yang-Mills theories.\footnote{See \cite{yangmills1954}.} Higher $n$-stacks, 
derived algebraic spaces, locales and condensed spaces are even more general 
space theories, which also play a role in mathematical physics and all of 
which have connections with infinity structures\index{Infinity category}.  
The concept of equivalence\index{Equivalence} plays a crucial role in 
them.\footnote{See \cite{anelcat2021,carchedi2020,maclane_moerdijk1992, 
scholze2021}.} Condensed spaces are particularly promising among these 
variants, as they unite many space concepts. They are defined as 
functors\index{Functor}
\[
T \colon \text{\bf Profinite Sets}^\mathrm{op} \xlongrightarrow{~~~~~~} 
\text{\bf Set},
\]
which satisfy the sheaf condition with respect to a suitable topology. Every 
topological space can be represented as a condensed space. Profinite sets are 
limits of so-called directed systems of finite sets. Examples include the 
$p$-adic integers $\mathbb{Z}_p=\varprojlim \mathbb{Z}/p^n\mathbb{Z}$, which 
are the limit of the finite sets $\mathbb{Z}/p^n\mathbb{Z}$ for $n \ge 1$, or 
the numbers $\hat{\mathbb{Z}}=\varprojlim \mathbb{Z}/m\mathbb{Z}$ as the limit 
of all sets $\mathbb{Z}/m\mathbb{Z}$ for $m \ge 1$. 

If we want to study topological spaces\index{Topology} $X$ not only up to 
homotopy\index{Homotopy theory} but up to homeomorphism\index{Homeomorphism}, 
it is advisable to replace the space $X$ with the elementary topos\index{Topos} 
$\mathbf{Sh}(X)$, because $X$ can be reconstructed with it, provided $X$ is a 
Hausdorff space.\footnote{Instead of Hausdorff space, the property 
\enquote{sober} is sufficient. The subobject classifier $\Omega$ satisfies 
$\Omega(U)=\{V \mid V \subset U \; \text{open}\}$ and thus the union of all 
open sets in $X$ can be reconstructed. If $X$ is \enquote{sober}, this allows 
the determination of $X$, see \cite[Ch. IX, Sec. 3]{maclane_moerdijk1992}.} 
Even more flexible is the infinity topos\index{Topos} $\mathbf{Sh}_\infty(X)$ 
in the form of an $(\infty,1)$-category\index{Infinity category} of 
infinity sheaves\index{Sheaf} of infinity groupoids\index{Infinity groupoid} on 
$X$. Such generalised spaces can, despite their degree of abstraction, be used 
as efficient replacements for classical topological spaces\index{Topology} 
and yet preserve the properties of traditional space concepts. 

All different types of generalised spaces can be assigned homology 
groups\index{Homology group} and other invariants\index{Invariant}. We have 
already encountered the singular homology groups\index{Homology group} 
$H_n(X,\mathbb{Z})$ with integer coefficients of a topological 
space\index{Topology} or a manifold $X$. Such homology 
groups\index{Homology group}, also those with coefficients other 
than $\mathbb{Z}$, contain primary invariants\index{Invariant} of $X$ and of 
geometric objects on $X$, such as the Chern classes of $X$ or of 
sheaves\index{Sheaf} $\mathcal{E}$ on $X$. These invariants\index{Invariant} 
express geometric curvature properties of $X$ or of $\mathcal{E}$ and thereby 
classify these objects in a certain way. Refinements of the homology 
groups\index{Homology group}, such as the differential 
homology\index{Homology group} or the 
Deligne homology\index{Homology group}\index{Deligne, Pierre}, allow the 
definition of secondary invariants\index{Invariant}, which refine the primary 
invariants\index{Invariant} and describe some 
properties of spaces and sheaves more precisely. In physics, 
Chern classes and their generalisations play an important role. The simplest 
example of this is the magnetic flux in Maxwell's theory\index{Maxwell, 
James Clerk}, which can be expressed as Chern class $c_1$ of a $U(1)$-bundle, 
as Paul Dirac\index{Dirac, Paul} discovered.

There are algebraic versions of the homology groups\index{Homology group} and 
the homotopy groups\index{Homotopy theory} of algebraic schemes, which can be 
best understood through the $\mathbb{A}^1$-homotopy theory discovered by 
Vladimir Voevodsky\index{Voevodsky, Vladimir} and others. A  
well-known example among these are the classical Chow groups $CH_i(X)$ of a 
scheme $X$. Such groups are also called motivic homology groups\index{Homology 
group} and homotopy groups\index{Homotopy theory}, because the motives imagined 
by Grothendieck\index{Grothendieck, Alexander} are abstractions of the 
homological properties of algebraic schemes. 
Voevodsky\index{Voevodsky, Vladimir} proved conjectures of arithmetic, such as 
the Milnor conjecture and the Bloch-Kato conjecture, using the 
$\mathbb{A}^1$-homotopy theory, and was awarded the 
Fields Medal for this in 2002. Such results connect periods, motives and 
other invariants\index{Invariant} of algebraic objects with 
number-theoretical questions, partly in generalisation of the 
Riemann hypothesis\index{Riemann hypothesis}.\footnote{See 
\cite{hms2017,ms2013,soule}.}

\chapter{Type Theory\index{Type theory}}

From the discovery of the antinomies\index{Antinomy} of set theory\index{Set 
theory}, Bertrand Russell\index{Russell, Bertrand} concluded that a hierarchy 
of mathematical objects is necessary. As a solution, he invented type 
theory\index{Type theory}. His ideas were incorporated into the book 
\enquote{Principia mathematica} with Alfred North 
Whitehead\index{Whitehead, Alfred North}. The presentation therein seemed 
complicated and initially did not catch on. Leon
Chwistek\index{Chwistek, Leon} and Frank P. Ramsey\index{Ramsey, Frank} 
attempted to rectify some of these problems and created the simple 
type theory\index{Type theory}. However, it was an article by 
Alonzo Church\index{Church, Alonzo} that made this approach widely known using 
the $\lambda$-calculus.\footnote{See 
\cite{church1940,lambekscott1986,russell}.} 

Per Martin-L\"of\index{Martin-L\"of, Per} has further developed the simple 
type theory\index{Type theory} from around 1971 and built the 
intuitionistic\index{Intuitionism} dependent 
type theory\index{Type theory} as a possible foundation of mathematics 
alongside set theory\index{Set theory} and 
category theory\index{Category theory}. He has also written 
philosophical works in which he linked Kant's theory of 
synthetic judgements\index{Synthetic judgement} with 
type theory\index{Type theory} and G\"odel's\index{Goed@G\"odel, Kurt} 
incompleteness theorems\index{Incompleteness theorem}.\footnote{See 
\cite{martin-loef1984,martin-loef1994}.} A first-order fragment of dependent 
type theory\index{Type theory} was developed by Mih\'aly Makkai\index{Makkai, 
Mih\'aly} around 1995. 

The concept of equality\index{Identity}\index{Equality} plays a significant 
role in dependent type theory\index{Type theory}. In the initial versions, this 
was extensionally\index{Extensionality}\footnote{A property is 
intensional when the content or the essence play a role and  
extensional when it depends on the external scope.} formulated 
and thus narrowly defined. Later versions included an 
intensional\index{Intensionality} identity type\index{Identity type}, which was 
Martin-L\"of's\index{Martin-L\"of, Per} most influential idea. The pair of 
concepts intensional versus extensional is related to 
Frege's\index{Frege, Gottlob} distinction between sense and meaning. 

Around 1994, Thomas Streicher\index{Streicher, Thomas} and Martin 
Hofmann\index{Hofmann, Martin} constructed remarkable models\index{Model 
theory} of Martin-L\"of's\index{Martin-L\"of, Per} theory in which the types 
form groupoids\index{Groupoid} and hence have interesting identity 
types\index{Identity type}. Vladimir Voevodsky\index{Voevodsky, Vladimir}, 
Steve Awodey\index{Awodey, Steve}, Michael Warren\index{Warren, Michael} and 
others later extended the theory of Martin-L\"of\index{Martin-L\"of, Per} to 
homotopy type theory\index{Homotopy type theory} (alias univalent 
foundations), which gives dependent type theory\index{Type theory} a 
homotopically\index{Homotopy theory} motivated 
topological\index{Topology} interpretation\index{Interpretation}. This fits 
well with the character of Martin-L\"of's\index{Martin-L\"of, Per} type 
theory\index{Type theory}, leads to a more intuitive view, and is an 
inspiration for the further development of type theory\index{Type theory}. In 
particular, the introduction of (higher) inductive types has opened new 
perspectives. Another consequence of Martin-L\"of's\index{Martin-L\"of, Per} 
type theory\index{Type theory} and its variants is the possibility of machine 
verification of proofs. This opens up interesting future perspectives for 
mathematics and its applications.\footnote{See \cite{awodeywarren}, 
\cite{hofmannstreicher}, \cite{voevodsky2006}, and \cite{voevodsky2013}.} 

\section{Types and Terms\index{Type theory}} 

A fundamental notation of type theory\index{Type theory} forms 
judgements\index{Judgement} of the form 
\[
a:A. 
\]
By this we understand a term $a$ that belongs to the type $A$. Two terms 
$a,b : A$ are called judgementally equal, denoted 
\[
a \equiv b, 
\]
if the designation $b$ is just an alternative designation of $a$. This form 
of judgemental equality is the closest possible form of 
equality\index{Equality}\index{Identity} in type theory\index{Type theory}. 

In contrast to set theory\index{Set theory}, the 
judgement\index{Judgement} $a:A$ avoids intensional\index{Intensionality} 
complications in type theory\index{Type theory}. If a term $a$ is given, then
is usually also the associated type $A$ that $a$ inhabits, i.e., each term $a$ 
has a well-defined type. In set theory\index{Set theory} this is quite 
different, because an element $a$ can be an element of several sets\index{Set 
theory}. The set-theoretical operations $\cap$ and $\cup$ are obviously  
intensional\index{Intensionality}. Therefore, they are not present in this form 
in type theory\index{Type theory}.

There are important elementary types, such as the empty type $\mathbf{0}$, 
which is not inhabited, or the single-element type $\mathbf{1}$, which is 
inhabited by exactly one term usually denoted as $\ast$. In addition, there 
is the two-element Boolean type $\mathbf{Bool}$, whose terms are the classical 
truth values\index{Truth} $\top$ (true) and $\bot$ (false).

\begin{figure}[ht!]
\begin{tikzpicture}
\draw[color=gray,ultra thick,domain=0:1.45] plot (1+\x,\x^2);
\draw[color=gray,ultra thick,domain=0:1.45] plot (1-\x,-\x^2);
\draw[color=gray,ultra thick,domain=0:1.21] plot (4+\x,\x^4);
\draw[color=gray,ultra thick,domain=0:1.21] plot (4-\x,-\x^4);
\draw[color=gray,ultra thick,domain=0:1.13] plot (7+\x,\x^6);
\draw[color=gray,ultra thick,domain=0:1.13] plot (7-\x,-\x^6);
\end{tikzpicture}
\caption{\label{fig:type_family}Fibers $B(x)$ of a family of types for three 
values of $x:A$.}
\end{figure}

In dependent type theories\index{Type theory}, types may additionally depend on 
parameters (i.e., term variables) and in this way form a family of types. Given 
a type $A$, a family $B$ is given by types $B(x)$ for all terms $x:A$. We 
imagine that a family is a geometric deformation in which the types $B(x)$ vary 
as fibers depending on terms $x:A$ (see Fig.~\ref{fig:type_family}). Such 
geometric imaginations help our intuition a lot but also may lead astray in 
certain situations.

\section{Universes\index{Universe} and Families of Types\index{Type theory}}

As in set theory\index{Set theory}, the impredicative 
assumption of a type of all types creates a paradox\index{Paradox} analogous to 
Russell's antinomy\index{Antinomy}, as noted by Jean-Yves Girard\index{Girard, 
Jean-Yves}. Martin-L\"of\index{Martin-L\"of, Per} and others have developed 
predicative and impredicative 
paradox-free versions over the past $50$ years using type-theoretical 
universes\index{Universe}. Grothendieck universes\index{Universe} were 
introduced by Alexander Grothendieck\index{Grothendieck, Alexander} to address 
set-theoretical problems in category theory\index{Category theory}. They are 
based on the assumption of strongly inaccessible cardinal 
numbers\index{Cardinal number}.\footnote{Universes were originally an idea of 
Paul Bernays\index{Bernays, Paul}, see \cite{grayson}. For type-theoretical 
universes see \cite[Sec. 1.3]{voevodsky2013}.}

The universes\index{Universe} $U$ postulated in type theory\index{Type theory} 
themselves form types that are inhabited by certain types $A:U$ as terms. 
Universes\index{Universe} are often iterated in a given type theory\index{Type 
theory} in ascending order of size in a finite or even infinite sequence
\[
U:U':U'':U''': \cdots 
\]
each representing a term in the next universe, i.e., $U$ is a type within the 
universe\index{Universe} $U'$, this in turn in $U''$ and so on. Such a sequence 
of universes\index{Universe} generates a high degree of flexibility in type 
theory\index{Type theory}. As a rule, at least one additional universe $U'$ is 
required with $U:U'$ to enable certain type constructions in $U$.

Useful properties of universes\index{Universe} $U$ are that they can first 
define families of types and secondly types become terms $A,B:U$ and can thus 
be compared within $U$. A family of types $B(x)$ with the parameter $x:A$ is 
given under these conditions by a function
\[
B \colon A  \xlongrightarrow{~~~~~~} U,
\]
where $U$ is a given universe\index{Universe} and the types $B(x):U$ are terms 
in $U$. We will introduce functions in type theory\index{Type theory} soon. 

What can be said about the existence\index{Existence} of type-theoretical 
universes\index{Universe} and thus of types? Type theory\index{Type theory} is 
a syntactic\index{Syntax} construct, so the concept of 
existence\index{Existence} for objects in type theories\index{Type theory} is 
-- at least from a from a nominalistic\index{Nominalism} position -- 
meaningless. The crucial question is consistency\index{Consistency}, i.e., the 
contradiction-free nature\index{Contradiction-freeness} of such deductive 
systems\index{Deductive system}. In this context, however, the second 
incompleteness\index{Incompleteness theorem} of  
G\"odel\index{Goed@G\"odel, Kurt} applies.

The existence\index{Existence} of universes\index{Universe} is often discussed 
in a semantic\index{Semantics} context in set-theoretical\index{Set theory} or 
categorical\index{Category theory} models. Axioms based on very large 
transfinite cardinal numbers\index{Cardinal number} offer a solution to the 
questions of 
consistency\index{Consistency}\index{Contradiction-freeness} and 
existence\index{Existence}. However, this answer shifts the problem to the 
question of consistency\index{Consistency}\index{Contradiction-freeness} in 
superordinate deductive systems\index{Deductive system} that can describe such 
numbers, and does not necessarily make it easier. However, in our opinion, even 
in models, existence\index{Existence} is a subordinate -- if not irrelevant --  
question compared to 
consistency\index{Consistency}\index{Contradiction-freeness}.

\section{Judgements and Inference Rules} 

Type theories\index{Type theory} can be understood as deductive 
systems\index{Deductive system} consisting of judgements\index{Judgement} and 
inference rules. These usually follow the formalism of the calculus of 
inductive constructions by Thierry Coquand\index{Coquand, Thierry}, G\'erard 
Huet\index{Huet, G\'erard} and Christine 
Paulin-Mohring\index{Paulin-Mohring, Christine}.\footnote{See 
\cite{coquand1986,coquand1990}.}

As we have seen, proofs in deductive systems\index{Deductive system} arise from 
trees using axioms, judgements and inference rules. The inference rules in type 
theories\index{Type theory} produce from a finite list $\mathcal{J}_1, \ldots, 
\mathcal{J}_n$ of given judgements as premises a new judgement $\mathcal{J}'$. 
An inference rule is thus expressed in the form 
\[
\frac{\mathcal{J}_1 \quad \mathcal{J}_2 \quad \cdots \quad \mathcal{J}_n}
{\mathcal{J}'}.
\]
The basis for judgements are contexts $\Gamma$ consisting of a series of terms
\[
x_0:A_0,x_1:A_1(x_0),\ldots,x_n:A_n(x_0,\ldots,x_{n-1}),  
\]
where the $A_i$ are dependent types and the $x_i$ are variables that stand for 
terms to be substituted. The contexts of a type theory\index{Type theory} $T$ 
form the syntactic (or contextual) category\index{Category theory} 
$\mathcal{S}_T$, which is naturally associated with $T$. Their morphisms arise 
from obvious mappings between contexts which go back to arrows between types.

Type theory\index{Type theory} is based on four fundamental 
judgements\index{Judgement} given a universe\index{Universe} $U$ and using 
contexts $\Gamma$ (see Fig.~\ref{fig:type_judgements}). 

\begin{figure}[ht!]
\begin{tabular}{l|l}
$\Gamma \vdash A:U$          & Type \\
$\Gamma \vdash a:A$          & Term \\
$\Gamma \vdash a \equiv_A b$ & Judgemental equality of terms\\
$\Gamma \vdash A \equiv_U B$ & Judgemental equality of types
\end{tabular}
\caption{\label{fig:type_judgements}The four basic judgements of type 
theory\index{Type theory}.}
\end{figure} 

There are a few elementary inference rules which are useful for the 
rewriting of judgements. Among the tautological rules are those that express 
that judgemental equality of terms and types is an equivalence relation, for 
example 
\[
\frac{\Gamma \vdash A}{\Gamma \vdash A \equiv A}, \quad \frac{\Gamma \vdash 
A\equiv B}{\Gamma \vdash B \equiv A}, \quad \text{and} \quad 
\frac{\Gamma \vdash A \equiv B \quad \Gamma \vdash B \equiv C}{\Gamma \vdash A 
\equiv C}.
\]
Other rules express that each term belongs to a type, for example 
\[
\frac{\Gamma \vdash x:A}{\Gamma \vdash A}.  
\]
In some sense, this rule can be inverted to obtain the following rule which 
introduces a new variable (resp. a constant family) and is often called 
\enquote{weakening rule}:  
\[
\frac{\Gamma \vdash A \quad \Gamma \vdash b:B}{\Gamma, \; 
x:A \vdash b:B} \quad  \text{ resp. } \quad 
\frac{\Gamma \vdash A \quad \Gamma \vdash b \equiv b':B}{\Gamma, \; 
x:A \vdash b \equiv b':B}. 
\]
The rule of the \enquote{generic element} is a further weakening rule 
which implies a conclusion of the form
\[
\frac{\Gamma \vdash A}{\Gamma, \; x:A \vdash x:A},
\]
and thus introduces the identity\index{Identity}\index{Equality} function.

\medskip
Another important rule is the substitution of term variables, i.e., switching 
from one term name to another one. It can be used to elegantly treat the fiber 
of a family $B(x)$ with parameter $x:A$ at a base point $a:A$. The 
inference rule looks like 
\[
\frac{\Gamma \vdash a:A \quad \Gamma, x:A \vdash B(x)}{\Gamma \vdash 
B[\frac{a}{x}]}. 
\]
Instead of $B[\frac{a}{x}]$, we occasionally write $B(a)$, even if this is 
something different. There are numerous additional rules which rewrite contexts 
in a useful way. All such elementary rules can be found in the 
literature.\footnote{See \cite[App. A]{voevodsky2013} and \cite[Ch. 1]{rijke}.} 

\medskip 
In addition to these elementary cases, inference rules can arise from four 
sorts of rules for every type species (see Fig.~\ref{fig:type_inference}). 

\begin{figure}[ht!]
\begin{tabular}{l|l}
A-Form & Formation rule for the 
generation of a type A \\
A-Intro & Introduction rule for terms of type A \\
A-Elim & Elimination rule for a type A \\
A-Comp & Computation rule for the application of the \\ 
& elimination rule on constructors of A 
\end{tabular}
\caption{\label{fig:type_inference}The four inference rules of type 
theory\index{Type theory}.} 
\end{figure}

We will list and explain them in each case. Some of the following rules and 
constructions are known in the older literature under the names 
$\lambda$-abstraction, $\alpha$-conversion (renaming of bound variables), 
$\eta$-conversion (function extensionality) and $\beta$-reduction (substitution 
of arguments).

\section{Function Type\index{Type theory}}

The most important basic type construction for any two types $A,B$ is the 
extensional\index{Extensionality} function type (or mapping type)
\[
(A \xlongrightarrow{~~~~~~} B)
\]
which consists of functions (or mappings) $f \colon A \longrightarrow B$. 
This is a primitive type, i.e., it is not definable from the remaining 
structures of the theory. In set theory\index{Set theory}, the 
function type is not primitive, because functions can be defined by their graph 
as a subset.

\medskip 
Each function $f$ corresponds to a rule, with which for a given $a:A$ a value 
$f(a):B$ can be obtained. One occasionally introduces functions by writing 
\[
\lambda x.f(x)
\]
using the notation from the $\lambda$-calculus of Church\index{Church, Alonzo}. 
Which function values $f(x)$ can occur is implicitly determined by the rules 
of the given dependent type theory\index{Type theory}. The four inference rules 
for the function type $(A \longrightarrow B)$ are described in 
Fig.~\ref{fig:functiontype} using $\lambda$-notation.

\begin{figure}[ht!] 
\fbox{
\begin{minipage}{0.8\textwidth}
\begin{align*}
&\text{F}-\text{Form:} && \frac{\Gamma \vdash A:U \quad  
\Gamma \vdash B:U}{\Gamma \vdash (A \longrightarrow B):U} \cr
&\text{F}-\text{Intro:} && 
\frac{\Gamma, x:A \vdash f(x):B(x)} 
{\Gamma \vdash \lambda x.f(x): (A \longrightarrow B)} \cr
&\text{F}-\text{Elim:} &&
\frac{\Gamma \vdash f:(A \longrightarrow B) \quad  
\Gamma \vdash a:A}{\Gamma \vdash f(a):B[\frac{a}{x}]} \cr
&\text{F}-\text{Comp:} && 
\frac{\Gamma, x:A \vdash f(x):B \quad \Gamma \vdash a:A} 
{\Gamma \vdash (\lambda x.f(x))(a) \equiv f[\frac{a}{x}]:B[\frac{a}{x}]}
\end{align*}
\end{minipage}
}
\caption{\label{fig:functiontype}The inference rules of the function type.}
\end{figure}

While the formation rule produces a function type for two types $A,B$, the 
introduction rule defines a function $f$ by giving its values $f(x)$ at 
arguments $x:A$. The elimination rule describes the evaluation of a function 
$f$ at an argument $a$. The corresponding computation rule computes the value 
of a function in $\lambda$-notation.

\medskip 
For every type $A$, the identical function  
\[
{\rm id}_A \colon A \xlongrightarrow{~~~~~~} A
\]
is a canonical term in the type $(A \longrightarrow A)$. It can be constructed 
using the formation rule, since the \enquote{weakening rule} implies the 
judgement
\[
\Gamma, \; x:A \vdash x:A. 
\]
The composition of two functions (see Fig.~\ref{fig:type_composition})

\begin{figure}[ht!]
\begin{tikzpicture}
\node (a) at (0,0) {$A$};
\node (b) at (3,0) {$B$};
\node (c) at (6,0) {$C$};
\draw[->] (a) to node[scale=.7] [above] {$f$} (b);
\draw[->] (b) to node[scale=.7] [above] {$g$} (c);
\draw[->] (a) to [bend right] node[scale=.7] [above] {$g \circ f$} (c);
\end{tikzpicture}
\caption{\label{fig:type_composition}The composition $g \circ f$ of 
concatenable functions $f,g$.}
\end{figure}

can be deduced from the elementary inference rules. One can also show that the 
composition of three or more functions is associative. One can therefore 
refrain from using brackets.\footnote{See \cite[Sec. 2.2]{rijke}.} 

With the help of the function type, types can be related to each other. For 
example, there is a sequence of mappings between number systems\footnote{In 
addition to $\mathbb{N}$, $\mathbb{Z}$, $\mathbb{Q}$ and $\mathbb{R}$, the 
complex numbers $\mathbb{C}$, the quaternions $\mathbb{H}$ and the Cayley 
octaves $\mathbb{O}$.}
\[
\mathbb{N} \longrightarrow \mathbb{Z} \longrightarrow \mathbb{Q} 
\longrightarrow \mathbb{R} \longrightarrow \mathbb{C}  \longrightarrow 
\mathbb{H} \longrightarrow \mathbb{O}, 
\]
under which the term $5:\mathbb{N}$ is sent to its respective images 
$5_\mathbb{Z}$, $5_\mathbb{Q}$, $5_\mathbb{R}$ etc., which we also denote by 
$5$ for simplicity.

\medskip 
The negation of a type $A$, denoted as $\neg A$, is defined as the function  
type 
\[
\neg A \equiv (A \longrightarrow \mathbf{0}). 
\]
This is a feasible definition, since  $(A \longrightarrow \mathbf{0})$ only has 
a term if $A$ does not have terms and vice versa. Using this, one can show that 
there is a function
\[
A \xlongrightarrow{~~~~~~} \neg \neg A  
\]
that usually does not have an inverse. This proves that 
intuitionistic\index{Intuitionism} logic fits more naturally into type 
theory\index{Type theory} than classical logic with its law of the excluded 
middle.

\section{Dependent Function Type\index{Type theory}}

The function type can be significantly generalised to the important dependent 
function type which plays a central role in dependent type theory\index{Type 
theory}. It is denoted by 
\[
\prod_{x:A} B(x)
\]
and arises from a family of types $B(x)$, i.e., by a function
\[
B \colon A  \xlongrightarrow{~~~~~~} U
\]
into a given universe\index{Universe} $U$ with $B(x):U$. A term   
\[
s:\prod_{x:A} B(x)
\]
in the dependent function type is a generalised mapping $s$ that assigns a 
$s(x):B(x)$ to each $x:A$, but otherwise has similar calculation rules as 
mappings. For convenience, it is occasionally also denoted as 
\[
s \equiv \lambda x.s(x)
\]
using Church's\index{Church, Alonzo} $\lambda$-notation. A dependent function 
can be viewed as a section of the family $B(x)$ (see 
Fig.~\ref{fig:type_section}).  

\begin{figure}[ht!]
\begin{tikzpicture}
\node (a) at (-1.5,0.3) {$s$};
\draw[color=gray,ultra thick,dashed, domain=-1:8.45] plot 
(\x,-0.2*\x+0.03*\x^2+0.15);
\draw[color=gray,ultra thick,domain=0:1.45] plot (1+\x,\x^2);
\draw[color=gray,ultra thick,domain=0:1.45] plot (1-\x,-\x^2);
\draw[color=gray,ultra thick,domain=0:1.21] plot (4+\x,\x^4);
\draw[color=gray,ultra thick,domain=0:1.21] plot (4-\x,-\x^4);
\draw[color=gray,ultra thick,domain=0:1.13] plot (7+\x,\x^6);
\draw[color=gray,ultra thick,domain=0:1.13] plot (7-\x,-\x^6);
\end{tikzpicture}
\caption{\label{fig:type_section}Dependent function $s: \prod_{x:A} B(x)$ as 
a section of a family $B(x)$.}
\end{figure}  

The four inference rules in the case of the dependent product type $\Pi_{x:A} 
B(x)$ are similar to the ordinary function type. We describe them in 
Fig.~\ref{fig:producttype} using the $\lambda$-notation.

\begin{figure}[ht!] 
\fbox{
\begin{minipage}{0.8\textwidth}
\begin{align*}
&\Pi-\text{Form:} && 
\frac{\Gamma \vdash A:U \quad  
\Gamma, x:A \vdash B(x):U}{\Gamma \vdash \prod_{x:A} B(x):U}  \cr
&\Pi-\text{Intro:} && 
\frac{\Gamma, x:A \vdash s(x):B(x)} 
{\Gamma \vdash \lambda x.s(x): \prod_{x:A} B(x)} \cr
&\Pi-\text{Elim:} &&
\frac{\Gamma \vdash s:\prod_{x:A} B(x) \quad  
\Gamma \vdash a:A}{\Gamma \vdash s(a):B[\frac{a}{x}]} \cr
&\Pi-\text{Comp:} && 
\frac{\Gamma, x:A \vdash s(x):B(x) \quad \Gamma \vdash a:A} 
{\Gamma \vdash (\lambda x.s(x))(a) \equiv s[\frac{a}{x}]:B[\frac{a}{x}]}
\end{align*}
\end{minipage}
}
\caption{\label{fig:producttype}The inference rules of the product type.}
\end{figure}

If $B(x)=B$ is a constant family, the special case of the mapping type 
\[
\prod_{x:A} B \equiv (A \longrightarrow B)
\]
arises. In the dependent case, there is only a map 
\[
\biggl( \prod_{x:A} B(x)\biggr) \xlongrightarrow{~~~~~~} \biggl( A 
\xlongrightarrow{~~~~~~} \sum_{x:A} B(x) \biggr), 
\]
since not every function on the right hand side is a section of the family 
$B(x)$.  

\section{Inductive Types\index{Type theory}}

Beyond (dependent) function types, the principles behind inductive types and 
their higher variants give rise to fascinating type constructions. The 
well-known type of natural numbers is one prominent example among these. Being 
an inductive type is, in general, not related to complete 
induction\index{Complete induction}, since the type $A$ does not have to be 
infinite. The terminology of \enquote{inductive} thinking is rather taken from 
the epistemological\index{Epistemology} method which is opposed to deductive 
thinking. An inductive type is called higher inductive type, if its 
constructors involve identity types\index{Identity type}. In the following 
sections we will also see examples of (higher) inductive types.

Inductive types have two distinguished features. First, they are defined by a 
finite list of constructors which express fundamental principles. A simple 
example for this is the unit type $\mathbf{1}$ with one constructor
\[
\ast:\mathbf{1} 
\]
defining a single term. The empty type $\mathbf{0}$ even does not require any 
constructor but is nevertheless a non-trivial type. The second feature is 
that they serve to produce sections 
\[
s: \prod_{x:A} B(x) 
\]
in type families $B(x)$ parametrised over an inductive type $A$ which are 
related to the constructors. Such sections are dependent functions on the type 
$A$. 

An existential statement\index{Existence} for inductive types is not intended 
by their definition with constructors. They are indeed free inductive 
constructions in the sense of Dedekind\index{Dedekind, Richard}. Only 
their consistency poses a crucial question. The constructors of inductive types 
suggest uniqueness due to their simplicity. The simple logical construction 
idea and the uniqueness are in line with logicism\index{Logicism} and the 
Platonic world of ideas\index{Platonic idealism}.\footnote{See \cite[Part 
II]{lambekscott1986}.}  

Inductive types satisfy the usual four sorts of inference rules. In this 
context, the elimination rule carries the main information of an induction 
principle. The rule tells us how to obtain a section of a family of 
dependent types which is parametrised over the inductive type. 

\section{Product and Sum}

The binary product type 
\[
A \times B
\]
and the binary disjoint sum type (or coproduct type) 
\[
A+B
\]
are quite similar to the corresponding set-theoretic\index{Set theory} 
constructions. In the former case, the terms are given by pairs $(a,b)$ with 
$a:A$ and $b:B$, in the latter case a term $c:A+B$ either fulfils $c:A$ or 
$c:B$ (but not both). The inference rules in the product case are based on both 
projections (see Fig.~\ref{fig:binary_product}). 

\begin{figure}[ht!] 
\fbox{
\begin{minipage}{0.8\textwidth}
\begin{align*}
&\times-\text{Form:} && 
\frac{\Gamma \vdash A:U \quad \Gamma \vdash B:U}{\Gamma \vdash A \times B:U} \cr
&\times-\text{Intro:} && 
\frac{\Gamma \vdash a:A \quad \Gamma \vdash b:B}{\Gamma \vdash (a,b):A \times 
B} \cr
&\times-\text{Elim:} && 
\frac{\Gamma \vdash z:A \times B}{\mathrm{pr}_1(z):A} \quad 
\frac{\Gamma \vdash z:A \times B}{\mathrm{pr}_2(z):B} \cr
&\times-\text{Comp}: &&  
\frac{\Gamma \vdash a:A \quad \Gamma \vdash b:B}{\Gamma \vdash 
\mathrm{pr}_1(a,b) \equiv a:A} \quad \frac{\Gamma \vdash a:A \quad \Gamma 
\vdash b:B}{\Gamma \vdash \mathrm{pr}_2(a,b) \equiv b:B}
\end{align*}
\end{minipage}
}
\caption{\label{fig:binary_product}The inference rules for the binary product 
type.}
\end{figure}

The elimination and the computation rule for the product type can also be 
formulated in a way that $A \times B$ becomes visible as an inductive type. 
They tell us how a function with values in a family $C$ dependent on $x:A$ and 
$y:B$ becomes a function dependent on $z=(x,y):A \times B$, i.e., they can be 
reduced to the tautological function  
\[
\mathrm{ind}_\times: \biggl( \prod_{x:A} \prod_{y:B} C(x,y)\biggr) 
\xlongrightarrow{~~~~~~} \biggl(\prod_{z:A \times B} C(z)\biggr).
\]
The inference rules of the sum type use the maps $\mathrm{inl} \colon 
A \rightarrow A+B$ and $\mathrm{inr} \colon B \rightarrow A+B$. We choose the 
formulation as an inductive type (see Fig.~\ref{fig:binary_sum}).

\begin{figure}[ht!] 
\fbox{
\begin{minipage}{0.8\textwidth}
\begin{align*}
&+-\text{Form:} && 
\frac{\Gamma \vdash A:U \quad \Gamma \vdash B:U}{\Gamma \vdash A+B:U} \cr
&+-\text{Intro:} && 
\frac{\Gamma \vdash a:A}{\Gamma \vdash \text{inl}(a): A+B} \quad \frac{\Gamma 
\vdash b:B}{\Gamma \vdash \text{inr}(b): A+B} \cr
&+-\text{Elim:} && 
\frac{  
\begin{array}{l}
\Gamma, \; z:A+B \vdash C(z) \quad \Gamma, \; x:A \vdash 
f(x):C(\text{inl}(x)) \\ 
\Gamma, \; y:B \vdash g(y):C(\text{inr}(y)) \quad  \Gamma \vdash e:A+B
\end{array} 
}{\Gamma \vdash 
\text{ind}_+(f,g)(e):C[\frac{e}{z}]} \cr &+-\text{Comp}: &&  
\frac{   
\begin{array}{l}
\Gamma, \; z:A+B \vdash C(z) \quad \Gamma, \; x:A \vdash 
f(x):C(\text{inl}(x)) \\  
\Gamma, \; y:B \vdash g(y):C(\text{inr}(y)) \quad \Gamma \vdash a:A
\end{array}
}{\Gamma \vdash 
\text{ind}_+(f,g)(\text{inl}(a)) \equiv f(a):C[\frac{\text{inl}(a)}{z}]} \cr 
&+-\text{Comp}: &&  
\frac{
\begin{array}{l}
\Gamma, \; z:A+B \vdash C(z) \quad \Gamma, \; x:A \vdash 
f(x):C(\text{inl}(x)) \\  \Gamma, \; y:B \vdash 
g(y):C(\text{inr}(y)) \quad \Gamma \vdash b:B
\end{array}
}{\Gamma \vdash 
\text{ind}_+(f,g)(\text{inr}(b)) \equiv g(b):C[\frac{\text{inr}(b)}{z}]}
\end{align*}
\end{minipage}
}
\caption{\label{fig:binary_sum}The inference rules for the binary sum type.}
\end{figure}

\medskip
The function $\mathrm{ind}_+$ can in this case be written as 
\[
\mathrm{ind}_+ \colon \biggl( \prod_{x:A} C(\mathrm{inl}(x)) \biggr)
\xlongrightarrow{~~~~~~} \biggl( \biggl(\prod_{y:B} C(\mathrm{inr}(y))\biggr) 
\xlongrightarrow{~~~~~~} \prod_{z:A+B} C(z) \biggr).
\]
Here, we used that functions $f \colon A \times B \longrightarrow C$ 
can also be written as 
\[
A \xlongrightarrow{~~~~~~} (B \xlongrightarrow{~~~~~~} C),
\]
since $f(a,-)$ can be viewed as a function $B \longrightarrow C$. This principle
from the area of combinatorial logic is called \enquote{currying}, after 
Haskell B. Curry\index{Curry, Haskell B.} and is very useful in type 
theory\index{Type theory}. 

\medskip
Alternative versions of the functions 
$\mathrm{ind}_+$ and $\mathrm{ind}_\times$ are 
\[
\mathrm{ind}_+ \colon (A \xlongrightarrow{~~~~~~} X) \xlongrightarrow{~~~~~~} 
\biggl( (B \xlongrightarrow{~~~~~~} X) \xlongrightarrow{~~~~~~} (A+B 
\xlongrightarrow{~~~~~~} X) \biggr),  
\]
and 
\[
\mathrm{ind}_\times \colon (X \xlongrightarrow{~~~~~~} A\times B) 
\xlongrightarrow{~~~~~~} (X \xlongrightarrow{~~~~~~} A) \times (X 
\xlongrightarrow{~~~~~~} B).
\]
The special case $X=\mathbf{0}$ is remarkable in the first case, since it 
expresses the logical equation 
\[
\neg A \wedge \neg B \Longrightarrow \neg (A \vee B).  
\] 

\section{Dependent Pair Type\index{Type theory}}

The dependent pair type 
\[
\sum_{x:A} B(x)
\]
arises from a family $B(x)$ as a generalisation of the product type 
$A \times B$. Terms in it are pairs $(x,y)$ with $y:B(x)$ for an $x:A$. There 
is a canonical projection mapping 
\[
\mathrm{pr}_1 \colon \biggl(\sum_{x:A} B(x)\biggr) \xlongrightarrow{~~~~~~} A 
\]
to the first factor. We can view $\sum_{x:A} B(x)$ as the total space of the 
family $B(x)$. The second projection produces a dependent function. 

\medskip
The four inference rules in the case of the dependent pair type $\sum_{x:A} 
B(x)$ can be formulated as in the case of $A \times B$ using both projections. 
For the inference rules as an inductive type see Fig.~\ref{fig:pair_type}.  

\newpage

\begin{figure}[ht!] 
\fbox{
\begin{minipage}{0.8\textwidth}
\begin{align*}
&\Sigma-\text{Form:} && 
\frac{\Gamma \vdash A:U \quad  
\Gamma, \; x:A \vdash B(x):U}{\Gamma \vdash \sum_{x:A} B(x):U}  \cr
&\Sigma-\text{Intro:} && 
\frac{\Gamma, \; x:A \vdash B(x) \quad \Gamma \vdash a:A \quad \Gamma \vdash 
b:B[\frac{a}{x}]} 
{\Gamma \vdash (a,b): \sum_{x:A} B(x)} \cr
&\Sigma-\text{Elim:} && \frac{
\begin{array}{l}
\Gamma, \; z:\sum_{x:A} B(x) \vdash C(z) \quad \Gamma \vdash p:\sum_{x:A} B(x) 
\\
\Gamma, \; x:A, \; y:B(x) \vdash g(x,y):C[\frac{(x,y)}{z}]
\end{array}
}
{\Gamma \vdash 
\mathrm{ind}_\Sigma(g)(p):C[\frac{p}{z}]} \cr
&\Sigma-\text{Comp:} && 
\frac{
\begin{array}{l}
\Gamma, z:\sum_{x:A} B(x) \vdash C(z) \quad \quad \Gamma \vdash a:A   \\ 
\Gamma \vdash b:B[\frac{a}{x}] \quad \Gamma, \; x:A, 
\; y:B(x) \vdash g:C[\frac{(x,y)}{z}]
\end{array}
}{\Gamma \vdash \mathrm{ind}_\Sigma (g)(\mathrm{pair}(a,b)) \equiv 
g[\frac{(a,b)}{(x,y)}]:C[\frac{(a,b)}{(x,y)}]}
\end{align*}
\end{minipage}
}
\caption{\label{fig:pair_type}The inference rules of the dependent pair type.}
\end{figure}

Here, the function $\mathrm{ind}_\Sigma$ is given by 
\[
\mathrm{ind}_\Sigma \colon \biggl(\prod_{x:A} \prod_{y:B(x)} 
C(\mathrm{pair}(x,y))\biggr) 
\xlongrightarrow{~~~~~~} \prod_{z:\Sigma_{x:A} B(x)} C(z).
\]
It uses the tautological pairing function 
\[
\mathrm{pair} \colon \prod_{x:A} \biggl( B(x) \xlongrightarrow{~~~~~~} 
\sum_{y:A} B(y) \biggr). 
\]
The binary product type 
\[
\sum_{x:A} B \equiv A \times B
\]
arises in the special case that $B(x) \equiv B$ is constant. 

\section{Identity Type\index{Identity type}}

An important new type is the identity type\index{Identity type}, denoted by
\[
\mathrm{Id}_A(a,b) \; \text{  or with }\;  (a=_A b).
\]
It is a primitive type like the mapping type and depends on the respective type 
theory\index{Type theory}. If it is inhabited by a term $p$, $p$ expresses a 
propositional equality\index{Equality}\index{Identity} between $a$ and $b$. 
This is a weakening of the judgemental equality $a \equiv b$, which expresses 
that $a$ and $b$ denote the same term. We imagine the terms 
\[
p:\mathrm{Id}_A(a,b)
\]
as paths $p$ from $a$ to $b$. The type $\mathrm{Id}_A(a,a)$ is always inhabited 
by a canonical term $\mathrm{refl}_a$, i.e., the constant path corresponding to 
the judgement $a \equiv a$.  

\medskip 
In the case of the inductive identity type\index{Identity type}, four inference 
rules apply (see Fig.~\ref{fig:id_type}). Here, $a:A$ is a base point and 
$C(x,p)$ is a type dependent on $x:A$ and $p:\mathrm{Id}_A(a,x)$ which 
contains terms of the form $c(x,p):C(x,p)$. 

\begin{figure}[ht!] 
\fbox{
\begin{minipage}{0.9\textwidth}
\begin{align*}
&\mathrm{Id}-\text{Form:} && 
\frac{\Gamma \vdash A:U \quad \Gamma \vdash a:A \quad  \Gamma \vdash 
b:A}{\Gamma \vdash \mathrm{Id}_A(a,b):U} \cr
&\mathrm{Id}-\text{Intro:} && 
\frac{\Gamma \vdash a:A} 
{\Gamma \vdash \mathrm{refl}_a:\mathrm{Id}_A(a,a) } \cr
&\mathrm{Id}-\text{Elim:} && 
\frac{\Gamma \vdash a,x:A \quad \Gamma \vdash p:\mathrm{Id}_A(a,x)
\quad \Gamma \vdash c_a:C(a,\mathrm{refl}_a)}
{\mathrm{ind}_\mathrm{Id}(c_a)(x,p):C(x,p)} \cr
&\mathrm{Id}-\text{Comp:} && 
\frac{\Gamma \vdash a:A \quad \Gamma \vdash c_a:C(a,\mathrm{refl}_a)}
{\Gamma \vdash \mathrm{ind}_\mathrm{Id}(c_a)(a,\mathrm{refl}_a) 
\equiv c_a:C(a,\mathrm{refl}_a)}
\end{align*}
\end{minipage}
}
\caption{\label{fig:id_type}The inference rules of the identity type.}
\end{figure}

Here, the function $\mathrm{ind}_\mathrm{Id}$ is given by  
\[
\mathrm{ind}_\mathrm{Id} \colon C(a,\mathrm{refl}_a) 
\xlongrightarrow{~~~~~~} \prod_{x:A} \prod_{p:\mathrm{Id}_A(a,x)} C(x,p). 
\]

The elimination rule by Martin-L\"of\index{Martin-L\"of, Per} is also 
called path induction rule\index{Path induction} and represents a 
far-reaching generalisation of Leibniz's\index{Leibniz, Gottfried Wilhelm} 
invariance rule\index{Leibniz's invariance rule}. One can formulate this rule 
in such a way that $a$ is not a base point but a variable.\footnote{See  
\cite[Rem. 5.1.4]{rijke} and \cite[Sec. 1.12]{voevodsky2013}.} 

The identity type\index{Identity type} represents the concept of polymorphism 
from theoretical computer science, because this type does indeed depend on the 
type $A$, but does so in a universal way. The identity type\index{Identity 
type} $\mathrm{Id}_A$ without a given pair $(a,b)$ can be considered as a 
dependent type
\[
\mathrm{Id}_A \equiv \biggl(\sum_{(a,b):A\times A} \mathrm{Id}_A(a,b) \biggr)
\xlongrightarrow{~~~~~~} A \times A
\]
over $A \times A$ and is called the path space of $A$. As in 
topology\index{Topology}, the constant path $\mathrm{refl}_a$, the 
concatenation of paths $q \circ p$ and the inverse path $p^{-1}$ in 
Martin-L\"of-type theory\index{Martin-L\"of, Per}\index{Type theory} 
fulfil the following relationships by means of suitable mappings
\begin{align*}
\mathbf{1} & \xlongrightarrow{~~~~~~} \mathrm{Id}_A(a,a) \\
\mathrm{Id}_A(a,b) \times \mathrm{Id}_A(b,c) & \xlongrightarrow{~~~~~~} 
\mathrm{Id}_A(a,c) \\
\mathrm{Id}_A(a,b) & \xlongrightarrow{~~~~~~} \mathrm{Id}_A(b,a).
\end{align*}
The formulas
\begin{align*}
p^{-1} \circ p & = \mathrm{refl}_a \\
p \circ p^{-1} & = \mathrm{refl}_b \\
(p \circ q) \circ r & = p \circ (q \circ r) 
\end{align*}
do not hold under judgemental equality\index{Equality}\index{Identity} but 
only when using the twice iterated identity type\index{Identity type}. This 
kind of higher conditions forms the structure of an infinity 
groupoid\index{Infinity groupoid}, which is associated with each type 
$A$.\footnote{See \cite{hofmannstreicher,KV91,lumsdaine2010}.} 

If two dependent functions $f,g \colon A \longrightarrow \prod_{x:A} B(x)$ 
have identical values at each term $x$ in the sense of the identity 
type\index{Identity type}, then they are not necessarily identical as 
dependent functions. Such a requirement of function 
extensionality\index{Extensionality} only follows with additional assumptions 
such as the univalence axiom\index{Univalence} of Voevodsky\index{Voevodsky, 
Vladimir}. We will discuss this later.

\medskip 
Let us give three more useful tools for identity types\index{Identity type}. 
For any function $f \colon A \longrightarrow B$, there is a function 
\[
\mathrm{ap}_f \colon \prod_{x,y:A} \biggl(\mathrm{Id}_A(x,y) 
\xlongrightarrow{~~~~~~} \mathrm{Id}_B(f(x),f(y))\biggr)
\]
which is called \enquote{action on paths}. It implies a sort of continuity for 
functions, since this rule would follow in the continuous case, since 
we have a composition of continuous functions with $f$. The construction of 
$\mathrm{ap}_f$ is defined by path induction\index{Path induction} and the 
definitorial requirement  
\[
\mathrm{ap}_f(\mathrm{refl}_x) \equiv \mathrm{refl}_{f(x)}.
\footnote{See \cite[Sec. 5.3]{rijke} and \cite[Sec. 2.2]{voevodsky2013}.}  
\]
There is also a more general function called \enquote{transport on paths}. If 
$B(x)$ is a type family over $A$, then\footnote{See \cite[Sec. 5.4]{rijke} and 
\cite[Sec. 2.3]{voevodsky2013}.} 
\[
\mathrm{trans}_B \colon \prod_{x,y:A} \biggl(\mathrm{Id}_A(x,y) 
\xlongrightarrow{~~~~~~} (B(x) \xlongrightarrow{~~~~~~} B(y))\biggr).
\]
Using the transport function, we can define a dependent version of 
\enquote{action on paths} using $f: \prod_{x:A} B(x)$ as a 
function\footnote{See 
\cite[Sec. 5.4]{rijke} and \cite[Sec. 2.3]{voevodsky2013}.} 
\[
\mathrm{apd}_f \colon \prod_{x,y:A} \biggl(\mathrm{Id}_A(x,y) 
\xlongrightarrow{~~~~~~} 
\mathrm{Id}_{B(y)}(\mathrm{trans}_B(p,f(x)), f(y))\biggr). 
\]
It is based on the definitorial equation\footnote{See \cite[Sec. 5.4]{rijke} 
and \cite[Sec. 2.3]{voevodsky2013}.} 
\[
\mathrm{apd}_f(\mathrm{refl}_x) \equiv \mathrm{refl}_{f(x)}. 
\]
With the help of these mappings we obtain a path lifting property for families
$B(x)$ which one can express by saying that every family $B(x)$ is a 
fibration (see Fig.~\ref{fig:type_lifting}) since topological fibrations are 
defined using exactly this property.\footnote{See \cite[Sec. 
2.3]{voevodsky2013}.} 

\begin{figure}[ht!]
\centering
\begin{tikzpicture}[yscale=.5,xscale=2]
    \draw (0,0) arc (-90:170:8ex) node[anchor=south east] {$A$} arc 
(170:270:8ex);
    \draw (0,6) arc (-90:170:8ex) node[anchor=south east] {$\sum_{x:A} B(x)$} 
arc (170:270:8ex);
    \draw[->] (0,5.8) -- node[auto] {$\mathrm{pr}$} (0,3.2);
    \node[circle,fill,inner sep=1pt,label=left:{$x$}] (b1) at (-.5,1.4) {};
    \node[circle,fill,inner sep=1pt,label=right:{$y$}] (b2) at (.5,1.4) {};
    \draw[decorate,decoration={snake,amplitude=1}] (b1) -- node[auto,swap] 
{$p$} 
(b2);
    \node[circle,fill,inner sep=1pt,label=left:{$f(x)$}] (b1) at (-.5,7.2) {};
    \node[circle,fill,inner sep=1pt,label=right:{$f(y)$}] (b2) at (.5,7.2) {};
    \draw[decorate,decoration={snake,amplitude=1}] (b1) -- node[auto] {$f(p)$} 
(b2);
\end{tikzpicture}
\caption{\label{fig:type_lifting}Families are fibrations since they have the 
path lifting property.\footnote{See \cite[Sec. 2.3]{voevodsky2013}. The figure 
also originates from this source.}}
\end{figure}

\section{Natural Numbers}

The type of natural numbers $\mathbb{N}$ is an inductive type, given by the 
constructor
\[
\mathbb{N}=
\begin{cases}
0 : \mathbb{N} \cr 
S \colon \mathbb {N} \longrightarrow \mathbb{N} \quad 
\text{(successor function)}. 
\end{cases}
\]
Similar to Dedekind\index{Dedekind, Richard}, this type only provides the 
initial element $0$ and the successor function $S$. The natural 
numbers are then given by the infinite sequence 
\[
0, 1 \equiv S(0), 2 \equiv S(S(0)), \ldots   
\]
The inductive type $\mathbb{N}$ of natural numbers has five inference rules, 
where $C(x)$ is a type dependent on $x:\mathbb{N}$ 
(see Fig.~\ref{fig:nat_type}). 

\begin{figure}[ht!] 
\fbox{
\begin{minipage}{0.9\textwidth}
\begin{align*}
&\mathbb{N}-\text{Form:} && 
\frac{}{\vdash \mathbb{N}:U} \cr
&\mathbb{N}-\text{Intro:} && 
\frac{}{\vdash 0:\mathbb{N}} \quad 
\frac{\Gamma \vdash n:\mathbb{N}} 
{\Gamma \vdash S(n):\mathbb{N}} \cr
&\mathbb{N}-\text{Elim:} && 
\frac{\Gamma \vdash c_0:C[\frac{0}{x}] \;\, \Gamma, x:\mathbb{N},c_x:C(x) 
\vdash c_{S(x)}:C[\frac{S(x)}{x}] \;\, \Gamma \vdash n:\mathbb{N}}{\Gamma 
\vdash \mathrm{ind}_\mathbb{N}(c_0,c_{S})(n):C[\frac{n}{x}]} \cr
&\mathbb{N}-\text{Comp}_0: &&  
\frac{\Gamma \vdash c_0:C[\frac{0}{x}] \quad \Gamma, x:\mathbb{N},c_x:C(x) 
\vdash c_{S(x)}:C[\frac{S(x)}{x}]}
{\Gamma \vdash \mathrm{ind}_\mathbb{N}(c_0,c_{S})(0)\equiv 
c_0:C[\frac{0}{x}]} 
\cr &\mathbb{N}-\text{Comp}_S: &&  
\frac{\Gamma \vdash c_0:C[\frac{0}{x}] \;\, \Gamma, x:\mathbb{N},c_x:C(x) 
\vdash c_{S(x)}:C[\frac{S(x)}{x}]\;\, \Gamma \vdash n:\mathbb{N}}
{\Gamma \vdash \mathrm{ind}_\mathbb{N}(c_0,c_{S})(S(n)) \equiv 
c_{S(n)}:C[\frac{S(n)}{x}]}
\end{align*}
\end{minipage}
}
\caption{\label{fig:nat_type}The inference rules for the type of natural 
numbers.}
\end{figure}

The last three rules can be expressed using the function 
\[
\mathrm{ind}_\mathbb{N} \colon C(0) \xlongrightarrow{~~~~~~} 
\biggl(\biggl(\prod_{n:\mathbb{N}} C(n) \xlongrightarrow{~~~c_S~~~} C(S(n))
\biggr) \xlongrightarrow{~~~~~~} \prod_{n:\mathbb{N}} C(n) \biggr).
\]

\medskip 
After we have defined the natural numbers, we can develop elementary 
arithmetic, i.e., the computation rules for addition, multiplication, and 
exponentiation, together with the commutativity, associativity, and 
distributivity laws. To begin, we look at simple equations. We define addition 
as  
\[
m+0 \equiv m \quad \text{and} \quad m+S(n) \equiv S(m+n).
\]
The equations hold judgementally and define a function 
\[
\mathrm{add}: \mathbb{N} \xlongrightarrow{~~~~~~} (\mathbb{N} 
\xlongrightarrow{~~~~~~} \mathbb{N}), 
\]
using path induction\index{Path induction}, of which commutativity and 
associativity still have to be proven. Now one can show equations like  
$1+1 \equiv 2$: 
\[
1+1 \equiv S(0)+S(0) \equiv S(S(0)+0) \equiv S(S(0)) \equiv 2.
\]
More generally, one can define all primitive recursive 
functions\index{Recursion} using this method.\footnote{See \cite[Sec. 
3.3]{rijke} and \cite[Sec. 1.9]{voevodsky2013}.} The multiplication $m \cdot n$ 
for $n \neq 0$ is defined as
\[
m \cdot 1 \equiv m \quad \text{and} \quad m \cdot S(n) \equiv m \cdot n + m,  
\]
and exponentiation via 
\[
a^0 \equiv 1 \quad \text{and} \quad a^{S(n)} \equiv a^n \cdot a.
\]
All the equations hold judgementally. 

\medskip 
For the addition function $\mathrm{add}$ four computation rules remain to be 
proven: 
\[
0+n=n \quad S(m)+n=S(m+n) \quad m+n=n+m \quad (m+n)+0=m+(n+0).
\]
On purpose, we formulate these equations using propositional equality 
\enquote{=}, since the proof can only be given using it. As an example, we want 
to explain the proof of $0+n=n$. It proceeds by induction on $n$. Obviously, 
the case $n=0$ is trivial, since by definition of $\mathrm{add}$ one has $0+0 
\equiv 0$. The induction step needs:
\[
0+n=0 \vdash 0+S(n)=S(n). 
\]
But by definition of $\mathrm{add}$ one has the equation $0+S(n) \equiv S(0+n)$.
Hence, it is sufficient to show: 
\[
0+n=0 \vdash S(0+n)=S(n). 
\]
This equations readily follows by applying \enquote{action on paths}, i.e., 
the function $\mathrm{ap}_S$, induced by the successor map $S$ at the path 
$p:0+n=n$ as argument. At this point, propositional equality \enquote{=} is 
needed. This completes the proof. 

The proof of $S(m)+n=S(m+n)$ is similar. We refer to the literature for all 
remaining proofs of the arithmetic of $\mathbb{N}$.\footnote{See \cite[Sec. 
2.13]{voevodsky2013} and \cite[Sec. 5.6]{rijke}.}

\medskip 
As an application of the coproduct type, we construct the integers 
$\mathbb{Z}$ via an iterated sum
\[
\mathbb{Z}=\mathbb{N}+(\mathbf{1} + \mathbb{N}). 
\]
Strictly positive and strictly negative numbers are given here by the 
two copies of  $\mathbb{N}$, and the integer $0$ corresponds to the term in 
the type $\mathbf{1}$.\footnote{See \cite[Sec. 6.10]{voevodsky2013} and 
\cite[Sec. 4.4]{rijke}.}

\section{Higher Inductive Types\index{Type theory}}

We will now introduce the examples $\mathbb{S}^1$, $\mathbb{S}^2$,  
$\mathbb{I}$, and $\mathbb{T}^2$ of higher inductive types. They are, so to 
speak, born out of pure logic and are determined by their calculation rules and 
their universal properties. A first example is the circle $\mathbb{S}^1$, which 
is introduced by 
\[
\mathbb{S}^1 \equiv 
\begin{cases}
0 : \mathbb{S}^1 \quad \text{(base point)} \cr 
\mathrm{loop} : \mathrm{Id}_{\mathbb{S}^1}(0,0) \quad 
\text{(loop generator)},  
\end{cases}
\]
i.e., by postulating a base point and a loop that defines the circle. It should 
be noted that this circle only contains one term $0$. Similarly, the interval  
$\mathbb{I}$ can be introduced by
\[
\mathbb{I} \equiv 
\begin{cases}
0,1 : \mathbb{I} \quad \text{ (start and end point)} \cr 
\mathrm{path} : \mathrm{Id}_\mathbb{I}(0,1)\quad 
\text{(interval generator)}.  
\end{cases}
\]
The order here is reversed compared to topology\index{Topology}, where the real 
interval $[0,1]$ defines the topological\index{Topology} path from $0$ to $1$. 
This is because paths are primitive objects in type theory\index{Type theory}. 
Both $\mathbb{S}^1$ and $\mathbb{I}$ correspond -- due to their few 
inhabitants -- as types not to a real circle or interval, but are synthetic 
constructions of their homotopy types\index{Homotopy theory}.

The $2$-sphere $\mathbb{S}^2$ has the two constructors
\[
\mathbb{S}^2 \equiv 
\begin{cases}
\mathrm{base}: \mathbb{S}^2 \quad \text{ (base point)}, \cr 
\mathrm{surf}:\mathrm{Id}_{\mathbb{S}^2}(\mathrm{refl}_\mathrm{base},\mathrm{
refl}_\mathrm{base}) \quad \text{(surface generator)}.
\end{cases}
\]
$\mathbb{S}^2$ can be constructed as a suspension of $\mathbb{S}^1$. 
A suspension $\Sigma A$ of any type $A$ is defined by three constructors
\[
\Sigma A \equiv 
\begin{cases}
\mathrm{N}: \Sigma A \quad \text{ (north pole)}, \cr 
\mathrm{S}: \Sigma A \quad \text{ (south pole)}, \cr 
\mathrm{mer} \colon A \longrightarrow \mathrm{Id}_{\Sigma A}(N,S) \quad 
\text{(meridian function)}.
\end{cases}
\]
The $2$-torus $\mathbb{T}^2$ has the four constructors
\[
\mathbb{T}^2 \equiv 
\begin{cases}
\mathrm{base}: \mathbb{T}^2 \quad \text{(base point)}, \cr 
\mathrm{p} : (\mathrm{base}=\mathrm{base}) \quad \text{(first loop)}, \cr 
\mathrm{q} : (\mathrm{base}=\mathrm{base}) \quad \text{(second loop)}, \cr 
\mathrm{surf}:\mathrm{Id}_{\mathbb{T}^2}(p \circ q,q \circ p) \quad 
\text{(surface generator)}.
\end{cases}
\]
This corresponds to the topological construction of a torus from a square by 
identifying opposite sides.

The inference rules for these types can be found in the  
literature.\footnote{See \cite[Ch. 6]{voevodsky2013}.} For the case of the 
circle, they are listed in Fig.~\ref{fig:circle_type}. 

\begin{figure}[ht!] 
\fbox{
\begin{minipage}{0.90\textwidth} 
\begin{align*}
&\mathbb{S}^1-\text{Form:} && 
\frac{}{\Gamma \vdash \mathbb{S}^1:U} \cr
&\mathbb{S}^1-\text{Intro:} && 
\frac{}{\Gamma \vdash \mathrm{base}:\mathbb{S}^1} \quad 
\frac{}{\Gamma \vdash \mathrm{loop}:(\mathrm{base}=\mathrm{base})} \cr
&\mathbb{S}^1-\text{Elim:} && 
\frac{
\begin{array}{l}
\Gamma,\, x:\mathbb{S}^1 \vdash C(x) \quad  
\Gamma \vdash c:C(\mathrm{base}) \\ 
\Gamma \vdash 
\ell:\mathrm{Id}_{C(\mathrm{base})}(\mathrm{trans}_C(\mathrm{loop},c),c) 
\quad \Gamma \vdash p:\mathbb{S}^1
\end{array}
}
{\mathrm{ind}_{\mathbb{S}^1}(c,\ell)(p):C[\frac{p}{x}]} \cr
&\mathbb{S}^1-\text{Comp}_1: &&  
\frac{\Gamma \vdash c:C(\mathrm{base}) \quad \Gamma \vdash \ell:
\mathrm{Id}_{C(\mathrm{base})}(\mathrm{trans}_C(\mathrm{loop},c),c)}
{\mathrm{ind}_{\mathbb{S}^1}(c,\ell)(\mathrm{base}) \equiv c} \cr
&\mathbb{S}^1-\text{Comp}_2: &&  
\frac{\Gamma \vdash c:C(\mathrm{base}) \quad \Gamma \vdash \ell:
\mathrm{Id}_{C(\mathrm{base})}(\mathrm{trans}_C(\mathrm{loop},c),c)}{
\mathrm{apd}_{\mathrm{ind}_{\mathbb{S}^1}(c,\ell)}(\mathrm{loop}) \equiv \ell} 
\end{align*}
\end{minipage}
}
\caption{\label{fig:circle_type}The inference rules for $\mathbb{S}^1$.}
\end{figure}

The function $\mathrm{ind}_{\mathbb{S}^1}$ is in this case given by  
\[
\mathrm{ind}_{\mathbb{S}^1} \colon \biggl(\sum_{c:C(\mathrm{base})} 
\mathrm{Id}_{C(\mathrm{base})}(\mathrm{trans}_C(\mathrm{loop},c),c)\biggr)
\xlongrightarrow{~~~~~~} \biggl(\prod_{x:\mathbb{S}^1} C(x)\biggr).
\]
The induction principle can be nicely illustrated as in 
Fig.~\ref{fig:circle_lifting}.  

\begin{figure}[ht!]
\centering
\begin{tikzpicture}
    \draw (0,0) ellipse (3 and .5);
    \draw (0,3) ellipse (3.5 and 1.5);
    \begin{scope}[yshift=4]
      \clip (-3,3) -- (-1.8,3) -- (-1.8,3.7) -- (1.8,3.7) -- (1.8,3) -- (3,3) 
-- 
(3,0) -- (-3,0) -- cycle;
      \draw[clip] (0,3.5) ellipse (2.25 and 1);
      \draw (0,2.5) ellipse (1.7 and .7);
    \end{scope}
    \node (P) at (4.5,3) {$\sum_{x:\mathbb{S}^1} C(x)$};
    \node (S1) at (4.5,0) {$\mathbb{S}^1$};
    \draw[->>,thick] (P) -- (S1);
    \node[fill,circle,inner sep=1pt,label={below right:$\mathrm{base}$}] at 
(0,-.5) {};
    \node at (-2.6,.6) {$\mathrm{loop}$};
    \node[fill,circle,inner sep=1pt] (b) at (0,2.3) {};
      \node at (-.3,2.3) {$c$};
      \node[fill,circle,inner sep=1pt] (tb) at (0,1.8) {};
to[out=180,in=180] (tb);
      \draw[dashed] (b) arc (-90:90:2.9 and 0.85) arc (90:270:2.8 and 
1.1);
      \begin{scope}
        \clip (b) -- ++(.1,0) -- (.1,1.8) -- ++(-.2,0) -- ++(0,-1) -- ++(3,2) 
-- 
++(-3,0) -- (-.1,2.3) -- cycle;
        \draw[dotted,thick] (.2,2.07) ellipse (.2 and .57);
        \begin{scope}
++(3,3) -| (b);
          \clip (.2,0) rectangle (-2,3);
          \draw[thick] (.2,2.07) ellipse (.2 and .57);
        \end{scope}
      \end{scope}
      \node at (1,1.2) {$\ell: \mathrm{trans}_C(\mathrm{loop},c)=c$};
\end{tikzpicture}
\caption{\label{fig:circle_lifting}The type-theoretic induction principle for 
$\mathbb{S}^1$ is illustrated by a lifting.\footnote{See \cite[Sec. 
6.3]{voevodsky2013}. The figure also originates from this source.}}
\end{figure}

\medskip 
By a theorem of Michael Shulman from 2011, the fundamental group of 
$\mathbb{S}^1$ is -- as in the classical topology -- given by\footnote{See 
\cite{licata}.} 
\[
\pi_1(\mathbb{S}^1,\mathrm{base}) \cong \mathbb{Z}.
\]
The proof is remarkable since it established the useful encode-decode method in 
homotopy type theory\index{Homotopy type theory}. It shows that the universal 
cover of $\mathbb{S}^1$ is given by an inductive type which is a 
\enquote{homotopical} variant of the reals. Its terms are the integers  
$\mathbb{Z}$ and there is exactly one path from $n$ to $n+1$ for all integers 
$n$. Some higher homotopy groups\index{Homotopy group} of spheres were computed 
with similar methods.\footnote{See \cite[Ch. 22]{rijke} and \cite{brunerie}.}

\section{Propositions and Contractibility} 

In this section, we look types\index{Type theory} $P$ which are called 
propositions (or mere propositions). These are defined via the 
requirement:\footnote{See \cite[Prop. 12.1.3]{rijke} and \cite[Sec. 
3.3]{voevodsky2013}.} 
\[
\text{For all } x,y:P \text{ one has } x= y \text{ in } P.  
\]
Hence, for any pair of terms in $P$ there is a path, if there are any terms. 
The property to be a (mere) proposition is therefore equivalent to the 
existence of a term in
\[
\text{is-prop}(P) \equiv \prod_{x,y:P} \mathrm{Id}_P(x,y). 
\]
This requirement is very restrictive. For example, it implies that $P$ is 
either not inhabited or contractible as we will see now.

\begin{figure}[ht!]
\begin{tikzpicture}[line cap=round,line join=round,scale=0.10] 
\draw[rotate around={0:(23.5,35)}] (23.5,35) ellipse (18cm and 12cm); 
\draw[] (8,36) node {$a$}; 
\draw[] (39,26) node {$x$}; 
\draw[-] (10,36) -- (38,28) ;
\draw[-] (10,36) -- (25,47) ;
\draw[-] (10,36) -- (14,45) ;
\draw[-] (10,36) -- (16,24) ;
\end{tikzpicture}
\caption{\label{fig:type_contractible}Point-like type centred at $a$.}
\end{figure}

A type\index{Type theory} is contractible or point-like (see 
Fig.~\ref{fig:type_contractible}) if there is a term $a:P$ with the 
condition:\footnote{See \cite[Sec. 10.1]{rijke} and \cite[Sec. 
3.11]{voevodsky2013}.}  
\[
\text{For all } x:P \text{ one has } x = a. 
\]
This is satisfied, if the type 
\[
\text{is-contr}(P) \equiv \sum_{a:A} \prod_{x:A} \mathrm{Id}_P(x,a)   
\]
is inhabited. Using the equivalence\index{Equivalence} sign  $\simeq$, our 
notation for this property is 
\[
P \simeq \mathbf{1}. 
\]
At first sight, this definition looks as if we had treated the notion of a 
path-connected topological space\index{Topology}.  However, the inherent 
continuity of type-theoretical\index{Type theory} structures implies that in 
this way indeed the topological\index{Topology} notion of contractibility 
is realised. This was a crucial insight of Vladimir Voevodsky\index{Voevodsky, 
Vladimir} et al. at the beginning of the year 2010 and led to a simple 
definition for the equivalence\index{Equivalence} of types\index{Type theory} 
as we will see. In addition, from the definitions it immediately follows that 
(mere) propositions which are inhabited are precisely the contractible 
types.\footnote{See \cite[Prop. 12.1.3]{rijke}, \cite{shulman2018}, and 
\cite[Lemma 3.11.3]{voevodsky2013}.} This implies that for (mere) propositions, 
at least in classical logic, one has: 
\[
\text{Either one has } P \simeq \mathbf{1} \text{ or } P \text{ is not 
inhabited}.
\]

\medskip
Now we assign a (mere) proposition to each type\index{Type theory} $A$. It 
is defined as a higher inductive type. The 
constructors of $||A||$ are: 

\begin{itemize}
\item For all $a:A$ there is an $|a|:||A||$. 
\item For all $x,y:||A||$ one has $x=y$ in $||A||$. 
\end{itemize}

The induction principle for $||A||$ is as follows:\footnote{See \cite[Sec. 
14.1 and Sec. 14.2]{rijke} and \cite[Sec. 3.7 and Sec. 6.9]{voevodsky2013}.} 

\begin{itemize}
\item If $P$ is a (mere) proposition and $f \colon A \longrightarrow P$ a  
function, then there is a function $g \colon ||A|| \longrightarrow P$, such that
$g(|a|)=f(a)$ for all $a:A$.
\end{itemize}

Consequently, for the type $||A||$ one has either $P \simeq \mathbf{1}$ or  
$P$ is not inhabited given the previous explanations. There is no formula 
for $||A||$, except in special situations.\footnote{See \cite[Exercises 
3.14 and 3.15]{voevodsky2013}.} From the recursion principle one 
concludes $P=||P||$ for every (mere) proposition. 

\medskip 
With the help of (mere) propositions we define subtypes of a type $X$ (resp. 
predicates) as functions 
\[
B \colon X \xlongrightarrow{~~~~~~} \mathrm{Prop}_U,
\]
i.e., as families $B(x)$ of (mere) propositions in $U$. The resulting subtype 
then consists of all terms $x:X$ such that $B(x)$ is inhabited. The power type 
of $X$ is then defined from this  as $\mathrm{Pot}(X) \equiv (X \longrightarrow 
\mathrm{Prop}_U)$. Often it is required that this condition does not 
essentially depend on the choice of universes\index{Universe} 
(\enquote{propositional resizing}). 

\section{Type Theory\index{Type theory} and 
Intuitionistic\index{Intuitionism} Higher-Order Logic}

Type theory\index{Type theory} is often referred to as a form of higher-order 
logic, as the internal logic of a given type theory\index{Type theory} 
generalises first-, second-, third-order logic and so on. In the usual 
second-order logic, the quantifiers, unlike in first-order 
logic, can range over predicate variables $\varphi$ and thus 
over subsets 
\[
A=\{x \mid \varphi(x)\}
\]
of domains of individuals, which are traversed by the variable 
$x$.\footnote{See \cite[Ch. I]{manin}.} Third- and higher-order logics 
generalise this to more general iterations of such situations. However, it is 
more precise to say that type theory\index{Type theory} includes higher logic 
through the Curry-Howard correspondence\index{Curry-Howard correspondence}. 
It is described with the expression \enquote{propositions as types}, because it 
forms a correspondence between logical propositions and dependent 
types (see Fig.~\ref{fig:curry_howard}).\footnote{See \cite[Ch. 
3]{voevodsky2013} and \cite{wadler} for 
the Curry-Howard correspondence. It is related to the 
Brouwer-Heyting-Kolmogorov and to the realizability interpretation.} 

\begin{figure}[ht!]
\begin{tabular}{c|c|c|c|c|c|c|c}
$\bot$ & $\top$ & $P \vee Q$ & 
$P \wedge Q$ & $P \Rightarrow Q$ & $\neg P$ & $ \exists a \, P(a)$ & 
$\forall a \, P(a)$ 
\\\hline
$\mathbf{0}$ & $\mathbf{1}$  & $P + Q$        & 
$P \times Q$ & $P \to Q$ & $P \to \mathbf{0}$  & $\sum_{a:A} P(a)$     & 
$\prod_{a:A} P(a)$ 
\end{tabular}
\caption{\label{fig:curry_howard}The 
Curry-Howard correspondence\index{Curry-Howard correspondence} 
provides a link between logical expressions and type\index{Type theory} 
constructions.}
\end{figure}

In this correspondence, each type $A$ is conversely assigned the proposition 
which expresses the information whether $A$ is inhabited or not. For example, 
the type $A+B$ corresponds to the logical proposition $A \vee B$, since $A+B$ 
is inhabited if and only if $A$ or $B$ (or both) are inhabited. 

The types $P+Q$ and $\sum_{a:A} P(a)$ are not mere propositions, even if $P$ 
and $Q$ are, since any term $r:P+Q$ retains the irrelevant information 
whether $r$ lives in $P$ or in $Q$. Therefore, a more precise correspondence 
is 
\[
P \vee Q \equiv ||P+Q|| \text{ and } \exists a \, P(a) \equiv ||\sum_{a:A} 
P(a)||.
\]
An illustrative example is given by the logical proposition 
\[
A \wedge B \Longrightarrow A \vee B. 
\]
It corresponds to the mapping type
\[
A \times B \xlongrightarrow{~~~~~~} ||A+B||. 
\]
A proposition is provable\index{Provability} if the associated type $P$ is 
inhabited by a term $p:P$. The Curry-Howard correspondence\index{Curry-Howard 
correspondence} was originally discovered by Haskell B. Curry\index{Curry, 
Haskell B.} in 1934 as a fact relating logical implications and functions. It 
was William Howard\index{Howard, William} in 1969 who went all the way to 
relating $\lambda$-calculus and natural deduction. Since the $\lambda$-calculus 
can express all computable\index{Computability} functions, any proof is 
understood as a program according to the maxim \enquote{proofs are programs}. In 
this sense, type theory is comparable to the code of a programming language.

In the $\lambda$-calculus, a proof of the above proposition can be expressed as 
$\lambda x \lambda y.\mathrm{inl}(x)$ or $\lambda x \lambda y.\mathrm{inr}(y)$, 
since there are two possible proofs given by composing interference rules for 
$\wedge$ and $\vee$:
\[
\cfrac{A \wedge B}{\cfrac{A}{A \vee B}} \quad \text{resp.} \quad 
\cfrac{A \wedge B}{\cfrac{B}{A \vee B}}
\]
The general picture is therefore that the proof of a logical proposition $P$ 
is a procedure (or program) to obtain a term in the corresponding type $P$. The 
foregoing in the proof tree produces an expression in the $\lambda$-calculus. 
This is obviously a constructive\index{Constructivism} concept of 
provability\index{Provability}. 

\medskip A further example is the implication 
\[
 \neg A \wedge \neg B \Longrightarrow \neg (A \vee B).  
\]
It corresponds to the type\footnote{In \cite[Ch. 1.11]{voevodsky2013} the 
corresponding $p:P$ is constructed.}  
\[
P \colon (A \xlongrightarrow{~~~~~~} \mathbf{0}) \times (B 
\xlongrightarrow{~~~~~~} \mathbf{0}) 
\xlongrightarrow{~~~~~~} (A+B \xlongrightarrow{~~~~~~} \mathbf{0}).
\]

\section{Proof Assistants and Famous Theorems}

Checking the correctness and completeness of proofs of mathematical theorems 
such that no doubts remain is usually hard. In the history of proofs of very 
famous theorems like the solution of Fermat's conjecture, the proof of the 
four-colour theorem\index{Four-colour theorem}, the proof of the 
Feit-Thompson theorem\index{Feit-Thompson theorem}, or the solution 
of the Kepler conjecture\index{Kepler conjecture}, this was an issue in each 
case. As a side remark, Vladimir Voevodsky\index{Voevodsky, Vladimir} also has 
developed dependent type theory\index{Type theory} further together with 
others because a hardly detectable error in one of his papers was pointed out 
to him.\footnote{See \cite{bordg,voevodsky2014}.}    

In recent years, proof assistance systems such as AGDA, ROCQ, ISABELLE/HOL, 
LEAN and others have been developed, with which checking of proofs and  
verification of software with machine support can be demonstrably carried out 
correctly in the framework of intuitionistic\index{Intuitionism}
Martin-L\"of type theory\index{Martin-L\"of, Per}\index{Type theory}. Such 
systems support the work with meaningful commands. There are games, like the 
\enquote{Natural Number Game} for LEAN, which help us to explore such software 
solutions and learn the commands necessary to try proofs 
ourselves.\footnote{See isabelle.in.tum.de, leanprover.github.io, 
rocq-prover.org, wiki.portal.chalmers.se/agda as well as the unimath library on 
github.com.} 

The proof of the Feit-Thompson theorem\index{Feit-Thompson theorem} was 
verified in 2012 by Georges Gonthier\index{Gonthier, Georges} using 
ROCQ.\footnote{See \cite{gonthier2012}.} The proof of the Kepler 
conjecture\index{Kepler conjecture} by Tom Hales\index{Hales, Tom} in 1998 
using the programming language OCAML was later validated by him and others 
using ISABELLE/HOL.\footnote{See \cite{hales2024}.} The proof of the 
four-colour theorem\index{Four-colour theorem} was checked by Georges 
Gonthier\index{Gonthier, Georges} in 2005 using ROCQ.\footnote{See 
\cite{gonthier2005}.}  The verification of a mathematical problem 
by Peter Scholze\index{Scholze, Peter}, called the liquid tensor 
experiment\index{Liquid tensor experiment}, was carried out in 2021 by Johan 
Commelin\index{Commelin, Johan} and others using LEAN.\footnote{See 
\cite{scholze2021}.} Finally, a formal proof of the ontological\index{Ontology} 
proof\index{Proof of God} for the existence of God, was given by Christoph 
Benzm\"uller\index{Benzm\"uller, Christoph} and Bruno Woltzenlogel 
Paleo\index{Woltzenlogel Paleo, Bruno} using ISABELLE/HOL.\footnote{See 
\cite{benzmueller}.}

\section{G\"odel's\index{Goed@G\"odel, Kurt} System $T$}

We have already seen the Ackermann-P\'eter function $A(m,n)$ which is given by
\begin{align*} 
A(0,n)&=n+1 \cr 
A(m+1,0)&=A(m,1) \cr 
A(m+1,n+1)&=A(m,A(m+1,n)).
\end{align*}
This function is not primitive recursive\index{Recursion}, because it can be 
shown that its growth is stronger than that of any primitive 
recursive\index{Recursion} function. There are two ways to define functions 
which are provably computable\index{Computability} in 
Dedekind-Peano arithmetic\index{Dedekind-Peano arithmetic}. On the one hand, 
one can proceed as in the case of the function $A(m,n)$ which is defined by a 
so-called double recursion\index{Recursion}. More generally, one can define 
(nested) multiple recursive\index{Recursion} functions.\footnote{See 
\cite{peter}.} On the other hand, such functions are strongly related to 
transfinite recursion\index{Recursion}. In the case of the function $A(m,n)$ 
the ordinal\index{Ordinal number} number $\omega^2$ is needed since the 
arguments $(m,n)$ may be ordered lexicographically as: 
\[
(0,0)<(0,1)< (0,2)< \cdots < (1,0)<(1,1)<(1,2)< \cdots  
\]
In general, all in Dedekind-Peano 
arithmetic\index{Dedekind-Peano arithmetic} provably computable 
functions\index{Computability} can be classified in a similar way using the 
L\"ob-Schwichtenberg-Wainer 
hierarchy\index{Loeb@L\"ob-Schwichtenberg-Wainer hierarchy} and ordinal 
numbers\index{Ordinal number} $<\epsilon_0$.\footnote{See 
\cite{rose,schwichtenberg}.}

\medskip
A minimal quantifier-free type theory\index{Type theory} that can describe  
such computable\index{Computability} functions is G\"odel's\index{Goed@G\"odel, 
Kurt} system $T$.\footnote{See \cite[Part III]{lambekscott1986}.} 
With the help of $T$, he demonstrated the relative 
consistency\index{Consistency}\index{Contradiction-freeness} of arithmetic 
in a 1958 article in the journal Dialectica.\footnote{See 
\cite[Vol. II, 1958]{goedel} and \cite[Thm. 4.8]{rose}. 
G\"odel's\index{Goed@G\"odel, Kurt} idea originated in 1941. It goes back to 
Hilbert's\index{Hilbert, David} attempt to prove the continuum 
hypothesis, see \cite{hilbert1926}.} G\"odel\index{Goed@G\"odel, Kurt} believed 
that $T$ realised the extension of the finite method hoped for by 
Hilbert\index{Hilbert, David} and himself and is more intuitively accessible 
than the well-ordering of ordinal numbers\index{Ordinal 
number} $\alpha < \varepsilon_0$, which implies the
consistency\index{Consistency}\index{Contradiction-freeness} of arithmetic. 

The system $T$ contains the inductive type of natural numbers $\mathbb{N}$ 
and is closed under mappings and binary products $A \times B$ including their 
two projections. Furthermore, $T$ contains a very general recursion 
rule\index{Recursion} given by a recursor 
\[
R \colon (\mathbb{N} \xlongrightarrow{~~~~~~} (A \xlongrightarrow{~~~~~~} 
(\mathbb{N} \xlongrightarrow{~~~~~~} (A \xlongrightarrow{~~~~~~} A )))) 
\xlongrightarrow{~~~~~~} A 
\]
such that the formula 
\[
R(0,s,t) \equiv s, \quad R(S(r),s,t) \equiv t(r)(R(r,s,t)) 
\]
is satisfied.\footnote{See \cite{tait}. The variables are $r:\mathbb{N}$, 
$s:A$, and $t \colon \mathbb{N} \rightarrow (A \rightarrow A)$.}

\medskip
In this way, $T$ represents all (provably) computable\index{Computability} 
functions $f \colon \mathbb{N}^n \longrightarrow \mathbb{N}$ that arise via  
transfinite recursion\index{Recursion} over ordinal numbers\index{Ordinal 
number} $\alpha < \varepsilon_0$. As conclusion rules, the rules of the 
calculus of inductive constructions are required, as long as they apply to $T$. 

\medskip 
The system $T$ allows the so-called Dialectica 
interpretation\index{Interpretation} of Heyting arithmetic\index{Heyting 
arithmetic}. To this end, G\"odel\index{Goed@G\"odel, Kurt} assigned to each 
arithmetic formula $\varphi(z)$ with free variables $z$ by an inductive 
process a formula 
\[
\varphi^D(z)=\exists x \forall y \, \varphi_D(x,y,z), 
\]
where $\varphi_D$ is a quantifier-free formula in $T$. In addition, 
it holds that for every intuitionistic\index{Intuitionism} proof of 
$\varphi(z)$ in Heyting arithmetic\index{Heyting arithmetic} there exists a 
term $t=t(z)$ in $T$ such that there is a proof of 
$\varphi_D(t,y,z)$ within $T$. With the help of double negation, this results 
in an interpretation\index{Interpretation} of 
Dedekind-Peano arithmetic\index{Dedekind-Peano arithmetic} in $T$.

\section{Topological Interpretation\index{Interpretation}}

Martin Hofmann\index{Hofmann, Martin} and Thomas 
Streicher\index{Streicher, Thomas} discovered around 1994 that Martin-L\"of's 
type theory\index{Martin-L\"of, Per}\index{Type theory} allows a 
homotopically\index{Topology} motivated topological\index{Topology} 
interpretation\index{Interpretation} with non-trivial identity 
types\index{Identity type}, by establishing a connection with 
groupoids\index{Groupoid}. A few years later, Steve Awodey\index{Awodey, 
Steve}, Michael Warren\index{Warren, Michael}, Vladimir 
Voevodsky\index{Voevodsky, Vladimir} and others contributed further building 
blocks to the development of type theory\index{Type theory}, resulting in the 
emergence of homotopy type theory\index{Homotopy type theory}. Types $A$ are 
interpreted as representatives of homotopy types\index{Homotopy theory} of 
topological spaces\index{Topology} and mappings $A \longrightarrow B$ as 
continuous 
mappings.\footnote{See \cite{awodeywarren,hofmannstreicher,voevodsky2013}.}  

For example, the topological interpretation\index{Interpretation} of the 
dependent pair type\index{Type theory} $\sum_{x:A} B(x)$ is a topological 
space\index{Topology} together with a continuous projection mapping 
\[
\mathrm{pr} \colon \sum_{x:A} B(x) \xlongrightarrow{~~~~~~} A,
\]
so that the fibres $\mathrm{pr}^{-1}(x)$ correspond to the types $B(x)$. The 
dependent function type\index{Type theory} $\prod_{x:A} B(x)$ can be 
interpreted as the space of continuous sections of this projection mapping. 

Such an interpretation makes only sense if homotopies\index{Homotopy theory}  
and identity types\index{Identity type} are closely related to each other.  
The identity type\index{Identity type} $\mathrm{Id}_A(a,b)$ is interpreted in 
homotopy type theory\index{Homotopy type theory} as the space of paths between 
$a$ and $b$ and the type 
\[
\mathrm{Id}_A \equiv \sum_{(a,b):A\times A} \mathrm{Id}_A(a,b)
\]
as the path space of $A$. Using constant paths and distinguishing between 
starting and ending points, we get a sequence of arrows
\[
A \xlongrightarrow{~~~~~~} \mathrm{Id}_A \xlongrightarrow{~~~~~~} A \times A
\]
given by $a \mapsto (a,a,\mathrm{refl}_a)$ and $p \mapsto (a,b)$ for 
$p:\mathrm{Id}_A(a,b)$. If a starting point $a$ is fixed, the 
path spaces $\mathrm{Id}_A(a,-)$ correspond to contractible topological 
spaces\index{Topology}, with the homotopy\index{Homotopy theory} 
to the point path $\mathrm{refl}_a$ running over the parametrisation of each 
path. This fact explains why the path induction rule\index{Path induction} is 
valid. The iterated identity type\index{Identity type}
\[
\mathrm{Id}_{\mathrm{Id}_A(a,b)}(p,q)
\]
between two paths $p,q: \mathrm{Id}_A(a,b)$ is interpreted as 
homotopy\index{Homotopy theory} between $p$ and $q$, i.e., as a path in the 
path space\index{Interpretation}. Further iterations result in the structure of 
an infinity groupoid\index{Infinity groupoid}, which is associated with the 
type $A$. This shows that homotopy type theory\index{Homotopy type theory} 
internally represents all features of homotopy types.

\medskip
The topological interpretation\index{Interpretation} of Martin-L\"of's 
type theory\index{Martin-L\"of, Per}\index{Type theory} has 
applications in homotopy theory\index{Homotopy theory} itself. This shows 
in a surprising way the unity of mathematics through the connection 
between homotopy theory\index{Homotopy theory} and the foundations. It is 
important to understand that this interpretation\index{Interpretation} 
is synthetic and only becomes apparent in topologically\index{Topology} 
oriented semantic\index{Semantics} models. Dependent type theory\index{Type 
theory} can also be defined and used without such notions. However, it seems to 
be the case that the nature of type theory\index{Type theory} comes into its 
own when types are provided with a topological\index{Topology} 
interpretation\index{Interpretation}.

\section{A Stratification by Homotopy\index{Homotopy theory} Levels}

Vladimir Voevodsky\index{Voevodsky, Vladimir} has defined a hierarchy for types 
that is ordered by homotopy levels\index{Homotopy theory}. With the help of 
iterated identity types\index{Identity type}, the (mere) propositions as well 
as sets\index{Set theory} and groupoids\index{Groupoid} are generalised. 

\medskip 
In the lowest level $n=-2$, the contractible (resp. point-like) types $A$  
are located with the condition 
\[
\text{is-contr}(A) \equiv \sum_{x:A} \prod_{y:A} \mathrm{Id}_A(x,y) 
\]
is inhabited, i.e., if there is a term $x:A$ such that for all $y:A$ already 
$x=y$ applies. We have used the notation $A\simeq {\mathbf 1}$ for this.

\medskip
In general, we define recursively\index{Recursion}
\[
A \text{ has level } n+1 \colon \text{ for all } x,y : 
A \text{ has } \mathrm{Id}_A(x,y) \text{ level }n.
\]

\medskip 
Thus, at level $n=-1$, the condition becomes 
\[
\prod_{x,y:A} \text{is-contr}(\mathrm{Id}_A(x,y)). 
\]
Surprisingly, one can show that this condition is equivalent to the definition 
\[
\text{is-prop}(A) \equiv \prod_{x,y:A} \mathrm{Id}_A(x,y)  
\]
of a (mere) proposition.\footnote{See \cite[Prop. 12.1.3]{rijke} and \cite[Lem. 
3.11.10]{voevodsky2013}.} Hence the level $n=-1$ consists precisely of (mere) 
propositions. 

\medskip
In this way, we get a stratification on all types, which for $n \ge 0$ 
corresponds to the stratification by homotopy types\index{Homotopy theory}. 
Types $A$ at level $n=0$ are called sets\index{Set theory} and satisfy the 
condition 
\[
\text{is-set}(A) \equiv \prod_{x,y:A} \text{is-prop}(\mathrm{Id}_A(x,y)).    
\]
Why are they called sets? By definition, in this 
case for any two terms $a,b : A$ the identity types\index{Identity type} 
$\mathrm{Id}_A(a,b)$ are (mere) propositions. In the topological 
interpretation\index{Interpretation}, all path spaces $\mathrm{Id}_A(a,b)$ 
are thus either empty or point-like, i.e., there are no non-trivial paths or 
other higher homotopic\index{Homotopy theory} information. The 
collection of all connected components thus forms a 
topological\index{Topology} $0$-type, i.e., a disjoint union 
of point-like spaces, which can be interpreted as a set\index{Set theory}. 

\medskip The natural numbers form a set\index{Set theory}, since its identity 
types are given by (mere) propositions. The proof also uses the encode-decode 
method and shows that the non-trivial identity types\index{Identity type} are 
inhabited by the elements $\mathrm{refl}_n$ with $n:\mathbb{N}$. The proof of 
this result also settles that the successor map $S \colon \mathbb{N} 
\longrightarrow \mathbb{N}$ is injective, as required by 
Dedekind\index{Dedekind, Richard} in his definition of 
$\mathbb{N}$.\footnote{See \cite[Sec. 11.3]{rijke} and 
\cite[Sec. 2.13]{voevodsky2013}.}  

\medskip
The types at level $n=1$ can be understood as groupoids\index{Groupoid} 
and not as categories\index{Category theory}, as might be suspected. In the 
topological interpretation\index{Interpretation}, this means that there are 
non-trivial paths, but paths between paths are trivial. The type 
$\mathbb{S}^1$ is of level $n=1$ and corresponds to the homotopy 
type\index{Homotopy theory} of the circle $S^1=\{(x,y) \in \mathbb{R}^2 \mid 
x^2+y^2=1\}$. However, both are different as types\index{Type theory}, as 
$\mathbb{S}^1$ is not a set\index{Set theory} in the 
type-theoretical\index{Type theory} sense. The circle $S^1$, on the contrary, 
can be understood as a set\index{Set theory} in level $n=0$. 

\medskip 
The higher types at level $n \ge 2$ are given by $n$-groupoids\index{Groupoid}, 
which represent topological $n$-homotopy types\index{Homotopy theory}. Examples 
of types which live in level $n$ but not in lower levels exist for all $n \ge 
2$, even for $n=\infty$. Furthermore, for any type there is an idempotent 
truncation operation $A \longrightarrow ||A||_n$ which maps $A$ to a type in 
level $n$. For $n=-1$ this coincides with the previously defined 
propositional truncation operator $||A||$. Compared with classical homotopy 
theory\index{Homotopy theory}, this is related to Postnikov 
towers\index{Postnikov tower} 
\[
\cdots \xlongrightarrow{~~~~~~} ||A||_{n} \xlongrightarrow{~~~~~~} ||A||_{n-1} 
\xlongrightarrow{~~~~~~} \cdots \xlongrightarrow{~~~~~~} ||A||_{0} 
\xlongrightarrow{~~~~~~} ||A||_{-1} \equiv ||A||.
\]
Truncation is defined as a higher inductive type.\footnote{See 
\cite[Sec. 7.3 and Sec. 8.8]{voevodsky2013} for truncations and examples.}

\medskip
In type theory\index{Type theory} categories\index{Category theory} can be 
defined synthetically. Besides numerous assumptions, which reflect the axioms 
of a category\index{Category theory} with objects in a corresponding type 
$A_0$, the essential condition for this is that the natural map for objects 
$a,b:A_0$ 
\[
\text{id-to-iso} \colon \mathrm{Id}_{A_0}(a,b) 
\xlongrightarrow{~~~~~~} \mathrm{Iso}_A(a,b)
\]
is an equivalence\index{Equivalence}. Here, on the left hand side, one has 
the identity type\index{Identity type} inside the object type $A_0$ of $A$, and 
on the right hand the notion of an isomorphism\index{Isomorphism} in the 
categorical\index{Category theory} sense, i.e., with mutually inverse 
morphisms $f,g$ and the type-theoretic\index{Type theory} propositional 
equalities\index{Equality}\index{Identity} concerning both compositions 
$f \circ g$ and $g \circ f$. An example for such a category\index{Category 
theory} is $\mathbf{Set}_U$, the category\index{Category theory} of sets 
\index{Set theory} inside some universe\index{Universe} $U$. With the help of 
truncation, the type of cardinal numbers\index{Cardinal number} in $U$ can be 
defined by\footnote{See \cite[Sec. 9.1]{voevodsky2013}.} 
\[
\mathbf{Card}_U \equiv ||\mathbf{Set}_U||_0. 
\]

\section{Isomorphism\index{Isomorphism} and Equivalence\index{Equivalence}} 

The concept of isomorphism\index{Isomorphism} in mathematics is extremely 
important and omnipresent. Two objects $A$ and $B$ in a given 
category\index{Category theory} are isomorphic if there are morphisms, i.e., 
structure-preserving arrows
\[
f \colon A \xlongrightarrow{~~~~~~} B 
\]
and 
\[
g \colon B \xlongrightarrow{~~~~~~} A,  
\]
such that the compositions $f \circ g$ and $g \circ f$ coincide with the 
identity mappings\index{Identity}\index{Equality} $\mathrm{id}_B$ and 
$\mathrm{id}_A$. If this is the case, we call $f$ and $g$ 
isomorphisms\index{Isomorphism}.

\medskip
To give a simple example, consider groups with $2$ elements within the 
category of groups\index{Category theory}. There is only one of these up to 
isomorphism\index{Isomorphism}. It is customary to write this group either as 
the additive group 
\[
\mathbb{Z}/2\mathbb{Z}=\{0,1\} \text{ with } 1+1=0
\]
or as the multiplicative group
\[
\mu_2=\{\pm 1\} \text{ with } (-1) \cdot (-1)=+1.
\]
The isomorphism\index{Isomorphism}
\[
f \colon \mathbb{Z}/2\mathbb{Z} \xlongrightarrow{~~~~~~} \mu_2
\]
is given by the assignment $f(0)=+1$ and $f(1)=-1$. 

Even though this example seems simple, many constructions in modern mathematics 
involve complications with isomorphisms\index{Isomorphism}. There, isomorphic 
objects should be identified as much as possible without concrete 
identifications being known or unique. A good example of this is the 
homeomorphism problem\index{Homeomorphism} in topology\index{Topology}. The 
question of whether two topological\index{Topology} spaces $X$ and $Y$ are 
isomorphic, or equivalently homeomorphic\index{Homeomorphism}, is an 
undecidable\index{Undecidability} problem according to a theorem by 
Markov\index{Markov, Andrey}. Even if $X$ and $Y$ are 
isomorphic\index{Isomorphism}, a concrete isomorphism\index{Isomorphism} is 
usually not easy to specify. On the other hand, the additional information 
contained in existing isomorphisms\index{Isomorphism} is often very important, 
as we have seen with the infinity groupoids\index{Infinity groupoid}. 
When studying such things, one can imagine that the concept of 
equality\index{Equality}\index{Identity} has a deep philosophical meaning 
beyond mathematics. In Martin-L\"of type theory\index{Martin-L\"of, 
Per}\index{Type theory}, this is expressed in the existence\index{Existence} of 
identity types\index{Identity type}. 

\medskip
In type theory\index{Type theory}, for two given types $A,B$ in a 
universe\index{Universe} $U$, the isomorphism type\index{Isomorphism}
\[
\mathrm{Iso}(A,B) \text{ or } (A \cong B)
\]
and the equivalence type\index{Equivalence}
\[
\mathrm{Eq}(A,B) \text{ or } (A \simeq B)
\]
are defined. 

\medskip
The isomorphism type\index{Isomorphism} $\mathrm{Iso}(A,B)$ is -- as we did at 
the beginning of this chapter for categories\index{Category theory} -- defined 
by the existence\index{Existence} of mutually inverse mappings $f \colon A 
\longrightarrow B$ and $g \colon B \longrightarrow A$. Here we use 
type-theoretic\index{Type theory} equality\index{Equality}\index{Identity} 
though, i.e., we require that 
\[
f \circ g = \mathrm{id}_B \text{ and } g \circ f = \mathrm{id}_A.  
\]
In the case where the types $A$ and $B$ are sets\index{Set theory} in the sense 
of type theory\index{Type theory} this definition makes sense, and the natural 
mapping
\[
\mathrm{Id}_U(A,B) \xlongrightarrow{~~~~~~} \mathrm{Iso}(A,B)
\]
between identity type\index{Identity type} and isomorphism 
type\index{Isomorphism} is an isomorphism\index{Isomorphism}. From level 
$n \ge 1$ on, the notion of isomorphism\index{Isomorphism} is usually not used 
anymore since the identity types\index{Identity type} are too large. This was 
basically already the discovery of Martin Hofmann\index{Hofmann, Martin} and 
Thomas Streicher\index{Streicher, Thomas}, which led to their groupoid 
model\index{Groupoid}\index{Model theory} of type theory\index{Type 
theory}.\footnote{See \cite[Sec. 5]{awodey2014}, \cite{hofmannstreicher} and 
\cite[Sec. 9.1]{voevodsky2013}.} 

\medskip
The notion of equivalence\index{Equivalence} in 
type theory\index{Type theory} is based on the notion of 
homotopy\index{Homotopy theory} from topology. Two dependent functions 
\[
f,g \colon \prod_{x:A} B(x)
\]
are homotopic, denoted by $f \sim g$, if there is a term in  
\[
(f \sim g) \equiv \prod_{x:A} \mathrm{Id}_{B(x)}(f(x),g(x)).
\]
This is a priori not equivalent to $f=g$ in the dependent function type. We 
will get to know the univalence axiom\index{Univalence} which forces this 
property.

\medskip
Using the notion of homotopy\index{Homotopy theory}, 
equivalences\index{Equivalence} between two types\index{Type theory} $A,B$ can 
be defined. However, we will not imitate the definition of 
isomorphisms\index{Isomorphism} and just use homotopies\index{Homotopy theory} 
instead of identities\index{Identity type}. Instead, we declare a function 
\[
f \colon A  \xlongrightarrow{~~~~~~} B 
\]
to be an equivalence\index{Equivalence}, if it has a section $g$ and a 
retraction $h$, i.e., one has 
\[
f \circ g \sim \mathrm{id}_B \text{ and } h \circ f \sim \mathrm{id}_A.  
\]
In other words, a function $f \colon A \longrightarrow B$ is called an 
equivalence\index{Equivalence} if the type 
\[
\text{is-equiv}(f) \equiv \biggl(\sum_{g:B \to A}  (f \circ g \sim 
\mathrm{id}_B) \times \sum_{h:B \to A} (h \circ f \sim \mathrm{id}_A)\biggr) 
\]
is inhabited. Alternatively, it can be shown that a mapping $f \colon A 
\longrightarrow B$ is an equivalence\index{Equivalence}, if 
all homotopy fibres\index{Homotopy theory} 
\[
\mathrm{hfiber}(f,b) \equiv \sum_{a:A} \mathrm{Id}_B(f(a),b)
\]
of $f$ are contractible for each $b$. Both definitions generate the same  
equivalence type\index{Equivalence} 
\[
\mathrm{Eq}(A,B) \equiv \sum_{f:A \to B} \prod_{b:B}  
\text{is-contr}(\mathrm{hfiber}(f,b)).
\]
Another notation for this type is $(A \simeq B)$.\footnote{See \cite[Sec. 
9]{rijke} and \cite[Sec. 4.4 and 4.5]{voevodsky2013}.}

\section{Univalence\index{Univalence}}

Vladimir Voevodsky\index{Voevodsky, Vladimir} has further developed the type 
theory\index{Martin-L\"of, Per}\index{Type theory} of Martin-L\"of with its 
central identity types\index{Identity type}. When considering simplicial 
sets\index{Simplicial set} as models\index{Model theory} of type 
theory\index{Type theory}, he discovered the univalence 
axiom\index{Univalence}.\footnote{See the introduction 
of \cite{voevodsky2013} and \cite{kapulkin2021}.} This axiom is not a rule 
built 
into type theory\index{Type theory} like the rest, but an axiom that must be 
imposed afterwards. It does not hold in all models\index{Model theory} of type 
theory\index{Type theory} and its significance is currently not entirely clear. 
Nevertheless, we will try to explain the univalence axiom\index{Univalence}.

For two types $A,B$ in a universe\index{Universe} $U$, there is always 
a mapping\footnote{See \cite[Sec. 2.10]{voevodsky2013}.}
\[
v \colon \mathrm{Id}_U(A,B) \xlongrightarrow{~~~~~~} \mathrm{Eq}(A,B).
\]
The univalence axiom\index{Univalence} states that the mapping 
$v$ itself is an equivalence\index{Equivalence} in a higher 
universe\index{Universe} $U$, i.e., 
\[
\mathrm{Id}_U(A,B) \simeq \mathrm{Eq}(A,B).
\] 
If $A$ and $B$ are sets in the sense of type theory\index{Type theory}, then 
the concept of equivalence\index{Equivalence} coincides with the concept of 
isomorphism\index{Isomorphism} of sets. 

For each property $P$, then with respect to equivalence\index{Equivalence}, one 
has a variant of Leibniz's\index{Leibniz, Gottfried Wilhelm} invariance 
rule\index{Leibniz's invariance rule}
\[
\frac{A \simeq B \quad P(A)}{P(B)},
\]
i.e., the property $P$ is preserved when replacing 
equivalent\index{Equivalence} types. This follows from the univalence 
axiom\index{Univalence}, as this invariance rule\index{Leibniz's invariance 
rule} for identity\index{Identity}\index{Equality} follows from 
path induction\index{Path induction} and under the mapping 
\[
v \colon \mathrm{Id}_U(A,B) \xlongrightarrow{~~~~~~} \mathrm{Eq}(A,B)
\]
is still preserved.\footnote{See \cite{awodey2014} for a proof.} 
Topologically, this principle is trivial, because the assignment $A \mapsto 
P(A)$ takes discrete values $\bot$ and $\top$ and is continuous. From this even 
follows the seemingly stronger statement
\[
A \simeq B \xRightarrow{~~} A = B, 
\]
which results from the inverse $v^{-1}$ of the mapping $v$ (up to 
homotopy\index{Homotopy theory}). This can be interpreted as 
equivalences\index{Equivalence} being able to be upgraded to 
identities\index{Identity}\index{Equality}. 

The hierarchy of identity\index{Identity type}\index{Identity}\index{Equality}, 
isomorphism\index{Isomorphism} and equivalence types\index{Equivalence} 
provides room for some open questions. For example, it is not clear in which 
categorical\index{Category theory} semantic\index{Semantics} 
interpretations\index{Interpretation} of type theory\index{Type theory} the 
univalence axiom\index{Univalence} is valid and how 
univalence\index{Univalence} can be added to categories\index{Category 
theory} as an additional property.\footnote{This is related to the so-called 
Rezk completion, see \cite[Sec. 9.9]{voevodsky2013}.}   

\medskip
The univalence axiom\index{Univalence} has consequences for the 
question of function extensionality\index{Extensionality}. If  
\[
f,g: \prod_{x:A} B(x) 
\]
are two dependent functions, then there is a natural function 
\[
\mathrm{Id}_{\prod_{x:A} B(x)} (f,g) \xlongrightarrow{~~~~~~} \prod_{x:A} 
\mathrm{Id}_{B(x)}(f(x),g(x)).
\]
If this function is an equivalence\index{Equivalence} then function 
extensionality\index{Extensionality} holds by definition.\footnote{See 
\cite[Axiom 2.9.3]{voevodsky2013}.} The univalence 
axiom\index{Univalence} implies the equivalence\index{Equivalence}  
\[
\biggl( \prod_{x:A} \text{is-contr} B(x) \biggr) \xlongrightarrow{~~~~~~}
\text{is-contr} \biggl(\prod_{x:A} B(x) \biggr),
\]
which resembles a weak form of function 
extensionality\index{Extensionality}.\footnote{See 
\cite[Thm. 4.9.4]{voevodsky2013}.} However, one can show that this implies the 
strong form of function extensionality\index{Extensionality}.\footnote{See 
\cite[Thm. 4.9.5]{voevodsky2013}.}   

\section{Categorical\index{Category theory} 
Interpretation\index{Interpretation} of Type Theory\index{Type theory}}

Every type theory\index{Type theory} $T$ has a natural 
interpretation\index{Interpretation} in its own syntactic   
category\index{Category theory} $\mathcal{S}_T$. It is widely believed that 
every interpretation\index{Interpretation} of $T$ in a model 
category\index{Model category} (or a path category\index{Path category})  
$\mathcal{C}$ can be factorised by the syntactic   
category\index{Category theory} $\mathcal{S}_T$, i.e., is given by a 
functor\index{Functor} of categories\index{Category theory}
\[
\mathcal{S}_T \xlongrightarrow{~~~~~~} \mathcal{C}.
\]
This initiality conjecture has not been completely settled, but the approach by 
Thomas Streicher\index{Streicher, Thomas} was significantly extended by Menno 
de Boer\index{Boer, Menno de}.\footnote{See \cite{boer,streicher1991}.} 

\medskip 
In such interpretations\index{Interpretation}, two essential problems occur. 
On the one hand, general homotopical\index{Homotopy theory} 
categories\index{Category theory} a priori do not have the necessary properties 
to interpret the structures in type theory\index{Type theory}, like the 
identity type\index{Identity type} and others. In addition, a fundamental 
coherence problem arises when substitutions in type theory\index{Type theory} 
(or any other deductive system\index{Deductive system}) are interpreted 
as pullbacks of fibrations, since in categories\index{Category theory} the 
concatenation of pullbacks is only unique up to isomorphism\index{Isomorphism}.
To solve both these issues, Vladimir Voevodsky\index{Voevodsky, Vladimir} has 
suggested the construction of contextual categories\index{Category theory} 
($C$-systems) arising from homotopical\index{Homotopy theory} 
categories\index{Category theory} with universes\index{Universe}. $C$-Systems 
are named after John Cartmell\index{Cartmell, John} who invented the notion of 
contextual (i.e., syntactic) categories in his framework of generalised 
algebraic theories.\footnote{For contextual categories see \cite{cartmell1986}, 
and \cite{streicher1991}. $C$-systems were introduced in 
\cite{voevodsky2015,voevodsky2023}. With the help of a construction by Jean 
B\'enabou\index{B\'enabou, Jean} \cite{benabou1985}, each fibration can be 
\enquote{split} in such a way that the coherence problem can be successfully 
tackled. Possible further solutions have been proposed by Andrew 
Pitts\index{Pitts, Andrew} (split-type categories), Peter Dybjer\index{Dybjer, 
Peter} (categories with families), Martin Hofmann\index{Hofmann, Martin} 
(categories with attributes), and Steve Awodey\index{Awodey, Steve} (natural 
models). See \cite{ahrens2018}, \cite{hofmann1994} for a survey.}  

\medskip The syntactic category\index{Category theory} $\mathcal{S}_T$ of a 
type theory $T$ is also called contextual category\index{Category theory}, 
since the definition of $\mathcal{S}_T$ leads to the notion of contextual 
categories\index{Category theory}. In a contextual category\index{Category 
theory} $\mathcal{C}$ all objects have a unique \enquote{length} $n \ge 0$. 
There is a unique, terminal object $1$ of length $0$ in $\mathrm{Ob}_0 
\mathcal{C}$. An object $X$ of length $n>0$ in $\mathrm{Ob}_n \mathcal{C}$ has 
a projection map 
\[
p_X \colon X \rightarrow \mathrm{ft} X,  
\]
where $\mathrm{ft} X$ is an assigned object of length $n-1$ in 
$\mathrm{Ob}_{n-1} \mathcal{C}$. In the case of syntactic 
categories\index{Category theory} such maps correspond to projections of 
contexts
\[
x_0:A_0,\ldots,x_{n}:A_{n}(x_0,\ldots,x_{n-1}) \mapsto
x_0:A_0,\ldots,x_{n-1}:A_{n-1}(x_0,\ldots,x_{n-2}).
\]
The most important property of contextual categories\index{Category theory} 
claims that for every morphism $f \colon Y \rightarrow \mathrm{ft} X$ there is 
an object $f^* X$ such that the diagram displayed in 
Fig.~\ref{fig:type_contextual} is a pullback diagram.\footnote{See 
\cite[Def. 1.2.1]{kapulkin2021}.}  

\begin{figure}[ht!]
\begin{tikzpicture}
  \matrix (m) [matrix of math nodes,row sep=3em,column sep=4em,minimum 
width=2em]
  {
     f^*X & X \\
      Y & \mathrm{ft} X \\};
  \path[-stealth]
    (m-1-1) edge node [left] {$p_{f^*X}$} (m-2-1)
            edge node [above] {$q(f,X)$} (m-1-2)
    (m-2-1.east|-m-2-2) edge node [below] {$f$}
            node [above] {} (m-2-2)
    (m-1-2) edge node [right] {$p_X$} (m-2-2);
\draw +(-.2,0) -- +(0,0)  -- +(0,.2);
\end{tikzpicture}
\caption{\label{fig:type_contextual}This cartesian diagram is fundamental for 
contextual categories\index{Category theory}.}
\end{figure}

\medskip Now we want to demonstrate how the additional structure of a 
universe\index{Universe} $U$ in a homotopical\index{Homotopy theory} 
category\index{Category theory} $\mathcal{C}$ leads to a new contextual 
category\index{Category theory} $\mathcal{C}_U$ which is very closely connected 
with $\mathcal{C}$. The philosophy behind this is that in some sense contexts 
correspond to families. A universe\index{Universe} in a category $\mathcal{E}$ 
consists of an object $U$ together with a function $p \colon \tilde U 
\longrightarrow U$ (universal family) such that for each map 
$f \colon X \longrightarrow U$ there is a choice of a pullback square as in the 
Fig.~\ref{fig:type_universe}.\footnote{See \cite[Def. 1.3.1]{kapulkin2021}.} 

\begin{figure}[ht!]
\begin{tikzpicture}
  \matrix (m) [matrix of math nodes,row sep=3em,column sep=4em,minimum 
width=2em]
  {
     (X;f) & \tilde U \\
      X & U \\};
  \path[-stealth]
    (m-1-1) edge node [left] {$p_{X;f}$} (m-2-1)
            edge node [above] {$Q(f)$} (m-1-2)
    (m-2-1.east|-m-2-2) edge node [below] {$f$}
            node [above] {} (m-2-2)
    (m-1-2) edge node [right] {$p$} (m-2-2);
\draw +(-.2,0) -- +(0,0)  -- +(0,.2);
\end{tikzpicture}
\caption{\label{fig:type_universe}In categories\index{Category theory} with 
universes\index{Universe} there are canonical pullbacks from the universal 
family $p \colon \tilde U \rightarrow U$.}
\end{figure}

\medskip The contextual category\index{Category theory} $\mathcal{C}_U$ arises 
as follows. The unique object of length $0$ is the terminal object $1$. Every 
object $X$ in $\mathcal{C}$ arises from a morphism $f_1 \colon 1 
\longrightarrow U$ such that $(1;f_1)=X$. By iteration of the construction in 
the diagram we can form objects $(1;f_1,\ldots,f_n)$ in $\mathcal{C}$. With this 
we define objects of length $n\ge 1$ in $\mathcal{C}_U$ as  
\[
\mathrm{Ob}_{n} \mathcal{C}_U := \{(f_1,\ldots,f_{n}) \in (\mathrm{Mor} 
\mathcal {C})^n \mid f_{i} \colon (1;f_1,\ldots,f_{i-1}) \rightarrow U \; (1 
\le i \le n)\}. 
\]
The resulting category\index{Category theory} $\mathcal{C}_U$ is indeed 
contextual and comes with a fully faithful, but not always injective, functor 
$\mathcal{C}_U \rightarrow \mathcal{C}$.\footnote{See \cite[Def. 1.2.1 and 
Prop. 1.3.3]{kapulkin2021} and \cite[Sec. 2]{voevodsky2015}.}

\medskip
Universes\index{Universe} help to solve the coherence problem, since they 
form a sort of fixed framing for families. More precisely, for morphisms 
\[
X'' \xlongrightarrow{~~~h~~~} X' \xlongrightarrow{~~~g~~~} X 
\xlongrightarrow{~~~f~~~} U
\]
one can show that the pullbacks of the family $(X;f) \longrightarrow X$ satisfy 
\[
(g \circ h)^* (X;f) = h^*(g^*(X;f)).
\]
This solves the coherence problem.\footnote{See \cite[Sec. 2]{voevodsky2015}.} 
 
\medskip 
The syntactic category\index{Category theory} $\mathcal{S}_T$ of a dependent 
type theory\index{Type theory} $T$ is a path category\index{Path category} 
under mild conditions.\footnote{See \cite{vdB}.} The fibrations correspond to 
certain projections of contexts. The path objects $\mathrm{Path}(A)$ are given 
by the contexts
\[
x:A, y:A, p:\mathrm{Id}_A(x,y) 
\]
and the maps in Fig.~\ref{fig:type_pathobject} are given by $r(a) \equiv 
(a,a,\mathrm{refl}_a)$ and $\mathrm{pr}(x,y,p) \equiv (x,y)$. 

\begin{figure}[ht!]
\begin{tikzpicture}
  \matrix (m) [matrix of math nodes,row sep=3em,column sep=4em,minimum 
width=2em]
  {
      & \mathrm{Path}(A) \\
      A & A \times A \\};
  \path[-stealth]
    (m-2-1.east|-m-2-2) edge node [below] {$(\mathrm{id},\mathrm{id}$)}
            node [above] {} (m-2-2)
    (m-1-2) edge node [right] {$\mathrm{pr}$} (m-2-2)
    (m-2-1) edge node [above] {$r$} (m-1-2);
\end{tikzpicture}
\caption{\label{fig:type_pathobject}Path object in type theory\index{Type 
theory} are defined in the same way as the corresponding objects in homotopy 
theory\index{Homotopy theory}.}
\end{figure}

Universes\index{Universe} allow to define identity structures (path spaces) 
on any category\index{Category theory} with universe\index{Universe}. To do 
this, it is sufficient to treat the universal case and define a map 
\[
 \mathrm{Id} \colon \tilde U \times_U \tilde U 
\xlongrightarrow{~~~~~~} U  
\]
together with a specified lift $r$ of the diagonal (see 
Fig.~\ref{fig:type_id}). Then all path spaces for objects $X$ can be obtained 
via pullback maps.

\begin{figure}[ht!]
\begin{tikzpicture}
  \matrix (m) [matrix of math nodes,row sep=3em,column sep=4em,minimum 
width=2em]
  {
    & \mathrm{Id}^* \tilde U & \tilde U \\
    \tilde U  &  \tilde U  \times_U \tilde U & U \\};
  \path[-stealth]
    (m-2-1) edge node [left] [below] {$(\mathrm{id},\mathrm{id})$} (m-2-2)
    (m-2-1) edge [dashed] node [above] {$r$} (m-1-2)
    (m-2-2) edge node [below] {$\mathrm{Id}$} (m-2-3)
    (m-1-2) edge node [right] {$\mathrm{Id}^*p$} (m-2-2)
    (m-1-2) edge node [above] {$Q(\mathrm{Id})$} (m-1-3)
    (m-1-3) edge node [right] {$p$} (m-2-3);
\draw +(1.3,0) -- +(1.5,0)  -- +(1.5,.2);
\end{tikzpicture}
\caption{\label{fig:type_id}Identity\index{Identity}\index{Equality} 
structures in categories\index{Category theory} with universes\index{Universe} 
can be defined by a pullback using the universal family $p \colon \tilde U 
\rightarrow U$.}
\end{figure}

\medskip A universe\index{Universe} $U$ in a homotopical\index{Homotopy 
theory} category\index{Category theory} is called univalent\index{Univalence} 
if the map $r$ is an equivalence\index{Equivalence}. Further logical structures 
are defined using universes\index{Universe} as well. Simplicial localisations 
preserve such structures.\footnote{See \cite[App.]{kapulkin2021}.}

\medskip Summarizing, we have seen how to construct from topologically relevant 
model (or path\index{Path category}) categories\index{Model category} suitable 
target categories\index{Category theory} for the 
interpretation\index{Interpretation} of contextual categories\index{Category 
theory} of dependent type theories\index{Type theory}. Simplicial localisations 
of such categories\index{Category theory} yield locally cartesian 
closed\index{Cartesian closed category} 
$(\infty,1)$-categories\index{Infinity category} resp. elementary 
$(\infty,1)$-topoi\index{Elementary topos}\index{Topos} as corresponding 
target categories\index{Category theory} in the world of 
$(\infty,1)$-categories\index{Infinity category} (see 
Fig.~\ref{fig:great_picture}).

\begin{figure}[ht!]
\begin{tikzpicture}
  \matrix (m) [matrix of math nodes,row sep=3em,column sep=4em,minimum 
width=2em]
  {
    \text{Contextual categories in HoTT} & \text{Elementary 
$(\infty,1)$-topoi} 
\\
    \text{Contextual categories in MLTT}& \text{LCC 
$(\infty,1)$-categories.} 
\\};
  \path[-stealth]
    (m-1-1) edge[>->] node [left] {} (m-2-1)
            edge node [above] {} (m-1-2)
    (m-2-1.east|-m-2-2) edge node [above] {}
            node [above] {} (m-2-2)
    (m-1-2) edge[>->] node [right] {} (m-2-2);
\end{tikzpicture}
\caption{\label{fig:great_picture}Locally cartesian closed (LCC)
$(\infty,1)$-categories\index{Infinity category}  and elementary 
$(\infty,1)$-topoi\index{Elementary topos}\index{Topos}, are possible target 
categories of contextual categories\index{Category theory} with logical rules 
from Martin-L\"of type theory\index{Type theory}\index{Martin-L\"of, Per} 
(MLTT) and from homotopy type theory\index{Homotopy type theory} (HoTT).}
\end{figure}

Elementary infinity topoi\index{Elementary topos}\index{Topos} satisfy the same 
conditions as elementary topoi\index{Elementary topos}\index{Topos} (with 
limits resp. colimits as homotopy (co)limits), and have in addition 
sufficiently many object classifiers (i.e., univalent\index{Univalence} 
universes\index{Universe}).\footnote{See \cite{shulman2018,shulman2019}.} 

\chapter{Semantics\index{Semantics} and Reality}

In today's societies, the process of truth finding\index{Truth} 
is fraught with difficulties. The truth\index{Truth} of 
scientific theses must be constantly questioned, because we only get closer to 
it in individual steps and even the scientific operation is not error-free. In 
the non-scientific area, complex or conditional statements are difficult to 
convey and some insights do not reach the public at all or slowly. In addition, 
there are different schools of thought, prejudices, influences on opinions 
and not least stupidity and ignorance. 

Leibniz's\index{Leibniz, Gottfried Wilhelm} dream, to be able to determine and 
prove the truth\index{Truth} of assertions in all areas, not only encounters 
such banal, but also systematic limits of undecidability\index{Undecidability} 
and incompleteness\index{Incompleteness theorem}. 
Truth finding\index{Truth} is nevertheless unavoidable. The search for 
a philosophical concept of truth\index{Truth} goes far before 
Leibniz\index{Leibniz, Gottfried Wilhelm}. The quality of this 
concept is described by Thomas Grundmann\index{Grundmann, Thomas} in reference 
to the correspondence theory\index{Correspondence theory of truth} in his book 
on epistemology\index{Epistemology}:
\begin{quote}
The concept of truth is strictly speaking not an epistemological 
concept. It picks out a relation between mind and world, namely the 
reference to something that is the case. It is therefore very similar to other 
basic semantic concepts such as reference or meaning. Although 
truth is a semantic concept, it plays a prominent role within 
epistemology.\footnote{Machine translation of \cite[Ch. 2]{grundmann}.}
\end{quote}
Over the centuries, an extremely rich literature has emerged, with partly 
controversial definitions of truth\index{Truth}, in which truth 
bearers\index{Truth bearer} and truth makers\index{Truth maker} are understood 
differently.\footnote{The correspondence theory has different variants in the 
literature. In addition, there are many attitudes of scepticism that reject 
aspects of reality or the concept of truth, or shift them to social consensus 
and speech acts. Descriptions of these theories can be found in 
\cite{grundmann,halbach1996,reclam,schrenk}.} We will completely neglect most 
of these theories and focus on the language-analytical 
semantic\index{Semantics} truth theories\index{Truth}.

The latter were introduced by Tarski\index{Tarski, Alfred} in the 20th century 
with the help of formal languages\index{Formal language}. His method led to the 
concept of semantics\index{Semantics}, which he could define via the concept 
of a metalanguage\index{Metalanguage} that is richer than the given object 
language\index{Object language}.\footnote{See \cite{brendel2015}, 
\cite{halbach1996}, \cite{kripke1959,kripke1965,kripke1975,kripke1980}, 
\cite{tarski1935,tarski1969}, and \cite{schrenk}.} 
Suitable semantics\index{Semantics} illustrate mathematics and thus simplify 
understanding. In mathematics, set theory\index{Set theory} with the 
Zermelo-Fraenkel axioms\index{Zermelo-Fraenkel axioms} or category 
theory\index{Category theory} is often used as a 
metalanguage\index{Metalanguage}. This gave rise to the field of model 
theory\index{Model theory}. However, this is only one of many options, as forms 
of type theory\index{Type theory}, category theory\index{Category theory} or 
set theory\index{Set theory} can be used both as object language \index{Object 
language} and as metalanguage\index{Metalanguage} and are mutually 
interpretable. The question of 
consistency\index{Consistency}\index{Contradiction-freeness} of a given 
theory does not become easier, because richer languages rather complicate this 
question than solve it. In this respect, semantic\index{Semantics} truth 
theories\index{Truth} are closely related to the coherence 
theory\index{Coherence theory of truth} of truth.

The concepts of equivalence\index{Equivalence} and 
equality\index{Equality}\index{Identity}, which we encountered in type 
theory\index{Type theory} and in higher categories\index{Category theory}, 
allows a new structural view of the foundations of mathematics. The associated 
semantic\index{Semantics} concept of truth\index{Truth} -- in conjunction with 
the concept of equivalence\index{Equivalence} -- enables insights that are 
relevant with regard to natural and philosophical questions of science.

\section{Finding the Truth\index{Truth}}

The pre-Socratic scholar Anaximander of Miletus\index{Anaximander of Miletus} 
lived around 600 BC. Although virtually nothing of his writings has survived, 
he is considered a significant natural scientist due to posthumous 
descriptions. Pliny the Elder\index{Pliny, the Elder} wrote that Anaximander 
had opened the door to nature:
\begin{quote}
Rerum fores aperuisse, Anaximander Miletus traditur primus.\footnote{See 
\cite[Book II]{plinius}.} 
\end{quote}
In one of his books, Carlo Rovelli\index{Rovelli, Carlo} vividly described how 
the development of natural sciences was influenced by 
Anaximander\index{Anaximander of Miletus} experienced a tremendous upswing 
because he replaced traditional flawed views with scientific ways of thinking. 
Rovelli\index{Rovelli, Carlo} writes: 
\begin{quote}
Anaximander opened the door to physics, geography, the study of meteorological 
phenomena, and biology. Beyond these important contributions, he set in motion 
the process that led to the rethinking of our world view: our way of gaining 
knowledge, which is based on the rebellion against the obvious. In this sense, 
Anaximander can undoubtedly be called one of the fathers of scientific 
thinking.\footnote{Machine translation of \cite[p. 12]{rovelli}.}  
\end{quote}
This was vividly documented in Anaximander's\index{Anaximander of Miletus} 
observations of celestial bodies. His logical arguments support the sphericity 
of the Earth's surface -- as opposed to the Earth as a disc -- and its 
free-floating position in space through observations and precise conclusions. 
Mathematics as a science had already developed outside the time of 
Anaximander\index{Anaximander of Miletus} earlier in Babylonian and Sumerian 
cultures and after this period within the influential school of 
Pythagoras\index{Pythagoras}. 

Rovelli's\index{Rovelli, Carlo} example wonderfully illustrates the complexity 
of the concept of truth\index{Truth} beyond Anaximander's\index{Anaximander of 
Miletus} considerations. Because when we say that the Earth's surface is a 
sphere, this statement is not literally true, not even the weakening that it 
is an ellipsoid. The poles of the Earth are flattened in different ways. If we 
look even more closely, we see that the Earth has an irregular surface and is  
filled with matter inside, which we do not fully understand in its innermost 
structure and whose composition has gaps in space. Therefore, the best form of 
this statement is that the Earth's surface is approximately spherical. The 
formulation by means of an approximation is not imprecise, because with a 
clever mathematical description, the deviation of the Earth's surface from an 
ideal sphere can be numerically limited by an inequality that a disc shape 
cannot fulfil. Thus, after a few steps, we have arrived at a formulation that 
excludes all claims that describe the Earth as a disc, which was the goal of 
all 
these considerations. Anaximander's\index{Anaximander of Miletus} arguments 
historically provided the first proof of this. In this example, we see very 
well 
that it can be difficult and sometimes requires several steps of thought to 
precisely grasp truth\index{Truth}. 

\section{The Opposite of Truth\index{Truth}}

There is also the opposite of truth\index{Truth}, untruth. 
Leibniz\index{Leibniz, Gottfried Wilhelm} had already thought about the origin 
of untrue statements and attributed this to the lack or ignoring of evidence 
and other errors. He wrote about this in the \enquote{Nouveaux essais sur 
l'entendment humain} from 1704: 
\begin{quote}
After we have spoken sufficiently of all the means which let us recognise or 
suspect the truth, we want to say something about our errors and incorrect 
judgements\index{Judgement}. People must often be mistaken because there are so 
many disagreements among them. The causes of this can be traced back to the 
following four: 1) the lack of evidence, 2) the little skill to use it, 3) the 
lack of good will to make use of it, 4) the wrong rules of 
probability.\footnote{Machine translation of \cite[Ser. VI, Vol. 6]{leibniz}. 
German translation by Ernst Cassirer\index{Cassirer, Ernst} (digitally 
available at www.project-gutenberg.org).}
\end{quote}
By false probability rules, Leibniz\index{Leibniz, Gottfried Wilhelm} 
understood unreliable guesses and estimates that people use in their 
judgements\index{Judgement}. Indeed, humans often fall prey to irrational 
prejudices.\footnote{Astonishing examples of this can be found in 
\cite{rosling}.} 

False statements can arise from errors that are oversights, 
if they are not consciously made. Closely related to the error is the concept 
of fallacy, which means a faulty implication in the sense of logic 
and in which the acting person feels justified. In colloquial language 
and the logic used there, fallacies are often made. Particularly 
noticeable is that the implication $A \Rightarrow B$ is often reversed without 
reason. Let's look at two examples and consider the implications: 
\begin{quote}
I am currently eating a soup $\Longrightarrow$ I am alive
\end{quote}
and 
\begin{quote}
My boss does not appreciate my performance $\Longrightarrow$ My salary 
will not be increased. 
\end{quote}
Both implications are valid in general. This is independent 
of whether I am currently eating a soup or my boss really 
does not appreciate my performance, because an implication $A \Rightarrow B$ is 
always true, if the statement $A$ is false. In the first implication, the 
reversal 
\begin{quote}
I am alive $\Longrightarrow$ I am currently eating a soup 
\end{quote}
is usually nonsensical. However, this type of reversal is constantly applied 
in the second implication, even though there can obviously be many other 
reasons why salaries are not increased. Another fallacy is the linking of 
correlations with causations. A frequently cited example of this is the 
correlation between birth rate and the occurrence of storks in certain 
residential areas. It is obviously nonsensical to conclude the causality that 
babies are brought by storks. The truth\index{Truth} is that there are other 
parameters that both promote or prevent, such as possibly the rural or urban 
location.  

A state that goes beyond errors and fallacies and lasts longer 
we call a mistake. This is a constellation in the 
thought world of a person, in which the truth\index{Truth} of certain 
statements is believed because the situation at the time of the mistake 
suggests apparent evidence.

Mistakes can arise in various ways. A special case are 
faulty or contradictory theories, in which any statement 
can be inferred. Russell's paradox\index{Paradox}, which appeared in 
Frege's\index{Frege, Gottlob} work, can be considered the fundamental 
mistake of naive set theory\index{Set theory}. Only the 
Zermelo-Fraenkel axioms\index{Zermelo-Fraenkel axioms} and Russell's 
type theory\index{Type theory} have shown a way out. In physics, 
Albert Einstein\index{Einstein, Albert} cleared up some of the mistakes of 
physics in his works in the annus mirabilis 1905, including the 
existence\index{Existence} of the ether as the basic substance of the universe. 
Mistakes lead to new insights in this way.

A second and very threatening form of mistakes arises through 
misinformation, conspiracy theories, targeted manipulation of 
information or filter bubbles in communication media, which put people in the 
situation that they no longer question statements themselves. Resilience 
against such influences and the training of independent 
thinking are therefore indispensable. 

\section{The Truth Theory of Tarski\index{Tarski, Alfred} and 
Kripke\index{Kripke, Saul}} 

Alfred Tarski\index{Tarski, Alfred} wanted a formally correct theory of 
truth\index{Truth} which further develops the Aristotelian\index{Aristotle} 
correspondence theory\index{Correspondence theory of truth} and 
avoids paradoxes\index{Paradox}. In a famous work from 1933, he himself wrote:
\begin{quote} The present work is almost entirely devoted to a single problem, 
namely that of the definition of truth; its essence consists in the fact that 
one has to construct -- with regard to this or that language -- a 
factually accurate and formally correct definition of the term \enquote{true 
statement}. This problem, which is counted among the classical 
questions of philosophy, causes significant 
difficulties.\footnote{Machine translation of \cite{tarski1935}. It first 
appeared in 1933 in Polish language in \cite{tarski1933}.} 
\end{quote}
Tarski\index{Tarski, Alfred} succeeded in this in the context of languages, 
which he called poorer languages. These include in particular the formal 
languages\index{Formal language} of logic and mathematics. He wrote about this: 
\begin{quote} 
In the further course of the treatise, I will only consider those 
built up by scientific methods, which are known today, i.e., the formalised 
languages of the deductive sciences; their characterisation is given at the 
beginning of \S 2 ... With regard to the \enquote{poorer} languages, the 
problem of the definition of truth finds a positive solution.\footnote{Machine 
translation of \cite{tarski1935}.}
\end{quote}
His main result Tarski\index{Tarski, Alfred} formulated in this work 
as follows: 
\begin{quote}
For each formalised language of finite order, we can construct in the 
metalanguage\index{Metalanguage} a formally correct and factually accurate 
definition of the true statement, by using exclusively 
expressions of general-logical character, further expressions 
of the language itself, as well as terms from the field of morphology of the 
language, i.e., the names of the language expressions and the structural 
relations existing between them.\footnote{Machine translation of 
\cite{tarski1935}.}
\end{quote}
In richer languages, which include the natural languages, this 
construction cannot be applied directly. Tarski's\index{Tarski, Alfred} 
method can be interpreted in such a way that he constructed a 
semantics\index{Semantics} for certain poorer languages using a richer 
metalanguage\index{Metalanguage}, which describes the truth\index{Truth} 
of statements. More concretely, Tarski used for his 
theory of truth a formal language\index{Formal language} $L$ and a 
superordinate metalanguage\index{Metalanguage} $M$ which contains $L$. He 
used for his new theory a 
truth predicate\index{Truth predicate} $T$, with which the 
truth\index{Truth} of statements from $L$ within $M$ can be checked. The letter 
$T$ stands for \enquote{true} or for Tarski\index{Tarski, Alfred}.

Each truth predicate\index{Truth predicate} $T$ must fulfil a famous
adequacy condition of Tarski\index{Tarski, Alfred}, which 
is colloquially formulated as follows:
\begin{quote}
The sentence \enquote{Snow is white} is true if and only if snow 
is white.\footnote{Machine translation of \cite{tarski1935,tarski1969}.}
\end{quote}
This equivalence is called Tarski's\index{Tarski, Alfred} 
biconditional\index{Biconditional} in the literature. 
Biconditionals\index{Biconditional} of this kind are sometimes more abstractly 
expressed as 
\begin{quote}
The sentence \enquote{p} is true if and only if p
\end{quote}
or in the formal way 
\[
T(p) \xLeftrightarrow{~~~} \tilde p.
\]
Here, the left side is the application of the truth predicate\index{Truth 
predicate} $T$ to $p$, or more precisely to a standard name for $p$, such as a 
G\"odel number\index{Goed@G\"odel, Kurt} $\ulcorner p \urcorner$. The 
right-hand side consists of the translation $\tilde p$ of $p$ into the 
metalanguage\index{Metalanguage} $M$. The biconditional is therefore a 
syntactic\index{Syntax} statement within $M$.\footnote{See 
\cite{brendel2015,halbach1996,tarski1969}.}   

Tarski\index{Tarski, Alfred} demanded in his theory of truth\index{Truth} 
such syntactic\index{Syntax} biconditionals\index{Biconditional} as 
necessary conditions of adequacy. These appear at first 
glance like tautologies. Indeed, biconditionals\index{Biconditional} are also 
called quotation deletions and have inspired the redundancy theory of truth, 
which -- similar to Frege's comment on the prime number property of the number 
$5$ -- considers the whole concept of truth\index{Truth} to be dispensable. 

For the actual definition of truth\index{Truth}, Tarski\index{Tarski, Alfred} 
used a large repertoire. As the associated metalanguage\index{Metalanguage}, he 
mostly used set theory\index{Set theory} with the Zermelo-Fraenkel 
axioms\index{Zermelo-Fraenkel axioms}. The truth 
predicate\index{Truth predicate} $T$ associated with an object 
language\index{Object language} $L$ is in this case 
given by the satisfiability of formulas in set-theoretical\index{Set theory} 
models. The truth\index{Truth} of statements in $L$ is verified in the 
richer language $M$ by a proof. Another method of his is to reduce the 
truth\index{Truth} of non-atomic statements which are linked by logical signs 
like $\wedge$ or $\vee$, by complete induction\index{Complete induction} over 
the length of the concatenation to the truth\index{Truth} of the atomic 
statements contained therein. A difficulty arose with statements that 
contain quantifiers and free variables at the same time and 
are neither true nor false in themselves. One way to circumvent this 
problem is the method of quantifier elimination.\footnote{See 
\cite{tarski1931} and the overview article \cite{tarski1969}.} 

Tarski's theory\index{Tarski, Alfred} avoids self-referential 
paradoxes\index{Paradox}, but cannot be applied to natural languages 
in which such are possible. Kripke\index{Kripke, Saul} 
made an important step in 1975 to extend the theory of truth\index{Truth} to 
more general cases.\footnote{See \cite{kripke1975}.} For this 
he used infinite hierarchies of formal languages\index{Formal language} 
\[
L_1 \subset L_2 \subset \cdots,
\]
so that $L_i$ each has a truth predicate\index{Truth predicate} 
$T_i$ in $L_{i+1}$ and limited himself to the fact that the $T_i$ can only 
be partially defined. In this way, he managed to reach a fixed point 
in the limit of the languages $L_i$. 

\section{The Semantics\index{Semantics} of Arithmetic}

The semantics\index{Semantics} of formal languages\index{Formal language} and 
deductive systems\index{Deductive system} serves to give 
syntactic\index{Syntax} concepts a meaning in the Fregean\index{Frege, 
Gottlob} sense. Semantics\index{Semantics} is traditionally mostly defined as a 
set-theoretical interpretation\index{Interpretation} of a given deductive 
system\index{Deductive system}, so that all inference rules remain fulfilled in 
the language of sets. The constant symbols, the function and relation symbols 
of the formal language\index{Formal language} are assigned with sets and 
mappings. For this, a base set $A$ is used, so that the constants are elements 
of $A$, the $n$-ary functions $f$ are interpreted by set-theoretic mappings 
\[
f \colon A^n \xlongrightarrow{~~~~~~} A 
\]
and the $n$-ary relations are interpreted as subsets of $A^n$. The resulting 
structure of sets and mappings is called a model. A formula $\varphi$ is said 
to be true in this model\index{Model theory} if it is satisfiable under the 
interpretation\index{Interpretation}.\footnote{An atomic formula $\varphi$ in 
$n$ variables corresponds to a relation and is mapped by the 
interpretation\index{Interpretation} to a subset of $A^n$. Then 
$\varphi(t_1,\ldots,t_n)$ is true if the interpretation\index{Interpretation} 
of the vector $(t_1,\ldots,t_n)$ lies in this subset. For composite formulas, 
truth is defined by reducing to the atomic components. See \cite[Ch. II, Sect. 
2]{manin}.} Satisfiability is checked using the calculation rules of set 
theory\index{Set theory}, i.e., usually with the Zermelo-Fraenkel 
axioms\index{Zermelo-Fraenkel axioms} and possible additional chosen axioms. 
The model relation\index{Model theory}  
\[
\models \varphi
\]
expresses that the formula $\varphi$ is satisfied in every conceivable 
model\index{Model theory} and can therefore be called universally valid. 

An illustrative case is the complete Dedekind-Peano 
arithmetic\index{Dedekind-Peano arithmetic} with first-order predicate 
logic as the object language\index{Object language}. The truth\index{Truth} of 
arithmetic statements is defined by their satisfiability in the standard 
model\index{Model theory} of natural numbers. Andrzej 
Mostowski\index{Mostowski, 
Andrzej} was able to show, using the transfinite hyperarithmetical hierarchy, 
that a truth predicate\index{Truth predicate} $T$ in the second-order 
Dedekind-Peano arithmetic\index{Dedekind-Peano arithmetic} is 
syntactically\index{Syntax} formulable.\footnote{See \cite{mostowski1951}.} 

The $4$-squares theorem is an example of a concrete arithmetic statement of 
the form 
\[
\forall n \, \exists a \, \exists b \, \exists c \, \exists d \quad 
n=a^2+b^2+c^2+d^2.
\]
This statement is true\index{Truth}, if -- colloquially speaking -- for every 
$n$ a tuple $(a,b,c,d)$ can be specified that satisfies the theorem. In the 
set-theoretic model, the satisfiability can be verified by the surjectivity of 
the mapping 
\[
\mathbb{N}^4 \xlongrightarrow{~~~~~~} \mathbb{N}, \; (a,b,c,d) \mapsto 
a^2+b^2+c^2+d^2. 
\]
The infinity of the set 
\[
\mathbb{N}=\{0,1,2,3,4,\ldots\}
\]
plays a significant role in this semantics\index{Semantics}. But where does 
this infinite set come from? The model\index{Model theory} by John von 
Neumann\index{Neumann, John von} is 
\[
0=\emptyset, 1=\{\emptyset\}, 2=\{\emptyset,\{\emptyset\}\}, \ldots, 
n+1=\{0,1,\ldots,n\}, \ldots 
\]
An alternative is the definition 
\[
0=\emptyset, 1=\{\emptyset\}, 2=\{\{ \emptyset \}\}, \ldots, n+1=\{ n \}, 
\ldots 
\]
In both cases, $\mathbb{N}$ results in complex sets. None of these possible 
definitions is to be preferred over another or particularly useful from a 
structural point of view. 

Dedekind\index{Dedekind, Richard} was aware of this ambiguity problem and he 
pursued a different approach. He considered an arbitrary set $X$ with an 
initial element $\ast$ and an injective self-mapping 
\[
S \colon X \xlongrightarrow{~~~~~~} X, 
\]
such that $\ast$ is not in the image of $S$. This results in a set 
of natural numbers $\mathbb{N}$ within $X$ by 
\[
0=\ast, 1=S(\ast), 2=S(S(\ast)), \ldots  
\]
He then showed in a uniqueness theorem that all these avatars of the natural 
numbers are isomorphic\index{Isomorphism} to each other. For this, he used 
second-order predicate logic. In the first stage, this theorem does not indeed 
hold and there exist deviating non-standard models\index{Non-standard model}. 

In older literature, interpretations\index{Interpretation} 
are usually based on naive, material set theory\index{Set theory}. 
A somewhat more general and often more suitable approach is 
Henkin semantics\index{Henkin semantics}\index{Semantics} in the form of an 
interpretation\index{Interpretation} in category theory\index{Category theory}. 
To this day, set theory\index{Set theory} is the preferred basis for 
mathematical investigations and is considered the standard 
semantics\index{Semantics} of mathematics. It conceals intrinsic questions, as 
set theory\index{Set theory} is an axiomatic theory that itself 
has different models. The continuum hypothesis\index{Continuum hypothesis} is 
an example of how the truth\index{Truth} of a statement like 
\begin{quote}
There exists a cardinal number $\kappa$ with $\aleph_0 < \kappa < 2^{\aleph_0}$ 
\end{quote}
can depend on the chosen model\index{Model theory} and can only be decided by 
additional axioms as the forcing construction\index{Forcing} shows.

\section{Completeness\index{Completeness theorem} and 
Model Theory\index{Model theory}}

Leopold L\"owenheim\index{Loew@L\"owenheim, Leopold} and Thoralf 
Skolem\index{Skolem, Thoralf} founded a mathematical theory of 
semantics\index{Semantics} in the 20th century, laying the foundation for 
model theory\index{Model theory}\footnote{See \cite[Ch. X]{manin} on 
model theory.} and the works of Tarski\index{Tarski, Alfred} and 
G\"odel\index{Goed@G\"odel, Kurt}. The L\"owenheim-Skolem theorem is 
named after them and was first proven by 
L\"owenheim\index{Loew@L\"owenheim, Leopold} in 1915. Around 1920, 
Skolem\index{Skolem, Thoralf} provided a further, more general proof. 

The L\"owenheim-Skolem theorem states that the existence\index{Existence} 
of an infinite model\index{Model theory} of a finitely or countably axiomatised 
mathematical theory, formulated in first-order predicate logic, implies for 
every (infinite) cardinal number\index{Cardinal number} $\kappa$ a 
model\index{Model theory} of cardinality $\kappa$. This theorem, which is 
usually derived from the completeness theorem\index{Completeness theorem} as a 
corollary, shows that deductive systems\index{Deductive system} always have 
models\index{Model theory} of different cardinality\index{Cardinal number} and 
thus models\index{Model theory} in first-order logic are never unique. This is 
a 
remarkably noteworthy result. 

As it turned out, Skolem's\index{Skolem, Thoralf} proof method was strong 
enough to derive G\"odel's\index{Goed@G\"odel, Kurt} completeness theorem, 
which he proved in his dissertation in 1929. G\"odel\index{Goed@G\"odel, Kurt} 
referred in his dissertation to the proof by Skolem\index{Skolem, Thoralf}:
\begin{quote}
A similar method was used by Th. Skolem to prove the well-known theorem named 
after him and L\"owenheim.\footnote{Machine translation of \cite[Vol. I, 
1930]{goedel}.} 
\end{quote}
The completeness theorem\index{Completeness theorem} states that any statement 
formulated in a first-order axiomatic theory of predicate 
logic is provable if and only if it is universally 
valid, i.e., fulfilled in every set-theoretical model. In other words:
\[
\models \varphi \text{ if and only if  } \vdash \varphi. 
\]
It is clear that every provable statement is universally valid. Only the 
converse is the difficult part of this theorem. Modern proofs for the 
completeness theorem\index{Completeness theorem} usually go back to the 
dissertation of Leon Henkin\index{Henkin, Leon}, who modified and 
generalised G\"odel's\index{Goed@G\"odel, Kurt} proof. His proof also yielded a 
version of this theorem in type theory\index{Type theory}, i.e., in 
higher-level logic. The most important idea in the proof is to demonstrate the 
existence\index{Existence} of sufficiently many suitable models, because only 
then is universality a sufficiently strong condition.\footnote{See 
\cite{henkin1949,henkin1950, manin}. The proof by contradiction of the theorem 
(in classical logic) assumes that $\varphi$ is universally valid, but not 
provable. Then the addition of $\neg \varphi$ remains consistent. By adjoining 
free symbols for constants to the underlying formal language, a model is formed 
in which $\neg \varphi$ is valid. This is a contradiction, since $\varphi$ was 
valid. A topos-theoretical proof goes back to Pierre Deligne\index{Deligne, 
Pierre}. For further aspects see \cite{awodey2020} and \cite[Ch. 
IX]{maclane_moerdijk1992}.} 

The L\"owenheim-Skolem theorem implies several seemingly 
paradoxical\index{Paradox} statements, such as the existence\index{Existence} 
of countable non-standard models\index{Non-standard model} of the real 
numbers and countable or uncountable non-standard models\index{Non-standard 
model} of the natural numbers. A non-standard model\index{Non-standard model} 
of $\mathbb{N}$ can be constructed by a kind of adjunction. For this, we choose 
a new element $c$, which is assumed to be different from all natural numbers, 
and demand infinitely many conditions  
\[
c>n \text{ for each } n \in \mathbb{N}
\]
which are expressible in first-order predicate logic. The L\"owenheim-Skolem 
theorem, under this assumption, yields a new model\index{Model theory} with 
an infinitely large number $c$. Every model\index{Model theory} of $\mathbb{N}$ 
contains at least the numerals, i.e., the ordinary natural numbers of the form 
$n=S^n(0)$, which are created by $n$-fold application of the successor function 
$S$. Only the standard model\index{Model theory} of the numerals is unique in 
all models\index{Model theory} up to isomorphism\index{Isomorphism}.

Such statements and constructions do not contradict the uniqueness theorems for 
the natural and real numbers found in textbooks, as these are usually 
formulated in second-order predicate logic. For 
example, Dedekind\index{Dedekind, Richard} showed the uniqueness of the natural 
numbers up to isomorphism\index{Isomorphism} using his recursion 
theorem\index{Recursion}. The L\"owenheim-Skolem theorem even implies 
countable models\index{Model theory} of the real numbers. The uncountability of 
the real numbers in such countable models\index{Model theory} is a -- in a 
certain sense -- correct statement, which is called Skolem's 
paradox\index{Paradox}.\footnote{Skolem\index{Skolem, 
Thoralf} rejected the existence of uncountably infinite sets.}

These two theorems form the beginning of the field of 
model theory\index{Model theory} within mathematics. It investigates 
which set-theoretical interpretations\index{Interpretation} exist for 
axiomatic mathematical structures and which statements about them follow from 
knowledge of the logical foundations. Among the most important results 
of this field, in addition to the mentioned theorems, is the 
compactness theorem, which was found by Kurt G\"odel\index{Goed@G\"odel, Kurt} 
in his proof of the completeness theorem\index{Completeness 
theorem}.\footnote{The compactness theorem states that every (potentially 
infinite) set of formulas has a model if and only if every finite 
subset has models.} 

\section{The Undefinability Theorem of Truth\index{Truth}}

In the formal language\index{Formal language} of Dedekind-Peano 
arithmetic\index{Dedekind-Peano arithmetic} with first-order predicate logic, 
arithmetical sentences $\varphi$ are in principle characterisable as true in 
the standard model\index{Model theory} by considering the truth 
predicate\index{Truth predicate} $T$ as a function of the G\"odel 
number\index{Goed@G\"odel, Kurt} of $\varphi$. A theorem by 
Tarski\index{Tarski, Alfred} states that the subset of natural numbers which 
consists of the G\"odel numbers\index{Goed@G\"odel, Kurt} of such true 
arithmetical sentences is not an arithmetical set. This implies that the 
Tarskian\index{Tarski, Alfred} truth predicate\index{Truth predicate} $T$ is 
not definable by a formula in the underlying formal language\index{Formal 
language} of first order. For this, a second- or higher-order arithmetic is 
needed, for which there are further assertions, the truth\index{Truth} of which 
is not definable. In particular, the set of arithmetical sentences true in the 
standard model\index{Model theory} in the union of all levels is not 
decidable\index{Decidability}. It is not even recursively 
enumerable\index{Recursion}, due to the undecidability\index{Undecidability} of 
Hilbert's 10th problem\index{Hilbert problems}, because diophantine sets and 
thus recursively enumerable\index{Recursion} sets are always arithmetical. 
Tarski's\index{Tarski, Alfred} theorem was known to G\"odel\index{Goed@G\"odel, 
Kurt} and is often referred to as the undefinability theorem of 
truth\index{Truth}. It should be noted that this result does not contradict 
Tarski's\index{Tarski, Alfred} definition of truth\index{Truth}. But it shows 
that strictly richer languages\index{Formal language} are needed to define a 
truth predicate\index{Truth predicate}. To prove this theorem, assume that 
the concept of truth\index{Truth} is given by a class sign 
\[
\text{True}(x)
\]
and apply G\"odel's diagonal lemma\index{Goed@G\"odel, Kurt} 
as in the proof of the first incompleteness theorem\index{Incompleteness 
theorem} to $\neg \text{True}(x)$ or colloquially
\begin{quote}
I am not true.
\end{quote}
The fixed point then results in a contradiction.\footnote{See 
\cite[p. 256]{manin} and \cite{smullyan2013}.} From the correct mathematical 
proof of Tarski's\index{Tarski, Alfred} theorem, it follows that the set of 
G\"odel numbers\index{Goed@G\"odel, Kurt} of provable arithmetical 
sentences does not coincide with the set of G\"odel 
numbers\index{Goed@G\"odel, Kurt} of sentences true in the standard model, 
from which a variant of G\"odel's first incompleteness 
theorem\index{Incompleteness theorem} follows.\footnote{See 
\cite[Ch. II, Sec. 11]{manin}.}

G\"odel\index{Goed@G\"odel, Kurt} first thought about a 
semantic\index{Semantics} proof after his engagement with the completeness 
theorem\index{Completeness theorem} and exchanged ideas about it with 
Tarski\index{Tarski, Alfred}.\footnote{See \cite{vonplato2020}.} In such 
semantic\index{Semantics} versions of the incompleteness 
theorem\index{Incompleteness theorem}, there is a G\"odel 
sentence\index{Goed@G\"odel, Kurt} $Q$ which is true in the standard 
model\index{Model theory} of natural numbers, but in the object 
language\index{Object language} is neither provable nor refutable. In the 
non-standard models\index{Non-standard model}, such theorems are generally not 
true, because otherwise the completeness theorem\index{Completeness theorem} 
would imply provability\index{Provability}. The concept of truth\index{Truth} 
thus corresponds to provability\index{Provability} in the respective 
model\index{Model theory}.

It should be noted that the first incompleteness 
theorem\index{Incompleteness theorem} is not the opposite of the completeness 
theorem\index{Completeness theorem}. The naming of these two theorems is 
confusing. The sentence $Q$ can be replaced by a specific theorem. A later 
found example for this is Goodstein's theorem, which holds in the standard 
model\index{Model theory} of natural numbers and can be proven by transfinite 
induction\index{Transfinite induction}.\footnote{See 
\cite[Ch. 4]{hoffmann2018}.} 

\section{Truth Theory\index{Truth} of Non-classical Logic\index{Intuitionism}}

There is also a concept of semantics\index{Semantics} for non-classical 
logics. In intuitionistic logic\index{Intuitionism}, the rules of the deductive 
system\index{Deductive system} differ from classical logic, so that the 
concept of truth\index{Truth} is closer to the concept of 
computability\index{Computability}. In models, the satisfiability with 
intuitionistic logic\index{Intuitionism} is to be demonstrated. Since the time 
of Brouwer\index{Brouwer, Luitzen Egbertus Jan}, there 
has been a close connection between intuitionism\index{Intuitionism} and 
topology\index{Topology}. Tarski\index{Tarski, Alfred} has constructed an 
interpretation\index{Interpretation} of intuitionistic 
logic\index{Intuitionism} which uses the 
Heyting algebra\index{Heyting algebra} $\mathbf{Off}(X)$ of open sets in 
a topological space\index{Topology} $X$. In it, the negation operator 
\[
\neg U = (X \setminus U)^\circ 
\]
applied to an open set $U$ is given by the open core $(X \setminus U)^\circ$ of 
the complement of $U$ in $X$. With this definition it is generally not the case 
that $\neg \neg U=U$, unless $U$ coincides with the open core of its closure in 
$X$, which is almost never the case. Thus, $\mathbf{Off}(X)$ is generally not a 
Boolean algebra\index{Boolean algebra}.\footnote{See 
\cite[Ch. VI]{maclane_moerdijk1992} and the original work \cite{tarski1938}.}

Even more general than intuitionistic logic\index{Intuitionism} is modal logic. 
This is very old and is connected with the theory of possible 
worlds\index{Theory of possible worlds} by Leibniz\index{Leibniz, Gottfried 
Wilhelm}. Possible worlds\index{Theory of possible worlds} are 
semantic\index{Semantics} realisations in the form of different 
models. Modal logic was in parts already before
Leibniz\index{Leibniz, Gottfried Wilhelm} is known, for example, in 
Jean Buridan\index{Buridan, Jean}, a medieval philosopher. It 
contains, in addition to the two usual truth values\index{Truth} true and 
false, possible and necessary statements:
\begin{align*}
\Diamond p \colon &  \text{ It is possible that } p \\
\Box p \colon &  \text{ It is necessary that } p.
\end{align*}
This allows for statements that are called contingent because they are 
possible, but not necessary. 

Ruth Barcan\index{Barcan, Ruth} laid the foundations of modal predicate 
logic in her 1946 dissertation, building on the 
modal logical calculus of Clarence Irving Lewis\index{Lewis, 
Clarence Irving}. Before Lewis\index{Lewis, Clarence Irving}, 
Hugh MacColl\index{MacColl, Hugh} had already established a modal logical 
system around 1900.\footnote{See \cite{barcan1946,maccoll1906}.} 
These systems consist of the calculus of classical logic and the modus 
ponens as a conclusion rule from the formula
\[
\Box(A \xRightarrow{~~} B) \xRightarrow{~~} (\Box A 
\xRightarrow{~~} \Box B)  
\]
and the additional conclusion rule 
\[
\frac{A}{\Box A}.
\]
All corresponding formulas for $\Diamond$ follow from the relationship
\[
\Diamond A= \neg \Box \neg A.
\]
Barcan\index{Barcan, Ruth} discovered formulas, which are now known as Barcan 
formulas 
\[
\Diamond \exists x Fx \xRightarrow{~~}\exists x \Diamond Fx
\]
and 
\[
\Box \forall x Fx \xRightarrow{~~} \forall x \Box Fx.  
\]
The second formula has gained slightly more acceptance than the first, 
as the first formula implies an existence statement\index{Existence} in the 
real world for all objects to which the formula can be applied. 
In addition, there is the Buridan formula\index{Buridan, Jean} 
\[
\Diamond \forall x Fx \xRightarrow{~~} \forall x \Diamond Fx 
\]
and its reversal 
\[
\forall x \Diamond Fx \xRightarrow{~~}  \Diamond \forall x Fx  
\]
which together imply the interchangeability of $\Diamond$ and $\forall$. 

Saul A. Kripke\index{Kripke, Saul} introduced semantics\index{Semantics} 
for non-classical logics\index{Intuitionism} 
that generalise classical semantics\index{Semantics}. This 
Kripke semantics\index{Kripke-Joyal semantics} generalises the  
Tarskian\index{Tarski, Alfred} semantics\index{Semantics}. At the age of $18$, 
Kripke\index{Kripke, Saul} proved a completeness theorem\index{Completeness 
theorem} for modal logic. Before him, Kurt 
G\"odel\index{Goed@G\"odel, Kurt}, Alfred Tarski\index{Tarski, Alfred}, Evert 
Willem Beth\index{Beth, Evert Willem} and Andrzej 
Grzegorczyk\index{Grzegorczyk, Andrzej} had studied 
the semantics\index{Semantics} of non-classical logic. These 
semantic\index{Semantics} interpretations\index{Interpretation} were 
later developed further by Andr\'e Joyal\index{Joyal, Andr\'e} and others, 
leading to the concept of a categorical semantics\index{Semantics} which 
connects type theory\index{Type theory} (also called higher-order logic) with 
category theory\index{Category theory} and is called 
Kripke-Joyal semantics\index{Kripke-Joyal semantics}.\footnote{See 
\cite{beth1956,grzegorczyk1964,kripke1965,tarski1938} on intuitionistic logic, 
\cite{kripke1959} on modal logic and \cite[Vol. I, p. 299]{goedel}.}  

Around Kripke\index{Kripke, Saul} and Barcan\index{Barcan, Ruth} 
there was a long-lasting and still not finally clarified
priority dispute over the introduction of some terms into philosophy, 
particularly about the ideas underlying the theory of names as treated in the 
book \enquote{Naming and necessity} by Kripke\index{Kripke, Saul}.\footnote{See 
\cite{kripke1980}. Many people were involved in this dispute.}

\section{Categorical\index{Category theory} Semantics\index{Semantics} of 
Type Theory\index{Type theory}}

So far we have only studied the set-theoretical semantics\index{Semantics} 
of axiomatic theories which are based on first (or second) order predicate 
logic. However, in the current book the interpretation\index{Interpretation} 
of type theory\index{Type theory} in category theory\index{Category theory} is 
the new paradigm which we want to investigate more closely. This approach is 
in some sense based on the work of Leon Henkin\index{Henkin, Leon}. He has 
proven a completeness theorem\index{Completeness theorem} for type 
theory\index{Type theory}, i.e., for higher-level logic.\footnote{See 
\cite{henkin1950}. See \cite[Part II]{lambekscott1986} for a different proof. 
In 
these generalisations, the compactness theorem and the L\"owenheim-Skolem 
theorem are abandoned in their original formulation.} In his proof, types are 
interpreted\index{Interpretation} using objects and morphisms in 
$\mathbf{Set}$, i.e., with sets and mappings, where not all objects and 
morphisms in $\mathbf{Set}$ appear as images. This results in 
subcategories\index{Category theory} of $\mathbf{Set}$. This was mistakenly 
sometimes referred to as a disadvantage in parts of the literature and for this 
reason first-order logic was preferred. But Henkin semantics\index{Henkin 
semantics}\index{Semantics}\index{Category theory} are in principle more 
suitable for considering mathematics structurally\index{Structuralism}. 
In the following, we want to get rid of Henkin's\index{Henkin, Leon} 
set-theoretic background and pass to a true categorical\index{Category theory} 
language. Such interpretations\index{Interpretation} have the additional 
advantage that they abstract away many technical aspects of type 
theories\index{Type theory}.

Semantic\index{Semantics} categorical\index{Category theory} models\index{Model 
theory} are given by interpretations\index{Interpretation} of a type 
theory\index{Type theory} $T$, i.e., functors 
\[
\mathcal{S}_T \xlongrightarrow{~~~~~~} \mathcal{C}
\]
from the syntactic category\index{Category theory} $\mathcal{S}_T$ to a 
suitable category $\mathcal{C}$. Voevodsky\index{Voevodsky, Vladimir} studied a 
model\index{Model theory} of homotopy type theory\index{Homotopy type theory} 
in the category\index{Category theory} $\mathbf{sSet}$ of simplicial 
sets\index{Simplicial set} to verify the univalence 
axiom\index{Univalence}.\footnote{See \cite{kapulkin2021}.} More generally, 
homotopical infinity topoi\index{Elementary topos}\index{Topos}, i.e. 
simplicial localisations of model categories\index{Model category}, path 
categories\index{Path category} or \enquote{tribes}, provide potential 
semantical\index{Semantics} categories for dependent homotopy type 
theory\index{Homotopy type theory} with the univalence\index{Univalence} axiom. 
They allow techniques of homotopy theory\index{Homotopy theory} from 
$\mathbf{Top}$ and $\mathbf{sSet}$ to be abstracted. We have already discussed 
the technical problems which come with this approach.

\medskip The assignment of the syntactic categories\index{Category theory} 
$\mathcal{S}_T$ to a type theory\index{Type theory} $T$ has a reversal. Each 
higher elementary topos\index{Elementary topos}\index{Topos} is based on a 
syntactic\index{Syntax} type theory\index{Type theory} called 
Mitchell-B\'enabou language\index{Mitchell-B\'enabou language}. This 
categorical\index{Category theory} logic is equivalent to a higher-level logic, 
which reflects the internal logic of arguing with the arrows in all diagrams of 
the category\index{Category theory} because of the Curry-Howard 
correspondence\index{Curry-Howard correspondence}.\footnote{See 
\cite[Ch. VI]{maclane_moerdijk1992}.} 

Surprisingly, the Mitchell-B\'enabou language\index{Mitchell-B\'enabou 
language} does not necessarily satisfy the law of the excluded middle\index{Law 
of excluded middle}. Since a typical example of an elementary 
topos\index{Elementary topos}\index{Topos} is the category\index{Category 
theory} $\mathbf{Sh}(X)$ of sheaves\index{Category theory} on a topological 
space\index{Topology} $X$, this can be well read from the Heyting 
algebra\index{Heyting algebra} $\mathbf{Off}(X)$ of the open 
sets\index{Category theory} in $X$, which we have already seen is usually not a 
Boolean algebra\index{Boolean algebra}. The category\index{Category theory} 
$\mathbf{Set}$, on the other hand, has classical logic with the law of the 
excluded middle as its internal logic.

The Mitchell-B\'enabou language\index{Mitchell-B\'enabou language} has a 
natural semantics\index{Semantics} in the underlying elementary infinity  
topos\index{Elementary topos}\index{Topos}. However, the validity of logical 
formulas there is sometimes difficult to grasp. The Kripke-Joyal 
semantics\index{Kripke-Joyal semantics} has proven to be a useful tool, which 
goes back to the works of Beth\index{Beth, Evert Willem}, 
Grzegorczyk\index{Grzegorczyk, Andrzej}, Kripke\index{Kripke, Saul} and 
Joyal\index{Joyal, Andr\'e} on the semantics\index{Semantics} of non-classical 
logic. 

The technique of forcing\index{Forcing} in set theory\index{Set theory} can be 
explained with the help of Kripke-Joyal semantics\index{Kripke-Joyal semantics} 
in the context of suitable sheaves\index{Sheaf}, as William 
Lawvere\index{Lawvere, William} and Myles Tierney\index{Tierney, Myles} had 
recognised from about 1970 onwards. The original semantics\index{Semantics} of 
Kripke\index{Kripke, Saul} similarly emerges as a special case by considering 
sheaves\index{Sheaf} on a suitable partially ordered set 
$\mathbf{P}$.\footnote{See \cite{kripke1965,tierney1972} and \cite[Ch. 
VI]{maclane_moerdijk1992}.} 

\section{The Symbolic Construction of Reality}

We initially discussed the correspondence\index{Correspondence theory of truth} 
and coherence theories\index{Coherence theory of truth} of truth\index{Truth}, 
which both seemed quite far apart to us. In the further course, we approached 
the coherence theory\index{Coherence theory of truth} in a 
syntactic\index{Syntax} way through type theory\index{Type theory}. Alfred 
Tarski\index{Tarski, Alfred}, on the other hand, has tried by his own account 
to further develop the correspondence theory\index{Correspondence theory of 
truth} of truth\index{Truth} with the help of suitable 
semantics\index{Semantics}, which we found for type theory\index{Type theory} 
in categorical semantics\index{Semantics}. In the original formulations of the 
correspondence theory\index{Correspondence theory of truth}, however, a 
relationship with reality was demanded. How does this fit together? 

Let's pause for a moment. The original correspondence 
theory\index{Correspondence theory of truth} is based on the assumption of a 
reality that raises many questions. Is there a reality at all, or is this a 
naive notion? What kind of existence\index{Existence} are we talking about when 
it comes to nature and physical reality? How much of it can we perceive as 
humans? Is there a reality beyond physical reality, or is it always reducible 
to physical foundations, as reductionism\index{Reductionism} postulates?

Physical objects have very abstract properties that are only sufficiently 
comprehensible with mathematical methods. Hermann Weyl\index{Weyl, Hermann} has 
referred to the mathematical approach to understanding physics and other 
natural sciences as a symbolic construction of reality. In the introduction to 
an essay on this topic, he wrote in vivid and clear words, referring to 
Democritus\index{Democritus}:
\begin{quote}
It would probably not be a bad choice to date constructive natural sciences and 
critical philosophy from the day Democritus proclaimed: \enquote{Sweet and 
bitter, cold and warm as well as the colours, all this exists only in opinion, 
not in reality 
($\nu \acute{o} \mu \omega, o \grave{\nu} \; \varphi \acute{\nu} \sigma 
\epsilon \iota$); what really exists are unchangeable substance particles and 
their movement in empty space.} Indeed, without questioning the standpoint of 
naive realism, there is no philosophy, and a theoretically constructive science 
of nature is impossible as long as one takes the phenomena as they are given 
with perception at face value.\footnote{Machine translation of 
\cite[Vol. IV, p. 289]{weyl1968}.}
\end{quote}
In the rest of the text, Weyl\index{Weyl, Hermann} explains how science is 
shaped by the constructing spirit of humans and how mathematics permeates the 
concepts of physics such as space and time. A well-known example of this is 
the modelling of space-time by manifolds. We have questioned this description 
-- 
following Riemann's\index{Riemann, Bernhard} habilitation thesis -- as it is 
not clear whether space-time can be modelled continuously on a small scale or 
whether a discrete mathematical model is more appropriate. 

Weyl\index{Weyl, Hermann} can be considered as the discoverer of gauge theories 
in physics. He translated the idea of covariance\index{Covariance} into 
the differential geometric concept of gauge theories. These are given 
by a Lie group $G$, so that the algebra
$\mathcal{C}^\infty(M,G)$ of $G$-valued differentiable functions on 
a space-time manifold $M$ operate on the physical gauge fields and the 
Lagrangian function $\mathcal{L}$ of the theory. Gauge fields are 
$\mathfrak{g}$-valued $1$-forms $A$ on $M$, where $\mathfrak{g}$ is 
the Lie algebra of $G$.\footnote{The setup includes the group $G$, a 
principal fibre bundle $P$ on $M$ with group $G$ and the associated 
vector bundle $E$ under the adjoint representation. The gauge fields $A$ are 
connections on $E$ and the Lagrangian function $\mathcal{L}$ includes the 
curvature $F$ of the respective connection $A$.} An element $g \in 
\mathcal{C}^\infty(M,G)$ operates on the gauge fields by 
\[
A \mapsto g^{-1} A g+ g^{-1} dg. 
\]
In his first approach around 1919, Weyl\index{Weyl, Hermann} dealt with a 
synthesis of electrodynamics and general relativity. Electrodynamics is a gauge 
theory with the gauge group $G=U(1)$, i.e., the complex numbers of absolute 
value $1$ and the Lie group $\mathfrak{g}=i \mathbb{R}$. An elegant and simple 
formulation for this are the Maxwell\index{Maxwell, James Clerk} equations in 
the compact notation
\begin{align*}
d F &= 0 \\
\ast d \ast F & = \mu_0 \, j.
\end{align*}
Here, the field strength $F=dA$ is an exact $2$-form, the magnetic vector 
potential $A$ is a $1$-form, the current density $j$ is a $1$-form, 
$\mu_0$ the permeability constant, and $\ast$ is the Hodge $\ast$-operator. The 
gauge transformations are given by 
\[
A \mapsto A + d \varphi, \quad \varphi \text{  a scalar function}. 
\]
In a work from 1929 titled \enquote{Electron and gravity},\footnote{See 
\cite{weyl1968}.} Weyl\index{Weyl, Hermann} described the formalism of gauge 
theory even more precisely in a quantum mechanical context, by introducing a 
complex phase factor $\exp(i\varphi)$ that left the physics invariant. The 
gauge group $G$ in this situation was also $U(1)$. 

In more general cases, the gauge group $G$ is not abelian. Such theories 
are called Yang-Mills theories, after a work by Chen N. 
Yang\index{Yang, Chen} and Robert L. Mills\index{Mills, Robert} from 
1954.\footnote{See \cite{yangmills1954}.} The standard model of 
particle physics contains two gauge theories, the electroweak interaction 
with the gauge group $SU(2) \times U(1)$ and the strong interaction with 
$G=SU(3)$. General relativity, i.e., gravitation, can 
be understood as a gauge theory with gauge group 
$G=SO(3,1)$ according to Ryoyu Utiyama\index{Utiyama, Ryoyu}. 

The operation of the gauge group on the gauge fields is analogous to 
homotopies\index{Homotopy theory} of paths in topology. Therefore, it is 
useful to consider the gauge groupoid $\mathcal{G}(M)$ of all gauge fields. 
Similar to the fundamental groupoid\index{Fundamental groupoid} $\Pi_1$, it 
carries more information than just the quotient space of the orbits of the 
gauge fields under the gauge group, as it also considers the automorphisms of 
each individual gauge field. This view has only established itself in recent 
years.

In addition to his considerations about the symbolic construction of 
reality, Weyl\index{Weyl, Hermann} was interested in 
metaphysics\index{Metaphysics} and phenomenology.\footnote{The name 
metaphysics goes back to \cite{aristoteles}, phenomenology to
Edmund Husserl\index{Husserl, Edmund}.} He was particularly influenced by 
Kant\index{Kant, Immanuel}, Fichte\index{Fichte, Johann Gottlieb}, 
Leibniz\index{Leibniz, Gottfried Wilhelm} and Husserl\index{Husserl, Edmund}. 
In metaphysics\index{Metaphysics}, attempts are made to conceptually approach 
the dimensions and limits of human existence\index{Existence}. Such 
philosophical considerations, despite their large questions, which are 
extremely thought-provoking, also only provide access through language. 
Leibniz\index{Leibniz, Gottfried Wilhelm} 
had already realised this and believed that man's gain in knowledge 
would be limited to a symbolic calculus and that the actual 
understanding of the world would be reserved for God alone. 

Much suggests that for the description of reality in the natural sciences a 
symbolic calculus with suitable semantics\index{Semantics} is sufficient in the 
considered science, which makes a priori assumptions in the respective theory. 
Such considerations touch on the field of natural philosophy.\footnote{See the 
essays on natural philosophy in \cite[Ch. V]{schrenk}. Kurt 
G\"odel\index{Goed@G\"odel, Kurt} was of the opinion that physics, like 
mathematics, was synthetic a priori in its theory assumptions, see \cite[Vol. 
III, p. 360]{goedel}.} The deep, unresolved question of whether reality itself 
has a Platonic\index{Platonic idealism} quality beyond the symbolic 
construction remains open. For it could be that our 
world, which most people imagine to be made purely of matter, 
is simply made up of facts that we can perceive at most in our minds through 
sensory experience or by recognising physical laws. Such sceptical, 
antirealist\index{Antirealism} attitudes 
have a tradition dating back to Protagoras\index{Protagoras} in antiquity and 
also appear in the immaterialism\index{Immaterialism} of George 
Berkeley\index{Berkeley, George}.\footnote{See \cite{berkeley} and 
\cite[Ch. IV]{schrenk}.} They are the complete opposite of 
materialism\index{Materialism}, which only recognises the 
existence\index{Existence} of the material world. Ludwig 
Wittgenstein\index{Wittgenstein, Ludwig} coined the well-suited famous 
sentences in his \enquote{Tractatus logico-philosophicus}:
\begin{quote}
The world is everything that is the case. \\
The world is the totality of facts, not of 
things.\footnote{Machine translation of \cite{wittgenstein1922}. 
Wittgenstein\index{Wittgenstein, Ludwig} can be seen as a precursor of 
post-structuralism in his later works.}
\end{quote}

\section{Equivalence\index{Equivalence} and Truth Scepticism}

Based on our considerations about the nature of reality, there are 
at least two corresponding ways to consider the concept of truth\index{Truth}. 
The correspondence theory\index{Correspondence theory of truth} of 
truth\index{Truth} in its original version makes sense in the experimental 
natural sciences, where the existence\index{Existence} of a reality is usually 
not doubted and the results of physical and other experiments can at first 
glance be evaluated as a correspondence between theory and reality. On the 
other hand, there are sceptical, antirealist\index{Antirealism} views, which 
regard reality and our perception of it only as a variant of a 
semantics\index{Semantics}. In such cases, as well as in mathematics, a 
coherence theory\index{Coherence theory of truth} of truth\index{Truth} is more 
applicable. 

We must acknowledge that the concept of truth\index{Truth} 
is not seen equally in every scientific discipline and 
it may be impossible to have a uniform concept of truth\index{Truth} that 
represents a minimal consensus. This 
observation has given rise to newer philosophical attitudes that view the 
concept of truth\index{Truth} sceptically.\footnote{This includes truth 
pluralism, which traces back to Crispin Wright\index{Wright, Crispin} and 
Michael Lynch\index{Lynch, Michael}, as well as the truth relativism of John 
MacFarlane\index{MacFarlane, John} and others. The latter is also widespread in 
post-structuralist philosophies and includes the idea of deconstruction. For 
mathematical aspects, see \cite{tasic2012}.}

The diverse positions of truth relativism\index{Truth relativism} assume that 
truth\index{Truth} is always dependent on a context or an evaluation 
perspective and cannot be absolute. Some of these views have the advantage of 
avoiding self-referentiality and are suitable for solving the difficulties 
in finding truth\index{Truth} in natural languages and thus explaining 
phenomena such as faultless disagreements or defining the validity of future 
statements. They therefore require interesting further developments of the 
concept of semantics\index{Semantics}.

On the other hand, provably true statements are prototypes for absolute 
truths\index{Truth} and form an important prerequisite for science and for the 
success of our coexistence. Only with larger changes in science, such as 
Kuhn's\index{Kuhn, Thomas}\footnote{See \cite{kuhn1962}.} paradigm shifts, are 
truths\index{Truth} possibly corrected. For example, the laws of 
Newtonian\index{Newton, Isaac} physics were modified by the general theory of 
relativity. In this process, the description of physical space changed from 
flat Euclidean geometry to curved manifolds, in which Euclidean spaces are only 
contained infinitesimally in the form of tangent spaces.

An extreme position of truth relativism\index{Truth relativism} asserts that 
there can be no absolute truths\index{Truth} in science because there are 
fundamentally no ultimate justifications and therefore all forms of 
truth\index{Truth} are only relative. The possibility of Kuhnian\index{Kuhn, 
Thomas} paradigm shifts plays a rather controversial role here, as they 
seemingly describe a disruptive transition. However, Thomas Kuhn\index{Kuhn, 
Thomas} has defended himself against claims that his theory supports truth 
relativism\index{Truth relativism}.

How should such claims be evaluated? Vittorio H\"osle\index{Hoes@H\"osle, 
Vittorio} emphasised that this form of truth relativism\index{Truth relativism} 
is affected by an antinomy\index{Antinomy}, as the lack of an ultimate 
justification can be applied to truth relativism\index{Truth relativism} itself 
and thus leads to a performative self-contradiction.\footnote{See 
\cite{hoesle1990}. Theodor Adorno\index{Adorno, Theodor} and Karl 
Popper\index{Popper, Karl} also rejected truth relativism for other reasons.} 
This alone is already a significant counter-argument. Further criticism is 
directed against tendencies of arbitrariness, which result from the ethical 
consequences of such a theory due to insufficient societal consensus.

What about in science? Mathematical truths\index{Truth}, especially theorems, 
are based on certain assumptions, denoted by $\Gamma$, and assert statements 
$A$. They are usually in all calculi of the form of a judgement\index{Judgement}
\[
\Gamma \vdash A. 
\]
The assumptions in $\Gamma$ can contain axioms. Through a correct proof of such 
a theorem, the absolute truth\index{Truth} of the judgement\index{Judgement} 
$\Gamma \vdash A$ is shown or alternatively the relative truth\index{Truth} of 
the statement $A$ under the premise $\Gamma$. We have seen further absolute 
mathematical truths\index{Truth} in the discussion of the 
G\"odel-Lucas-Penrose\index{Penrose, Roger}\index{Lucas, 
John}\index{Goed@G\"odel, Kurt} argument. 

In summary, it can be said that in mathematics there are both relative 
and absolute truths\index{Truth}, which -- given the constant 
premise $\Gamma$ and in the same deductive system\index{Deductive system} 
-- apply under all conceivable circumstances and in all possible 
worlds\index{Theory of possible worlds}. In any case, there are absolute 
truths\index{Truth} in mathematics that are associated with the concept of 
provability\index{Provability}.

With theories of physics, it is not much different. Sometimes it is said, 
that truth\index{Truth} in the natural sciences is fundamentally based on 
measurements and agreement with nature. This ignores the 
hypothesis- and theory-building and falls short. Physical theories 
depend on underlying assumptions, like mathematical theories  
on their axioms. Kuhnian\index{Kuhn, Thomas} paradigm shifts change these 
and enable new insights, but derivations within the old theory 
remain valid in the classical limit. Similar to mathematics, 
there are therefore relative and absolute truths\index{Truth} in physics.

It is worth revisiting the example of the spherical shape of the Earth's 
surface. Using principles of equivalence\index{Equivalence} or invariance, -- 
similar to our discussion about the existence\index{Existence} of ideal circles 
-- neither an infinitely long descending chain of ultimate justifications nor a 
Kuhnian\index{Kuhn, Thomas} paradigm shift will ever change the approximate 
spherical shape of the Earth' surface. Such truths\index{Truth} are called 
stable or absolute. Already Peirce\index{Peirce, Charles} pointed out the 
necessity of correcting insights through constant doubt, falsifications and the 
reorientation of theories when finding absolute truths\index{Truth} in 
science.\footnote{See \cite{peirce} and the article on pragmatism in 
\cite{reclam}.}

In generalisation of this example, there is a convincing approach to 
refute truth relativism\index{Truth relativism} by classifying 
truths\index{Truth} with the help of the concept of 
equivalence\index{Equivalence} in a stable or invariant sense. 
The idea of such an antirelativistic criterion of objectivity for 
truth\index{Truth} goes back to Hermann Weyl\index{Weyl, Hermann}, who liked to 
use the concept of equivalence\index{Equivalence} and saw himself in this 
respect in the tradition of Felix Klein\index{Klein, Felix}.\footnote{See 
\cite{deppert1988}. Hermann Weyl\index{Weyl, Hermann} saw 
open problems in the application of his ideas to quantum mechanics, 
as the measurement process introduces discontinuities. Even from today's 
perspective, the understanding of the measurement process in quantum mechanics 
still appears incomplete. Attempts to explain it through non-linear theories 
with hidden variables were not successful.} 

How can Weyl's\index{Weyl, Hermann} idea be formulated in the language that we 
have developed in this book? The truth\index{Truth} of a statement $\varphi$ in 
an object language\index{Object language} $L$ is described by 
interpretations\index{Interpretation} $L \longrightarrow M$ in a suitable 
semantic\index{Semantics} metalanguage\index{Metalanguage} $M$ in which 
$\varphi$ can be proven. It has to be shown that -- as in the example of the 
spherical shape of the Earth's surface -- every conceivable deconstruction of 
truth\index{Truth} is based on a refined interpretation\index{Interpretation} 
$L' \longrightarrow M'$, so that the statement $\varphi$ corresponds to an 
equivalent statement $\varphi'$ in $L'$. The iterative process of 
deconstruction terminates with this strategy in appropriate 
equivalence\index{Equivalence} classes of statements, which reinforces the 
meaningfulness of the demand for an absoluteness of truth\index{Truth}. It is 
noteworthy that the concept of equivalence\index{Equivalence} itself in a 
certain way represents a deconstruction of the concept of 
identity\index{Identity}\index{Equality}, which does not contradict our 
argument against truth relativism\index{Truth relativism} though. 

\section{The Circle Closes}

Which questions have we answered? The concepts of 
equivalence\index{Equivalence} and univalence\index{Univalence} in 
dependent type theory\index{Type theory} explain the relation of 
equivalent\index{Equivalence} objects and are related to new forms of 
equality\index{Equality}\index{Identity}. Such concepts acquire a logical 
character, as equivalences\index{Equivalence} become exactly the 
transformations of abstract mathematical objects that preserve all relevant 
structures. This  is related to the non-uniqueness problem\index{Non-uniqueness 
problem} and therefore with the Platonic world of ideas\index{Platonic 
idealism} and is seen as the epitome of structural\index{Structuralism} 
thinking. This entire edifice of thought is connected with ideas from 
Carnap\index{Carnap, Rudolf}, Frege\index{Frege, Gottlob}, 
Grothendieck\index{Grothendieck, Alexander},  Leibniz\index{Leibniz, Gottfried 
Wilhelm}, and Tarski\index{Tarski, Alfred} on the invariance of logical 
truths\index{Truth}.\footnote{See \cite{carnap1928,awodey2017,awodey2018}. It 
is worthwhile to revisit Hume's principle in this context.} 

What did we learn about truth\index{Truth} and semantics\index{Semantics}?
The most popular semantics\index{Semantics} for dependent type 
theory\index{Type theory} is categorical semantics\index{Semantics}. But this 
is only one possibility, because type theory\index{Type theory}, category 
theory\index{Category theory} and set theory\index{Set theory} can be mutually 
interpreted in the form of object languages\index{Object language} or 
metalanguages\index{Metalanguage}.\footnote{See \cite{awodey2011}.} Through 
this perspective, the concepts of truth\index{Truth} and 
semantics\index{Semantics} in mathematics are somewhat demystified, as they 
mean nothing more and nothing less than provability\index{Provability} in 
richer deductive systems\index{Deductive system}. 

\begin{figure}[ht!]
\begin{onion}*{0.50in}
     \annulus*[gray!20]{3}{90}{210}[Type theory]
     \annulus*[gray!20]{3}{-30}{90}[Category theory]
     \annulus*[gray!20]{3}{210}{330}[Set theory]
     \annulus*[gray!30]{2}{90}{270}[syntax and deductive systems]     
     \annulus*[gray!30]{2}{-90}{90}[semantics and truth]
     \annulus*[gray!40]{1}{90}{210}[type]
     \annulus*[gray!40]{1}{-30}{90}[object]
     \annulus*[gray!40]{1}{210}{330}[set]
\end{onion}
\caption{\label{fig:circle_picture}A survey of the three foundations of 
mathematics.}\index{Syntax}\index{Deductive system}\index{Set 
theory}\index{Type theory}\index{Semantics}\index{Truth}
\end{figure} 

The circular image (see Fig.~\ref{fig:circle_picture}) shows the three 
foundations set theory\index{Set theory}, type theory\index{Type theory} and 
(higher) category theory\index{Category theory} and their essential properties 
in an overall view. On the outer ring they are depicted as equals. The second 
ring describes the relationship of the respective theory in relation to 
syntactic\index{Syntax} or semantic\index{Semantics} aspects, because this -- 
despite the possibility of mutual interpretations\index{Interpretation} -- 
plays a role. In the innermost ring, corresponding basic objects are found. In 
the table in Fig.~\ref{fig:comparison}, some relationships between the three 
foundations are further specified.\footnote{See \enquote{Relation between type 
theory and category theory} on ncatlab.org.}

\begin{figure}[ht!]
\resizebox{\textwidth}{!}{%
\begin{tabular}{|l|l|l|}\hline
{\bf Type theory}\index{Type theory} & 
{\bf Category theory}\index{Category theory} & 
{\bf Set theory, logic}\index{Set theory} \\
& & {\bf and topology}\index{Topology} \\\hline
type &  object &  set, space  \\\hline
term $a:A$ & $1 \xlongrightarrow{a} A$ & element $a \in A$ \\\hline 
universe\index{Universe} & object classifier\index{Object classifier}  & 
universe\index{Universe} \\\hline 
type $\mathbf{1}$ & terminal object & $\top$, true \\\hline
type $\mathbf{0}$ & initial object & $\bot$, false \\\hline
type of propositions & subobject classifier\index{Subobject classifier} 
$\Omega$ & truth value set \\\hline
power type & $\mathrm{Hom}(-,\Omega)$ & power set \\\hline 
product type & product & $\wedge$ \\\hline
sum type  & coproduct & $\vee$ \\\hline 
function type & exponential object\index{Exponential object} & 
$\Rightarrow$ \\\hline
$(A \to \mathbf{0})$ & $\mathrm{Hom}(A,0)$ & $\neg A$ \\\hline 
dependent type & morphism & family, fibration \\\hline
$\prod B(x)$ & sections of a fibration &  $\forall$ \\\hline
$\sum B(x)$ &  total space of a fibration &  $\exists$ \\\hline
identity type\index{Identity type} & path object & $=$, path space \\\hline
inductive type & colimit & induction, recursion\index{Recursion} \\\hline
coinductive type & limit & coinduction \\\hline 
equivalence\index{Equivalence} & isomorphism\index{Isomorphism}, higher & 
bijection, homotopy-\index{Homotopy theory}\\
& equivalence\index{Equivalence}  & equivalence\index{Equivalence}  \\\hline 
substitution & composition & cut rule \\\hline
program & generalised element & proof\index{Proof theory} \\\hline
Curry-Howard\index{Curry-Howard correspondence} & 
Mitchell-B\'enabou language\index{Mitchell-B\'enabou language} & 
intuitionistic\index{Intuitionism} logic \\
correspondence &  & of higher order (IHOL) \\\hline
Martin-L\"of type theory\index{Type theory} & locally cartesian closed 
& IHOL with intensional\index{Intensionality} \\
& $(\infty,1)$-category\index{Category theory} &  \\\hline
homotopy type theory\index{Homotopy type theory} & 
elementary $(\infty,1)$-topos\index{Elementary topos}\index{Topos} &  IHOL 
with intensional\index{Intensionality} \\ 
&& equality\index{Equality}\index{Identity} and 
univalence\index{Univalence}\\\hline 
\end{tabular}}
\caption{\label{fig:comparison}A detailed comparison of the three foundations 
of mathematics.}
\end{figure} 

What questions remain open? The 
consistency question\index{Consistency}\index{Contradiction-freeness} in 
mathematics is due to G\"odel's\index{Goed@G\"odel, Kurt} 
incompleteness theorem\index{Incompleteness theorem} in all three 
foundations in principle unsolved. Only through additional axioms of transfinite 
nature can the consistency\index{Consistency}\index{Contradiction-freeness} of 
smaller parts of mathematics be secured. This approach shifts 
the consistency problem\index{Consistency}\index{Contradiction-freeness} into 
other deductive systems\index{Deductive system} at the price of further 
incompleteness\index{Incompleteness theorem}. 

In dependent type theory\index{Type theory}, it is open which axioms besides  
univalence\index{Univalence} are sensible and necessary. What conditions must 
be placed on $(\infty,1)$-categories\index{Infinity category} to allow 
type-theoretical\index{Type theory} interpretations\index{Interpretation} is 
also still an open question. This raises the fundamental question of whether 
(higher) category theory\index{Category theory} and dependent type 
theory\index{Type theory} represent complete foundations of mathematics  
independently of set theory\index{Set theory}, as predicted by William 
Lawvere\index{Lawvere, William} and others.

Which consequences arise from the mathematics described in this book?  
There will most likely be a change of teaching, research and publication 
practices through digital assistance systems for the verification of proofs and 
algorithms which support us -- possibly using artificial intelligence 
techniques\index{Artificial intelligence} -- in the creative search for new 
mathematics and in achieving more transparency in the spirit of the open 
science idea. In a way, this approach is similar to that of the Bourbaki 
group\index{Bourbaki}, who tried to do this in the form of 
texts.\footnote{Similar open science projects are MathOverflow, nLab, 
n-Category Caf\'e, Stacks Project and Stack Exchange.} 
The possibilities that arise from such developments can be seen as 
a partial realisation of Leibniz's\index{Leibniz, Gottfried Wilhelm} dream 
of a lingua universalis\index{Lingua universalis}.

\printendnotes

\backmatter 
\renewcommand{\refname}{Bibliography}

\renewcommand{\indexname}{Name and Subject Index} 
\printindex


\begin{thebibliography}{Wittgenstein19xxx}
\bibitem[Ahrens2018]{ahrens2018} Benedikt Ahrens, Peter LeFanu Lumsdaine, 
Vladimir Voevodsky: Categorical structures for type theory in univalent 
foundations, Logical Methods in Computer Science, Vol. 14(3), 1-16 (2018). 
\bibitem[AnelCat2021]{anelcat2021} Mathieu Anel, Gabriel Catren: New spaces in 
mathematics, Cambridge University Press (2021).
\bibitem[vonAquin1259]{aquin} Thomas von Aquin: Quaestiones disputatae de 
veritate, Quaestio I, Philo\-sophische Bibliothek, Band 384, Felix Meiner 
Verlag (1986).
\bibitem[Aristoteles2009]{aristoteles} Aristoteles: Metaphysik, Philosophische 
Bibliothek, Band 308, Felix Meiner Verlag (2009).
\bibitem[Artemov2019]{artemov2019} Sergei Artemov: The provability of 
consistency, Preprint arXiv (2019).
\bibitem[ArtinMazur1969]{artinmazur1969} Michael Artin, Barry Mazur: \'Etale 
homotopy, Lecture Notes in Mathematics, Vol. 100 (1969).
\bibitem[AwodWarr2009]{awodeywarren} Steve Awodey, Michael A. Warren: 
Homotopy theoretic models of identity types, Mathematical Proceedings 
of the Cambridge Philosophical Society, Vol. 146, 45-55 (2009).
\bibitem[Awodey2011]{awodey2011} Steve Awodey: From sets to types, to 
categories, to sets, in: Foundational Theories of Classical and Constructive 
Mathematics (Ed. G. Sommaruga), Springer Verlag, 113-125 (2011).
\bibitem[Awodey2014]{awodey2014} Steve Awodey: Structuralism, invariance, and 
univalence, Philosophia Mathematica, Vol. 22, 1-11 (2014).
\bibitem[Awodey2017]{awodey2017} Steve Awodey: Carnap and the invariance of 
logical truth, Synthese, Vol. 194, 65-78 (2017). 
\bibitem[Awodey2018]{awodey2018} Steve Awodey: Univalence as a principle of 
logic, Indagationes Mathematicae, Vol. 29(6), 1497-1510 (2018).
\bibitem[Awodey2020]{awodey2020} Steve Awodey: Sheaf representations and 
duality in logic, Outstanding Contributions to Logic, Vol. 20, Springer Verlag, 
39-57 (2021).
\bibitem[Balaguer1998]{balaguer1998} Michael Balaguer: Platonism and 
anti-Platonism in mathematics, Oxford University Press (1998). 	
\bibitem[Barcan1946]{barcan1946} Ruth Barcan: A functional calculus of first 
order based on strict implication, The Journal of Symbolic Logic, Vol. 11(1), 
1-16 (1946).
\bibitem[Beeson1985]{beeson1985} Michael Beeson: Foundations of constructive 
mathematics, Ergebnisse der Mathematik und ihrer Grenzgebiete, 3. Folge, Vol. 
6, Springer Verlag (1985).
\bibitem[BellCook1992]{bellantonicook1992} Stephen Bellantoni, Stephen 
Cook: A new recursion theoretic characterization of the polytime functions, 
Computational Complexity, Vol. 2(2), 97-110 (1992).
\bibitem[B\'enabou1985]{benabou1985} Jean B\'enabou: Fibred categories and the 
foundations of naive category theory, Journal of Symbolic Logic, Vol. 50, 
10-37 (1985). 
\bibitem[Benacerraf1965]{benacerraf1965} Paul Benacerraf: What numbers could 
not be, The Philosophical Review, Vol. 74, 47-73 (1965).
\bibitem[Benacerraf1973]{benacerraf1973} Paul Benacerraf: Mathematical truth,  
Journal of Philosophy, Vol. 70, 661-680 (1973).
\bibitem[Ben-David2019]{ben-david2019} Shai Ben-David, Pavel Hrube\v{s}, Shay 
Moran, Amir Shpilka, Amir Yehudayoff: Learnability can be undecidable, Nature 
Machine Intelligence, Vol. 1, 44-48 (2019).
\bibitem[Benzm\"uller2023]{benzmueller} Christoph Benzm\"uller: A simplified 
variant of G\"odel's ontological argument, Preprint arXiv (2023). 
\bibitem[Berg2018]{vdB} Benno van den Berg: Path categories and propositional 
identity types, ACM Transactions on computational logic, Vol. 19(2), 1-32 
(2018).
\bibitem[BergMoer2018]{vdBM} Benno van den Berg, Ieke Moerdijk: Exact 
completion of path categories and algebraic set theory: Part I, Journal of Pure 
and Applied Algebra, Vol. 222(10), 3137-3181 (2018).
\bibitem[Berkeley1710]{berkeley} George Berkeley: A treatise concerning the 
principles of human knowledge, Aaron Rhames Printer (1710).
\bibitem[Beth1956]{beth1956} Evert Willem Beth: Semantic construction of 
intuitionistic logic, Koninklijke Nederlandse Akademie van Wetenschappen, 
Vol. 19(11), 357-388 (1956).
\bibitem[Blanke1996]{blanke1996} Detlev Blanke: Leibniz und die Lingua 
Universalis, Sitzungsberichte der Leibniz-Soziet\"at, Band 13, Heft 5, 27-35 
(1996).
\bibitem[Blom2011]{blom2011} Philipp Blom: B\"ose Philosophen: Ein Salon in 
Paris und das vergessene Erbe der Aufkl\"arung, Hanser Verlag (2011).
\bibitem[Boardman73]{boardman} Michael Boardman, Rainer Vogt: Homotopy 
invariant algebraic structures on topological spaces, Lecture Notes in 
Mathematics, Vol. 347, Springer Verlag (1973).
\bibitem[Boehm2024]{boehm} Omri Boehm, Daniel Kehlmann: Der bestirnte Himmel 
\"uber mir, Propyl\"aen Verlag (2024).
\bibitem[Boer2020]{boer} Menno de Boer: A proof and formalization of the 
initiality conjecture of dependent type theory, Licentiate Thesis Stockholm 
University (2020).
\bibitem[Bolzano1837]{bolzano1837} Bernard Bolzano: Wissenschaftslehre, 
Seidel Verlag (1837).
\bibitem[Bolzano1851]{bolzano1851} Bernard Bolzano: Paradoxien des 
Unendlichen, Reclam Verlag (1851).
\bibitem[Boole1847]{boole1847} George Boole: The mathematical analysis of 
logic: being an essay towards a calculus of deductive reasoning, 
Macmillan (1847).
\bibitem[Bordg2021]{bordg} Anthony Bordg: A replication crisis in 
mathematics?, The Mathematical Intelligencer, Vol. 43(4), 48-52 (2021).
\bibitem[Brendel2015]{brendel2015} Elke Brendel: Die Wahrheit \"uber den  
L\"ugner, De Gruyter Verlag (2015). 
\bibitem[Bridges1994]{bridges1994} Douglas Bridges: Computability, Graduate 
Texts in Mathematics, Vol. 146, Springer Verlag (1994).
\bibitem[BridgMin1984]{bridges1984} Douglas Bridges, Ray Mines: What is 
constructive mathematics?, The Mathematical Intelligencer, Vol. 6(4), 32-38 
(1984). 
\bibitem[Brown1973]{brown1973} Kenneth S. Brown: Abstract homotopy theory and 
generalized sheaf cohomology, Transaction of the American Mathematical Society, 
Vol. 186, 419-458 (1973).
\bibitem[Brunerie2016]{brunerie} Guillaume Brunerie: On the homotopy groups of 
spheres in homotopy type theory, Preprint arXiv (2016).
\bibitem[Cantor1932]{cantor} Georg Cantor: Gesammelte Werke, Springer Verlag 
(1932).
\bibitem[Carchedi2020]{carchedi2020} David J. Carchedi: Higher orbifolds and 
Deligne-Mumford stacks as structured infinity topoi, Memoirs of the American 
Mathematical Society, Vol. 264 (2020).
\bibitem[Carnap1928]{carnap1928} Rudolf Carnap: Der logische Aufbau der 
Welt, Weltkreis Verlag (1928). 
\bibitem[Carnap1934]{carnap1934} Rudolf Carnap: Logische Syntax der 
Sprache, Springer Verlag (1934).
\bibitem[Carnap1950]{carnap1950} Rudolf Carnap: Empiricism, semantics, and 
ontology, Revue Internationale de Philosophie, Vol. 4(11), 20-40 (1950). 
\bibitem[Cartmell1986]{cartmell1986} John Cartmell: Generalised 
algebraic theories and contextual categories, Annals of Pure and Applied 
Logic, Vol. 32, 209-243 (1986). 
\bibitem[Cassou2011]{cassou2011} Pierre Cassou-Nogu\`es: On 
G\"odel's \enquote{Platonism}, Philosophia Scientiae, Vol. 15(2), 137-171 
(2011).
\bibitem[CavNoeth1937]{cavailles} Jean Cavaill\`es, Emmy Noether (Hrsg.): 
Briefwechsel Cantor-Dedekind, Hermann \'Editeurs (1937).
\bibitem[Cheng2021]{cheng2021} Yong Cheng: Current research on G\"odel's
incompleteness theorem, The Bulletin of Symbolic Logic, Vol. 27(2), 113-167 
(2021).
\bibitem[Church1936]{church1936} Alonzo Church, An unsolvable problem of 
elementary number theory, American Journal of Mathematics, Vol. 58, 345-363 
(1936).
\bibitem[Church1940]{church1940} Alonzo Church: A formulation of the simple 
theory of types, The Journal of Symbolic Logic, Vol. 5, 56-68 (1940).
\bibitem[Cisinski2015]{cisinski2015} Denis-Charles Cisinski: Cat\'egories 
sup\'erieurs et th\'eorie des topos, S\'eminaire Bourbaki, 67\`eme ann\'ee 
(2014-2015), Expos\'e No. 1097, 1-57 (2015).
\bibitem[Cisinski2025]{cisinski2025} Denis-Charles Cisinski, Bastiaan Cnossen, 
Kim Nguyen, Tashi Walde: Formalization of higher categories, Preprint (2025).
\bibitem[Clarke2000]{clarke2000} Samuel Clarke, Gottfried Wilhelm Leibniz: 
Correspondence (Ed. R. Ariew), Hackett Publ. Co. (2000). 
\bibitem[Clarke2020]{clarke2020} Justin Clarke-Doane: Morality and mathematics, 
Oxford University Press (2020).
\bibitem[Cohen1963]{cohen1963} Paul J. Cohen: The independence of the continuum 
hypothesis, Part I, Proceedings of the National Academy of Sciences of the 
United States of America,  Vol. 50(6), 1143-1148 (1963).
\bibitem[Cohen1964]{cohen1964} Paul J. Cohen: The independence of the continuum 
hypothesis, Part II, Proceedings of the National Academy of Sciences of the 
United States of America, Vol. 50(6), 105-110 (1964).
\bibitem[Connes1995]{connes1995} Alain Connes: Non-commutative geometry, 
Academic Press (1995).
\bibitem[CoqHuet1986]{coquand1986} Thierry Coquand, G\'erard Huet: 
The calculus of constructions, INRIA rapport de recherche 530 (1986).
\bibitem[CoqPaulin1990]{coquand1990} Thierry Coquand, Christine 
Paulin-Mohring: Inductively defined types, Lecture Notes in Computer Science, 
Vol. 417, Springer Verlag, 50-66 (1990).
\bibitem[Coquand2014]{coquand2014} Thierry Coquand: Th\'eorie des types 
d\'ependants et axiome d'univalence, S\'eminaire Bourbaki, 66\`eme ann\'ee 
(2013-2014), Expos\'e No. 1085, 1-18 (2014).
\bibitem[Couturat1901]{couturat1901} Louis Couturat: La loqique de Leibniz, 
F\'elix Alcan \'Editeur (1901). 
\bibitem[CoutLeau1903]{couturat1903} Louis Couturat, L\'eopold Leau: 
Histoire de la langue universelle, Hachette \'Editeurs (1903). 
\bibitem[Davis1965]{davis1965} Martin Davis (Ed.): The undecidable, Raven 
Press (1964).
\bibitem[Davis1993]{davis1993} Martin Davis: How subtle is G\"odel's theorem? 
More on Roger Penrose, Behavioral and Brain Science, Vol. 16, 611-612 (1993).
\bibitem[Dedekind1930]{dedekind} Richard Dedekind: Gesammelte 
mathematische Werke, Vieweg Verlag (1930).
\bibitem[Deppert1988]{deppert1988} Wolfgang Deppert, Kurt H\"ubner, Arnold 
Oberschelp, Volker Weidemann (Hrsg.): Exakte Wissenschaften und ihre 
philosophische Grundlegung, Peter Lang Verlag, 445-467 (1988).
\bibitem[DeRisi2007]{derisi2007} Vincenzo De Risi: Geometry and monadology, 
Leibniz's analysis situs and philosophy of space, Birkh\"auser Verlag (2007).
\bibitem[Dershowitz1993]{dershowitz} Nachum Dershowitz: Trees, ordinals and
termination, Lecture Notes in Computer Science, Vol. 668, Springer Verlag, 
243-250 (1993).
\bibitem[Diophant]{diophant} Diophant: Arithmetica (about 250 AD). 
\bibitem[DwyKan1980]{dwyerkan} William Dwyer, Daniel Kan: Simplicial 
localizations of categories, Journal of Pure and Applied Algebra, Vol. 18(1), 
17-35 (1980).
\bibitem[Eco1993]{eco1993} Umberto Eco: La ricerca della lingua perfetta nella 
cultura europea, Editori Laterza (1993).
\bibitem[Edmonds2021]{edmonds2021} David Edmonds: Die Ermordung des 
Professor Schlick, C. H. Beck Verlag (2021).
\bibitem[Einstein1916]{einstein1916} Albert Einstein: Die Grundlage der 
allgemeinen Relativit\"atstheorie, Annalen der Physik, Band 49, 769-822 
(1916). 
\bibitem[Feferman1962]{feferman1962} Solomon Feferman: Transfinite recursive 
progressions of axiomatic theories, Journal of Symbolic Logic, Vol. 27(3), 
259-316 (1962).
\bibitem[Feferman1978]{feferman1978} Solomon Feferman: Constructive theories of
functions and classes, Logic Colloquium 1978, North Holland, 159-224 (1978).
\bibitem[Feferman1996]{feferman1996} Solomon Feferman: Penrose's G\"odelian 
argument, Psyche, Vol. 2(7), 21-32 (1996).
\bibitem[Feige2014]{feige} Daniel M. Feige: Philosophie des Jazz, Suhrkamp 
Verlag (2014). 
\bibitem[Fichte1794]{fichte1794} Johann Gottlieb Fichte: Grundlage der 
gesamten Wissenschaftslehre, Gabler Verlag (1794).
\bibitem[Field1980]{field1980} Hartry Field: Science without numbers, Oxford 
University Press (1980).
\bibitem[Field1989]{field1989} Hartry Field: Realism, mathematics, and 
modality, Blackwell Publ. (1989).
\bibitem[Frege1879]{frege1879} Gottlob Frege: Begriffsschrift, Louis Nebert 
Verlag (1879).
\bibitem[Frege1884]{frege1884} Gottlob Frege: Die Grundlagen der Arithmetik, 
Koebner Verlag (1884).
\bibitem[Frege1893]{frege1893} Gottlob Frege: Grundgesetze der Arithmetik, 
Pohle Verlag (1893).
\bibitem[Frege1962]{frege1962} Gottlob Frege: Funktion, Begriff, Bedeutung, 
Vandenhoeck und Ruprecht (1962).
\bibitem[Frege1976]{frege1976} Gottlob Frege: Wissenschaftlicher Briefwechsel, 
Felix Meiner Verlag (1976).
\bibitem[Frege1983]{frege1983} Gottlob Frege: Nachgelassene Schriften, 
2. Auflage, Felix Meiner Verlag (1983).
\bibitem[Frege2003]{frege2003} Gottlob Frege: Logische Untersuchungen, 5. 
Auflage, Vandenhoeck und Ruprecht Verlag (2003).
\bibitem[Friedlander1982]{friedlander1982} Eric Friedlander: \'Etale homotopy 
of simplicial schemes, Annals of Mathematics Studies, Vol. 104, Princeton 
University Press (1982).
\bibitem[Gentzen1935]{gentzen} Gerhard Gentzen: Untersuchungen über das 
logische Schlie{\ss}en I, Mathematische Zeitschrift, Band 39, 176-210 (1935).
\bibitem[G\"odel2003]{goedel} Kurt G\"odel: Collected works, Vols. I-V, Oxford 
University Press (1986-2003).
\bibitem[GoerssJard2009]{goerssjardine} Paul G. Goerss, John F. Jardine: 
Simplicial homotopy theory, Birkh\"auser Verlag (2009).   
\bibitem[Gonthier2005]{gonthier2005} Georges Gonthier: A computer-checked 
proof of the four-color theorem, Preprint (2005).
\bibitem[Gonthier2012]{gonthier2012} Georges Gonthier: A machine-checked proof 
of the odd-order theorem (2012). 
\bibitem[Goodfellow2016]{goodfellow2016} Ian Goodfellow, Yoshua Bengio, Aaron 
Courville: Deep learning, MIT Press (2016). 
\bibitem[Grayson2018]{grayson} Daniel Grayson: An introduction to univalent 
foundations for mathematicians, Bulletin of the American Mathematical Society, 
Vol. 55, 427-450 (2018).
\bibitem[Grothen1983]{grothendieck} Alexander Grothendieck: A la 
poursuite des champs, Preprint (1983).
\bibitem[Groth2015]{groth2015} Moritz Groth: A short course on 
$\infty$-categories, in: Handbook of homotopy theory, Chapman and Hall, 
549-617 (2020). 
\bibitem[Grundmann2008]{grundmann} Thomas Grundmann: Analytische Einf\"uhrung 
in die Erkenntnistheorie, De Gruyter Verlag (2008).
\bibitem[Grzeg1964]{grzegorczyk1964} Andrzej Grzegorczyk: A 
philosophically plausible interpretation of intuitionistic logic, Indagationes 
Mathematicae, Vol. 26, 596-601 (1964).
\bibitem[Halbach1996]{halbach1996} Volker Halbach: Axiomatische 
Wahrheitstheorien, Akademie Verlag (1996).
\bibitem[Hales2024]{hales2024} Tom Hales: The formal proof of the Kepler 
conjecture -- a critical retrospective, Preprint arXiv (2018). 
\bibitem[Hardy1940]{hardy} Godfrey Harold Hardy: A mathematician's apology, 
Cambridge University Press (1940).
\bibitem[Hawking2021]{hawking} Stephen Hawking, Roger Penrose: Was sind Raum 
und Zeit?, Klett-Cotta Verlag (2021).
\bibitem[Hebb1949]{hebb1949} Donald O. Hebb: The organization of behavior - a 
neuropsychological theory, Wiley (1949).
\bibitem[Heidegger1957]{heidegger} Martin Heidegger: Identit\"at und Differenz, 
Klett-Cotta Verlag (1957).
\bibitem[Henkin1949]{henkin1949} Leon A. Henkin: The completeness of the 
first-order functional calculus, The Journal of Symbolic Logic, Vol. 14(3), 
159-166 (1949).
\bibitem[Henkin1950]{henkin1950} Leon A. Henkin: Completeness in the theory of 
types, The Journal of Symbolic Logic, Vol. 15(2), 81-91 (1950). 
\bibitem[Hilbert1926]{hilbert1926} David Hilbert: \"Uber das Unendliche, 
Mathematische Annalen, Vol. 95, 161-190 (1926).
\bibitem[Hilbert1934]{hilbert1934} David Hilbert, Paul Bernays: Grundlagen der 
Mathematik I, Grundlehren der mathematischen Wissenschaften, Band 40, 
Springer Verlag (1934). 
\bibitem[Hilbert1939]{hilbert1939} David Hilbert, Paul Bernays: Grundlagen der 
Mathematik II, Grundlehren der mathematischen Wissenschaften, Band 41, Springer 
Verlag (1939). 
\bibitem[Hilbert2013]{ewaldsieg} William Ewald, Wilfried Sieg (Eds.): 
David Hilbert's lectures on the foundations of arithmetic and logic 1917-1933, 
Springer Verlag (2013).
\bibitem[Hobbes1651]{leviathan} Thomas Hobbes: Leviathan, Andrew Crooke (1651).
\bibitem[Hobbes1655]{hobbes1655} Thomas Hobbes: De corpore, Andrew Crooke 
(1655).
\bibitem[H\"osle1990]{hoesle1990} Vittorio H\"osle: Die Krise der Gegenwart und 
die Verantwortung der Philosophie, C. H. Beck Verlag (1990). 
\bibitem[Hoffmann2017]{hoffmann2017} Dirk W. Hoffmann:  Die G\"odelschen 
Unvollst\"andigkeitss\"atze, 2. Auflage, Springer Spektrum (2017).
\bibitem[Hoffmann2018]{hoffmann2018} Dirk W. Hoffmann: Grenzen der Mathematik, 
3. Auflage, Springer Spektrum (2018). 
\bibitem[Hofmann1994]{hofmann1994} Martin Hofmann: On the interpretation of 
type theory in locally cartesian closed categories, Lecture Notes in Computer 
Science, Vol. 933, 427-441 (1994).
\bibitem[HofmannStr1998]{hofmannstreicher} Martin Hofmann, Thomas Streicher: 
The groupoid interpretation of type theory, in: Twenty-five years of 
constructive type theory (Venice, 1995), Oxford Logic Guides, Vol. 36, 
83-111 (1998).
\bibitem[HuberMS2017]{hms2017} Annette Huber, Stefan 
M\"uller-Stach: Periods and Nori motives, Ergebnisse der Mathematik und ihrer 
Grenzgebiete, 3. Folge, Vol. 65, Springer Verlag (2017).
\bibitem[Hume1739]{hume} David Hume: A treatise of human nature, John Noon 
(1739). 
\bibitem[Illusie2004]{illusie2004} Luc Illusie: What is a topos?, Notices of 
the American Mathematical Society, Vol. 51(9), 1060-1061 (2004).
\bibitem[JorNim2019]{reclam} Stefan Jordan, Christian Nimtz: Grundbegriffe der 
Philosophie, 2. Auflage, Reclam Verlag (2019). 
\bibitem[Jost2017]{jost2017} J\"urgen Jost: Object oriented models vs. data 
analysis -- is this the right alternative?, Boston studies in the philosophy 
and history of science, Vol. 327, Springer Verlag, 253-286 (2017). 
\bibitem[Jost2019]{jost2019} J\"urgen Jost: Leibniz und die moderne 
Naturwissenschaft, Springer Verlag (2019).
\bibitem[Joyal2002]{joyal2002} Andr\'e Joyal: Quasicategories and Kan 
complexes, Journal of Pure and Applied Algebra, Vol. 175(13), 207-222 (2002).
\bibitem[Joyal2017]{joyal2017} Andr\'e Joyal: Notes on clans and tribes, 
Preprint arXiv (2017).
\bibitem[Kant1781]{kant} Immanuel Kant: Kritik der reinen Vernunft, 
Johann Friedrich Hartknoch Verlag (1781).  
\bibitem[KaprVoev1991]{KV91} Mikhail Kapranov, Vladimir Voevodsky: 
$\infty$-groupoids and homotopy types, Cahiers de topologie et g\'eometrie 
differ\'entielle, Vol. 32(1), 29-46 (1991).
\bibitem[Kapulkin2021]{kapulkin2021} Krzysztof Kapuklin, Peter LeFanu Lumsdaine: 
The simplicial model of univalent foundations (after Voevdosky), Journal of the 
European Mathematical Society, Vol. 23, 2071-2126 (2021).
\bibitem[Kleene1952]{kleene1952} Stephen C. Kleene: Metamathematics, North 
Holland Publishing Company (1952). 
\bibitem[Kripke1959]{kripke1959} Saul A. Kripke: A completeness theorem in 
modal logic, The Journal of Symbolic Logic, Vol. 24(1), 1-14 (1959). 
\bibitem[Kripke1965]{kripke1965} Saul A. Kripke: Semantical analysis of 
intuitionistic logic I, Studies in Logic and the Foundations of Mathematics, 
Vol. 40, 92-130 (1965).  
\bibitem[Kripke1975]{kripke1975} Saul A. Kripke: Outline of a theory of 
truth, The Journal of Philosophy, Vol. 72, 690-716 (1975). 
\bibitem[Kripke1980]{kripke1980} Saul A. Kripke: Naming and necessity, 
Blackwell (1980). 
\bibitem[Kuhn1962]{kuhn1962} Thomas S. Kuhn: The structure of scientific 
revolutions, University of Chicago Press (1962).
\bibitem[Kutyniok2023]{kutyniok2023} Gitta Kutyniok: The mathematics of 
artificial intelligence, Proceedings of the ICM 2022, EMS Press, 5118-5139 
(2023). 
\bibitem[LambScott1986]{lambekscott1986} Joachim Lambek, Philip J. Scott: 
Introduction to higher-order categorical logic, Studies in Advanced 
Mathematics, Vol. 7, Cambridge University Press (1988).
\bibitem[Lawvere1964]{lawvere1964} William Lawvere: An elementary theory of 
the category of sets, Proceedings of the National Academy of Sciences of the 
United States of America, Vol. 52, 1506-1511 (1964).
\bibitem[Leibniz]{leibniz} Gottfried Wilhelm Leibniz: Akademieausgabe der 
Leibniz Edition (digitally available at www.leibnizedition.de).
\bibitem[Leibniz1989]{leibniz1989} Gottfried Wilhelm Leibniz: Philosophical 
papers and letters (Ed. L. E. Roemker), Kluwer Acad. Publ. (1989).
\bibitem[Licata2013]{licata} Daniel R. Licata, Michael Shulman: Calculating the 
fundamental group of the circle in homotopy type theory, in: ACM/IEEE Computer 
Symposium on Logic in Computer Science (LICS 2013), 223-232 (2013).
\bibitem[Linnebo2017]{linnebo} {\O}ystein Linnebo: Philosopy of 
mathematics, Princeton University Press (2017). 
\bibitem[Llull1985]{lullus1985} Ram\'on Llull: Die neue Logik, Philosophische 
Bibliothek, Band 379, Felix Meiner Verlag (1985).
\bibitem[Llull1290]{lullus2001} Ram\'on Llull: Ars brevis, Philosophische 
Bibliothek, Band 518, Felix Meiner Verlag (2001).
\bibitem[Locke1690]{locke1690} John Locke: An essay concerning humane 
understanding, Thomas Basset and Edward Mory (1690).
\bibitem[Lolli2017]{lolli2017} Gabriele Lolli: Ambiguit\`a, un viaggio fra 
letteratura e matematica, Il Mulino (2017).
\bibitem[Lucas1961]{lucas1961} John R. Lucas: Minds, machines and G\"odel, 
Philosophy, Vol. 36, 112-127 (1961).
\bibitem[Lumsdaine2010]{lumsdaine2010} Peter LeFanu Lumsdaine: Weak 
$\omega$-categories from intensional type theory, Logical Methods in Computer 
Science, Vol. 6, 1-19 (2010).
\bibitem[Lurie2008]{lurie2008} Jacob Lurie: What is an $\infty$-category?, 
Notices of the American Mathematical Society, Vol. 55(8), 949-950 (2008).
\bibitem[Lurie2009]{lurie2009} Jacob Lurie: Higher topos theory, Annals of 
Mathematical Studies, Vol. 170, Princeton University Press (2009).  
\bibitem[MacColl1906]{maccoll1906} Hugh MacColl: Symbolic logic and its 
applications, Longmans, Green, and Co. (1906). 
\bibitem[MacLaneMo1992]{maclane_moerdijk1992} Saunders MacLane, 
Ieke Moerdijk: Sheaves in geometry and logic, Springer Universitext (1992). 
\bibitem[Manin2010]{manin} Yuri I. Manin: A course in mathematical logic for 
mathematicians, 2nd edition, Graduate Texts in Mathematics, Vol. 53, Springer 
Verlag (2010).
\bibitem[MartinL\"of1984]{martin-loef1984} Per Martin-L\"of: Intuitionistic 
type theory, Bibliopolis (1984).
\bibitem[MartinL\"of1994]{martin-loef1994} Per Martin-L\"of: Analytic and 
synthetic judgements in type theory, in: Kant and contemporary epistemology, 
Kluwer Academic Publishers, 87-99 (1994).
\bibitem[McCullPitts1943]{cullochpitts} Warren McCulloch, William Pitts: 
A logical calculus of the ideas immanent in nervous activity, Bulletin of 
Mathematical Biophysics, Vol. 5, 115-133 (1942).
\bibitem[McDaniel2017]{mcdaniel} Kris McDaniel: The fragmentation of being, 
Oxford University Press (2017). 
\bibitem[deMorgan1847]{demorgan1847} Augustus De Morgan: Formal logic, Taylor 
and Walton (1847).
\bibitem[Mostowski1951]{mostowski1951} Andrzej Mostowski: A classification of 
logical systems, Studia Philosophica, Vol. 4, 237-274 (1951). 
\bibitem[MSPiontkow2011]{mspiontkow2011} Stefan M\"uller-Stach, Jens 
Piontkowski: Elementare und algebraische Zahlentheorie, 2. Auflage, 
Vieweg+Teubner Verlag (2011).
\bibitem[MS2013]{ms2013} Stefan M\"uller-Stach: What is a period?, Notices of 
the American Mathematical Society, Vol. 61(8), 898-899 (2013). 
\bibitem[MS2023]{ms2023} Stefan M\"uller-Stach: Richard Dedekind: 
,,Was sind und was sollen die Zahlen?'' und ,,Stetigkeit und 
Irrationalzahlen'', 2. Auflage, Springer Verlag (2023).
\bibitem[MS2024]{ms2024} Stefan M\"uller-Stach: Max Dehn, Axel Thue, and the 
undecidable, in: Max Dehn -- Polyphonic Portrait (Jemma Lorenat, John McCleary, 
Volker Remmert, David Rowe, Marjorie Senechal eds.), History of 
Mathematics, Vol. 46, American Mathematical Society (2024). 
\bibitem[MS2027]{ms2027} Stefan M\"uller-Stach: Emmy Noether und 
ihre Bedeutung f\"ur die moderne Mathematik, in: Wie kommt das Neue 
in die Welt? (Hrsg. Mechthild Koreuber) (2027).
\bibitem[Nagel2023]{nagel} Thomas Nagel: Moral feelings, moral reality, and 
moral progress, Oxford University Press (2023).
\bibitem[Noether1983]{noether} Emmy Noether: Gesammelte Abhandlungen, Springer 
Verlag (1983). 
\bibitem[Peano1957]{peano} Giuseppe Peano: Opere scelte, Edizioni Cremonese 
(1957).
\bibitem[Peirce1878]{peirce} Charles S. Peirce: How to make our ideas clear, 
Popular Science Monthly, Vol. 12, 286-302 (1878).
\bibitem[Penrose1989]{penrose1989} Roger Penrose: The emperor's new mind, 
Oxford University Press (1989). 
\bibitem[Penrose1994]{penrose1994} Roger Penrose: Shadows of the mind, 
Oxford University Press (1994). 
\bibitem[Penrose2006]{penrose2006} Roger Penrose: A road to reality: a complete 
guide to the laws of the universe, Vintage Books (2006).
\bibitem[P\'eter1936]{peter} R\'ozsa P\'eter: \"Uber die mehrfache Rekursion, 
Math. Annalen, Vol. 113, 489-527 (1936).
\bibitem[Plato2020]{vonplato2020} Jan von Plato: Can mathematics be proved 
consistent?, Sources and Studies in the History of Mathematics and Physical 
Sciences, Springer Verlag (2020).
\bibitem[Plinius]{plinius} Plinius der \"Altere: Naturalis historia (circa 77 
A.D).
\bibitem[Poincar\'e1902]{poincare1902} Henri Poincar\'e: La science et 
l'hypoth\'ese, \'Editions Flammarion (1902).
\bibitem[Polya1945]{polya1945} George P\'olya: How to solve it, Princeton 
University Press (1945). 
\bibitem[Post1936]{post1936} Emil Post: Finite combinatory processes, The 
Journal of Symbolic Logic, Vol. 1, 103-105 (1936).
\bibitem[Post1947]{post1947} Emil Post: Recursive unsolvability of a problem 
of Thue, The Journal of Symbolic Logic, Vol. 12, 1-11 (1947). 
\bibitem[Quillen1967]{quillen} Daniel Quillen: Homotopical algebra, 
Lecture Notes in Mathematics, Vol. 43, Springer Verlag (1967).
\bibitem[Quine1951]{quine1951} Willard van Orman Quine: The two dogmas of 
empiricism, The Philosophical Review, Vol. 60, 20-43 (1951).
\bibitem[RadeToep1930]{rademachertoeplitz} Hans Rademacher, Otto 
Toeplitz: Von Zahlen und Figuren, Springer Verlag (1930).
\bibitem[RiehlVer2022]{riehlverity} Emily Riehl, Dominic Verity: Elements 
of $\infty$-category theory, Studies in Advanced Mathematics, Vol. 194, 
Cambridge University Press (2022).
\bibitem[Riemann1990]{riemann} Bernhard Riemann: Gesammelte Werke, Springer 
Verlag (1990). 
\bibitem[Rijke2025]{rijke} Egbert Rijke: Introduction to homotopy type theory, 
Studies in Advanced Mathematics, Vol. 219, Cambridge University Press (2025).
\bibitem[Rosenblatt1958]{rosenblatt1958} Frank Rosenblatt: The perceptron -- a  
probabilistic model for information storage and organization in the brain, 
Psychological Review, Vol. 65(6), 386-408 (1958). 
\bibitem[Rose1984]{rose} Harvey E. Rose: Subrecursion -- functions and 
hierarchies, Oxford Logic Guides, Vol. 9 (1984).
\bibitem[Rosling2019]{rosling} Hans Rosling: Factfulness, Ullstein Taschenbuch 
(2019).
\bibitem[Rosser1936]{rosser1936} Barkley Rosser: Extensions of some theorems of 
Gödel and Church, Journal of Symbolic Logic, Vol. 1(3), 87-91 (1936).
\bibitem[Rovelli2019]{rovelli} Carlo Rovelli: Die Geburt der 
Wissenschaft -- Anaximander und sein Erbe, Rowohlt Verlag (2019).
\bibitem[RussWhite1910]{russell} Bertrand Russell, Alfred North Whitehead: 
Principia mathematica, Cambridge University Press (1910).
\bibitem[Scholze2021]{scholze2021} Peter Scholze: Liquid tensor 
experiment, Experimental Mathematics, 1-6 (2021).
\bibitem[Schrenk2017]{schrenk} Markus Schrenk (Hrsg.): Handbuch Metaphysik, 
Metzler Verlag (2017).
\bibitem[Schr\"oder1890]{schroeder1890} Ernst Schr\"oder: Algebra der Logik, 
drei B\"ande, Teubner Verlag (1890-1905).
\bibitem[Sch\"utte2020]{schuette2020} Kurt Sch\"utte, Helmut 
Schwichtenberg: Mathematische Logik, in: The Legacy of Kurt Sch\"utte, Springer 
Verlag, 71-91 (2020).
\bibitem[Schwichten1971]{schwichtenberg} Helmut Schwichtenberg: Eine 
Klassifikation der $\epsilon_0$-rekursiven Funktionen, Zeitschrift f\"ur 
mathematische Logik und Grundlagen der Mathematik, Band 17, 61-74 (1971).
\bibitem[Searle1992]{searle1992} John Searle: The rediscovery of the mind, 
MIT Press (1992).
\bibitem[Shor1997]{shor} Peter W. Shor: Polynomial-time algorithms for 
prime factorization and discrete logarithms on a quantum computer, SIAM Journal 
on Computing, Vol. 26, 1484-1509 (1997).
\bibitem[Shulman2018]{shulman2018} Mike Shulman: Towards elementary  
$\infty$-toposes, Talk at Voevodsky Memorial Conference (2018).  
\bibitem[Shulman2019]{shulman2019} Mike Shulman: All $(\infty,1)$-toposes
have strict univalent universes, Preprint arXiv (2019).
\bibitem[Sieg2013]{sieg} Wilfried Sieg: Hilbert's programs and beyond, Oxford 
University Press (2013).
\bibitem[Siegel1966]{siegel} Carl L. Siegel: Gesammelte Abhandlungen, Springer 
Verlag (1966). 
\bibitem[Siegelmann1999]{siegelmann1999} Hava T. Siegelmann: Neural networks 
and analog computation, Birkh\"auser Verlag (1999).   
\bibitem[Skolem1970]{skolem} Thoralf Skolem: Selected works in logic, 
Universitetsforlaget Oslo (1970).
\bibitem[Smullyan2013]{smullyan2013} Raymond M. Smullyan: Truth and 
provability, The Mathematical Intelligencer, Vol. 35(1), 21-24 (2013).
\bibitem[Soul\'e2002]{soule} Christophe Soul\'e: The work of Vladimir 
Voevodsky, Proceedings of the ICM 2002, Vol. 1, 99-103 (2002).
\bibitem[Stepien2017]{stepien2017} {\L}ukasz T. St\k{e}pie\'n, Tomasz J.  
St\k{e}pie\'n: On the consistency of the arithmetic system, Journal of 
Mathematics and Systems Science, Vol. 7, 43-55 (2017).
\bibitem[Streicher1991]{streicher1991} Thomas Streicher: Semantics of type 
theory, Birkh\"auser Verlag (1991).
\bibitem[Struik1976]{struik} Dirk J. Struik: Abriss einer Geschichte der 
Mathematik, Deutscher Verlag der Wissenschaften (1976).
\bibitem[Tait2006]{tait} William Walker Tait: G\"odel's interpretation 
of intuitionism, Philosophia Mathematica Vol. 14(2), 308-338 (2006).
\bibitem[Takeuti2013]{takeuti2013} Gaisi Takeuti: Proof theory, 2nd edition, 
Dover (2013).
\bibitem[Tapp2013]{tapp} Christian Tapp: An den Grenzen des Endlichen, 
Springer Verlag (2013).
\bibitem[Tarski1931]{tarski1931} Alfred Tarski: Sur les ensembles 
d\'efinissables de nombres r\'eels I, Fundamenta Mathematicae, Vol. 17, 
210-239 (1931).
\bibitem[Tarski1933]{tarski1933} Alfred Tarski: Poj\k{e}cie prawdy w 
j\k{e}zykach nauk dedukcyjnych, Nak\l adem Towarzystwa Naukowego 
Warszawskiego (1933).
\bibitem[Tarski1935]{tarski1935} Alfred Tarski: Der Wahrheitsbegriff in den 
formalisierten Sprachen, Studia Philosophica, Vol. 1, 261-405 (1935). 
\bibitem[Tarski1938]{tarski1938} Alfred Tarski: Der Aussagenkalk\"ul und die 
Topologie, Fundamenta Mathematicae, Vol. 31, 103-134 (1938).
\bibitem[Tarski1969]{tarski1969} Alfred Tarski: Truth and proof, Scientific 
American, 63-77 (1969).
\bibitem[Tasi\'c2012]{tasic2012} Vladimir Tasi\'c: Poststructuralism and 
deconstruction: a mathematical history, Cosmos and History: The Journal of 
Natural and Social Philosophy, Vol. 8(1), 177-189 (2012). 
\bibitem[Thue1977]{thue} Axel Thue: Selected mathematical papers, 
Universitetsforlaget Oslo (1977).
\bibitem[Tierney1972]{tierney1972} Myles Tierney: Sheaf theory and the 
continuum hypothesis, Lecture Notes in Mathematics, Vol. 274, 13-42 (1972).
\bibitem[Turing1936]{turing1936} Alan M. Turing: On computable numbers with an 
application to the Entscheidungsproblem, Proceedings of the London Mathematical 
Society, Vol. 42, 230-265 (1936).
\bibitem[Turing1939]{turing1939} Alan M. Turing: Systems of logic based on 
ordinals, Proceedings of the London Mathematical Society, Vol. 45(2), 161-228 
(1939).
\bibitem[Turing1950]{turing1950} Alan M. Turing: Computing machinery and 
intelligence, Mind, Vol. LIX(236), 433-460 (1950).
\bibitem[Valiant2013]{valiant2013} Leslie Valiant: Probably approximately 
correct, Basic Books (2013).
\bibitem[Voevodsky2006]{voevodsky2006} Vladimir Voevodsky: A very short note 
on the homotopy $\lambda$-calculus, Preprint (2006). 
\bibitem[Voevodsky2013]{voevodsky2013} Vladimir Voevodsky et al.: Homotopy 
type theory: univalent foundations of mathematics, Institute for Advanced Study 
Princeton (2013).
\bibitem[Voevodsky2014]{voevodsky2014} Vladimir Voevodsky: The origins and 
motivations of univalent foundations, Institute Letter, Institute for Advanced 
Study (2014).
\bibitem[Voevodsky2015]{voevodsky2015} Vladimir Voevodsky: A $C$-system 
defined by a universe category, Theory and Applications of Categories, Vol. 
30, 1181-1214 (2015). 
\bibitem[Voevodsky2023]{voevodsky2023} Vladimir Voevodsky: Martin-L\"of 
identity types in $C$-systems, Publications math\'ematiques de l'IH\'ES, Vol. 
138(1), 1-67 (2023).
\bibitem[Vossenkuhl2013]{vossenkuhl} Wilhelm Vossenkuhl: Frege -- Das Problem 
der Wahrheit, \enquote{SWR2 Wissen: Aula} radio contribution (2013). 
\bibitem[Wadler2015]{wadler} Philip Wadler: Propositions as types, 
Communications of the ACM, Vol. 58(12), 75-84 (2015).
\bibitem[Wardhaugh2022]{wardhaugh} Benjamin Wardhaugh: Begegnungen mit Euklid, 
Harper Collins (2022).
\bibitem[Weyl1968]{weyl1968} Hermann Weyl: Gesammelte Abhandlungen, Springer 
Verlag (1968).
\bibitem[Wiener1950]{wiener} Norbert Wiener: Human use of human beings -- 
cybernetics and society, The Riverside Press (1950).
\bibitem[Wigderson2019]{wigderson2019} Avi Wigderson: Mathematics and 
computation, Princeton University Press (2019).
\bibitem[Wigner1960]{wigner} Eugene Wigner: The unreasonable effectiveness
of mathematics in the natural sciences, Communications in Pure and Applied 
Mathematics, Vol. 13(1), 1-14 (1960).
\bibitem[Wilkins1668]{wilkins1668} John Wilkins: An essay towards a real 
character und a philosophical language, Royal Society (1668).
\bibitem[Willaschek2015]{willaschek2015} Marcus Willaschek: Der mentale Zugang 
zur Welt, Klostermann Verlag (2015). 
\bibitem[Willaschek2023]{willaschek2023} Marcus Willaschek: Kant -- Die 
Revolution des Denkens, C. H. Beck Verlag (2023) 
\bibitem[Wittgenstein1922]{wittgenstein1922} Ludwig Wittgenstein: 
Tractatus logico-philosophicus, Paul Kegan (1922).
\bibitem[YangMills1954]{yangmills1954} Chen N. Yang, Robert L. Mills: 
Conservation of isotopic spin and isotopic gauge invariance, Physical 
Review, Vol. 96, 191-195 (1954). 
\end{thebibliography}
\end{document}